%% file: Projekt.tex
\documentclass[twoside,11pt,a4paper]{article}

\input{MeinStyle/MeinStyle}
\newcommand{\Schalter}[1]{#1  %Schalter ein oder ausschalten!
}%Schalter
\newcommand{\SchalterZ}[2]{#1 %
}%Schalter

\input{StyleFiles/Style}

\begin{document}

%\color[rgb]{0.0.1,0.01,0.01} %Falls die Farbe Schwarz leer ist.

%Für Doppelseitiger Ausdruck bitte einschalten:
%~\thispagestyle{empty} \newpage \pagenumbering{arabic}

\newcommand{\MitRechnungen}[1]{#1 %Schaltet die Rechnungen ein.
}%MitRechnugen

\newcommand{\MitText}[1]{#1 %Schaltet den Textein ein.
%\newpage
}%MitText

\newcommand{\Frage}[1]{
%\Bem{Frage:}{\begin{flushright} \bfseries{ #1 ??? } \end{flushright}} 
}

\Schalter{
\Frage{Wie geht Tiny im Mathemodus? Siehe Stylefile:StyleReduktion  sricptstile sricptstile}
\Frage{Wie schreibe ich schiefe Bruchstriche :$\backslash$ Bruchstrich rizebox?}
}

\newcommand{\Bruchstrich}{\Big/}

%Hier kommt das Document

% \input{EDA-Beweis}

%\input{Titel/Titel-Pruefungsamt}
\input{Titel/Titel}

\SchalterZ{ 
  %\addtolength{\oddsidemargin}{-10mm}
  %\addtolength{\evensidemargin}{-20mm}
\input{Haupt/Reduktionssysteme/Reduktionssysteme} \newpage 

}{}%SchalterZ

\Schalter{

 \chapter{Konvergentes Reduktionssystem für $\fQG(n)$}
 \label{sec:Groebnerbasis}

  \newpage
  \input{Haupt/RechnungGBA} \newpage 
  \input{Haupt/EDA-Beweis}
 %\chapter{Freie Auflösung für $A^n$}
 %\label{sec:FreienAufloesung} \newpage

\input{Haupt/BiAufloesung}

\input{Haupt/HochschildAufloesung/HochschildAufloesung}

\input{Haupt/Ext}
}

\Schalter{
}

\SchalterZ{}{
 %Kapitel die Gerade Bearbeitet werden.

\input{Bezeichner}

}%SchalterZ

\newpage
\Frage{
   Fehlende Literatur:
   Bir35,  (Bichon, Julien) 
 }

   \nocite{MR2274824, MR1129886,Avenhaus,MR1213453,Newman,MR1227656,MR1287608,MR1714602,MR1189133,MR1714602,   				 ThomCollins,Rudin,  Magma1997,GLS03,GPS05,GAP4,Lueck, MR0255472}

   %Begleitende Literatur 
   \nocite{MR1034428,MR1422274,MR1868810,MR1053466,MR1089685,MR0038958,MR776465,MR533813,MR1511645}

%\Schalter{
  \bibliographystyle{alpha}

  \bibliography{MeinBib/MeinBib}
%}
\end{document}

%% file: MeinStyle/MeinStyle.tex
\usepackage{ifpdf}
\usepackage{ifthen}
\newboolean{mpdf}\setboolean{mpdf}{true}
\newboolean{deutsch}\setboolean{deutsch}{true}

\usepackage{MeinStyle/myarticlestyle}
\usepackage{MeinStyle/meinMathe}
\usepackage{lscape}

\usepackage{amscd}

\usepackage[arrow, matrix, all]{xy}

\newcommand{\kleinerText}[1]{ \scriptsize #1 \normalsize}

%% file: StyleFiles/Style.tex
\newcommand{\bs}[1]{\boldsymbol{#1}}

\input{StyleFiles/StyleRechnungen}

\input{StyleFiles/StyleReduktion}

\input{StyleFiles/StyleTorExt}

\newcommand{\Aev}{\fQG^{e}}

\newcommand{\fQG}{\mathcal{A}}

\newcommand{\Det}[1]{Det\left(#1 \right)}

\newcommand{\BlockIoffen}[1]{\begin{array}{|c|} #1 \\ \end{array} }

\newcommand{\BlockI}[1]{\begin{array}{|c|}\hline #1 \\ \hline \end{array} }
\newcommand{\BlockII}[1]{\begin{array}{|cc|}\hline #1 \\ \hline \end{array} }
\newcommand{\BlockIV}[1]{\begin{array}{|cccc|}\hline #1 \\ \hline \end{array} }

\newcommand{\MatheAbbildung}[4]{\mathe{\Matrix{#1 & \rightarrow& #3 \\ #2 & \mapsto &#4}}}
\newcommand{\Tor}[1]{Tor^{\fQG}\left( #1 \right)}
\newcommand{\Ext}[1]{Ext_{\fQG}\left( #1 \right)}

\newcommand{\T}[1]{\textnormal{#1} }

\newcommand{\MonomOrdnung}{<_{\Monome}}

\newcommand{\AufzaehlungZ}[1]{\begin{enumerate} #1 \end{enumerate} }
\newcommand{\AufzaehlungP}[1]{\begin{itemize}	#1 \end{itemize} }

\newcommand{\LMatrix}[1]{~ \begin{array}{llllllllllllll}  #1   \end{array}~ }

\newcommand{\Menge}[2]{ \left\{  \begin{matrix} #1 \end{matrix} ~\left|~ \begin{matrix}  #2   \end{matrix}  \right. \right\} }

\newcommand{\Sum}[2]{\sum\limits_{\textnormal{ \tiny $\begin{matrix}   #1    \end{matrix}$ }}^{\textnormal{ \tiny $\begin{matrix}  #2   \end{matrix}$ }}}
\newcommand{\Sumn}[1]{\Sum{#1}{n}}
\newcommand{\RSum}[2]{\sum\limits_{\textnormal{ \tiny $\begin{matrix}   #1    \end{matrix}$ }}^{\textnormal{ \tiny $\begin{matrix}  #2   \end{matrix}$ }}{'}}

\newcommand{\HH}{\mathcal{HH}}

%% file: StyleFiles/StyleRechnungen.tex
\newcommand{\DELTA}{\delta}
\newcommand{\Erzeuger}[2]{a_{#1,#2}}
\newcommand{\Zeile}{Z}
\newcommand{\Spalte}{S}

\newcommand{\Zi}[2]{-\sum_{i=2}^n \Erzeuger{#1}{i}\Erzeuger{#2}{i} + \DELTA_{#1,#2}}

\newcommand{\mZi}[2]{\sum_{i=2}^n \Erzeuger{#1}{i}\Erzeuger{#2}{i} - \DELTA_{#1,#2}}
\newcommand{\mZj}[2]{\sum_{j=2}^n \Erzeuger{#1}{j}\Erzeuger{#2}{j} - \DELTA_{#1,#2}}

\newcommand{\Zii}[3]{-\Erzeuger{#1}{1}\sum_{i=3}^n \Erzeuger{#2}{i}\Erzeuger{#3}{i} + \Erzeuger{#1}{1}\DELTA_{#2,#3}}

\newcommand{\mZii}[3]{\Erzeuger{#1}{1}\sum_{i=3}^n \Erzeuger{#2}{i}\Erzeuger{#3}{i} - \Erzeuger{#1}{1}\DELTA_{#2,#3}}

\newcommand{\Si}[2]{-\sum_{i=2}^n \Erzeuger{i}{#1}\Erzeuger{i}{#2} + \DELTA_{#1,#2}}
\newcommand{\mSi}[2]{\sum_{i=2}^n \Erzeuger{i}{#1}\Erzeuger{i}{#2} - \DELTA_{#1,#2}}

\newcommand{\mSj}[2]{\sum_{j=2}^n \Erzeuger{j}{#1}\Erzeuger{j}{#2} - \DELTA_{#1,#2}}

\newcommand{\Sii}[3]{-\Erzeuger{1}{#1}\sum_{i=3}^n \Erzeuger{i}{#2}\Erzeuger{i}{#3} + \Erzeuger{1}{#1}\DELTA_{#2,#3}}

\newcommand{\mSii}[3]{\Erzeuger{1}{#1}\sum_{i=3}^n \Erzeuger{i}{#2}\Erzeuger{i}{#3} - \Erzeuger{1}{#1}\DELTA_{#2,#3}}
\newcommand{\mSjj}[3]{\Erzeuger{1}{#1}\sum_{j=3}^n \Erzeuger{j}{#2}\Erzeuger{j}{#3} - \Erzeuger{1}{#1}\DELTA_{#2,#3}}

%% file: StyleFiles/StyleReduktion.tex
\newcommand{\rot}[1]{
{#1} }

\newcommand{\MengeE}{\rot{\mathcal{E}}}
\newcommand{\Red}{\rot{\mathcal{R}}}
\newcommand{\Nf}{\rot{\mathcal{N}} }

\newcommand{\Weg}{\rot{\Diagramm{  \ar@{.>}[r] & }}}
\newcommand{\RedWeg}{\rot{\Diagramm{  \ar@{.>}[r]_>{\MatheTiny{\Red}} & }}}
\newcommand{\RRedWeg}{\rot{\Diagramm{  \ar@{.>}[r]_>{\MatheTiny{\RRed}} & }}}
\newcommand{\RRedWegr}{\rot{\Diagramm{  \ar@{.>}[r]_>{\MatheTiny{\Rel}} & }}}
\newcommand{\RRedWegF}{\rot{\Diagramm{  \ar@{.>}[r]_>{\MatheTiny{\RelF}} & }}}
\newcommand{\RedpWeg}{\rot{\Diagramm{  \ar@{.>}[r]_>{\MatheTiny{\Red_{\P}}} & }}}
\newcommand{\RRedpWeg}{\rot{\Diagramm{  \ar@{.>}[r]_>{\MatheTiny{\RRed_{\P}}} & }}}

\newcommand{\AqRedWeg}{\rot{\Diagramm{  \ar@{.>}[r]_>{\MatheTiny{\Aq{\Red}}} & }}}
\newcommand{\regel}{\rot{V}}
\newcommand{\Regel}{\rot{\Diagramm{  \ar@{->}[r] & }}}
\newcommand{\redRegel}{\rot{\Diagramm{  \ar@{->}[r]_>{\MatheTiny{\red}} & }}}
\newcommand{\RedRegel}{\rot{\Diagramm{  \ar@{->}[r]_>{\MatheTiny{\Red}} & }}}
\newcommand{\RRedRegel}{\rot{\Diagramm{  \ar@{->}[r]_>{\MatheTiny{\RRed}} & }}}
\newcommand{\RRedRegelF}{\rot{\Diagramm{  \ar@{->}[r]_>{\MatheTiny{\RelF}} & }}}
\newcommand{\RRedRegelr}{\rot{\Diagramm{  \ar@{->}[r]_>{\MatheTiny{\Rel}} & }}}
\newcommand{\RedpRegel}{\rot{\Diagramm{  \ar@{->}[r]_>{\MatheTiny{\Red_{\P}}} & }}}
\newcommand{\RRedpRegel}{\rot{\Diagramm{  \ar@{->}[r]_>{\MatheTiny{\RRed_{\P}}} & }}}

\newcommand{\RedAq}{\rot{\Diagramm{   \ar@{<.>}[r]_>>>>{\MatheTiny{\Red}} & }}}

\newcommand{\RRedfAq}{\rot{\Diagramm{   \ar@{<.>}[r]_>>>>{\MatheTiny{\RRed_f}} & }}}
\newcommand{\RRedaqAq}{\rot{\Diagramm{   \ar@{<.>}[r]_>>>>{\MatheTiny{\Aq{\Red_f}}} & }}}

\newcommand{\red}{\rot{r}}

\newcommand{\Alphabet}{\rot{\mathcal{A}}}

\newcommand{\Monome}{\rot{\mathcal{M}}}
\newcommand{\MonomeM}{\Monome}
\newcommand{\Terme}{\rot{\mathcal{T}}}

\newcommand{\Algebra}{\rot{\ensuremath{\boldsymbol{A} } }}

\newcommand{\Rig}{ \rot{ \bs{R} }}
\newcommand{\RigF}{ \rot{ \widetilde{\mathcal{F}} }}
\newcommand{\RRed}{ \rot{ \mathcal{\widetilde{R}} }}

\newcommand{\FreieAM}{\rot{ \mathcal{F}_{\AnzA + \AnzM} } }
\newcommand{\FreieAlgebra}{\rot{ {\mathcal{F}} }}
\newcommand{\FreieA}{\FreieAlgebra }
\newcommand{\FreieAA}{\rot{ \mathcal{F}_{\AnzA} }}

\newcommand{\Modul}{\ensuremath{\rot{\boldsymbol{M}}} \xspace }
\newcommand{\AnzA}{\ensuremath{\rot{\boldsymbol{n_{\Algebra}}}}\xspace}
\newcommand{\AnzM}{\ensuremath{\rot{\boldsymbol{n_{\Modul}} }}\xspace}
\newcommand{\AnzUM}{\ensuremath{\rot{\boldsymbol{n_{\UModul}} }}\xspace}

\newcommand{\RSm}{\rot{\Regelsystem_{\Modul} }}

\newcommand{\UModul}{\ensuremath{\rot{\boldsymbol{\overline{M}}}\xspace }}

\newcommand{\AoA}{\ensuremath{\rot{\Algebra \tensor \Algebra^{op}}} }

\newcommand{\Rel}{\rot{\textnormal{Rel}}}
\newcommand{\RelF}{\widetilde{\Rel}}

\renewcommand{\P}{\rot{\mathcal{P}}}

\newcommand{\wahr}{\textnormal{\bfseries{wahr}}}
\newcommand{\falsch}{\textnormal{\bfseries{falsch}}}

\newcommand{\Rp}{\rot{\Red_{\P}}}
\newcommand{\Aq}[1]{\rot{ [ #1 ]}}

\newcommand{\RpAq}{\rot{\sim}}

\newcommand{\F}{\ensuremath{\mathcal{F}} }

\newcommand{\A}{\ensuremath{A} }
\newcommand{\p}{\ensuremath{p} }

\newcommand{\MatheTiny}[1]{#1}

\newcommand{\nachAqR}{\Diagramm{ \ar@{->}[r]_>{\MatheTiny{\Aq{\Red}}}&} }

\newcommand{\WegnachAqR}{\Diagramm{ \ar@{.>}[r]_>{\MatheTiny{\Aq{\Red}}}&} }
\newcommand{\WegnachR}{\Diagramm{ \ar@{.>}[r]_>{\MatheTiny{{\Red}}}&} }

\newcommand{\ErzeugendenSystem}[1]{S_{#1}}
\newcommand{\ESa}{\ErzeugendenSystem{\Algebra}}

\newcommand{\Regelsystem}{\red}
\newcommand{\RS}{\Regelsystem }
\newcommand{\vRS}{\widetilde{\Regelsystem} }
\newcommand{\RSa}{\Regelsystem_{\Algebra} }
\newcommand{\RSam}{\Regelsystem_{\Algebra, \Modul} }

\newcommand{\vRSa}{\widetilde{\Regelsystem}_{\Algebra} }

\newcommand{\RSkern}{\Regelsystem_{\textnormal{kern}} }

%% file: StyleFiles/StyleTorExt.tex
\newcommand{\Rang}[1]{Rang\left(\Phi_{#1 *} \right)}
\newcommand{\RanG}[1]{Rang \left( \Phi_{#1}^* \right) } 
\newcommand{\Dim}[1]{Dim\left( #1 \right)}

%% file: Titel/Titel.tex
%\input{Titel/Deckblatt} \newpage
%~ \thispagestyle{empty} \pagebreak \newpage 
\title{Reduktionssysteme zur Berechnung einer Auflösung\\{\tiny ~\\} der orthogonalen freien Quantengruppen $\boldsymbol{\fQG_o(n)}$	}

\author{Johannes Härtel%, Ivan Yudin
\\ jhaertel@math.uni-goettingen.de}
\date{%2.Advent
\today
}
\maketitle

%\pagenumbering{roman}

%%%%%%%%%%%%%%%%Abstract
%\vfill
%\begin{abstract}

~\thispagestyle{empty}
%\chapter*{Vorwort}
\label{sec:Einleitung}

\input{Titel/Einleitung}

%\input{Titel/Danksagung}
%\end{abstract}

%\tableofcontents

%
%\newpage 
%~\thispagestyle{empty}%\newpage

%\newpage
%\newpage~\thispagestyle{empty}
%\newpage~\thispagestyle{empty}

%% file: Titel/Einleitung.tex
In diesem Artikel werden Reduktionssysteme vorgestelltmöchte ich eine Methode zur Berechnung eines Reduktionssystems für eine Algebra

Schon Mitte der sechziger Jahre des 20. Jahrhunderts hatte  Kac \cite{MR0156921} eine neue Klasse ("`Ring\-gruppen"') von mathe\-matischen Objekten zur Analyse der Pontryagindualität eingeführt, die später unter dem Begriff "`Quantengruppen"' klassifiziert wurde.

Etwa zwanzig Jahre später beschäftigte sich Woronowicz  \cite{MR890482} \cite{MR943923} mit "`Pseudogruppen"', als er eine wichtige Familie von Matrixquantengruppen einführte. Sie lieferte neuartige Deformationen der zuvor noch als rigide eingestuften Liegruppe $SU(2)$. Diese Quanten-$SU(2)$-Gruppen lieferten wichtige Beispiele für die Quantenphysik, in der $SU(2)$ beispielsweise als Isospingruppe in der schwachen Wechselwirkung auftritt.

Woronowicz erklärte dabei die Quanten-$SU(2)$-Gruppen über die Charaktere einer axiomatisch vorgegebenen Hopfalgebra. Entsprechend stellen die Quanten-$SU(2)$-Gruppen Beispiele für Räume in der nicht-kommutativen Geometrie von Connes \cite{MR823176} dar. Diese Konstruktionen entsprechen denen von Drinfeld \cite{MR869575} und Jimbo \cite{MR797001} (vgl. \cite{MR889754}).

Später, im Jahr 1995, definierten Wang \cite{MR1316765} und Van Daele \cite{MR1382726} die orthogonalen und die unitären freien Quantengruppen. Diese stellten sich als universell heraus, in dem Sinn, dass jede kompakte Matrixquantengruppe isomorph ist zu einem direkten Produkt von Unteralgebren der freien Quantengruppe. Dieses Paper zeichnet sich durch einen sehr konkreten, direkten Zugang zu Quantengruppen aus. Meine Arbeit bezieht sich auf Wangs Beschreibung von freien Quantengruppen.

Entsprechend der Universalität von Wang und Van Daeles Konstruktion ist die Homologie der beiden freien Quantengruppen von besonderem Interesse. Unter der Voraussetzung, dass eine bestimmte Sequenz exakt ist, haben Collins und Thom in \cite{ThomCollins} Schlüsse über die Homologien der orthogonalen und der unitären freien Quantengruppen gezogen, mit dem Ziel, ihre $\ell^2$-Bettizahlen zu bestimmen. 
Eine Übersicht zu orthogonalen Quantengruppen findet sich bei Collins und Banica \cite{MR2341011}.

Das wesentliche Problem, auf das Collins und Thom in ihrer Arbeit stoßen, ist, dass sie zwar explizit Erzeuger und Relationen der freien Quantengruppe, jedoch keine Basis kennen. In dieser Arbeit wird der Formalismus von Reduktionssystemen genutzt, um exemplarisch eine solche Basis zu berechnen und die Exaktheit von Collins und Thoms projektiver Auflösung zu beweisen.

Reduktionssysteme (synonym Gröbnerbasen oder Standardbasen) wurden 1964 und 1965 von Buchberger \cite{MR2202562}, Hironaka \cite{MR0199184} sowie Knuth und Bendix \cite{MR0255472} unabhängig voneinander beschrieben. 1978 zeigte Bergman \cite{Bergman}, dass %mit Hilfe von Newmans Diamantenlemma \cite{Newman}
 auch Gröbnerbasen bezüglich Moduln berechnet werden können.

Reduktionssysteme entwickelten sich zu einem wichtigen Werkzeug in der theoretischen Informatik, für Theorembeweise und funktionale Programmierung. Es gibt inzwischen sehr gute Software zur Berechnung von Reduktionssystemen im Kontext kommutativer Algebren (z.B. SINGULAR \cite{GPS05}), und auch in der Literatur werden vornehmlich kommutative Algebren behandelt. Auch für nicht-kommutative Algebren gibt es Software (z.B. PLURAL \cite{GLS03}, MAGMA \cite{Magma1997}, GAP \cite{GAP4}), allerdings beschränkt diese sich auf Spezialfälle. Zur Berechnung werden Vervollständigungsalgorithmen (siehe \cite{MR0463136} und \cite{MR0255472}) benutzt. Die mir bekannten Implementierungen reduzieren in jedem Schritt alle bekannten Regeln vollständig und nutzen  Symmetrien nicht aus, daher terminieren sie für die freien Quantengruppen nur für $n<10$.

In Kapitel 1 trage ich die Ergebnisse für nicht-kommutative Algebren zusammen und formuliere sie direkt unter Verwendung des Diamantenlemmas. Insbesondere wird gezeigt, dass sich die Verifizierung eines Reduktionssystems für den Kern eines Modulhomomorphismus vereinfachen lässt, falls für eine Algebra ein Reduktionssystem bekannt ist.

Im zweiten Kapitel wird ein vollständiges Reduktionssystem bezüglich der freien orthogonalen Quantengruppen an und verifiziert. Für den Fall $n=2$ sind es nur sehr wenige Regeln und es kann  ein endlicher Automat angegeben werden, der sämtliche Basiselemente darstellt:

\input{EDA}

In Kapitel 3 werden die Ergebnisse aus Kapitel 1 und 2 genutzt, um eine Basis für die Kerne einer freien Auflösung der freien Quantengruppe als Bimodul zu beweisen.

%Für $n\geq 3$ findet sich eine Ordnung 

Insbesondere im Fall von Collins und Thoms Sequenz liefert diese Methode die korrekten Kerne.

Abschließend wird im letzten Kapitel mit der nun verifizierten Auflösung die Homologie der orthogonalen freien Quantengruppen explizit berechnet.

%% file: EDA.tex
\begin{center}$
\entrymodifiers={++[o][F-]}
\SelectTips{cm}{}
\xymatrix @-1pc {
*\txt{} 
& *\txt{}
& *\txt{}
& \ar `d^r[rr] `^u[rr] _{a_{1,2}} [rr] 
& *\txt{}
& \ar `u^l[ll] `^d[ll] _{a_{2,1}} [ll] 
\\
*\txt{}
\\
*\txt{}
\\
*\txt{}
\\
\ar `d^r[r] `^u[r] _{a_{1,2}} [r] 
\ar `dl^r[ddr]  [ddrrrr] ^{a_{1,1}}
& \ar `u^l[l] `^d[l] _{a_{2,2}} [l] 
	\ar[uuuurr] _{a_{2,1}}
& *\txt{}
& *\txt{}
& S \ar@(ur,ul) _{a_{2,2}}
		\ar[lll] _{a_{1,2}}
		\ar[rrr] ^{a_{2,1}}
		\ar[dd] _{a_{1,1}}
& *\txt{}
& *\txt{}
& \ar `d^r[r] `^u[r] _{a_{2,2}} [r] 
	\ar[uuuull] ^{a_{1,2}}
	& \ar `u^l[l] `^d[l] _{a_{2,1}} [l] 
	\ar `dr_l[ddl]  [ddllll] _{a_{1,1}}
	&*\txt{}
\\
*\txt{}
\\
*\txt{}
&*\txt{}
&*\txt{}
&*\txt{}
& \ar  `dr[dd] [dd]  ^>>>>>>>>>>{a_{2,2}}
&*\txt{}
& *\txt{}
	&*\txt{}
		&*\txt{}
			&*\txt{}
				&*\txt{}	&*\txt{}	&*\txt{}
\\
*\txt{}
\\
*\txt{}
&*\txt{}
&*\txt{}
&*\txt{}
& \ar `ul[uu]^>{a_{1,1}} [uu]  
}
$
\end{center}

%% file: Haupt/Reduktionssysteme/Reduktionssysteme.tex
\chapter{Reduktionssysteme}
\label{chapter:Reduktionssysteme}

\input{Haupt/Reduktionssysteme/EinleitungReduktionssysteme}

\input{Haupt/Reduktionssysteme/abstrakteReduktionssysteme}

% alt \input{Haupt/Reduktionssysteme/DefPraedikat}
% alt \input{Haupt/Reduktionssysteme/DefPradikatPaar}

\input{Haupt/Reduktionssysteme/ReduktionssystemeAlgebraDef}

\input{Haupt/Reduktionssysteme/ReduktionssystemeAlgebraKonfluenz}

\input{Haupt/Reduktionssysteme/ReduktionssystemeAlgebraOrdnungen}

\input{Haupt/Reduktionssysteme/ReduktionssystemeAlgebraNoethersch}

\input{Haupt/Reduktionssysteme/ReduktionssystemeModul}

\input{Haupt/Reduktionssysteme/P-Vollstaendig}

\input{Haupt/Reduktionssysteme/ModulPraedikat-Neu}
\input{Haupt/Reduktionssysteme/kern-Neu}

\newpage

% alt \input{Haupt/Reduktionssysteme/ReduktionssystemeCokern}\newpage

%% file: Haupt/Reduktionssysteme/EinleitungReduktionssysteme.tex
Zum Rechnen in durch Gleichungen definierten Strukturen können Reduktionssysteme genutzt werden. Beispiele für solche Strukturen, die in dieser Arbeit betrachtet werden, sind Halbgruppen, Gruppen, Ringe, Algebren und Moduln.
In den folgenden Kapiteln werden für Wörter über einem festen Alphabet Reduktionssysteme, man spricht auch von Wortersetzungssystemen oder Semi-Thue-Systemen, genutzt, um Gröbnerbasen zu Idealen zu berechnen.

~\\
Weitere Anwendungsbereiche für Reduktionssysteme, die hier nicht weiter betrachtet werden, sind Termersetzungssysteme für das Rechnen in logischen Strukturen. Sie finden Anwendung bei effizienten Verfahren zum automatischen Theorem-Beweisen in der Prädikatenlogik und zur Beschreibung von abstrakten Datentypen, insbesondere für Korrektheitsbeweise für Programme in einer hohen Programmiersprache.
\Frage{Hier noch die Literatur angeben}

~\\
Das Problem in durch Gleichungen definierten Strukturen ist, dass Gleichungen in beide Richtungen angewendet werden können und dass dadurch der Suchraum auch für eine Computer gestützte Berechnung zu groß wird. Der Trick bei Reduktionssystemen besteht darin die Gleichungen nur in eine Richtung anzuwenden. Der Suchraum wird dadurch stark eingeschränkt und man erhält eine  mächtige Simplifikation.
%\Frage{was ist eine mächtige Simplifikation}

~\\
Ein Beispiel dafür ist das Wortproblem.
Sei $\mathcal{E}$ eine Menge und sei $\sim$ eine Äquivalenzrelation auf $\mathcal{E}$. Seien $v,w \in \mathcal{E}$. Die Frage, ob $v$ und $w$ äquivalent, also $v \sim w$, sind,  heißt Wortproblem.
~\\
Das Problem ist schon für sehr einfache Gruppen nicht auf den ersten Blick lösbar.

\Bsp{ Betrachten wir die Gruppe $G= < a,b ~|~ aba, b^2>$. Gilt hier $a^2\stackrel{}{\sim}b?$\footnote{$ab\sim ababa \sim ba$, also ist $ b\sim baba \sim bbaa \sim a^2$.}
}
~\\
Um das Wortproblem zu lösen wählt man zunächst eine wohldefinierte Partialordnung "`>"' und Reduktionsrelationen $\mapsto$ auf $\mathcal{E}$, so dass gilt:
\AufzaehlungP{
\item{Reduktion:} Jede Relation verkleinert, also $\forall w_1, w_2 \in \mathcal{E} $ mit $w_1 \mapsto w_2$ gilt: $w_1 \sim w_2$ und $w_1 > w_2$.
\item{Konvergenz:} Zu jedem Element gibt es genau eine Normalform, also für jede Äquivalenzklasse $[w]$ existiert ein Element $\hat{w} \in [w]$, so dass es für jedes $w \in [w]$ eine Folge von Relationen gibt mit $w \mapsto \dots \mapsto \hat{w}$.
}

Die Reduktionsrelationen erhält man, indem man die Äquivalenzrelationen durch ein Gleichungssystem beschreibt und dann die Relationen nimmt, die man durch Anwenden der Gleichungen von rechts nach links erhält. Dann startet man die Vervollständigung nach Knuth-Bendix (\cite{MR0255472}). Dieser Algorithmus versucht zu einem gerichteten Gleichungssystem ein äquivalentes konfluentes Reduktionssystem zu finden.

~\\

Dieses aus drei Abschnitten bestehende Kapitel basiert auf \cite{Avenhaus}, \cite{Newman}, \cite{MR2032182} und \cite{MR506890}.
Zunächst werden abstrakte Reduktionssysteme vorgestellt, danach Ersetzungssysteme für Algebren und Moduln. Im letzten Abschnitt wird insbesondere eine effiziente Methode zur Berechnung des Kerns einer A-Modul Abbildung vorgestellt.

%% file: Haupt/Reduktionssysteme/abstrakteReduktionssysteme.tex
\section{Abstrakte Reduktionssysteme}
\label{sec:Reduktionssysteme}

In diesem Abschnitt werden die für Reduktionssysteme nötigen Begriffe definiert und es wird ein Beweis für das Diamantenlemma angegeben.

~\\
Damit man leicht mit einem Reduktionssystem arbeiten kann, muss es die so genannte Church-Rosser Eigenschaft erfüllen. Der Nachweis dieser Eigenschaft ist sehr mühsam. Leichter ist Konfluenz und noch leichter lokale Konfluenz zu beweisen. Am Ende dieses Abschnittes werden wir mit dem Diamantenlemma zeigen, dass diese drei Eigenschaften unter gewissen Voraussetzungen gleich sind.

\subsection{Definitionen}
\label{sec:Definitionen}

Zunächst definieren wir für eine beliebige Menge $\MengeE$ ein Reduktionssystem $\Red$.
\label{DefMenge}

\Def{Reduktionssystem}{
\label{DefReduktionssystem}
\label{DefReduktionsregel}
Eine Teilmenge $\Red \subset \MengeE \times \MengeE$ heißt Reduktionssystem auf $\MengeE$.
\\Ein Element $(x,y) \in \Red$ heißt Reduktionsregel und wir schreiben dann:
\mathe{ x \RedRegel y.
			}
Wir sagen: "`$x$ wird zu $y$ reduziert"'. Wenn eindeutig ist, welches Reduktionssystem gemeint ist, schreiben wir statt $x \RedRegel y$ oft  nur $x \Regel y$.
}

Trotz der Schreibweise ist ein Reduktionssystem keine Abbildung zwischen Mengen. Es ist durchaus erlaubt, dass $x$ zu zwei verschiedenen $y_1, y_2$ reduziert werden kann.

~\\

\Def{Unreduzierbar}{
\label{DefUnreduzierbar}
Ein Element $z \in  \MengeE$ heißt unreduzierbar, falls es keine Reduktionsregel gibt, bei der $z$ auf der linken Seite steht, also:
\mathe{ \forall w \in \MengeE \textnormal{ gilt: } (z,w) \notin \Red . }
}

\Def{Reduktionsweg \RedWeg }{
\label{DefReduktionsweg}
Ein Reduktionsweg ist eine Folge von Reduktionsregeln $(x_i,y_i)_{i\in I} \in \Red$, so dass für alle $i\in I$ gilt:
\mathe{ y_i = x_{i+1} .}
Es gibt also einen Graphen:
\mathe{ x_1 \RedRegel y_1 \RedRegel \dots \RedRegel y_n  \RedRegel \dots}
Wir schreiben einen endlichen Reduktionsweg als $ x_1   \RedWeg y_n$.
Wenn eindeutig ist, welches Reduktionssystem gemeint ist, schreiben wir statt $x \RedWeg y$ oft  nur $x \Weg y$.
Ein Reduktionsweg kann leer, endlich oder unendlich sein.
}

~\\
Zu einem Element gibt es Reduktionswege, die von ihm ausgehen. Manche dieser Wege enden in einem unreduzierbaren Element. Dies motiviert folgende Definition.

\Def{Normalform}{
\label{DefNormalform}
Eine Normalform eines Elements $x \in \MengeE$ ist ein unreduzierbares Element $\Nf(x) \in  \MengeE$, so dass es einen endlichen (oder leeren)  Reduktionsweg gibt, der in $x$ beginnt und in $\Nf(x) $ endet, also :

\mathe{\Diagramm{
 x 	
 \ar@{.>}[dr] \\
   &\Nf(x)
}}%Diagramm

}

\Bem{Normalform}{Eine Normalform ist im Allgemeinen nicht eindeutig. In einem abstrakten Reduktionssystem kann es für ein Element verschiedene  Reduktionswege geben, die in verschiedenen unreduzierbaren Elementen enden.}

\Def{Church-Rosser-äquivalent}{
\label{DefChurchRosserAq}
Zwei Elemente $y_1, y_n$ heißen Church-Rosser-äquivalent bezüglich $\Red$, falls sie über ungerichtete Reduktionsregeln aus $\Red$ miteinander verbunden sind; es also eine endliche Folge $(y_i,y_{i+1})_{i=1\dots n-1}$ gibt, wobei  $(y_i,y_{i+1})\in \Red$ oder  $(y_{i+1},y_{i}) \in \Red$.

Für zwei Church-Rosser-äquivalente Elemente $x, \bar{x}$ schreiben wir 
$x \RedAq \bar{x}$.
}

\subsection{Konfluenz}
\label{sec:Eigenschaften}
Ziel ist es möglichst einfach zu erkennen, ob die Normalformen eindeutig sind. Dazu definieren wir einige Eigenschaften für Reduktionssysteme.

\Def{noethersch}{
Ein Reduktionssystem $\Red$ heißt noethersch oder terminierend, wenn jeder Reduktionsweg endlich ist.
}
\Bem{Normalform}{
Falls $\Red$ noethersch ist, dann hat jedes Element wenigstens eine Normalform.}

\Def{lokal konfluent}{Ein Reduktionssystem $\Red$ heißt lokal konfluent, falls es für je zwei Reduktionsregeln $(x,y_1), (x,\bar{y}_1)$, die vom selben Element $x$ starten, zwei endliche Wege $(y_i,y_{i+1})_{i=1 \dots n}$, $(\bar{y}_j,\bar{y}_{j+1})_{j=1\dots \bar{n}}$ gibt, die mit $y_1$ bzw. $\bar{y}_1$ beginnen und im selben Element $z$ enden, also zu jedem:
\mathe{
\Diagramm{
 &&&x 	
  \ar@{->}[dll]^>{\MatheTiny{}} \ar@{->}[drr] \\
 &  y   &&  && \bar{y} \\
   \textnormal{existiert:}  & y  \ar@{.>}[dddrr] &&&&\bar{y}  \ar@{.>}[dddll]^>{\MatheTiny{}} \\
   ~\\
   ~ \\
 & &&  z
}%Diagramm
}
}

\Def{total konfluent}{Ein Reduktionssystem $\Red$ heißt total konfluent, falls es zu je zwei endlichen Wegen, die vom selben Element $x$ starten, zwei endliche Wege gibt, die diese Wege so fortsetzen, dass sie im selben Element $z$ enden, also zu jedem:
\mathe{
\Diagramm{
 &&&x 	
  \ar@{.>}[dll]^>{\MatheTiny{}} \ar@{.>}[drr]  \\
 & y   && && \bar{y} \\
   \textnormal{existiert:}  & y  \ar@{.>}[dddrr]  &&&&\bar{y}  \ar@{.>}[dddll]^>{\MatheTiny{}} \\
  \\
  \\
& &&  z
}%Diagramm
}
}

\Def{Church-Rosser Eigenschaft}{Ein Reduktionssystem $\Red$ erfüllt die Church-Rosser Eigenschaft, falls sich je zwei Church-Rosser-äquivalente Elemente $y, \bar{y}$ über jeweils nur absteigende Kanten zusammenführen lassen, also zu jedem:
\mathe{\Diagramm{
& y   \ar@{<.>}[rrrr]_>>>>{\MatheTiny{\Red}} &&&& \bar{x}  \\    \textnormal{existiert:}  & y  \ar@{.>}[dddrr]  &&&&\bar{y}  \ar@{.>}[dddll] \\
   ~\\   
 \\
 & &&  z
}}%Diagramm

}

\Def{Konvergenz}{
Ein Reduktionssystem $\Red$ heißt konvergent, wenn es noethersch ist und es  zu  jedem $x \in \MengeE$ genau eine Normalform gibt.
}

In der Literatur findet man anstelle von Konvergenz auch die Begriffe kanonisch oder vollständig. Wir werden den Begriff "`vollständig"' für einen Spezialfall von Konvergenz gebrauchen.

\Bem{totale Konfluenz $\neq$ lokale Konfluenz}{
Aus totaler Konfluenz folgt lokale Konfluenz, die Umkehrung gilt allerdings im Allgemeinen nicht. Dazu betrachten wir folgendes Gegenbeispiel:
}

\Bsp{Seien folgende Reduktionsregeln gegeben:
\AufzaehlungP{
\item $x_1 \Regel z_1$
\item $x_1 \Regel x_2$
\item $x_2 \Regel x_1$
\item $x_2 \Regel z_2$
}

Das zugehörige Reduktionssystem können wir darstellen durch:
\mathe{\Diagramm{
& x_1 \ar[r] \ar@/^/[d]&z_1\\
z_2 &x_2  \ar@/^/[u]\ar[l] 
}}

Von keinem Element aus gibt es zwei Reduktionsregeln, die nach $z_1$ bzw. $z_2$ reduzieren; alle anderen Paare von zwei Elementen können wir sogar auf die gleiche Normalform reduzieren. Also ist dieses Reduktionssystem lokal konfluent.
\mathe{\begin{array}{|c|c|c|c|c|} \hline
 &x_1 & x_2 & z_1 & z_2\\ \hline
x_1 & \Regel z_1 & \Regel z_1 & \Regel z_1 &\Regel z_2\\ \hline
x_2 & \Regel z_1 & \Regel z_1 & \Regel z_1 &\Regel z_2\\ \hline
z_1 & \Regel z_1 & \Regel z_1 & \Regel z_1 & \textnormal{nicht reduzierbar}\\ \hline
z_2 & \Regel z_1 & \Regel z_1 & \textnormal{nicht reduzierbar} &\Regel z_2\\ \hline
\end{array}
}
Die Elemente $z_1$ und $z_2$ sind unreduzierbar, aber es gibt in $x_1$ startende Reduktionswege, die nach $z_1$ und $z_2$ reduzieren. Also ist dieses Reduktionssystem nicht total konfluent.
\\
Da es einen unendlich langen Redukionsweg gibt, nämlich \mathe{x_1 \Regel x_2 \Regel x_1 \Regel \dots,} ist dieses Reduktionssystem nicht noethersch.
}

\Lemma{total konfluent $\Leftrightarrow$ Church-Rosser Eigenschaft}{
\label{totalkonfluent}
Ein Reduktionssystem erfüllt genau dann die Church-Rosser Eigenschaft, wenn es total konfluent ist.
\\
\Beweis{Lemma \ref{totalkonfluent}}{
\\
"`$\Leftarrow$"':\\
Sei für ein Reduktionssystem die Church-Rosser Eigenschaft erfüllt, dann sind alle Paare, die über ungerichtete Regeln miteinander verbunden sind, über ausschließlich absteigende Kanten zusammenführbar. Also ist das Reduktionssystem total konfluent.\\
"`$\folgt$"':\\
Für die Umkehrung wollen wir zeigen, dass sich in einem total konfluenten Reduktionssystem, je zwei Elemente, die über ungerichtete Regeln miteinander verbunden sind, auf das gleiche Element reduzieren lassen. Für den Beweis nutzen wir eine vollständige Induktion:\\
Seien $x_1$ und $x_2$ über $n$ Church-Rosser-äquivalent.
\AufzaehlungP{
\item{Induktionsbeginn $n=0$: } Für den Weg der Länge $0$ gilt $x_1=x_2$.
\item{Induktionsannahme: } Alle Elemente, die über $n$ ungerichtete Regeln miteinander verbunden sind, lassen sich auf das gleiche Element reduzieren.
\item{Induktionsschritt $n \folgt n+1$: } Dann gibt es ein $y$ mit:\\ 
\mathe{
\Diagramm{ & y \ar[dl] \ar@{<.>}[rr]^{n} && x_2 &\textnormal{ oder } &x_1 \ar[dr]\\
x_1 &&&&&&  y \ar@{<.>}[rr]^{n} && x_2 .
}%Diagramm
}

Nach Induktionsannahme lassen sich $y$ und $x_2$ über einen Weg \Diagramm{\ar@{-->}[r]&}zu $z$ reduzieren. \\
Im ersten Fall gibt es wegen totaler Konfluenz zu $x_1$ und $z$ Wege \Diagramm{\ar@{:>}[r]&}, die sich in $z'$ wieder zusammenführen lassen. Der Weg von $x_2$ nach $z$ lässt sich mit dem Weg $z$ nach $z'$ fortsetzen. \\
Im zweiten Fall lässt sich der Weg von $x_1$ nach $y$ durch den Weg $y$ nach $z$ fortsetzen. \\
\mathe{\Diagramm{ & y 
\ar@{-->}[dr]
 \ar[dl] \ar@{<.>}[rr]^{n} && x_2 \ar@{-->}[dl] &\textnormal{ oder } &x_1 \ar[dr]\\
x_1 \ar@{:>}[dr] &&z \ar@{:>}[dl]&&&&  y  \ar@{-->}[dr]\ar@{<.>}[rr]^{n} && x_2 \ar@{-->}[dl]\\
&z'&&&&&&z
}.}%Diagramm
\\
So ist in beiden Fällen die Church-Rosser Eigenschaft für $x_1$ und $x_2$ und deshalb auch für Wege der Länge $n+1$ erfüllt.
}
}}

\input{Haupt/Reduktionssysteme/Bsp-Church-Rosser-Aequi}

\subsection{Diamantenlemma}
\label{sec:DiamantenLemma}

Im Allgemeinen sind totale Konfluenz und die Church-Rosser Eigenschaft sehr schwer nachzuweisen. 
In diesem Abschnitt wollen wir zeigen, dass für noethersche Reduktionssysteme
diese beiden Eigenschaften und lokale Konfluenz äquivalent sind.

Den folgenden Satz findet man in seiner Grundform in \cite{Newman}. Sein Name kommt daher, dass sich die Definitionen von lokaler Konfluenz, totaler Konfluenz und der Church-Rosser Eigenschaft leicht an einem Diamanten veranschaulichen lassen.
Lokale Konfluenz bedeutet, dass je zwei Ecken, die von der Spitze über eine Kante erreichbar sind, über jeweils nur absteigende Kanten wieder verbunden werden können. Totale Konfluenz bedeutet, dass zwei beliebige Ecken, die über absteigenden Kanten von der Spitze erreichbar sind, über jeweils nur absteigende Kanten wieder verbunden werden können. Die Church-Rosser Eigenschaft bedeutet, dass zwei beliebige Ecken auf einem Diamanten über jeweils nur absteigende Kanten wieder verbunden werden können.
\mathe{\Diagramm{
&&&&x \ar@{->}[dlll]  \ar@{->}[drrr] \\
& y   \ar@{<.>}[drr]_>>>>>>>{\MatheTiny{\Red}} \ar@{.>}[ddddddrrr] &&&&&& \bar{y}  \ar@{.>}[ddddddlll]  \ar@{<.>}[dll]^>>>>>>>{\MatheTiny{\Red}} &&&&&&&\\
&&& \dots \ar@{<.>}[rr]_>>>>>>>{\MatheTiny{\Red}} && \dots&&&&&\\
    \\
   ~\\   ~\\
&&&&&& \\
& & & & z &&&&&&
}}%Diagramm %Mathe

\Satz{Newmans Diamantenlemma}{
\label{Newman}
Sei $\Red$ noethersch, dann ist äquivalent:

\AufzaehlungZ{
	\item Zu jedem $x \in \MengeE$ gibt es genau eine Normalform.
	\item Das Reduktionssystem $\Red$ ist total konfluent.
	\item Das Reduktionssystem $\Red$ ist lokal konfluent.
	\item Das Reduktionssystem $\Red$ erfüllt die Church-Rosser Eigenschaft.
	} 

}

\Beweis{Satz \ref{Newman}}{\label{DiamantenLemma}
Nehmen wir also an, dass $\Red$ noethersch ist, dann gibt es zu jedem Element aus $\MengeE$ wenigstens eine Normalform.
\AufzaehlungP{
	\item{$1 \folgt 2:$ } Sei die Normalform eindeutig. Wir betrachten zwei Reduktionswege, die im selben $x$ beginnen und zu $y$ bzw. $\bar{y}$ reduzieren:
	\mathe{
	\Diagramm{
 &x 	
  \ar@{.>}[dl] \ar@{.>}[dr] \\
  y   \ar@{.>}[d]&& \bar{y}\ar@{.>}[d]\\
  \Nf(y) &&\Nf(\bar{y})
}%Diagramm\\
}
Da das Reduktionssystem noethersch ist, haben $y$ und $\bar{y}$ eine Normalform. Diese Normalformen sind auch Normalformen für $x$. Da die Normalform eindeutig ist, gilt $\Nf(y)=\Nf(\bar{y})$. Also ist das Reduktionssystem total konfluent.

	\item{$2 \folgt 3:$} Zu zeigen: Aus total konfluent folgt lokal konfluent. Da Wege der Länge eins gerade Reduktionsregeln sind, ist jedes total konfluente auch ein lokal konfluentes Reduktionssystem.
	
	\item{$3 \folgt 1:$} Sei $\Red$ lokal konfluent. Wir werden mit einem Widerspruchsbeweis zeigen, dass die Normalform eindeutig ist.
	Nehmen wir also an, dass es zwei unterschiedliche Normalformen $y_n$ und $\bar{y}_m$ zu einem $x$ gibt.\mathe{
 \Diagramm{
 	 &&  \ar[dl] x \ar[dr]\\
 	 &\ar@{.>}[ddl] y_{1} && \bar{y}_1 \ar@{.>}[ddr]\\
 	 \\
	 y_n & && & \bar{y}_m\\ 
	 }
	 }
 
	 Da $\Red$ lokal konfluent ist, lassen sich  $y_1$ und $\bar{y}_1$ zu demselben Element $z_1$ reduzieren. Da $\Red$ noethersch ist, gibt es eine Normalform $\Nf(z_1)$. Da sich $y_n$ und $\bar{y}_m$ unterscheiden, gilt  $\Nf(z_1) \neq y_n$ oder $\Nf(z_1)\neq \bar{y}_m$. OBdA sei $\Nf(z_1)\neq y_n$. Dann hat auch $y_1$ zwei unterschiedliche Normalformen, nämlich $y_n$ und $\Nf(z_1)$. 

	 \mathe{
	  \Diagramm{
 	 &&  \ar[dl] x \ar[dr]\\
 	 &\ar@{.>}[ddl] y_{1}\ar@{.>}[dr] && \ar@{.>}[dl] \bar{y}_1 \ar@{.>}[ddr]\\
 	 &&z_1 \ar@{.>}[d]\\
	 y_n & & \Nf(z_1)& & \bar{y}_m\\ 
	 }
	 }	 
	 Nun nennen wir $y_1$ in $x_1$ um. Dann gilt: $x \Regel x_1$, insbesondere $x \neq x_1$ und $x_1$ hat zwei Normalformen. Genauso wie für $x$ finden wir zu $x_1$ ein $x_2$, für welches gilt:  $x_1 \RedRegel x_2$, insbesondere $x_1 \neq x_2$ und $x_2$ hat zwei Normalformen. Durch Wiederholung erhalten wir einen unendlich langen Reduktionsweg:
	 \mathe{x \Regel x_1 \Regel x_2 \Regel \dots \Regel x_n \Regel \dots ,}
	 der ein Widerspruch zu noethersch ist.

\item{$3 \Leftrightarrow 4 :$} Mit Lemma \ref{totalkonfluent} folgt die Behauptung.
}
}

\Bem{Konvergenz}{
Ein Reduktionssystem $\Red$ heißt konvergent, wenn es noethersch ist und eine (und damit alle) der folgenden Bedingungen erfüllt:

\AufzaehlungP{
	\item Zu jedem $x \in \MengeE$ gibt es genau eine Normalform.
	\item Das Reduktionssystem $\Red$ ist total konfluent.
	\item Das Reduktionssystem $\Red$ ist lokal konfluent.
	\item Das Reduktionssystem $\Red$ erfüllt die Church-Rosser Eigenschaft.
	} 
}

%% file: Haupt/Reduktionssysteme/Bsp-Church-Rosser-Aequi.tex
\Bsp{
Folgendes ist ein noethersches Reduktionssystem auf $\MengeE$, so dass es zu einem $x\in \MengeE$ unendlich viele Church-Rosser-äquivalente Elemente gibt:\\
Sei
\mathe{
\MengeE:=\Menge{x_i}{i\in \N} \cup \Menge{z_i}{i \in \N_0} \textnormal{ und } \Red:=\Menge{(x_i,z_i)}{i\in \N} \cup \Menge{(x_i,z_{i-1}}{i\in \N}, 
}
dann ist \Diagramm{ x_1   \ar@{<.>}[r]_>{\MatheTiny{\Red}} & x_i} für jedes $i\in \N$, da:\\
\Diagramm{ 
 & \ar@{->}[dl] x_1   \ar@{->}[dr] 
&& \ar@{->}[dl] x_2   \ar@{->}[dr] 
&& \ar@{->}[dl] \dots \ar@{->}[dr]
&& \ar@{->}[dl] x_n   \ar@{->}[dr]
&& \ar@{->}[dl] \dots \\
z_0 && z_1  && z_2  && z_{n-1}  && z_n
}
Wenn wir zusätzlich noch $z_0=z_1=\dots=z_n=\dots$ setzen, ist das Reduktionssystem sogar konvergent.
}

%% file: Haupt/Reduktionssysteme/ReduktionssystemeAlgebraDef.tex
\section{Reduktionssystem bezüglich einer Algebra}
\label{sec:ReduktionssystemFbezüglichEinerAlgebra}
In diesem Kapitel wollen wir eine Algebra $\Algebra$  als Quotient einer freien Algebra $\FreieA$ auffassen. So können wir jedes $a \in \Algebra$ als Äquivalenzklasse in $\FreieAlgebra$ betrachten. Wir werden  $\Red$ so wählen, dass $f_1, f_2 \in \FreieAlgebra$ genau dann in derselben Äquivalenzklasse liegen,  wenn sie dieselbe Normalform haben. 
Die lokale Konfluenz  dieses Reduktionssystems werden wir mittels minimaler Überschneidungen von Silben beweisen.

%\Frage{$\Red \subset ?$, $\MonomeM=?$}
%\Frage{Mehr Text vgl. Avenhaus 2.1 mit Algebra}

Die in dieser Arbeit betrachteten Spezialfälle von Reduktionssystemen
%, für durch Gleichungen definierten Strukturen (Gruppen, Algebren, Moduln, etc) 
nennt man auch Wortersetzungssysteme.
% Man hat man in der Regel eine Vielzahl an möglichen Reduktionssystemen, für die Konvergenz unterschiedlich schwer nach zuweisen ist.

~\\
%\Frage{Definition Algebra, freie Algebra}
%In diesem Abschnitt wird eine Methode vorsgestellt, wie für eine Algebra $\Algebra$, die durch Erzeuger und Relationen definiert ist, ein konvergentes Reduktionssystem gefunden wird.Die Algebra $\Algebra$ wird aufgefasst als Quotient der freien Algebra $\FreieAlgebra$ durch ein über endlich viele Relationen definiertes Ideal $\Ideal$:\mathe{ \Algebra \iso \FreieAlgebra / \Ideal .}

%~\\
%Um noethersch zu beweisen findet man eine noethersche Teilordnung $\geq$ auf $\MengeE \times \MengeE$, die $R$ als Teilmenge enthält, also:
%\mathe{R ~ \subseteq ~ \geq}
%
%~\\
%Um lokal konfluent zu beweisen sucht man sich eine Teilmenge $r$ der Reduzierungsregeln $R$. Die Teilmenge $r$ wird so gewählt, dass die lokale konfluentz leicht nachzuweisen ist und sich die Reduzierungsregeln aus $R$ von denen aus $r$ erzeugen läßt. Die Teilmenge $r$ sind in gewisser Weise die "`Atome"' von $R$.
%
%~\\

\subsection{Definitionen}
\label{sec:DefinitionenWort}
\label{secDefAlgebra}

Sei $\K$ ein Körper.
Sei $\Alphabet$ ein Alphabet.\\
%\Frage{$\K$ vertauscht mit $\Algebra$}
\Def{Terme $\Terme$}{
\label{DefTerm}
Ein (endliches) Produkt aus Elementen $a \in \Alphabet$ und einem Koeffizient $\lambda \in \K$ heißt Term. Die Menge der Terme bezeichnen wir mit $\Terme$, also: 
\mathe{ \Terme:= \Menge{ \lambda a_1 \cdots a_n}{a_1, \dots , a_n \in \Alphabet, \lambda \in \K^*, n \in \N}.
	}
}

\Def{Monome $\MonomeM$}{
\label{DefMonom}
Terme ohne Koeffizient aus $\K$ heißen Monome oder Worte. Die Menge der Monome bezeichnen wir mit $\MonomeM$, also:
\mathe{ \MonomeM := \Menge{ a_{1} \cdots a_{n} }{ a_1, \dots , a_n \in \Alphabet }. }
Es ist auch das leere Monom zugelassen.
}
\Ver{Im Folgenden werden wir nur noch Reduktionssysteme auf Mengen $\MengeE$ betrachten, die $\Terme$ enthalten.
\mathe{\MonomeM \subset \Terme \subset \MengeE,}
wobei wir $\MengeE$ später genauer festlegen werden.
}

\Def{Wortersetzungssystem}{
\label{DefWortersetzungssystem}
Ein Reduktionssystem $\red$, in dem bei jeder Reduktionsregel auf der linken Seite ein Monom steht, heißt Wortersetzungssystem, also:
\mathe{\red \subset \MonomeM \times \MengeE.}
}

%\subsection{Definitionen}
\label{sec:DefinitionenAlgebra}

\Def{Algebra $\Algebra$}{
\label{DefAlgebra}
Ein  $\K$-Vektorraum $\Algebra$ mit einer $\K$-bilinearen Verknüpfung

\mathe{\Algebra \times \Algebra \rightarrow \Algebra, }
für die zusätzlich gilt:
\begin{itemize}
	\item Assoziativität: $a * (b * c) = (a * b) * c~~~~ \forall a,b,c \in \Algebra$,
	\item Eins: $\exists 1 \in \Algebra: 1 * a = a*1 = a ~~~~\forall a \in \Algebra$,
\end{itemize}

heißt Algebra $\Algebra$. 
}

\Def{freie Algebra $\FreieA$}{ 
\label{DefFreieAlgebra}
Sei $\Alphabet$ ein Alphabet und sei $\MonomeM$ die Menge von Monomen über $\Alphabet$. Wir bezeichnen als freie Algebra  $\FreieA$ den Ring:

 \mathe{ \mathcal{F}:= \Menge{ \Sum{m \in \mathcal{M} }{}\lambda_m m }{ \# \{\lambda_m \neq 0 \} < \infty \\ \lambda_m \in \K \\ m \in \MonomeM}
	.}

Für ein $f \in \FreieA$ bezeichnen wir mit $\Traeger{f}$ die Menge von Monomen,  deren Koeffizient in $f$ ungleich Null ist. 
}

Also sind die Monome über dem Alphabet $\Alphabet$ eine lineare Basis für $\FreieA$. Manchmal schreiben wir für eine freie Algebra $\FreieA_n$ um deutlich zu machen, dass es sich um die freie Algebra über einem Alphabet mit $n$ Buchstaben handelt.

Aus einem Wortersetzungssystem $\red \subset \Monome \times \FreieA$ lässt sich ein größeres Reduktionssystem auf der freien Algebra definieren.

\Def{induziertes Reduktionssystem $\Red$}{
\label{DefInduziertesRedAlgebra}
Sei $\red \subset \Monome \times \FreieA$ ein Wortersetzungssystem, dann definieren wir:
\mathe{ \Red:= \Menge{\left( \lambda p x s + q, \lambda p (y) s + q\right) }{\lambda \in \K\\ (x,y) \in r\\p,s\in \Monome\\ q \in \FreieA  \textnormal{ mit }p x s \notin \Traeger{q} }
	.}
Wir sagen, $\Red$ ist das durch $\red$ auf $\FreieA$ induzierte Reduktionssystem. Um deutlich zu machen, dass $y$ ein Polynom und nicht nur ein Monom wie $x$ ist, setzen wir Klammern.
}
In keinem Summanden, der in $q$ vorkommt, darf $pxs$ als Wort vorkommen, da sonst dieser Summand mit $\lambda pxs$ zusammengefasst werden könnte.
Für Beweise ist es umständlich jeweils $pxs \notin \Traeger{q}$ zu formulieren (vergl. \cite{MR506890}). Im Folgenden  übertragen wir daher  $\Red$ in ein größeres Reduktionssystem $\RRed$, das auf einer größeren Menge definiert ist, in der mit $+$ nicht zusammengefasst werden kann.
Das Reduktionssystem $\RRed$ soll so gewählt werden, dass sich die Menge der unreduzierbaren Elemente nicht verändert. Es sollen alle Regeln aus $\Red$ enthalten sein, also müssen wir zusätzlich Regeln finden, die dem Zusammenfassen von Summanden in $\FreieA$ entsprechen. Es verlängern sich die Reduktionswege, wodurch die Argumente deutlicher werden.

%in $\RigF$ bezüglich $\mathcal{\widetilde{R}}$ gleich der Menge der unredzuierbaren Elemente in $\FreieA$ bezüglich $\Red$ ist.In $\mathcal{\widetilde{R}}$ müssen alle Regeln aus $\Red$ enthalten sein und zusätzlich Regeln gefunden werden für Summanden die in $\FreieA$ zusammen gefasst werden können.
\Def{Rig $\Rig$ }{
\label{DefRig}
Ein Rig ist  eine Menge $\Rig$,  die mit zwei Verknüpfungen Addition $+'$ und Multiplikation $\cdot$ versehen ist, so dass gilt:

\AufzaehlungP{
\item{$(\Rig, +')$ ist ein kommutativer Monoid mit neutralem Element 0:}
\Gleichung{
(a +'  b) +'  c &= a +'  (b +'  c),\\
0 +'  a &= a +'  0 = a,\\
a +'  b &= b +'  a.
}
\item{$(\Rig, \cdot)$ ist ein Monoid mit neutralem Element 1: }
\Gleichung{
(a\cdot b)\cdot c &= a\cdot (b\cdot c),\\
1\cdot a &= a\cdot 1 = a.
}
\item{Distributiv:}
\Gleichung{
a\cdot(b +'  c) &= (a\cdot b) +'  (a\cdot c),\\
(a +'  b)\cdot c &= (a\cdot c) +'  (b\cdot c).}
\item{Null ist ein Annulator: }
\Gleichung{
0\cdot a = a\cdot 0 = 0.
}

}%Aufzaehlung
}%Def
Wir wollen im Folgenden einen speziellen Rig betrachten. Dieser Rig soll die freie Algebra $\FreieA$ enthalten. Dazu betrachten wir zunächst folgendes Beispiel.

\Bsp{\label{BspRig}
Sei $(M,1,\cdot)$ ein Monoid und sei $(H(M),0,+')$ ein freier
kommutativer durch  $M$ erzeugter Monoid. Also ist $H(M)$ die
Menge von formalen Summen
\mathe{
m_1 +' m_2 +' \dots +' m_k
}
von Elementen aus $M$.

Die Menge $(H(M),1,0,+',\cdot)$ hat die Struktur eines Rigs,
wobei die Multiplikation aus
$M$ auf $H(M)$ erweitert wird:
 
\mathe{
(m_1 +' m_2 +' \dots +' m_k) (n_1 +' n_2 +' \dots +'
n_l) := \RSum{i=1}{k} \RSum{j=1}{l}  m_in_j .
}
}

\Bsp{Sei $H$ wie in Beispiel \ref{BspRig} definiert.
Die Menge von Termen $\Terme:= \{\lambda_w w ~|~ w \in \Monome \}$ ist über \mathe{ \lambda_w w \cdot \lambda_{\bar{w}} \bar{w} = (\lambda_w \lambda_{\bar{w}}) w \bar{w}} ein Monoid, daher ist $H(\Terme)$ ein Rig.
}

\Def{Rig $\RigF$}{
\label{DefRigF} 
Sei $\Terme$ die Menge von Termen über einem Alphabet $\Alphabet$. Wir definieren  den Rig über $\Alphabet$ durch:
  \mathe{\RigF:=H(\Terme),}
  wobei $H$ wie in Beispiel \ref{BspRig} definiert ist.
Im Folgenden sprechen wir auch von dem Rig $\RigF$, wenn eindeutig ist, welches $\Alphabet$ gemeint ist. 
} 
 
%Wir können die Elementen von $\tilde{F}$ als Polynome betrachten, in denen Terme mit gleiche Monome nicht zusammen gefasst werden könnnen.

Wir können den Rig  $\RigF$ auch auffassen als:
\mathe{ \RigF := \Menge{ \RSum{w \in \mathcal{M}}{} \RSum{i=1}{k_w} \lambda_{w,i} w  }{ \#\{k_w \neq 0 \} < \infty \\ \lambda_{w,i} \in \K^* \\ w \in \MonomeM }  .}

Ein Element $\widetilde{f} \in \RigF$ ist eine Menge von Familien\footnote{In einer Familie kann ein Element mehrmals vorkommen.} über $\MonomeM$  mit Elementen in $\K^*$.  Also ist $\tilde{f}$ ein Ausdruck von Summen mit Termen, die mit "`$+'$"' anstelle von "`$+$"' verknüpft sind. Im Gegensatz zur Summe in der freien Algebra können die Koeffizienten vor gleichen Monomen nicht zusammengefasst werden. 
Es ist also zu beachten, dass diese Summe nicht distributiv mit der Multiplikation des Körpers ist, also:
\mathe{ \lambda w +' \mu w \neq (\lambda + \mu)w .}
Zwei Ausdrücke sind dann gleich, wenn sie die gleichen Summanden in höchstens unterschiedlicher Reihenfolge besitzen.

\Ver{
\label{VerEinbettung}
Wir betrachten $\FreieAlgebra \subset \RigF$ in folgender Weise:
\mathe{ \textnormal{Sei } f= \Sum{w\in \mathcal{M}}{} \alpha_w w \in \mathcal{F} \textnormal{, dann ist }
f= \RSum{w \in \Monome}{} \RSum{i=1}{k_w} \alpha_{w} w \in \widetilde{\mathcal{F}},}
wobei $k_w = \left\{ \Matrix{0 &\textnormal{ für } \alpha_w=0\\ 1 & \textnormal{sonst}} \right.$  . Bei dieser Einbettung ist zu beachten, dass  sie mit den Algebra-Verknüpfungen nicht verträglich ist.

}

%\Ver{Multiplikation auf $\RigF$}{Seien $a= a_1 +'a_2 ,b \in \RigF$, dann ist die Multiplikation definiert durch \mathe{ a \cdot b = a_1b +' a_2b.}Zur Erinnerung: Summanden, die sich nur im Koeffizienten unterscheiden, werden nicht zusammengeführt.}

Wir wollen nun ein größeres Reduktionssystem definieren, das $\Red$ enthält.

\Def{induziertes Reduktionssystem $\RRed$}{
\label{DefInduziertesRedRig}
\label{Rtilde}
Sei $\red \Monome \times \FreieA$ ein Wortersetzungssystem und sei:
\Gleichung{ 
\Rel:=& \Menge{ (\lambda p x s +' q~,~ \lambda p (y) s +' q)}{ \lambda \in \K^*\\ p,s\in \Monome \\ q \in \RigF\\(x,y) \in \red },
\\
 \RelF :=& \Menge{ \left(\lambda w +' \mu w +' q, (\lambda + \mu)w +' q\right) }{\lambda, \mu \in \K^*\\ w \in \Monome,\\ q \in \RigF }. 
 }
Das durch
\mathe{\RRed := \Rel \cup  \RelF }
definierte Reduktionssystem heißt das durch $\red$ auf $\RigF$ induzierte.

%%In $\Rel$ sind deutlich mehr Reduktionsregeln als in $\Red$. Da anders als in $\Red$ nun ein $w$ im Träger von $q$ liegen darf. 
}

%% file: Haupt/Reduktionssysteme/ReduktionssystemeAlgebraKonfluenz.tex
\subsection{Konfluenz}
\label{sec:Konfluenz}
Sei $\Red$ bzw. $\RRed$ das durch $\red$ auf $\FreieA$ bzw. $\RigF$ induzierte Reduktionssystem.
In diesem Abschnitt zeigen wir, dass die Mengen von unreduzierbaren Elementen aus $\Red$ und  aus $\RRed$ übereinstimmen und dass sich Konfluenz von $\RRed$ nach $\Red$ überträgt, falls $\RRed$ noethersch ist.

\Satz{$\Nf = \widetilde{\Nf}$}{
\label{NgleichTN}
Sei $\Nf$ die Menge von unreduzierbaren Elementen in $\FreieA$ bezüglich $\Red$ und sei $\widetilde{\Nf}$ die Menge der unreduzierbaren Elemente in $\RigF$ bezüglich $\RRed$. 
Es gilt:
\mathe{\widetilde{\Nf} \subseteq \FreieA}
und
\mathe{\Nf = \widetilde{\Nf}.}

\Beweis{Satz \ref{NgleichTN}}{Sei $f \in \widetilde{\Nf}$. In $f$ können keine zwei Terme mit gleichen Monomen vorkommen, da diese mit Regeln aus $\RelF$ reduziert werden könnten. Das bedeutet, dass  $f$  in $\FreieA$ liegt.

Für den zweiten Teil des Satzes bemerken wir zunächst, dass $z \in \Nf$ genau dann gilt, wenn $z$ ein Polynom in $\FreieA$ ist, in dem kein Summand eine Silbe enthält, die sich mit einer Regel aus $\red$ reduzieren lässt. Mit \ref{VerEinbettung} können wir $z$ als ein Element aus $\RigF$ auffassen, für welches keine zwei Terme mit gleichen Monomen vorkommen. Es kann also keine Regel aus $\RRed$ angewendet werden. Also ist $z$ auch in $\RigF$ ein unreduzierbares Element. 

Für die umgekehrte Richtung sei $z \in \widetilde{\Nf}$. Das ist genau dann der Fall, wenn weder eine Regel aus $\RelF$ noch aus $\Rel$ angewendet werden kann. 
Da keine Regel aus $\RelF$ angewendet werden kann, können in $z$ keine zwei Terme mit gleichen Monomen vorkommen, wir können $z$ also auch als Element in $\FreieA$ auffassen. Da nach Voraussetzung auch keine Regel aus $\Rel$ angewendet werden kann, können wir auch keine Regel aus $\Red$ anwenden, also gilt $z \in \Nf$.

}%Beweis

}%Satz

%\Frage{Nun wollen wir die konvergenz von $\Red$ durch die von $\red$ weisen. Dazu vergrößern wir zunächst das Reduktionssystem auf $\RigF$ und fügen dann Regel hinzu. Die zugefügten Regel sorgen dafür, dass die Normalformen übereinstimmen.}

Ein Reduktionsweg in $\RRed$, der in $\FreieA$ beginnt und endet, impliziert jedoch keinen entsprechenden Reduktionsweg 
$\Red$. Dazu betrachten wir folgendes Beispiel.

\Bsp{
Sei $\red$ ein Wortersetzungssystem, das nur aus folgenden zwei Regeln besteht:
\mathe{x \redRegel z + y \textnormal{ und }y \redRegel \bar{z}.}
Dann gibt es einen Reduktionsweg in $\RRed$:
\mathe{x +' (-y) \RRedRegel (z +' y ) +' (-y) \RRedRegel (z+'\bar{z}) +' (-y),}
der in $\FreieA$ beginnt und endet; aber es gibt keinen entsprechenden Reduktionsweg in $\Red$, da es für $x-y$ in $\Red$ nur folgende Reduktionswege gibt:
 \mathe{
 \Diagramm{
 					& z + y - y= z \\
 x-y \ar@{->}[ur]_>{\MatheTiny{\Red}} \ar@{->}[dr]^>{\MatheTiny{\Red}}\\
 					& x-\bar{z} \ar@{->}[r]_>{\MatheTiny{\Red}} & z + y -\bar{z} \ar@{->}[r]_>{\MatheTiny{\Red}} & z + \bar{z} -\bar{z}=z 
 .}
 }

 %\mathe{  x-y &\Matrix{  					\RedRegel & z + y - y= z \\  					\RedRegel &x-\bar{z} &\RedRegel z + y -\bar{z} \RedRegel z + \bar{z} -\bar{z}=z .} }
}

Die umgekehrte Richtung gilt jedoch. Dazu betrachten wir folgendes Lemma.
\Lemma{$\RedWeg \subset\RRedWeg$}{\label{RsubsetTR}
Zu jedem Reduktionsweg in $\Red$ gibt es einen Reduktionsweg 
in $\RRed$.\\
\Beweis{Lemma \ref{RsubsetTR}}{
Direkt aus der Definition \ref{DefInduziertesRedRig} folgt, dass jede Regel aus $\Red$ sich durch eine Regel aus $\Rel$ und mehrere Regeln aus $\RelF$ schreiben lässt. Also lässt sich auch jeder Weg in $\Red$ in einen (eventuell längeren) Weg in $\RRed$ übersetzen.
}
}

\Satz{Konvergenz}{
\label{konfluentnoethersch} 
Sei $\RRed$ konvergent, dann ist auch $\Red$ konvergent .

\Beweis{Satz \ref{konfluentnoethersch}}{~\\
Sei $\RRed$ noethersch. Wegen Lemma \ref{RsubsetTR} lässt sich jede Reduktionsregel in $\Red$ in eine Reduktionsregel aus $\Rel$ und mehrere Reduktionsregeln aus $\RelF$ übersetzen. Ein Reduktionsweg $\RedWeg$ in $\FreieA$ kann demnach in einen längeren Reduktionsweg $\RRedWeg$in $\RigF$ überführt werden. Wenn dieser längere Reduktionsweg endlich ist, dann ist es auch der in $\FreieA$. Also ist $\Red$ noethersch.
Daraus folgt: Jedes $f \in \FreieA$ hat wenigstens eine Normalform bezüglich $\Red$. 

Für die Konfluenz betrachten wir folgenden Widerspruchsbeweis:\\
Sei nun $\Red$ nicht konfluent, dann gibt es ein $f \in \FreieA$ mit zwei verschiedenen Normalformen $z_1$ und $z_2$ bezüglich $\Red$. Es gibt also zwei Reduktionswege: $f \RedWeg z_1$ und $ f \RedWeg z_2$. Zu diesen Reduktionswegen in $\Red$ können wir mit \ref{RsubsetTR} Reduktionswege in $\RRed$ finden. Nun sind $z_1$ und $z_1$ nach \ref{NgleichTN} unreduzierbar bezüglich $\RRed$. Dies ist ein Widerspruch zur Konfluenz von $\RRed$.

}%Beweis

}%Satz

\Bem{Induzierter Reduktionsweg}{
\label{axb} 
Sei  $g_1 \RRedWeg g_2$ ein Reduktionsweg in $\RRed$ für Elemente aus dem Rig $\RigF$, dann ist auch 
\mathe{ a g_1 b +' q \RRedWeg a g_2 b +'q}
ein Reduktionsweg, wobei $a, b, q$ beliebige Elemente aus dem Rig $\RigF$ sein können.
}

\subsection{Überschneidungen}
\label{sec:Überschneidungen}
Im letzten Abschnitt hatten wir gesehen, dass es reicht die Konvergenz von $\RRed$ zu zeigen um die Konvergenz von $\Red$ zu beweisen. Jedoch ist es viel Arbeit alle möglichen Reduktionswege von jedem Element aus $\RigF$ zu überprüfen. Daher wollen wir 
nun zeigen, dass es für den Beweis der Konfluenz ausreicht sich auf eine Teilmenge zu beschränken.
Wie auch in den vorhergehenden Abschnitten seien $\Red$ und $\RRed$ durch ein Wortersetzungssystem $\red$ induziert.

\Def{Überschneidung $(w,\regel,\bar{\regel})$}{
Eine Überschneidung $(w,\regel,\bar{\regel})$ ist ein Monom $w\in \RigF$, zusammen mit zwei unterschiedlichen Reduzierungsregeln $\regel=(x,y), \bar{\regel}=(\bar{x},\bar{y}) \in \red $ , so dass es $p,\bar{p}, s, \bar{s} \in \MonomeM$ und $\lambda, \bar{\lambda} \in \K^*$ und $q, \bar{q} \in \RigF$ gibt mit:
\mathe{ w= \lambda p x s +' q = \bar{\lambda} \bar{p} \bar{x} \bar{s} +' \bar{q}.}
}

\Bem{leeres Monom}{Die Präfixe $p,\bar{p}$ und die Suffixe $s,\bar{s}$ können auch das leere Wort sein.}

Manche Überschneidungen lassen sich durch kürzere ersetzen. Betrachtet man zum Beispiel die Überschneidung $ p x  = \bar{p} \bar{x} $, dann beginnen die Monome $p$, $\bar{p}$ gleich, sie lassen sich also zerlegen in $p= q \cdot p_1$ und $ \bar{p} = q \cdot \bar{p}_1$, wobei $p_1$ oder $\bar{p}_1$ leer ist. Ähnliches gilt für Suffixe. Daher ist folgende Definition nahe liegend.

%oder kritisches Paar $\kP$
\Def{minimale Überschneidung, kritisches Paar}{ 
Eine minimale Überschneidung $(w,\regel,\bar{\regel})$ ist ein Monom $w \in \Monome$, zusammen mit zwei unterschiedlichen Reduzierungsregeln $\regel=(x,y), \bar{\regel}=(\bar{x},\bar{y}) \in \red$, so dass es $p, s \in \MonomeM$ gibt mit:
\mathe{ w= p x  =  \bar{x} s ~~~\textnormal{ oder }~~~ w = p x s =  \bar{x} . }
Eine minimale Überschneidung besteht also maximal aus einem Präfix und einem Suffix.
Manchmal sagen wir zu einer minimalen Überschneidung auch kritisches Paar.
}

\Def{zusammenführbar, behebbar}{
Eine minimale Überschneidung $(w,\regel_1,\bar{\regel}_1)$ heißt zusammenführbar oder behebbar bezüglich $\Red$ bzw. $\RRed$, wenn es ein $z \in \FreieA$ bzw. $\in \RigF$ zusammen mit zwei endlichen Reduktionswegen $(\regel_i)_{i=1\dots n}$ und $(\bar{\regel_i})_{i=1\dots \bar{n}}$ in $\Red$ bzw. $\RRed$ gibt, die im selben Element $ z$ enden und es somit Wege gibt, die die beiden Regeln wieder zusammenführen.
}
%}

\Def{vollständig}{Ein Ersetzungssystem $\red$ heißt vollständig bezüglich $\Red$ bzw. $\RRed$, falls jede minimale Überschneidung mit Reduktionsregeln aus $\Red$ bzw. $\RRed$ behebbar ist.
}

\Satz{$\red$ vollständig $\folgt \RRed$ lokal konfluent}{
\label{vollstaendigFolgtlokalkonfluent}
Sei $\RRed$ noethersch. Falls $\red$ vollständig bezüglich $\RRed$ ist, dann gilt:
\mathe{ \RRed \textnormal{ ist lokal konfluent.} }

\Beweis{Satz \ref{vollstaendigFolgtlokalkonfluent}}{Seien alle minimalen Überschneidungen aus $\red$ behebbar.\\
Es muss gezeigt werden, dass es zu jedem $g_1 \in \RigF$, zu dem es in $\RRed$ zwei verschiedene Reduktionsregeln $g_1\Matrix{\RRedRegel g_2 \\ \RRedRegel \bar{g}_2}$ gibt, zwei Reduktionswege $\Matrix{g_2 \RRedRegel \dots \RRedRegel\\ \bar{g}_2 \RRedRegel \dots \RRedRegel} z$ gibt, die im gleichen $z$ enden.

Die Reihenfolge der Summanden spielt keine Rolle. Deshalb schreiben wir die zu betrachtenden Summanden immer an den Anfang. Sei $(x,y)$ eine Reduktionsregel aus $\red$.  Für die Reduktionsregel in $\Rel$ schreiben wir dann:
\mathe{ pxs +' q \RRedRegel p(y)s +' q.}
Sei $(\bar{x}, \bar{y})$ eine weitere Reduktionsregel aus $\red$.
%Seien $a,b$ beliebige Terme in denen kein $x$ vorkommt, es also keine Reduktionsregel in $\Rel$ gibt.

Wir müssen folgende Arten von Überschneidungen betrachten:
\AufzaehlungP{
\item{Beginnend mit Reduktionsregeln jeweils aus $\RelF$:}~\\
Eine Reduktionsregel aus $\RelF$ überführt ein $+'$ in die übliche Addition. Die Reihenfolge, in der Reduktionsregeln aus $\RelF$ angewendet werden, spielt also keine Rolle.
\item{Beginnend mit einer Reduktionsregel aus $\RelF$ und einer aus $\Rel$:}\\ 
Wir weisen für die unterschiedlichen Fälle nach, dass sie zum gleichen Element reduziert werden können.\\
~\\Fall "`Getrennt"': Sei $g_1=\lambda pxs+'\mu \bar{x}+'\bar{\mu} \bar{x}+' q$.\\
Hier betrachten wir:
\Gleichung{ \lambda pxs+'\mu \bar{x}+'\bar{\mu} \bar{x}+' q \Matrix{ &\RRedRegelF& \lambda pxs+'(\mu+\bar{\mu}) \bar{x} +' q\\ &\RRedRegelr& \lambda p(y)s+'\mu \bar{x}+'\bar{\mu} \bar{x}+' q.}}
Diese Elemente lassen sich mit Regeln aus $\Rel$ bzw. $\RelF$ weiter reduzieren zu jeweils $ \lambda p(y)s+'(\mu +\bar{\mu}) \bar{x}+' q$.
~\\~\\Fall "`Inverses"': Sei $g_1=\lambda pxs +' (-\lambda pxs) +' q$ .\\
Hier betrachten wir:
\Gleichung{ \lambda pxs+' (-\lambda pxs) +' q \Matrix{ &\RRedRegelF& q\\ &\RRedRegelr& \lambda p(y)s+' (- \lambda pxs) +' q.}}
Der untere Teil lässt sich mit Regeln aus $\Rel$ weiter reduzieren zu $\lambda p(y)s+' (-\lambda p(y)s) +' q$ und dann mit Regeln aus $\RelF$ zu $q$ reduzieren.
\\~
\\Fall "`Überlagert"': Sei $g_1=\lambda pxs +' \mu pxs +' q$ .\\
Hier betrachten wir:
\Gleichung{ \lambda pxs +' \mu pxs +' q \Matrix{ &\RRedRegelF& (\lambda + \mu) pxs +' q\\ &\RRedRegelr& \lambda p(y)s +' \mu pxs +' q, }}
da gilt:
\Gleichung{\Matrix{ (\lambda + \mu)x &\RRedWeg&\\ \lambda(y) +' \mu (x) &\RRedWeg&  } (\lambda + \mu )(y)}

und wegen Bemerkung \ref{axb} lässt sich diese Überschneidung zum selben Element zusammenführen.

%brauchen wir nur die Summe der Silben $\lambda x +' \mu x$ betrachten:
%\mathe{ \lambda x +' \mu x \Matrix{ &\RRedRegelF& (\lambda + \mu)x \\ &\RRedRegelr& \lambda(y) +' \mu (x)} \textnormal{ lässt sich fortsetzen zu: } \Matrix{ (\lambda + \mu)x &\RRedRegelr&\\ \lambda(y) +' \mu (x) &\nachTR{\RelF \& \Rel}&  } (\lambda + \mu )(y).}

\item{Beginnend mit zwei Reduktionsregeln aus  $\Rel$:}\\ 

Fall "`Getrennt"': Sei $g_1=\lambda pxs +' \mu \bar{p}\bar{x}\bar{s} +'q$.\\
Hier betrachten wir:
\Gleichung{ \lambda pxs +' \mu \bar{p}\bar{x}\bar{s} +'q \Matrix{ &\RRedRegelr& \lambda p(y)s +' \mu \bar{p}\bar{x}\bar{s} +'q\\ &\RRedRegelr& \lambda pxs +' \mu \bar{p}(\bar{y})\bar{s} +'q,}}

da gilt:

\Gleichung{ \Matrix{ \lambda (y) +'  \mu \bar{x}  &\RRedWegr&  \\  \lambda x+' \mu (\bar{y})&\RRedWegr& } \lambda (y) +'  \mu ( \bar{y}) }
und wegen Bemerkung \ref{axb} lässt sich diese Überschneidung zum selben Element zusammenführen.

%\Frage{Statt der Monome $x$ und $\bar{x}$ kann man auch Polynome aus $\RigF$ einsetzen und so müssen wir nur die folgenden Fälle betrachten:\AufzaehlungP{\item $\lambda p x s +' \mu pxs +' \nu pxs$\item $\lambda p x s +' \mu pxs +' \lambda\bar{p} \bar{x}\bar{s}$\item $\lambda p x s +' \mu pxs +' \lambda \bar{p} \bar{x}\bar{s} +' \mu \bar{p}\bar{x}\bar{s}$}}

Fall: "`Überschneidung"': Sei $g_1=\lambda pxs\bar{x}\bar{s} +'q$.\\
Hier betrachten wir:
\Gleichung{ \lambda pxs\bar{x}\bar{s} +'q \Matrix{ &\RRedRegelr& \lambda p(y)s\bar{x}\bar{s} +'q\\ &\RRedRegelr& \lambda pxs(\bar{y})\bar{s} +'q,}}
da gilt:
\Gleichung{ \Matrix{  (y) s \bar{x}   &\RRedRegelr&  \\   x s (\bar{y}) &\RRedRegelr& }  y s (\bar{y})  }
und wegen Bemerkung \ref{axb} lässt sich diese Überschneidung zum selben Element zusammenführen.

}%Aufzählung
Dies waren die Fälle, die unabhängig davon sind, ob die minimalen Überschneidungen in $\red$ behebbar sind. Sei nun jede minimale Überschneidung aus $\red$ mit Regeln aus $\Red$ behebbar. Mit Satz \ref{RsubsetTR} gibt es dann auch einen Reduktionsweg in $\RRed$, der die minimale Überschneidung behebt. Sei $(x,y),(\bar{x},\bar{y})\in \red$. Wir betrachten noch folgende Fälle:
\AufzaehlungP{
\item{ minimale Teilüberschneidung: Sei $f=px= \bar{x}s$.}\\
Hier betrachten wir:
\Gleichung{ \lambda px b +' q \Matrix{ &\RRedRegelr& \lambda p(y) b +' q\\ &\RRedRegelr& \lambda a (\bar{y})sb +'q,}}
da nach Voraussetzung
\Gleichung{ f \Matrix{ &\RRedRegelr& \lambda p (y)  \\ &\RRedRegelr& \lambda (\bar{y}) s } \textnormal{ behebbar ist} }
und wegen Bemerkung \ref{axb} lässt sich diese Überschneidung zum selben Element zusammenführen.

\item{ minimale Totalüberschneidung: Sei $f=p x s=\bar{x}$.}\\
Hier betrachten wir:
\Gleichung{ \lambda a p x s b +' q \Matrix{ &\RRedRegelr& \lambda a p (y) s b +' q\\ &\RRedRegelr& \lambda a (\bar{y})b +'q,}}
da nach Voraussetzung:
\Gleichung{ f \Matrix{ &\RRedRegelr& \lambda p (y)s  \\ &\RRedRegelr& \lambda (\bar{y})  } \textnormal{behebbar ist } }
und wegen Bemerkung \ref{axb} lässt sich diese Überschneidung zum selben Element zusammenführen.

}%Aufzaehlung
Wir haben nun nachgewiesen, dass sich das Element für alle Überschneidungen, egal mit welcher Regel begonnen wird, auf dasselbe Element reduzieren lässt.

}%Beweis
}%Satz

%% file: Haupt/Reduktionssysteme/ReduktionssystemeAlgebraOrdnungen.tex
\newcommand{\beq}{\begin{eqnarray*}}
\newcommand{\eeq}{\end{eqnarray*}}

\subsection{Ordnungen}
\label{sec:AllgemeineBemerkungen}
Im vorhergehenden Kapitel haben wir gesehen, wie wir lokale Konfluenz nachweisen können. Als nächstes wollen wir eine Methode entwickeln, wie wir noethersch nachweisen können.

Dazu werden wir für den Rig $\RigF$ eine Ordnung definieren. 
%In einem Wortersetzungssystem $\Red$ definiert man zunächst eine Ordnung auf der Menge der Monome. 
Die Kunst besteht darin eine geeignete Ordnung auszuwählen, damit die Konvergenz von $\RRed$ leicht zu beweisen ist.

\subsubsection{Definitionen}
\label{sec:DefinitionenOrdnung}

\input{ordnungen}

\subsection{Konstruktion von $\red$ }
\label{sec:OrdnungAufDerAlgebra}
Sei $\Algebra$ eine endlich erzeugte Algebra gegeben durch die Erzeuger $a_1, \dots a_{\AnzA}$ und die Relationen $\left\{ s_1, \dots, s_{n_{I_{\Algebra}}}  \right\}$. Dann existiert eine exakte Sequenz von Vektorräumen:
\mathe{ 0 \nach  {I_{\Algebra}} \nach {\FreieA}_{\AnzA} \surj \Algebra \nach 0,}

wobei $I_{\Algebra}$, das durch $\left\{ s_1, \dots, s_{n_{I_{\Algebra}}}  \right\}$ erzeugte Ideal ist.

Sei auf dem Alphabet $\Alphabet = \{ a_1, \dots a_n \}$ eine Totalordnung gegeben und sei $\MonomOrdnung$ die kanonische Ordnung auf den Monomen über $\Alphabet$.
%Wir haben auf der Menge der Monome eine multiplikative noethersche Totalordnung gewählt.
Als Partialordnung auf $\RigF$ wählen wir die folgende:
\Def{Ordnung auf $\RigF$}{
\label{DefOrdnungRigF}
Solange die Terme, in denen das jeweils größte Monom vorkommt, für beide zu vergleichenden Elemente übereinstimmen, ignorieren wir diese Terme. 
\\
Wenn die größten Monome sich unterscheiden, dann vergleichen wir nach der kanonischen Ordnung $\MonomOrdnung$.
%Falls sie übereinstimmen, ist das Element größer mit den meisten Termen in dem das größte Monom vorkommt.
 }

\Bsp{
Seien $w_1, w_2, w_3\in \Monome$ mit $w_1\MonomOrdnung w_2 \MonomOrdnung w_3$ und $\lambda_1, \lambda_2, \lambda_3, \bar{\lambda}_3\in \K$ mit $\lambda_3\neq \bar{\lambda}_3$, dann gilt bezüglich der Ordnung auf $\RigF$:
\AufzaehlungP{
\item $\lambda_1 w_1 +' \lambda_2 w_2 < \lambda_3 w_3$, da $ w_2<w_3$.
\item $\lambda_1 w_1 +' \lambda_3 w_3 > \lambda_3 w_3$, da $w_1>0$.
\item $\lambda_1 w_1 +' \lambda_3 w_3 < \lambda_2 w_2 +'\lambda_3 w_3$, da $w_1<w_2$.
\item Aber: $\lambda_1 w_1 +' \bar{\lambda}_3 w_3$ und $ \lambda_2 w_2 +'\lambda_3 w_3$ sind unvergleichbar, da die größten Monome gleich sind, sich die Koeffizienten aber unterscheiden.
}
}

Wir schreiben jede Relation $s=0$ der Algebra so um, dass das jeweils größte Monom $w$ auf der linken Seite und der Rest, eine Summe von Termen, auf der rechten Seite des Gleichheitszeichens steht. Die Gleichung von links nach rechts gelesen ergibt dann eine Ersetzungsregel. Sie ist von der Form:
\mathe{ w \nach \sum_{i=1}^{n} \lambda_{i} w_{i} .}

Jede Reduktionsregel hat also die Eigenschaft, dass jedes Monom auf der rechten Seite kleiner ist als das auf der linken. 
%Da auf den Monomen eine noethersche Totalordnug definiert war, gibt es nur endlich viele Monome, die kleiner sind als $w$, also ist das Reduktionssystem endlich.

Die Menge dieser Reduktionsregeln bezeichnen wir mit $\red$. Sie induzieren eine Menge von Reduktionsregeln auf $\FreieA$ bzw. $\RigF$, die wir mit $\Red$ bzw. $\RRed$ (siehe: \ref{DefInduziertesRedAlgebra} bzw. \ref{DefInduziertesRedRig}) bezeichnen.

\Def{Gröbnerbasis}{
Falls $\red$ vollständig ist,  nennen wir die Menge $\left\{ s_1, \dots, s_{n_{I_{\Algebra}}}  \right\}$  eine Gröbnerbasis von $I_{\Algebra}$ bezüglich $\MonomOrdnung$.}

% Um deutlich zu machen, dass ein Reduktionssystem vollständig ist, schreiben wir statt $\red$ oft auch $\Vred$.
 
\Bem{lineare Basis}{
\label{BemLineareBasis}
Falls $\red$ vollständig ist, dann bilden die Wörter in $\FreieA$ die sich nicht reduzieren lassen eine lineare Basis für $\Algebra$.
}

Da $\RigF$ eine Rig ist, induzieren die 
 Verknüpfungen  für jede Ersetzungsregel eine Vielzahl an weiteren Regeln:
\mathe{  p \cdot w \cdot s +' q \RRedRegel \RSum{i=1}{n}  p \cdot \lambda_{i} w_{i} \cdot s +' q ,}
wobei $p,w, w_i s \in \Monome$ und $q \in \RigF$.
Da wir die Ordnung geeignet gewählt haben, verkleinert auch jede so erzeugte Regel.

\Bem{Wahl von $S$}{
\label{BemWahlSAlgebra}
Wir wollen hier noch einmal besonders darauf hinweisen, dass es zu einer Algebra verschiedene Mengen von Relationen gibt, die sie definieren. Von der Wahl der Relationen hängt es ab, ob das induzierte Reduktionssystem vollständig ist oder nicht.
}

%% file: ordnungen.tex
\newcommand{\BB}{\mathcal{B}}

\newcommand{\my}{\mu}
\newcommand{\ny}{\nu}

\newcommand{\maybeunderline}{}

%\newtheorem{Proposition}{Proposition}
%newtheorem{Beweis}{Beweis}

%\newtheorem{Beispiel}{Beispiel}

\newcommand{\SchalterMitKuerzen}[1]{%#1
}

\newcommand{\simless}{<\!\!\!\!\!\!_{_\sim}\,\,}

\Def{Striktordnung}{
Eine strikte Teilordnung, oder kurz {  Striktordnung}, ist eine Relation $<$
auf einer Menge $M$, die irreflexiv ($a\not<a$ für alle $a\in M$) und
transitiv ($a<b$ und $b<c$, dann $a<c$ für alle $a,b,c\in M$) ist.
}

\Def{Quasiordnung}{
Eine Quasiordnung ist eine Relation $\simless$
auf einer Menge $M$, die reflexiv ($a\simless a$ für alle $a\in M$) und
transitiv ist.
}

Es ist stets möglich, aus einer Quasiordnung eine Äquivalenzrelation zu
extrahieren, mittels:
\Gleichung{
a\approx b &:\Leftrightarrow& a\simless b \textnormal{ und } b\simless a.
}
Der verbleibende, strikte Anteil ist dann definiert durch:
\Gleichung{
a < b &:\Leftrightarrow& a\simless b \textnormal{ und nicht } b\simless a.
}
Äquivalenzrelation und strikter Anteil geben gemeinsam die Quasiordnung
voll\-ständig wieder. Offensichtlich ist eine Quasiordnung genau dann eine
Teilordnung, wenn ihre Äquivalenzrelation die Gleichheit ist.

\Def{Wohlordnung}{
Eine Striktordnung $<$ auf $M$ heißt {  (strikte) Wohlordnung}, falls
folgende Bedingungen erfüllt sind:
\AufzaehlungP{
\item Strikte Totalordnung: Für jedes Paar $a,b\in M$ ist entweder $a=b$ oder
$a<b$ oder $b<a$.
\item Noethersch: Es gibt keine unendliche, strikt absteigende Folge in $M$.
Äquivalent dazu: Jede absteigende Folge wird stationär, d.h.
\Gleichung{
(a_j)_{j\in\N} \subset M \textnormal{ mit } a_j\,\leq\, a_{j+1}\; \forall j\in
\N &\Rightarrow& \exists N\in N\; \forall j \geq N:\; a_j \,=\, a_N.
}
}
Eine Quasiordnung heißt noethersch, wenn ihr strikter Anteil noethersch
ist.
}

\Def{Multiplikativität}{
Eine Quasiordnung $\simless$ auf einer Halbgruppe $M$ (bspw. die Menge der
Monome über einem Alphabet) heißt {  multiplikativ}, wenn für alle
$p, a, b, s\in M$  gilt:

\Gleichung{
a \;<\; b &\Rightarrow pa \;<\; pb \qquad\textnormal{(Präfix-Invarianz)}\\
a \;\approx\; b &\Rightarrow pa \;\approx pb \\
\textnormal{und }\quad a \;<\; b &\Rightarrow as \;<\; bs\\
a \;\approx\; b &\Rightarrow as \;\approx bs 
\qquad\textnormal{(Suffix-Invarianz)}.
}

\SchalterMitKuerzen{
Wenn eine multiplikative Ordnung auf $M$ zudem die Umkehrungen
\Gleichung{
pa \;<\; pb &\Rightarrow& a \;<\; b\\
pa \;\approx\; pb &\Rightarrow& a \;\approx\; b\\
as \;<\; bs &\Rightarrow& a \;<\; b\\
\textnormal{und }\quad as \;\approx\; bs &\Rightarrow& a \;\approx\; b
}
erfüllt, so nennen wir sie {  multiplikativ \SchalterMitKuerzen{mit Kürzen}}.
}

}

\Def{kombinierte Ordnung}{
Seien $\simless_1$ und $\simless_2$ zwei Quasiordnungen auf $M$. Dann
definieren wir die {\bf kombinierte Ordnung} $\simless_{1,2}$ durch:
\Gleichung{
a\,\simless_{1,2}\,b &:\Leftrightarrow& a \,<_1\, b\;\textnormal{ oder }\;
(a\approx_1 b\; \textnormal{ und }\; a\,\simless_2\, b)
}
.}

In der kombinierten Ordnung $\simless_{1,2}$ wird zuerst nach der Ordnung
$\simless_1$ sortiert. Erst wenn zwei Elemente bezüglich $\simless_1$ nicht
angeordnet sind, wird $\simless_2$ zu Rate gezogen (s.a. \cite{Avenhaus},
Kapitel 1.3). Zwei Elemente sind bezüglich $\simless_{1,2}$ äquivalent,
wenn sie sowohl in $\approx_1$ als auch in $\approx_2$ äquivalent sind.

\Satz{kombinierte Ordnung}{
\label{SatzKomOrd}
Die kombinierte Ordnung $\simless_{1,2}$ ist eine Quasiordnung.
}
\Beweis{Satz \ref{SatzKomOrd}}{
Reflexivität: Folgt aus der Reflexivität von $\simless_1$ und $\simless_2$.

Transitivität: Sei $a\simless_{1,2}b\simless_{1,2}c$ mit $a,b,c\in M$. Dann
gibt es vier Fälle:
\AufzaehlungP{
\item
$a<_1 b<_1 c$:\quad $a<_{1,2} c$ folgt direkt.
\item
$a\approx_1 b<_1 c$:\quad Es gilt $a\,\simless_1\, c$, aber nicht
$a\,\approx_1\, c$, also $a \,<_1\, c$. %Folglich $a\,\simless_{1,2} c$.
\item
$a<_1 b \approx_1 c$:\quad Wie Fall (2).
\item
$a\approx_1 b\approx_1 c$:\quad Es gilt $a\approx_1 c$ und $a\simless_2 c$.
}
}

Sind $\simless_1$ und $\simless_2$ noethersch, so ist es auch $\simless_{1,2}$:\\ Sei
$(a_j)\subset \Alphabet$ eine in $\simless_{1,2}$ strikt absteigende Folge. Es kann nur
eine endliche Unterfolge geben, die auch in $\simless_1$ strikt abfällt. Alle
darauffolgenden Folgenglieder müssen in $\simless_2$ eine absteigende Folge
bilden. Diese muss stationär werden, da $\simless_2$ noethersch ist.

\Satz{multiplikativ \SchalterMitKuerzen{mit Kürzen}}{
\label{SatzMultiplikativ}
Sind $\simless_1$ und $\simless_2$ multiplikativ
\SchalterMitKuerzen{ mit Kürzen}, so ist es auch
$\simless_{1,2}$.

\Beweis{Satz \ref{SatzMultiplikativ}}{

Sei $p \in M$. Für die Präfix-Invarianz müssen  wir die folgenden Fälle betrachten:
\AufzaehlungP{
\item $a <_{1,2} b$\\
Falls $a <_{1} b$, dann folgt $pa <_1 pb$. Also gilt $pa <_{1,2} pb$.\\
Falls $a \approx_{1} b$ und $a <_2 b$, dann ist $pa \approx_{1} pb$ und $pa <_2 pb$, folglich gilt $pa <_{1,2} pb$.
\item $a \approx_{1,2} b$, dann ist $a \approx_{1} b$ und $a \approx_{2} b$.\\ Es folgt, dass $pa \approx_{1} pb$ und $pa \approx_{2} pb$. Also gilt $pa \approx_{1,2} pb$.
}
Die Suffix-Invarianz beweist man analog.
}

%Alt
%Seien $a,b,m\in \Alphabet$ beliebig. Sei zunächst $a \simless_{1,2} b$. Dann gilt
%entweder $a <_1 b$, und infolgedessen $ma \simless_{1,2} mb$ und $am <_{1,2}
%bm$, oder es gilt $a\approx_1 b$ und $a\simless_2 b$. Aufgrund der
%Kürzungseigenschaft von $<_1$ folgt $ma\not<_1 mb$, $am\not<_1 bm$, $mb\not<_1
%ma$ und $bm\not<_1 am$. Zusammen mit der Multiplikativität von $<_2$ ist damit
%$<_{1,2}$ multiplikativ.
%
%Sei nun $ma <_{1,2} mb$. Im Fall $ma <_1 mb$ folgt direkt $a <_1 b$. Sei also
%$ma \not<_1 mb$, $mb\not<_1 ma$ und $ma <_2 mb$. Aus der Multiplikativität von
%$<_1$ folgt $a\not<_1 b$ und $b\not<_1 a$ und aus der Kürzbarkeit von $<_2$
%folgt $a<_2 b$, zusammmen $a<_{1,2} b$.
%
%Analog zeigt man $am <_{1,2} bm\;\Rightarrow\; a <_{1,2} b$.
}

Sei nun $(\Alphabet, <)$ eine partiell geordnete Menge.

\Def{protolexikografische Ordnung}{
Die {\bf protolexikografische Ordnung} auf den Monomen über $\Alphabet$ ist
definiert durch:
\Gleichung{
&a_1\ldots a_n \,<_\textnormal{plex}\, a'_1\ldots a'_m\\
:\Leftrightarrow &n\,=\,m\, \textnormal{ und } \left\{\begin{array}{ll}
                 &a_1 \,<\, a'_1\\
\textnormal{ oder }&a_1 \,=\, a'_1,\, a_2 \,<\, a'_2\\
\textnormal{ oder }&\ldots\\
\textnormal{ oder }&a_1 \,=\, a'_1,\,\, \ldots,\, a_{n-1}\,=\,a'_{n-1},\,
a_n\,<\,a'_n\\
\end{array}\right.
}
für beliebige $a_j, a'_j \in \Alphabet, j = 1\ldots n$.
}

Diese Ordnung ist nur eine strikte Teilordnung, da Worte verschiedener Länge
nicht verglichen werden können. Sie ist multiplikativ \SchalterMitKuerzen{mit Kürzen} nach
Definition. Allerdings ist $<_\textnormal{plex}$ eine Totalordnung, wenn man
sich auf Monome fester Länge beschränkt, und das Alphabet $\Alphabet$ total
geordnet ist.

Ist $\Alphabet$ noethersch (z.B. weil $\Alphabet$ endlich ist), so ist auch
$<_\textnormal{plex}$ noethersch: Da wir nur Wörter gleicher und endlicher
Länge vergleichen, muss jede absteigende Folge letztlich auf eine absteigende
Folge von Buchstaben aus $\Alphabet$ an einer stationären Stelle in den Wörtern
hinauslaufen.

\Def{lexikografische Ordnung}{
Die {\bf lexikografische Ordnung} erhält man aus der protolexikografischen,
indem man $\Alphabet_e = \Alphabet \cup \{e\}$ setzt, mit einem neuen kleinsten Element $e<a$
für alle $a\in\Alphabet$. Die lexikografische Ordnung auf $\Alphabet$ ist dann
\Gleichung{
a \,<_\textnormal{lex}\, a' &:\Leftrightarrow& a e^{|a'|}
\,<_\textnormal{plex}\,a'e^{|a|},
}
wobei $<_\textnormal{plex}$ die protolexikografische Ordnung auf $\Alphabet_e$ ist.
D.h., es werden lediglich genügend viele $e$ an die Enden der Monome
angehängt, um ihre Längen anzugleichen.
}

Falls $\Alphabet$ total geordnet ist, ist auch die lexikografische Ordnung eine totale
Ordnung, aber nicht noethersch: Ein Gegenbeispiel ist die Folge $b > ab > aab >
aaab > \ldots$ mit $a<b$. Sie ist auch nicht multiplikativ: Mit $a<b$ ist zwar
$a < aa$, aber $ab > aab$.

\Def{Knuth-Bendix-Ordnung}{
Die {\bf Knuth-Bendix-Ordnung} auf den Monomen über $\Alphabet$ zu einer
Gewichtsfunktion $g:\Alphabet \rightarrow \N_+$ ist definiert durch:
\Gleichung{
a\,\simless\,a' &:\Leftrightarrow& G(a) \,\leq\, G(a'),
}
wobei $G(a_1\ldots a_n) := \sum_{j=1}^n g(a_j)$.
}

Die Knuth-Bendix-Ordnung ist eine Quasiordnung (mit reell-positiven
Koeffizienten ist sie genau dann eine strikte Totalordnung, wenn die Menge
der Gewichte $\{G(a_j):a_j\in \Alphabet \}$ linear unabhängig über $\Q$ ist). Sie ist
offensichtlich multiplikativ\SchalterMitKuerzen{ mit Kürzen}, und auch noethersch, da die
natürliche Ordnung auf $\N_+$ noethersch ist. 

\Def{Längenordnung}{
Die {\bf Längenordnung} $\simless_\textnormal{Länge}$ auf den Monomen über
$\Alphabet$ ist gegeben durch die Knuth-Bendix-Ordnung zur Gewichtsfunktion $G(a)=1
\,\forall a\in\Alphabet$.
}

Die Längenordnung zählt nur die Zahl der Buchstaben in einem Monom. Wie die
allgemeine Knuth-Bendix-Ordnung ist auch die Längenordnung im nicht-trivialen
Fall eine strikte Teilordnung, multiplikativ\SchalterMitKuerzen{mit Kürzen} und noethersch. In
der Längenordnung sind zwei Monome genau dann äquivalent, wenn sie die gleiche
Länge besitzen.

\Def{kanonische Ordnung}{
Die {\bf kanonische Ordnung} oder {\bf entgegengesetzt lexikographische Ordnung}
$\simless_\textnormal{kan}$ auf den Monomen über einer total geordneten Menge
$\Alphabet$ ist $\simless_\textnormal{Länge, plex}$.
}

Die kanonische Ordnung ist eine Teilordnung: Zwei Monome sind in ihr
äquivalent, wenn sie gleiche Länge haben und aus den gleichen Buchstaben
bestehen. Also ist $\approx_\textnormal{kan}$ die Gleichheit. Sie ist sogar
eine Totalordnung (wenn $\Alphabet$ total geordnet ist), da $<_\textnormal{plex}$ eine
Totalordnung auf den Mengen gleich langer Monome ist.

Sie ist zudem multiplikativ \SchalterMitKuerzen{mit Kürzen} und für noethersches $\Alphabet$ auch
noethersch, da sie eine Kombination von kürzbar-multiplikativen und
noetherschen Ordnungen ist. Sie lässt sich auch analog zur lexikographischen
Ordnung schreiben als
\Gleichung{
a \,<_\textnormal{kan}\, a' &:\Leftrightarrow e^{|a'|} a
\,<_\textnormal{plex}\,e^{|a|} a',
}
da kürzere Monome hier automatisch als kleiner eingeordnet werden.
Sofern nicht anders vereinbart, wird auf den Monomen über einem gegebenen
geordneten Alphabet stets die kanonische Ordnung angenommen.

\Def{Silbenordnung}{
Seien nun $(M, <_M)$ und $(N, <_N)$ strikt total geordnete Mengen (beispielsweise
Monome über einem Alphabet mit einer der obigen Ordnungen), und
\Gleichung{
\Monome(M, N) &:=& \{\my_1\ny_1\my_2\ny_2\ldots \my_n\ny_n\my_{n+1}:
\my_j \in M, \ny_j \in N, n\in \N_0\}
}
die Menge der Worte mit Silben abwechselnd aus $M$ und $N$.
Auf dieser Menge sind die Striktordnungen
\Gleichung{
\my_1\ny_1\ldots \my_{n+1} \,<_{M'}\, \my'_1\ny'_1\ldots \my'_{m+1}
&:\Leftrightarrow \my_1\my_2\ldots \my_{n+1} \,<_\textnormal{kan($M$)}\,
\my'_1\my'_2\ldots \my'_{m+1}\\
\my_1\ny_1\ldots \my_{n+1} \,<_{N'}\, \my'_1\ny'_1\ldots \my'_{m+1}
&:\Leftrightarrow \ny_1\ny_2\ldots \ny_n \,<_\textnormal{kan($N$)}\,
\ny'_1\ny'_2\ldots \ny'_m
}
durch die Projektionen von $\Monome(M,N)$ auf die Monome über $M$ und $N$
definiert. Besteht $N$ nur aus einem Element $c\in \Alphabet$, und ist $M$ seinerseits
die Menge der Monome über $\Alphabet \setminus\{c\}$ mit kanonischer Ordnung, so
nennen wir die kombinierte Ordnung $<_{N', M'}$ auf der Menge $\Monome(M,N)$ der
Monome über $\Alphabet$ {\bf Silbenordnung zum Trennbuchstaben $c$}. Allgemein, wenn
$N = \BB\subset \Alphabet$ mit induzierter Ordnung von $\Alphabet$ ist, und $M$ die Menge der
Monome über $\Alphabet \setminus \BB$ mit kanonischer Ordnung, dann nennen wir
$<_{N',M'}$ {\bf Silbenordnung zu den Trennbuchstaben $\BB$}.
}

Bei der Silbenordnung ist zu beachten, dass die kanonische Ordnung auf zwei
getrennten Ebenen stattfindet, die folgenden Beispiele mögen dies
erläutern:

\beispiel{
Sei $\Alphabet = \{a_1, a_2, a_3, b_1, b_2\}$ mit $a_1<a_2<a_3$ und $b_1<b_2$. Dann
gilt für Monome über $\Alphabet$ in kanonischer Ordnung:
\Gleichung{
a_3 &\,<\, a_1a_2  \textnormal{, da das linke Monom kürzer ist,}\\
a_1a_2 &\,<\, a_3a_1  \textnormal{, da an der ersten Unterscheidungsstelle gilt:
}a_1 < a_3.
}
Ist $\BB = \{b_1, b_2\} \subset \Alphabet$ die Menge der Trennbuchstaben, dann gilt
für die Silbenordnung zu den Trennbuchstaben $\BB$:
%(die wir unterstrichen hervorheben):
\mathe{\begin{array}{rcll}
a_3{\maybeunderline b_1}a_3 										 &\,<\,& a_1{\maybeunderline b_2}a_1a_2  & \textnormal{, da $b_1<b_2$,}\\
a_3{\maybeunderline b_1}a_3a_1a_2a_1a_2a_3a_1a_2 &\,<\,& a_3{\maybeunderline b_1}{\maybeunderline b_1} & \textnormal{, da $b_1<b_1b_1$,}\\
a_1{\maybeunderline b_1}a_1a_2{\maybeunderline b_2} &\,<\,& a_3{\maybeunderline
b_2}a_1a_2{\maybeunderline b_1} & \textnormal{, da $b_1b_2 < b_2b_1$,}\\
a_3a_2a_3 &\,<\, & {\maybeunderline b_1}
  & \textnormal{, da $ \emptyset < b_1$,}\\
a_1{\maybeunderline b_1}a_3 &\,<\,& a_1{\maybeunderline b_1}a_2a_1
  & \textnormal{, da $a_3 < a_2a_1$,}\\
a_2{\maybeunderline b_1}a_3a_2 &\,<\ & a_1a_3{\maybeunderline b_1}
  & \textnormal{, da $a_2 < a_1a_3$.}
\end{array}}
Im letzten Beispiel ist zu beachten, dass $ \emptyset$ durchaus ein Monom über
$\Alphabet$ ist und damit der linke wie der rechte Ausdruck je zwei Silben enthält,
der linke $a_2$ und $a_3a_2$, der rechte $a_1 a_3$ und $ \emptyset$. Nach
Längenordnung auf Ebene der Monome von Monomen sind sie gleich, also
$|(a_2)(a_3a_2)| = |(a_1a_3)( \emptyset)| = 2$, und es entscheiden die ersten
Silben $a_2$ und $a_1a_3$, die ihrerseits Monome über  $\Alphabet$ sind. Auf
dieser unteren Ebene ist dann $a_2 < a_1a_3$ wie üblich.
}

\Satz{Silbenordnung}{
\label{SatzSilbenordnung}
Sei $\Alphabet$ noethersch (z.B. endlich) und $\BB \subset \Alphabet$. Dann ist die
Silbenordnung zu den Trennbuchstaben $\BB$ auf der Menge $\Monome_{\Alphabet}$ der Monome über
$\Alphabet$ noethersch und multiplikativ\SchalterMitKuerzen{mit Kürzen}.
}
\Beweis{Satz \ref{SatzSilbenordnung}}{
Wir brauchen nur zu zeigen, dass die Projektionsordnungen $<_M$ und $<_N$ mit
$M=\Alphabet\setminus \BB$ und $N=\BB$ noethersch und multiplikativ\SchalterMitKuerzen{mit Kürzen} sind,
dann folgt die Behauptung aus den Eigenschaften der kombinierten Ordnung.

{\bf Noethersch}: Jede in $<_M$ strikt abfallende Folge $(a_j)\subset \Monome_{\Alphabet}$
induziert eine strikt absteigende Folge $(a'_j)\subset \Monome_M$ in den Monomen über
$M$. Die kanonische Ordnung auf $\Monome_M$ ist aber noethersch, da sie von der
kanonischen Ordnung auf $\Monome_{\Alphabet}$ induziert wird. Folglich gibt es keine strikt
absteigende Folge in $\Monome_M$, und damit ist $<_M$ noethersch. Analog gilt das für
$<_N$.

{\bf Multiplikativität\SchalterMitKuerzen{mit Kürzen}}: Da $<_{\textnormal{kan}}$ multiplikativ\SchalterMitKuerzen{ mit Kürzen} ist, gilt für $a\in M$.
\Gleichung{
\my_1\ny_1\ldots \my_{n+1} &<_M \my'_1\ny'_1\ldots \my'_{m+1}\\
\Leftrightarrow \;\,\qquad
\my_1\ldots \my_{n+1} &<_\textnormal{kan} \my'_1\ldots \my'_{m+1}\\
\Rightarrow \qquad
\my_1\ldots \my_{n+1}a &<_\textnormal{kan} \my'_1\ldots \my'_{m+1}a\\
\Leftrightarrow \quad
\my_1\ny_1\ldots \my_{n+1}a &<_M \my'_1\ny'_1\ldots \my'_{m+1}a.
}
Ist $a$ dagegen in $N$, so bleibt die Projektion nach $\Monome_M$ unverändert, und
es gilt nach wie vor $\my_1\ny_1\ldots \my_{m+1}a\;<_M\;
\my_1\ny_1\ldots\my_{m+1}a$. \\
Analog gilt das für Linksmultiplikation und $<_N$.
}

%% file: Haupt/Reduktionssysteme/ReduktionssystemeAlgebraNoethersch.tex
\subsection{Noethersch}
\label{sec:Noethersch}

\Satz{Noethersch}{\label{noetherschRig}Die Ordnung auf $\RigF$ ist noethersch, und für jedes $(x,y) \in \RRed$ gilt $x>y$, insbesondere ist $\RRed$ noethersch.
\label{Noethersch}

\Beweis{}{
Da die Ordnung $\MonomOrdnung$ noethersch und multiplikativ ist, ist es auch die induzierte Ordnung auf $\RigF$. Dass jede Reduktionsregel verkleinert, folgt direkt aus der Definition (\ref{Rtilde}).

}

}

\Satz{Vollständig}{
\label{vollstaendig}
Falls die Ordnung auf Monomen $\MonomOrdnung$ noethersch und multiplikativ ist und falls $\red$ vollständig ist, dann gilt:
 \AufzaehlungZ{
 \item
Die induzierten Reduktionssysteme $\Red$ und $\RRed$  auf $\FreieA$ bzw. $\RigF$ sind konvergent.
\item
Zwei Elemente aus $\FreieA$ haben genau dann dieselbe Normalform, wenn sie durch die kanonische Projektion $\pi_{\Algebra}: \FreieA \rightarrow \Algebra $ auf dasselbe Element abgebildet werden.
}

\Beweis{}{~\\
Zu $1.$: Da die Ordnung multiplikativ ist, gilt nach \ref{Noethersch}, dass auch  $\RRed$ noethersch ist. Da $\red$ vollständig ist, folgt aus \ref{vollstaendigFolgtlokalkonfluent}, dass $\RRed$ lokal konfluent ist. Mit \ref{Newman} ist $\RRed$ total konfluent und nach \ref{konfluentnoethersch} ist auch $\Red$ total konfluent.

Zu $2.$: "'$\folgt$"' Seien $f, \bar{f} \in \FreieA$ mit derselben Normalform. Zuerst überlegen wir  uns, dass $f -\Nf(f) \in I_{\Algebra}$ gilt:
Sei $(f~,~f_1)$ eine beliebige Reduktionsregel in $\Red$, die mit $f$ beginnt. Sei $s=0$ die Relation, die die Reduktionsregel $(w~,\Sumn{i} \lambda_{i} w)$ in $\red$ erzeugt, die $(f~,~f_1)$ induziert. Dann ist 
\mathe{ f- f_1 =  p (w- \Sumn{i} \lambda_{i} w) s \in I_{\Algebra}.}
Durch mehrfaches Anwenden sehen wir, dass auch $f- \Nf(f)\in I_{\Algebra}$.

Für die Differenz von $f$ und $\bar{f}$ gilt dann also:
\mathe{ 
f - \bar{f} = f -\Nf(f) + (\Nf(\bar{f}) - \bar{f}) \in I_{\Algebra}
;}
somit werden sie auf dasselbe Element in $\Algebra$ projiziert.

"'$\Leftarrow$"' Sei $\pi_{\Algebra}(f) = \pi_{\Algebra}(\bar{f})$,  dann gilt:
 
\Gleichung{
f - \bar{f} \in I_{\Algebra} 
&\folgt f - \bar{f} = \sum_{k=1}^m x_{i_k}
s_{i_k} y_{i_k}\\
& \folgt \textnormal{$f$ und $\bar{f}$ sind Church-Rosser-äquivalent bezüglich $\Red$}\\
&\stackrel{\ref{totalkonfluent}}{\folgt} \Nf(f) = \Nf(\bar{f}).
}

}
}

%% file: Haupt/Reduktionssysteme/ReduktionssystemeModul.tex
\section{Reduktionssystem bezüglich eines Moduls}
\label{sec:Modul}

Reduktionssysteme kann man auch für Moduln und für Modulabbildungen nutzen, vgl. \cite{MR506890}. 
Wir wollen in diesem Kapitel aber keine  neue Theorie entwickeln, sondern einen Modul als Teilmenge einer Algebra beschreiben. Genauer, wir beschreiben $\Modul$ als eine Untermenge der freien Algebra $\FreieA_{\AnzA + \AnzM}$, damit wir die in \ref{sec:ReduktionssystemFbezüglichEinerAlgebra} beschriebenen Methoden für ein Reduktionssystem für Algebren nutzen können.

%% file: Haupt/Reduktionssysteme/P-Vollstaendig.tex
\subsection{$\P$-Vollständigkeit}
\label{sec:PrädikatP}

In \ref{vollstaendigFolgtlokalkonfluent} hatten wir gesehen, dass wir Konfluenz nachweisen können, indem wir die minimalen  Überschneidungen untersuchen. Nun werden wir zeigen, dass wir nicht alle minimalen Überschneidungen betrachten müssen.

\Def{$\Red$-Prädikat $\P$}{
\label{DefRPraedikat}
 Sei $\P:\FreieA \nach \{ \wahr, \falsch \}$ ein Prädikat, das folgende Bedingung erfüllt:
\AufzaehlungP{
\item Falls $(x,y)\in \Red$ und $\P(x)=\wahr$, dann ist auch $\P(y)=\wahr$.
}
Ein solches Prädikat nennen wir $\Red$-Prädikat.
}

\Def{$\Rp$}{
\label{DefRp}
Für ein $\Red$-Prädikat $\P$ bezeichnen wir die Regeln aus $\Red$, die mit einem $x$ mit $\P(x)=\wahr$ beginnen, mit $\Rp$, also
\mathe{\Rp:=\Menge{(x,y)\in \Red}{\P(x)=\wahr}.}
Die Untermenge in $\MengeE$, die Urbild von $\wahr$ unter $\P$ ist, bezeichnen wir mit $\MengeE_{\P}$, also
\mathe{\MengeE_{\P}:=\Menge{x\in \MengeE}{\P(x)=\wahr}.}
}

\Bem{Folgerung}{
\label{BemP}
Sei $\P$ ein $\Red$-Prädikat. 
Folgende Eigenschaften ergeben sich für $\Rp$ direkt aus der Definition:
\AufzaehlungP{
\item Es gilt: \mathe{\Rp=\Menge{(x,y)\in \Red}{\P(x)=\wahr \textnormal{ und } \P(y)=\wahr}. }
\item Sei $y_1 \RedWeg y_n$  ein Reduktionsweg in $\Red$ und sei $\P(y_1)=\wahr$, dann gilt für jedes $i=1\dots n$:\mathe{\P(y_i)=\wahr.}
\item Sei $y_1 \RedWeg y_n$ ein Reduktionsweg in $\Red$ und sei $\P(y_n)=\falsch$, dann gilt für jedes $i=1\dots n$:\mathe{\P(y_i)=\falsch.}
\item Ist $\Red$ konvergent, dann ist auch die Teilmenge $\Rp$ auf $\MengeE_{\P}$ ein konvergentes Reduktionssystem.
}%Aufzaehlung
}

\Def{Silbenverträglichkeit}{
\label{DefSilbenvertraeglich}
Ein Prädikat $\P:\RigF \nach \{ \wahr, \falsch\}$ heißt silbenverträglich, falls für jedes $p,s\in \Monome$ folgende zwei Axiome erfüllt sind:
\AufzaehlungP{
\item {Summandenverträglichkeit: } \mathe{ \P(p+'s)=\wahr \Leftrightarrow \P(p)=\wahr=\P(s).}
\item {Faktorenverträglichkeit: }  \mathe{ \P(ps)=\wahr \folgt \P(p)=\wahr=\P(s).}
}
}

\Bem{Alternative Definition}{
Alternativ könnten wir Silbenverträglichkeit auch definieren durch das Axiom.\\
 Sei $q:=\RSum{w \in \Monome}{} \RSum{i=1}{N_w}  \lambda_{w,i} w$, dann gilt:
\mathe{ \P(q)=\wahr \folgt \P(w)=\wahr \textnormal{ für alle $w$ mit } N_w>0.}
}

\Satz{$\RedpRegel \subset \RRedpWeg$}{
\label{SatzEigenschaftenRRp2}
\label{pxsq}
 Sei $\P$ ein silbenverträgliches  $\RRed$-Prädikat, dann 
induziert jede Regel aus $\Red_{\P}$ einen Weg in $\RRed_{\P}$.

\Beweis{Satz \ref{SatzEigenschaftenRRp2}}{

Sei $(g_1,g_2)\in \Red_{\P}$, dann gibt es $(x,y)\in \red$ mit:
\mathe{ g_1= \lambda p x s + q \textnormal{ und } g_2= \lambda p (y) s + q .}
Mit Lemma \ref{RsubsetTR} lässt sich $(g_1, g_2)$ als Weg in $\RRed$ schreiben:
\mathe{ \lambda p x s +' q \RRedWeg \lambda p (y) s +' q.}
Nach Voraussetzung ist $\P(\lambda p x s +' q)=\wahr$ und mit Bemerkung \ref{BemP} liegt dann der gesamte Weg in $\RRed_{\P}$.

%
%Da $\RRed$ durch $\red$ induziert ist, gilt $(pxs+'q,pys+'q) \in \RRed$. Mit Bemerkung \ref{BemP} gilt für ein $\RRed$-Prädikat, dass aus $\P(pxs+'q)=\wahr$ folgt $\P(pys+'q)=\wahr$. Da $\P$ silbenverträglich ist, gilt auch $\P(x)=\P(y)=\wahr$. Eine Reduktionsregel aus $\Red_{\P}$ lässt sich also übersetzen in einen Weg in  $\RRed_{\P}$. 
}
}

\Satz{Konvergenz}{
\label{SatzKonvergenz2}
Sei $\P$ ein silbenverträgliches  $\RRed$-Prädikat. Sei $\RRed_{\P}$ konvergent, dann ist auch $\Red_{\P}$ konvergent.

\Beweis{Satz \ref{SatzKonvergenz2}}{
Durch mehrfaches Anwenden von Satz \ref{SatzEigenschaftenRRp2} folgt, dass alle Wege aus $\Red_{\P}$ einen Weg in $\RRed_{\P}$ bilden. Nun ist dieser Beweis  nach Ersetzung von $\FreieA, \RigF, \Red, \RRed$ durch  $\FreieA_P, \RigF_P, \Red_P, \RRed_P$ identisch mit dem für Satz \ref{konfluentnoethersch}. 
}%Beweis
}%Satz

\Def{$\P$-minimale Überschneidung}{
Sei $\P$ ein silbenverträgliches $\RRed$-Prädikat.
Eine minimale Überschneidung $(w, \regel, \bar{\regel})$ mit  $\P(w)=\wahr$, heißt $\P$-minimale Überschneidung.}

\Def{$\P$-vollständig}{
Ein Wortersetzungssystem $\red$ heißt $\P$-vollständig, falls jede $\P$-minimale Überschneidung behebbar bezüglich $\RRed$ ist.
}

\Satz{$\P$-vollständig $\folgt$ Konfluenz}{
\label{konfluentModul}
\label{konvergentModul}
Sei $\red$ ein $\P$-vollständiges Wortersetzungssystem, dann ist $(\RigF_{\P}, \RRed_{\P})$ ein lokal konfluentes Reduktionssystem.

\Beweis{Satz \ref{konfluentModul}}{~\\
Sei $f \in \RigF_{\P}$ und $(f,g_1), (f,g_2) \in \RRed_{\P}$, dann gilt auch $(f,g_1), (f,g_2) \in \RRed$.
Falls $(f,g_1), (f,g_2) \notin \Rel$ gilt, dann liegt wenigstens einer von beiden in $\widetilde{\Rel}$. Der Beweis ist nun genauso wie der für Satz \ref{vollstaendigFolgtlokalkonfluent}, nur dass wir anstelle von Bemerkung \ref{axb} den Satz \ref{pxsq} nutzen.

Falls $(f,g_1), (f,g_2) \in \Rel$ gilt, dann lassen sich $f$ und $g_i$ schreiben als:
\Gleichung{
f   &= \lambda_1 p_1 x_1 s_1 +' q_1 &= \lambda_2 p_2 x_2 s_2 +' q_2, \\
g_i &= \lambda_i p_i (y_i) s_i +' q_i , 
}
wobei $\lambda_i \in \K^*$ und  $(x_i,y_i) \in \red$, sowie $p_i, s_i \in \Monome$ und $q_i \in \RigF$.\\

Falls $q_1 \neq q_2$, dann muss $q_1= \lambda_2 p_2 x_2 s_2$ und $q_2= \lambda_1 p_1 x_1 s_1$ sein. In diesem Fall lassen sich $g_1$ und $g_2$ zu $\lambda_1 p_1 x_1 s_1+' \lambda_2 p_2 x_2 s_2$ reduzieren.

Sei also $q_1 = q_2$.\\
Falls $\lambda_1 p_1 x_1 s_1 = \lambda_2 p_2 x_2 s_2 $, dann ist der Beweis genauso wie der für Satz \ref{vollstaendigFolgtlokalkonfluent}.

Falls $\lambda_1 p_1 x_1 s_1 \neq \lambda_2 p_2 x_2 s_2 $, dann gibt es folgende Möglichkeiten:
\AufzaehlungZ{
\item{getrennte Silben: } Es gibt ein $\bs{w} \in \Monome$ mit:  
\mathe{ p_1 x_1 s_1 = p_1 x_1 \bs{w} x_2 s_2 = p_2 x_2 s_2.}
\label{getrennteSilben1}
\item{getrennte Silben vertauscht: } Es gibt ein $\bs{w} \in \Monome$ mit:  
\mathe{ p_1 x_1 s_1 = p_2 x_2 \bs{w} x_1 s_1 = p_2 x_2 s_2.}
\label{getrennteSilben2}
\item{volle Überschneidung: } Es gibt $\bs{w_1}, \bs{w_2} \in \Monome$ mit: 
\mathe{ x_1 = \bs{w_1} x_2 \bs{w_2}.}
\label{volleUberschneidung1}
\item{volle Überschneidung vertauscht: } Es gibt $\bs{w_1}, \bs{w_2} \in \Monome$ mit: 
\mathe{ \bs{w_1} x_1 \bs{w_2} =  x_2 .}
\label{volleUberschneidung2}
\item{Teil-Überschneidung: }Es gibt $\bs{w_1}, \bs{w_2} \in \Monome$ mit: 
\mathe{ \bs{w_1} x_1  =  x_2 \bs{w_2} .}
\label{teilUberschneidung1}
\item{Teil-Überschneidung vertauscht: } Es gibt $\bs{w_1}, \bs{w_2} \in \Monome$ mit: 
\mathe{ x_1 \bs{w_2}=  \bs{w_1}  x_2  .}
\label{teilUberschneidung2}
}

Die Fälle \ref{getrennteSilben1} und \ref{getrennteSilben2} können wir wieder wie in Satz \ref{vollstaendigFolgtlokalkonfluent} behandeln.

Die Fälle \ref{volleUberschneidung1} und \ref{volleUberschneidung2} sind symmetrisch, wir zeigen \ref{volleUberschneidung1}:\\
Hier gilt $\P(f)=\wahr$. Da $f=p_1 x_1 s_1+'q$ gilt,  ist wegen der Silbenverträglichkeit von $\P$  auch $\P(p_1 x_1 s_1)=\wahr$, also ist auch $\P(x_1)=\wahr$.\\
Demnach ist die minimale Überschneidung $\left( x_1 = \bs{w_1} x_2 \bs{w_2}~,~ (x_1,y_1)~,~(x_2,y_2)\right)$ eine $\P$-minimale Überschneidung. Nach Voraussetzung lässt sie sich beheben. Es gibt also Reduktionswege:
\mathe{\Matrix{ y_1 &\RRedpRegel \dots \RRedpRegel \\ \bs{w_1} y_2 \bs{w_2}&\RRedpRegel \dots \RRedpRegel} z.} Wegen Satz \ref{pxsq} gibt es dann auch Reduktionswege:
\mathe{\Matrix{ g_1 &\RedpRegel \dots \RedpRegel \\ g_2 &\RedpRegel \dots \RedpRegel} \hat{z}.}

Die Fälle \ref{teilUberschneidung1} und \ref{teilUberschneidung2} sind symmetrisch, wir zeigen \ref{teilUberschneidung1}:\\
Hier ist $\bs{w_1}$  ein Teilwort aus $p_1$ und
$\P(f)=\wahr$, also ist auch $\P(p_1 x_1 s_1)=\wahr$ und auch $\P(\bs{w_1} x_1)=\wahr$.\\
Demnach ist die minimale Überschneidung $\left( \bs{w_1}x_1 =  x_2 \bs{w_2}~,~ (x_1,y_1)~,~(x_2,y_2)\right)$ eine $\P$-minimale Überschneidung. Nach Voraussetzung lässt sie sich beheben. Es gibt also Reduktionswege:
\mathe{\Matrix{ 
\bs{w_1} y_1                 &\RRedpRegel \dots \RRedpRegel \\ 
         y_2\bs{w_2} \bs{w_2}&\RRedpRegel \dots \RRedpRegel} z.} Wegen Satz \ref{pxsq} gibt es dann auch Reduktionswege:
\mathe{\Matrix{ g_1 &\RedpRegel \dots \RedpRegel \\ g_2 &\RedpRegel \dots \RedpRegel} \hat{z}.}

~\\
Es lassen sich also alle Überschneidungen beheben, somit ist $\RRed_{\P}$  lokal konfluent.

}
}

%% file: Haupt/Reduktionssysteme/ModulPraedikat-Neu.tex
\subsection{$\AoA$-Moduln}
\label{sec:AnwendungFürBimodule}

%Im Abschnitt \ref{sec:Modul} hatten wir gesehen, dass wenn wir einen links Modul in eine Tensoralgebra einbetten, die gleichen Reduktionsmethoden wie für eine Algebra nutzen können.
%
%Statt für Bimodule einen entsprechenden Beweis zuführen nutzen wir Prädikate.

Wie in den vorhergehenden Abschnitten sei $\Alphabet$ das Alphabet mit den Buchstaben $a_1, \dots, a_{\AnzA}$ und sei $\Monome$ die Menge der Monome über $\Alphabet$, weiter sei $\Algebra$ eine Algebra und $\FreieAA$ die freie Algebra mit Erzeugern aus $\Alphabet$.
Bezüglich $\Algebra$ sei $\RSa$ ein vollständiges Reduktionssystem.

Nun sei $\Modul$ ein durch $ \bs{e}_1, \dots , \bs{e}_{\AnzM} $ endlich erzeugter $\AoA$-Modul. Mit  $\Monome_{\Algebra,\Modul}$ bezeichnen wir die Menge der Monome über dem Alphabet mit den Buchstaben $\Alphabet \cup \{ \bs{\epsilon}_1, \dots , \bs{\epsilon}_{\AnzM} \}$.
Sei $\FreieAM$ die freie Algebra über $\Monome_{\Algebra,\Modul}$ mit der Silbenordnung zu den Trennungsbuchstaben $\{ \bs{\epsilon}_1, \dots , \bs{\epsilon}_{\AnzM} \}$.

Nehmen wir an, dass es $\AnzUM\in \N$ und $S:=\{ \bs{{\bar{e}}}_1, \dots , \bs{{\bar{e}}}_{\AnzUM} \}$ gibt,  so dass folgende Sequenz exakt ist:
\mathe{ 0 \nach \UModul \inj \left(\AoA\right)^{\AnzM} \stackrel{\varpi_{\Modul}}{\surj} \Modul \nach 0, }
wobei $\UModul$ der durch  durch S erzeugte $\AoA$-Modul ist.

Ein Element $c \in \left(\AoA \right)^{\AnzM}$ schreiben wir als 
\mathe{ c= \Sum{\sigma=1}{N} p_{{\sigma}} \bs{\epsilon}_{i_{\sigma}} s_{{\sigma}}, \textnormal{ wobei $p_{{\sigma}}, s_{{\sigma}} \in \Algebra$ und $\bs{\epsilon}_{i_{\sigma}}$ Erzeuger von $\Modul$ sind. }}

Da $\UModul \inj \left(\AoA\right)^{\AnzM}$ injektiv ist, können wir jeden Erzeuger $\bs{\bar{e}}_{\iota}$ von $\UModul$ mit seinem Bild identifizieren. Es  gibt also  $\iota_{\sigma}\in \{1 \dots \AnzM \}$ und $\bar{p}_{\sigma},\bar{s}_{\sigma}\in \Algebra$, so dass

\mathe{\bar{e}_{\iota} = \Sum{\sigma=1}{\bar{n}_{\iota}} \bar{p}_{{\sigma}} \bs{\epsilon}_{\iota_{\sigma}} \bar{s}_{{\sigma}}.
}

Für jedes $\bar{p}_{\sigma},\bar{s}_{\sigma}$ wählen wir ein Urbild unter der kanonischen Projektion $\pi_{\Algebra}:\FreieAA \surj \Algebra$ und erhalten so ein Element in $\left(\FreieAA \tensor \FreieAA \right)^{\AnzM} \subset \FreieAM$:

\mathe{\tilde{e}_{\iota} = \Sum{\sigma=1}{\bar{n}_{\iota}} \tilde{p}_{{\sigma}} \bs{\epsilon}_{\iota_{\sigma}} \tilde{s}_{{\sigma}}.
}

Aus den Erzeugern von $\UModul$ werden wir nun ein Reduktionssystem $\RSm \subset \Monome_{\Algebra,\Modul} \times \FreieAM$ konstruieren. Dazu stellen wir die obige Summe so um, dass der größte Summand, wir bezeichnen ihn mit $x_{\iota}$, ohne seinen Koeffizienten $\lambda_{\iota}$ auf der linken Seite steht:
\mathe{\RSm := \Menge{ ( {x}_{\iota} ~,~ {x}_{\iota}- {\lambda}_{\iota}^{-1} \tilde{e}_{\iota})}{ \iota =  1 \dots \AnzUM  }.}

Um zu verdeutlichen, dass dieses Wortersetzungssystem vom Erzeugersystem $S$ des Moduls $\UModul$ abhängt, nennen wir $\RSm$ auch das durch $S$ induzierte Reduktionssystem.

\label{BemRedX}
Die Reduktionsregeln in $\RSm$ beginnen jeweils mit einem Wort, in dem genau ein Erzeuger des Moduls vorkommt.

\Bem{Wahl des Erzeugersystems}{
\label{WahlSModul}
Wenn man das Erzeugersystem $S$ für $\UModul$ geeignet gewählt hat, wird das induzierte Reduktionssystem $\red_{\Algebra, Q}$ vollständig. Diese Wahl ist oft jedoch sehr schwer. Es gibt Algorithmen, die aus einem gegebenen Erzeugersystem ein neues Erzeugersystem $S'$ konstruieren, so dass das induzierte Wortersetzungssystem vollständig ist. Leider terminieren diese Algorithmen nicht in jedem Fall.
}

Zur Erinnerung, in den vorhergehenden Kapiteln haben wir ein Wortersetzungssystem $\RSa$ bezüglich einer Algebra $\Algebra$ konstruiert. Wir haben $\RSa$ aus einem Erzeugersystem $\ESa$ für das Ideal $I_{\Algebra}$ konstruiert, wobei $I_{\Algebra}$  definiert war durch:
\mathe{ 0 \nach I_{\Algebra} \inj \FreieA_{n} \surj \Algebra \nach 0 .}

Sei $\RSm$ wie oben konstruiert, dann setzen wir:
\mathe{\RSam:=\RSa \cup \RSm.}

Wir wollen nun ein Prädikat auf $\FreieAM$ definieren, mit dem wir die Elemente des Moduls markieren können. Dazu betrachten wir zunächst folgende Definition.

\Def{$\FreieA(0), \FreieA(1)$}{
Die Menge der Elemente aus $\FreieAM$, deren Summanden keinen bzw. genau einen Erzeuger von $\Modul$ enthalten,  bezeichnen wir mit $\FreieA(0)$ bzw. $\FreieA(1)$.
}

\Lemma{$\FreieA(1)$}{
\label{FEins}
Die kanonische Projektion $\pi_{\Modul} : \FreieA(1) \nach \Modul$ ist surjektiv.

\Beweis{Lemma \ref{FEins} }{Sei $c \in \Modul$, dann gibt es $\iota_{\sigma}\in \{1 \dots \AnzM \}$ und $\bar{p}_{\sigma},\bar{s}_{\sigma}\in \Algebra$ mit:
\mathe{c=\Sum{\sigma=1}{N} p_{i_{\sigma}} \bs{\epsilon}_{i_{\sigma}} s_{i_{\sigma}}.}
 Da die  kanonische Projektion $\pi_{\Algebra}: \FreieAA \nach \Algebra$ surjektiv ist, können wir für jedes $p_i,s_i$ ein Urbild $\tilde{p}_i,\tilde{s}_i \in \FreieAA$ wählen. Dann gilt:
\mathe{\pi_{\Modul}\left(\Sum{\sigma=1}{N} \tilde{p}_{i_{\sigma}} \bs{\epsilon}_{i_{\sigma}} \tilde{s}_{i_{\sigma}} \right) = c.}
}
}

%\Schalter{

\Bem{$\pi_{\Modul}=  \varpi_{\Modul}   \circ \varpi_{\Algebra}$}{Mit $\left( \FreieAA \otimes \FreieAA \right)^{\AnzM} \subset \FreieA(1) \subset \FreieAM$ können wir die kanonische Projektion $\pi_{\Modul}: \FreieA(1) \nach \Modul$ als Verknüpfung der folgenden Abbildungen schreiben: 
\Gleichung{\varpi_{\Algebra}:&  \FreieA(1) \nach \left( \AoA \right)^{ \AnzM} \textnormal{ und }\\ 
\varpi_{\Modul}: &  \left( \AoA \right)^{ \AnzM} \nach \Modul.}
}

\Lemma{Zerlegung von $\Kern{\varpi_{\Algebra}}$}{
\label{LemmaKernZerlegung}

Sei $I:=\Kern{\pi_{\Algebra}: \FreieAA \surj {\Algebra}}$, dann lässt sich $\Kern{\varpi_{\Algebra}}$ wie folgt zerlegen:
\mathe{ \Kern{\varpi_{\Algebra}}=\sum_{i=1}^{\AnzM} I \epsilon_i \FreieAA  +  \sum_{i=1}^{\AnzM} \FreieAA \epsilon_i I.}
}

\Beweis{Lemma \ref{LemmaKernZerlegung}}{
Wir zeigen zunächst "`$\supseteq$"':\\
Alle Summanden liegen schon im Kern, da
\mathe{\varpi_{\Algebra}\left( I \epsilon_i \FreieAA \right) = \varpi_{\Algebra}\left( I \right)\epsilon_i \varpi_{\Algebra}\left(\FreieAA \right) = 0}
und auch
\mathe{\varpi_{\Algebra}\left( \FreieAA  \epsilon_i I\right) = \varpi_{\Algebra}\left( \FreieAA \right)\epsilon_i \varpi_{\Algebra}\left(I \right) = 0.}

Als nächstes zeigen wir "`$\subseteq$"':\\
Sei 
\mathe{ f = \Sum{i=1}{\AnzM}\Sum{j=1}{m_i} \lambda_{i,j} p_{i,j} \epsilon_{i} s_{i,j} \in \Kern{\varpi_{\Algebra}} \subset \FreieA(1),}
wobei $\lambda_{i,j} \in \K$ und $p_{i,j}, s_{i,j}\in \Monome$.

Es gilt also:
\mathe{ 
0=\varpi_{\Algebra}(f)= \Sum{i=1}{\AnzM}\Sum{j=1}{m_i} \lambda_{i,j} \pi_{\Algebra}(p_{i,j}) \epsilon_{i} \pi_{\Algebra}(s_{i,j}) 
.}
Nun ist $\left\{\epsilon_1, \dots, \epsilon_{\AnzM}\right\}$ eine Basis des freien $\Algebra$-Bimoduls $\left(\AoA\right)^{\AnzM}$, also ist für jedes $i=1, \dots {\AnzM}$:
\mathe{0=\Sum{j=1}{m_i} \pi_{\Algebra}(p_{i,j}) \otimes \pi_{\Algebra}(s_{i,j}).}
Anders gesagt, ist 
\mathe{\Sum{j=1}{m_i} (p_{i,j}) \otimes (s_{i,j}) \in \Kern{\pi_{\Algebra} \otimes \pi_{\Algebra} : \FreieAA \otimes \FreieAA \nach \Algebra \otimes \Algebra}.
}
Für diesen Kern gilt jedoch:
\mathe{ \Kern{\pi_{\Algebra} \otimes \pi_{\Algebra} : \FreieAA \otimes \FreieAA \nach \Algebra \otimes \Algebra} = I \otimes \FreieAA + \FreieAA \otimes I.}

Also ist 
\mathe{ \varpi_{\Algebra}(f) \subset \sum_{i=1}^{\AnzM} \left(I \epsilon_i \FreieAA  +   \FreieAA \epsilon_i I \right).}
}

\Satz{Normalform}{
\label{SatzRSChurchRosserAq}
Sei $\RSam$ wie oben konstruiert und sei $\Red$ durch $\RSam$ induziert. Seien $f,\bar{f} \in \FreieA(1)$, dann sind die Aussagen  
\AufzaehlungZ{
\item $ \pi_{\Modul}(f)=\pi_{\Modul}(\bar{f})$
\item $ f \RedAq \bar{f}$
}
äquivalent. Also stimmen, falls $\Red$ noethersch ist, die Menge der Normalformen und das Bild von $\pi_{\Modul}$ überein.

\Beweis{Satz \ref{SatzRSChurchRosserAq}}{~\\
$1 \folgt 2:$\\
Sei $ \pi_{\Modul}(f)=\pi_{\Modul}(\bar{f})$, also ist $ \pi_{\Modul}(f-\bar{f})=0$ und $\varpi_{\Algebra}(f-\bar{f}) \in \Kern{\varpi_{\Modul}}=\UModul$.
Es  gibt also  $\iota_{\sigma}\in \{1 \dots \AnzM \}$ und $\bar{p}_{\sigma},\bar{s}_{\sigma}\in \Algebra$, so dass

\mathe{\varpi_{\Algebra}(f-\bar{f}) = \Sum{\sigma=1}{n} \bar{p}_{{\sigma}} \bs{\epsilon}_{\iota_{\sigma}} \bar{s}_{{\sigma}}.
}

Für jedes $\bar{p}_{\sigma},\bar{s}_{\sigma}$ wählen wir ein Urbild unter der kanonischen Projektion $\pi_{\Algebra}:\FreieAA \surj \Algebra$ und erhalten so ein Urbild unter $\varpi_{\Algebra}$:

\mathe{ \Sum{\sigma=1}{n} \bar{p}_{{\sigma}} \bs{\epsilon}_{\iota_{\sigma}} \bar{s}_{{\sigma}} = \varpi_{\Algebra}\left( \Sum{\sigma=1}{n} \tilde{p}_{{\sigma}} \bs{\epsilon}_{\iota_{\sigma}} \tilde{s}_{{\sigma}}\right).
}

Nun gilt mit Lemma \ref{LemmaKernZerlegung}:
\mathe{
\underbrace{f-\bar{f}- \Sum{\sigma=1}{n} \tilde{p}_{{\sigma}} \bs{\epsilon}_{\iota_{\sigma}} \tilde{s}_{{\sigma}}}_{\in \Kern{\varpi_{\Algebra}}} 
=
\underbrace{g_1}_{\in \sum I \epsilon_i \FreieAA }  + \underbrace{g_2}_{ \sum \FreieAA \epsilon_i I}
}
Der Anteil in $I$ für jeden Summanden aus $g_1$ und $g_2$  lässt sich schreiben  als $p(x-y)s$, wobei $p,s \in \FreieAA$ und $(x,y)\in \RSa$; ähnlich gilt $\epsilon_{i}=(x_i,y_i)$ mit $(x_i,y_i)\in \RSm$. Zusammen sehen wir, dass $f \RedAq \bar{f}$.

~\\
$2 \folgt 1$:
Sei $g_1 \RedAq g_n$, dann gibt es eine endliche Folge $(g_i,g_{i+1})$, wobei  $(g_i,g_{i+1})\in \Red$ oder  $(g_{i+1},g_{i}) \in \Red$. Für die Behauptung  reicht es also zu zeigen: Wenn $g,\bar{g} \in \Red$ ist, dann gilt auch $\pi_{\Modul}(g)=\pi_{\Modul}(\bar{g})$.\\

Eine Regel $(g,\bar{g})\in \Red)$ lässt sich schreiben als:
\mathe{ (\lambda pxs+q, \lambda pys+q),}
 wobei $p,s \in \Monome_{\Algebra, \Modul}$, $q\in \FreieAM$ und $(x,y)\in \RSa \cup \RSm$.

Also gilt: $g-\bar{g}=p(x-y)s$ und somit auch
\mathe{ \varpi_{\Algebra}(g-\bar{g})=\varpi_{\Algebra}(p(x-y)s) = \varpi_{\Algebra}(p)\varpi_{\Algebra}(x-y)\varpi_{\Algebra}(s).}
Wir betrachten folgende zwei Fälle:
\AufzaehlungP{
\item Falls $(x,y)\in \RSa$, dann gilt $\varpi_{\Algebra}(x-y)=0$. Dann ist schon  $\varpi_{\Algebra}(g-\bar{g})=0$.
\item Falls $(x,y)\in \RSm$, dann gilt $x-y=\lambda^{-1}\tilde{e}$, wobei $\tilde{e}$ wie zu Beginn des Abschnitts aus einem Erzeuger $\bar{e}$ von $\UModul$ konstruiert ist. Es folgt $\varpi_{\Algebra}(x-y)=\bar{e}\in \UModul$, also $\varpi_{\Modul}(\bar{e})=0$.
}
Da $\pi_{\Modul}=  \varpi_{\Modul}   \circ \varpi_{\Algebra}$, folgt in beiden Fällen $\pi_{\Modul}(g)=\pi_{\Modul}(\bar{g})$.

}

}

\Def{Modul Prädikat $\P_{\Modul}$}{
Sei $\RRed$ das Reduktionssystem auf $\FreieAM$, das durch 
\mathe{ \RSam:= \RSa \cup \RSm}
induziert ist. Das Modul Prädikat $\P_{\Modul}$ ist folgendermaßen definiert. Sei $c \in \FreieAM$, dann gilt:
\mathe{ \P_{\Modul}(c) = \wahr \Leftrightarrow c \in \FreieA(0) \cup \FreieA(1)
.}
}

\Satz{Konvergenz}{
\label{noetherschModul}
Sei die Ordnung auf den Monomen $\Monome_{\Algebra,\Modul}$ so gewählt, dass sie noethersch und multiplikativ ist und sei $\RSam$ ein $\P_{\Modul}$-vollständiges Reduktionssystem, dann gilt:

\AufzaehlungZ{
\item Die Reduktionssysteme $\Red_{\P_{\Modul}}$ und $\RRed_{\P_{\Modul}}$ auf $\FreieA_{\P_{\Modul}}$ bzw. $\RigF_{\P_{\Modul}}$ sind konvergent.

\item Zwei Elemente aus $\FreieA(1) \subset \FreieA_{\P_{\Modul}}$ haben genau dann dieselbe Normalform, wenn sie durch die kanonische Projektion $\pi_{\Modul}: \FreieA(1) \nach \Modul$ auf dasselbe Element abgebildet werden.
}
\Beweis{Satz \ref{noetherschModul}
}{
Zu $1.$: Da die Ordnung multiplikativ ist, folgt aus Satz \ref{noetherschRig}, dass $\RRed_{\P_{\Modul}} \subset \RRed$ auch noethersch ist.
Da $\RSam$ $\P_{\Modul}$-vollständig ist, folgt aus Satz \ref{konfluentModul}, dass auch $\RRed_{\P_{\Modul}}$ lokal konfluent ist. Mit Satz \ref{Newman} ist $\RRed_{\P_{\Modul}}$ total konfluent und nach Satz \ref{konvergentModul}
ist auch $\Red_{\P_{\Modul}}$ total konfluent.

Zu $2.$: 
\\
"`$\folgt$"':\\
Seien $f, \bar{f} \in \FreieA(1)$ mit derselben Normalform, dann gilt auch $f \RedAq \bar{f}$ und mit Satz \ref{SatzRSChurchRosserAq} gilt \mathe{\pi(f)=\pi(\bar{f}).}

"`$\Leftarrow$"':\\
Sei $\pi(f)=\pi(\bar{f})$, dann gilt nach Satz \ref{SatzRSChurchRosserAq} auch 
$f \RedAq \bar{f}$. Da $\Red$ nach $1.$ konvergent ist, stimmen nach Satz \ref{Newman} die Normalformen von $f$ und $\bar{f}$ überein.

}
}

\Def{Gröbnerbasis}{
Sei $\RSam$ wie oben konstruiert. Falls $\RSam$ ein  $\P_{\Modul}$-vollständiges Reduktionssystem ist, dann heißt $\{\tilde{\bs{e}}_1, \dots, \tilde{\bs{e}}_{\AnzM}\}$  Gröbnerbasis von dem Untermodul $\UModul \subset \left(\AoA\right)^{\AnzM}$.
}

Um zu zeigen, dass $\Red$ konvergent ist, brauchen wir also nur nachzuweisen, dass die minimalen Überschneidungen in $\RSa \cup \RSm$ mit höchstens einem Modulerzeuger behebbar sind.

%% file: Haupt/Reduktionssysteme/kern-Neu.tex
\subsection{Kern eines $\AoA$-Homomorphismus}
\label{sec:KernEinesAoAHomomorphismus}

Sei $\Algebra=\FreieA \Bruchstrich I$ eine Algebra und $\RSa \subset \Monome \times \FreieA$ ein  Reduktionssystem für $I$.
Wir wollen nun eine Methode angeben, wie wir den Kern einer Abbildung zwischen freien $\AoA$-Moduln berechnen können.

Sei $\phi: \left(\AoA \right)^{n_1} \nach \left(\AoA \right)^{n_2}$ ein $\AoA$-Homomorphismus.
Den Graphen der Abbildung bezeichnen wir mit 
\mathe{\Gamma(\phi):=\Menge{(\phi(g),g) }{g \in \left( \AoA \right)^{n_1}}.}
Sei $n_Q:=n_1+n_2$. Mit $\left(\AoA \right)^{n_1} \oplus \left(\AoA \right)^{n_2} \iso \left(\AoA \right)^{n_Q}$ ist $\Gamma(\phi)$ Untermodul eines freien Moduls.
Wir setzen:
\mathe{Q:= \left(\AoA \right)^{n_Q} \Bruchstrich \Gamma(\phi)}
und erhalten ähnlich wie im Abschnitt \ref{sec:AnwendungFürBimodule} eine exakte Sequenz:

\mathe{ 0 \nach \bar{Q} \inj \left(\AoA\right)^{n_Q} \surj Q \nach 0,}
wobei $\bar{Q}=\Gamma(\phi)$.
Sei die Ordnung auf $\FreieAlgebra_{\AnzA+n_Q}$ die Silbenordnung zu den Trennungsbuchstaben $\bs{f}_1<\dots \bs{f}_{n_1} < \bs{e}_1 < \dots < \bs{e}_{n_2}$, wobei $\bs{f}_1, \dots, \bs{f}_{n_1}$ die Erzeuger des Urbildmoduls und  $\bs{e}_1 , \dots,  \bs{e}_{n_2}$ die Erzeuger des Bildmoduls sind.

Wie in Abschnitt \ref{sec:AnwendungFürBimodule} konstruieren wir aus einem Erzeugersystem für $\bar{Q}$ ein Reduktionssystem $\red_Q$.

\Bem{Wahl des Erzeugersystems}{
\label{WahlSKern}
Von der Wahl des Erzeugersystems für $\bar{Q}$ hängt es ab, ob das induzierte $\red_Q$ geeignet ist (vgl. \ref{WahlSModul}).
}

\Lemma{Erzeugersystem}{
\label{LemmaErzeugersystemKern}
Seien $\bs{f}_1, \dots, \bs{f}_{n_1}$ die Erzeuger des Urbildmoduls , dann ist 
\mathe{ \Menge{ \phi(\bs{f}_i), \bs{f}_i }{ i = 1 \dots n_1}}
ein Erzeugersystem für $\Gamma(\phi)=\bar{Q}$.

\Beweis{Lemma \ref{LemmaErzeugersystemKern}}{
Sei $g \in \left(\AoA \right)^{n_1}$, dann lässt es sich schreiben als:
\mathe{ g= \sum_{\sigma=1}^{N} b_{\sigma} \bs{f}_{i_{\sigma}} c_{\sigma}, }
wobei $b_{\sigma}, c_{\sigma} \in \AoA$ und $i_{\sigma} \in \{1, \dots, n_1\}$. Da $\phi$ ein $\AoA$-Homomorphismus ist, gilt:
\mathe{ \phi(g)= \sum_{\sigma=1}^{N} b_{\sigma} \phi( \bs{f}_{i_{\sigma}}) c_{\sigma}. }
Also lässt sich jedes $(\phi(g),g)\in \Gamma(\phi)$ durch $(\phi(\bs{f}_i), \bs{f}_i)$ erzeugen.
}
}

Die Reduktionsregeln in $\red_{Q}$ lassen sich wegen Bemerkung \ref{BemRedX} zerlegen in:
\Gleichung{
\red_{f}:=&\Menge{(x,y)\in \red_Q}{ i \in \{1,\dots, n_1\}, p,s \in \Monome \textnormal{ mit } x=p \bs{f}_i s } \textnormal{ und }\\
\red_{e}:=&\Menge{(x,y)\in \red_Q}{ i \in \{1,\dots, n_2\}, p,s \in \Monome \textnormal{ mit } x=p \bs{e}_i s }
.}

\Bem{Ordnung}{
\label{BemOrdnungf}
Nach der gewählten Ordnung ist jeder Erzeuger des Urbild-Moduls kleiner als ein Erzeuger des Bild-Moduls. Da eine Regel $(x,y)\in \red_f$ verkleinert, enthält auch jeder Summand, der in $y$ vorkommt, jeweils genau einen Erzeuger des Urbild-Moduls.
}

\newcommand{\fAq}{\Diagramm{   \ar@{<.>}[r]_>>>>{\MatheTiny{\Red_{\xi}}} & }}
\newcommand{\AQAq}{\Diagramm{   \ar@{<.>}[r]_>>>>{\MatheTiny{\RRed_{\Algebra, Q}}} & }}
\newcommand{\ProjA}{\pi_{1}}

Sei im Folgenden $\RRed_{\Algebra,Q}$ das durch  $\red_{\Algebra,Q}$ auf $\RigF_{\AnzA+ n_Q}$ induzierte Reduktionssystem, wobei $n_Q:=n_1+n_2$.

\Bem{$\Rel$}{
\label{RelAoA}
Eine Reduktionsregel  $(g,\bar{g})\in \RRed_{\Algebra,Q}$ ist von genau einer der folgenden Formen:

\newcommand{\EintragZ}[4]{#1: &#2 &\nach & #3 &\Matrix{#4}\\ \hline}

\mathe{\begin{array}{|lrcll|} \hline 
\EintragZ{\textnormal{Name}}{g}{\bar{g}}{\textnormal{ mit }}
\EintragZ{\RelF}{ \lambda_1 w +' \lambda_2 w +' q }{ (\lambda_1 + \lambda_2) w +' q }{~\\~\\~}
\EintragZ{\Rel_{\Algebra}}{ \lambda w x \bar{w} +' q}{\lambda w (y) \bar{w} +'q}{~\\(x,y)\in \red_{\Algebra}\\~}
\EintragZ{\Rel_f}{ \lambda p x s +' q}{\lambda p (y) s +'q}{~\\ (x,y)\in \red_f\\~}
\EintragZ{\Rel_e}{ \lambda p x s +' q}{\lambda p (y) s +'q}{~\\ (x,y)\in \red_e\\~}
\end{array}
}
 Dabei sei $\lambda, \bar{\lambda} \in \K^*$, sowie $  p,s \in \Monome$ und $w, \bar{w} \in \Monome_{\Algebra,Q}$ und $q \in \RigF_{\AnzA+n_Q}$.
}

Ähnlich wie $\varpi_{\Algebra}$ im Abschnitt \ref{sec:AnwendungFürBimodule} nutzen wir im Folgenden oft die beiden  Projektionen
\Gleichung{
\pi_0:& \RigF(0) \nach \Algebra,\\
\ProjA:& \RigF(1) \nach \left( \AoA \right)^{n_Q},
}
wobei $\RigF(0)$ bzw. $\RigF(1)$ die Teilmengen aus $\RigF_{\AnzA+n_Q}$ bezeichnen, deren Elemente nur Summanden, mit keinem bzw. genau einem Erzeuger enthalten. 

\newcommand{\Urbild}{1_f}
%\Def{$\RigF(\Urbild)$}{
Die Menge der Elemente aus $\RigF(1)$, deren Summanden genau einen Erzeuger aus dem Urbildmodul enthalten, bezeichnen wir mit $\RigF(\Urbild)$.
%}

\newcommand{\Killf}{\xi}%_{\Algebra,\textnormal{Bild}}}
Sei $\Killf: \RigF_{\AnzA+n_Q} \nach  \RigF_{\AnzA+n_Q}$ der Homomorphismus, der dadurch definiert ist, dass er die Erzeuger des Urbildmoduls $\bs{f}_i$ auf Null abbildet und für alle anderen Erzeuger die Identität ist. 
\Bem{$\Kern{\Killf}$}{
\label{BemKillf}
Jeder Summand eines Elements aus $\RigF(\Urbild)$ enthält einen Erzeuger des Urbildmoduls, also gilt $\RigF(\Urbild) \subset \Kern{\Killf}$.
}

Wir zerlegen  $\RRed_{\Algebra,Q}$ in die Teilmengen:
\Gleichung{ \Red_{\xi}:&= \Menge{(g,\bar{g})\in  \RRed_{\Algebra,Q} }{\Killf(g)=\Killf(\bar{g}) \textnormal{ und } x \in \RigF(1)},\\
\Red_{\Xi}:&=\RRed_{\Algebra,Q} - \Red_{\xi}.
}

Wir wollen uns zunächst überlegen, wie die Regeln in $\Red_{\xi}$ aussehen.
\Bem{Regeln in $\Red_{\xi}$}{
Sei $(g,\bar{g})\in \Red_{\xi}$, dann gilt mit den Bezeichnungen aus Bemerkung \ref{RelAoA}:
\AufzaehlungP{
\item Falls $(g,\bar{g})\in \RelF$, dann ist $w \in \RigF({\Urbild})$.
\item Falls $(g,\bar{g})\in \Rel_{\Algebra}$, dann ist $wx\bar{w} \in \RigF(\Urbild)$. 
\item Die Regel $(g,\bar{g})$ liegt nicht in $\Rel_e$. 
\item Falls $(g,\bar{g})\in \Rel_f$, dann ist $pxs, pys \in \RigF{\Urbild}$.
}
Es wird also niemals eine Regel aus $\red_e$ genutzt, aber es kann jeweils vorkommen, dass $q$ einen Erzeuger $\bs{e}_i$ des Bildmoduls hat.
}

%\Frage{also $x \in \RigF(1_e}\cup \RigF(1_{\Algebra})$}
%\Lemma{Eigenschaften von $\Red_{\xi}$}{
%\label{LemmaEigenschafteRRedf}
%
%
%
%Sei $\Rel_{f}(\Urbild)$ bzw. $\Rel_{\Algebra}(\Urbild)$ die Einschränkungen von $\Rel_{f}$ bzw. $\Rel_{\Algebra}$ auf $\RigF(\Urbild))$, dann gilt $\Red_{\xi}=\Rel_{f}(\Urbild) \cup \Rel_{\Algebra}(\Urbild)$.
%
%
%\Frage{Formulierung ist noch falsch, vielleicht $\Red_{\xi} \subset \RRed_{\Algebra,Q} ohne \Rel_e$ }
%
%\Beweis{Lemma \ref{LemmaEigenschafteRRedf}}{
%\\
%$"'\subset$"':\\
%Sei $(g,\bar{g})\in \Red_{\xi}$
%\AufzaehlungP{
%\item Falls $(g,\bar{g})\in \RelF$, also von der Form $(\lambda_1 w +' \lambda_2 w +' q ~,~ (\lambda_1 + \lambda_2) w +' q)$ ist und es gilt  $\Killf(\lambda_1 w +' \lambda_2 w +' q)= \Killf((\lambda_1 + \lambda_2) w +' q)$. Da $\Killf$ ein Homomorphismus ist folgt $\Killf(w)=0$ und mit Bemerkung \ref{BemKillf} auch $w\in \RigF(\Urbild)$.
%%
%\item Falls $(g,\bar{g})\in \Rel_{\Algebra}$, also von der Form $(\lambda w x \bar{w} +' q~,~\lambda w (y) \bar{w} +'q)$ ist und es gilt $\Killf(\lambda w x \bar{w} +' q)=\Killf(\lambda w (y) \bar{w} +'q)$. Da $\Killf$ ein Homomorphismus ist folgt $\Killf(\lambda w\Killf( x)\Killf( \bar{w})$
%}
%
%}
%\Frage{Versteh ich nicht}
%}

\Def{induzierter Modul $L$}{
Sei $\varpi_{\Algebra}:\RigF(1) \nach \left(\AoA\right)^{n_Q}$ wie im Abschnitt \ref{sec:AnwendungFürBimodule} konstruiert. Sei $L$ der Untermodul des Urbildmoduls, der folgendermaßen definiert ist:
\mathe{ L:= \left< \ProjA(x) - \ProjA(y) ~|~ (x,y)\in \red_f \right>.}
Wir nennen $L$ auch den durch $\red_f$ induzierten Modul.
}

\Lemma{Modul $L$}{
\label{LemmaL}
Seien $g_1,g_n \in \RigF(1)$ mit $g \fAq q_n$, dann gilt:
\mathe{ \ProjA(g)-\ProjA(g_n) \in L
.}

\Beweis{Lemma \ref{LemmaL}}{
Sei $g_1 \fAq q_n$, dann gibt es eine endliche Folge $(g_i,g_{i+1})$, wobei  $(g_i,g_{i+1})\in \Red_{\xi}$ oder  $(g_{i+1},g_{i}) \in \Red_{\xi}$. Für die Behauptung  reicht es also zu zeigen: Wenn $(g,\bar{g}) \in \Red_{\xi}$, dann auch $ \ProjA(g)-\ProjA(\bar{g}) \in L$. \\
Sei $(g,\bar{g})\in \Red_{\xi}$. Mit den Bezeichnern aus Bemerkung \ref{RelAoA} unterscheiden wir folgende Fälle:
\AufzaehlungP{
\item Falls $(g,\bar{g})\in \RelF$, dann gilt:
\Gleichung{
\ProjA(g)=&\ProjA(\lambda_1 w +' \lambda_2 w +' q ) = \lambda_1 \ProjA(w) + \lambda_2 \ProjA(w) + \ProjA(q) \\
=&   \ProjA((\lambda_1 + \lambda_2) w) + \ProjA( q )= \ProjA((\lambda_1 + \lambda_2) w +' q )\\
=& \ProjA(\bar{g}).}
Also gilt  $\ProjA(g)-\ProjA(\bar{g})=0 \in L$.
\item Falls $(g,\bar{g})\in \Rel_{\Algebra}$, dann liegt entweder $w$ oder $\bar{w}$ in $\RigF(1)$. Sei $\bar{w} \in \RigF(1)$. Es ist $\pi_0(x)=\pi_0(y)$, da $(x,y)\in \RSa$, und daher gilt:
\Gleichung{
\ProjA(g)=& \ProjA(\lambda w x \bar{w} +' q) = \lambda \pi_0(w) \pi_0(x) \ProjA(\bar{w}) + \ProjA(q) \\
=& \lambda \pi_0(w) \pi_0(y) \ProjA(\bar{w}) + \ProjA(q) =\ProjA(\lambda w y \bar{w} +' q) \\
=&\ProjA(\bar{g})
.}
Also gilt  $\ProjA(g)-\ProjA(\bar{g})=0 \in L$. Den Fall ${w} \in \RigF(1)$ behandeln wir analog.
\item Falls $(g,\bar{g})\in \Rel_{f}$, dann gilt:
\Gleichung{
\ProjA(g)=& \ProjA(\lambda pxs+'q) = \lambda \pi_0(p) \ProjA(x)\pi_0( s)+  \ProjA(q),\\
\ProjA(\bar{g})=& \ProjA(\lambda pys+'q) = \lambda \pi_0(p) \ProjA(y)\pi_0( s)+  \ProjA(q).
}
Also gilt wegen  $\ProjA(x)-\ProjA(y) \in L$ auch $\ProjA(g)-\ProjA(\bar{g})\in L$.
}%Aufzaehlung
}

}

Sei im Folgenden $\Aq{\RigF}:= \RigF_{\AnzA+ n_Q} \Bruchstrich \Diagramm{   \ar@{<.>}[r]_>>>>{\MatheTiny{\Red_{\xi}}} & }$ und sei
\mathe{ \Aq{\RRed_{\Algebra,Q}}:=\Menge{\left(\Aq{g},\Aq{\bar{g}} \right) \in \Aq{\RigF} \times \Aq{\RigF}}{(g,\bar{g}) \in \RRed_{\Algebra, Q} \textnormal{ und } \Aq{g} \neq \Aq{\bar{g}}}
.}

\Lemma{Äquivalenzklassen sind verträglich mit $\RigF(\Urbild)$}{
\label{LemmaII}
Sei $g \in \RigF(\Urbild)$, dann gilt das auch für jeden Repräsentanten in $\Aq{g}$.

\Beweis{Lemma \ref{LemmaII}}{
Sei $g \in \RigF(\Urbild)$ und sei $g_n \in \Aq{g}$, dann ist $g=g_1 \fAq q_n$, also gibt es eine endliche Folge $(g_i,g_{i+1})$, wobei  $(g_i,g_{i+1})\in \Red_{\xi}$ oder  $(g_{i+1},g_{i}) \in \Red_{\xi}$. Für die Behauptung  reicht es also zu zeigen: Wenn $(g,\bar{g}) \in \Red_{\xi}$ und $g$ oder $\bar{g} \in \FreieA(\Urbild)$, dann liegen beide in   $\RigF(\Urbild)$.

Sei $(g,\bar{g})\in \Red_{\xi}$. Mit den Bezeichnern aus Bemerkung \ref{RelAoA} unterscheiden wir folgende Fälle:
\AufzaehlungP{
\item Falls $(g,\bar{g})\in \RelF$, dann führt die Regel nur zwei Terme mit gleichen Monomen zusammen; es kommen in $\bar{g}$ die gleichen Monome wie in $g$ vor. \\
Also gilt $g\in \RigF(\Urbild) \Leftrightarrow \bar{g}\in \RigF(\Urbild)$.
\item Falls $(g,\bar{g})\in \Rel_{\Algebra}$, dann ist $x,y \in \RigF(0)$. Da entweder gilt $g$ oder $\bar{g} \in \RigF(\Urbild)$, muss  $q \in \RigF(\Urbild)$ und entweder $w$ oder $\bar{w}\in \RigF(\Urbild)$ sein.\\
Also gilt $g\in \RigF(\Urbild) \Leftrightarrow \bar{g}\in \RigF(\Urbild)$.
\item Falls $(g,\bar{g})\in \Rel_{f}$, dann betrachten wir die folgenden zwei Fälle:
\AufzaehlungP{
 \item	 
 				Sei $g \in  \RigF(\Urbild)$, dann sind auch $pxs, q \in \RigF(\Urbild)$. Mit Bemerkung \ref{BemOrdnungf} 
 				ist  $y\in \RigF(\Urbild)$ und da $p,s \in \Monome$, gilt auch $pys \in \RigF(\Urbild)$.
 \item 
 				Sei $\bar{g} \in \RigF(\Urbild)$, dann ist auch $ q \in \RigF(\Urbild)$. Da $(x,y)\in \red_f$, 
 				ist  $x\in \RigF(\Urbild)$ und da $p,s \in \Monome$, gilt auch $pxs \in \RigF(\Urbild)$.
 }
}%Aufzaehlung
 
 }%Beweis
 }%Lemma
 
%\Lemma{3allerneus}{
%Sei $(g,\bar{g})\in \RRed_{\Algebra,Q}$ mit $g \in \RigF(\Urbild)$, dann gilt:
%\mathe{(g,\bar{g})\in \Red_{\xi}.}
%
%\Beweis{}{
%Wir zeigen es durch einen Widerspruchsbeweis. Sei $(g,\bar{g})\in \RRed_{e}$ mit $g \in \RigF(\Urbild)$, dann lässt sie sich den Bezeichner aus Bemerkung \ref{RelAoA} schreiben als:
%\mathe{ (pxs+'q ~,~ pys +'q) \textnormal{, mit $(x,y)\in \red_e$.}}
%Da jedoch $g \in \RigF(\Urbild)$ kann insbesondere in $x$ kein Erzeuger des Bildmoduls vorkommen. $\Widerspruch$.
%}
%}
 
\Lemma{Reduktionsregeln sind verträglich mit $\RigF(\Urbild)$}{
\label{LemmaIII}
\label{LemmaIV}
Sei $(g,\bar{g})\in \RRed_{\Algebra,Q}$ mit $g \in \RigF(\Urbild)$, dann gilt:
\AufzaehlungZ{
\item $(g,\bar{g})\in \Red_{\xi}.$
\item $\bar{g} \in \RigF(\Urbild).$
}

\Beweis{Lemma \ref{LemmaIII}}{
Zu 1.:\\
Sei $(g,\bar{g})\in \RRed_{\Algebra,Q}$ und $g \in \RigF(\Urbild)$. Mit den Bezeichnern aus Bemerkung \ref{RelAoA} unterscheiden wir folgende Fälle:
\AufzaehlungP{
\item Falls $(g,\bar{g})\in \RelF$, dann führt die Regel nur zwei Terme mit gleichen Monomen zusammen; es kommen in $\bar{g}$ die gleichen Monome wie in $g$ vor. \\
Also folgt aus $g\in \RigF(\Urbild)$ auch $\bar{g}\in \RigF(\Urbild)$ und deshalb gilt $\xi(g)=0=\xi(\bar{g})$.
\item Falls $(g,\bar{g})\in \Rel_{\Algebra}$, dann ist $x,y \in \RigF(0)$. Da entweder gilt $g$ oder $\bar{g} \in \RigF(\Urbild)$, muss entweder $w$ oder $\bar{w}\in \RigF(\Urbild)$ sein.\\
Also folgt $g\in \RigF(\Urbild)$ und auch $\bar{g}\in \RigF(\Urbild)$ und deshalb gilt $\xi(g)=0=\xi(\bar{g})$.
\item Falls $(g,\bar{g})\in \Rel_{f}$, dann folgt aus Bemerkung \ref{BemOrdnungf},  dass auch $\bar{g}\in \RigF(\Urbild)$
}%Aufzaehlung

Zu 2.:\\
Mit 1. folgt aus $(g,\bar{g})\in \RRed_{\Algebra, Q}$, so dass gilt: $\Aq{g}=\Aq{\bar{g}}$. Mit Lemma \ref{LemmaII} folgt dann die Behauptung.

}
}

\Satz{Kern}{
\label{SatzKern}
Sei $\Aq{\RRed_{\Algebra,Q}}$ konvergent, dann gilt:
\mathe{ \left< \ProjA(x-y) ~|~ (x,y) \in \red_f \right> = \Kern{\phi}.
}

\Beweis{Satz \ref{SatzKern}}{~\\
"`$\subseteq$"':\\
Sei $(x,y)\in \red_f $, dann ist wegen Bemerkung \ref{BemOrdnungf} $\left(0,\ProjA(x-y)\right)\in \bar{Q}= \Gamma(\phi)$. Also ist $\phi(\ProjA(x-y))=0$.
\\
"`$\supseteq$"':\\
Sei $h \in \Kern{\phi}$, dann ist $(0,h) \in \Gamma(\phi)$. Wegen Satz \ref{SatzRSChurchRosserAq} gilt $(0,h) \Diagramm{   \ar@{<.>}[r]_>>>>{\MatheTiny{\Red_{\Algebra,Q}}} & } (0,0)$. Da wir jede Regel aus $\Red$ auch als einen Weg in $\RRed$ schreiben können (vgl. \ref{RsubsetTR}), gilt:
\mathe{(0,h) \Diagramm{   \ar@{<.>}[r]_>>>>{\MatheTiny{\Red_{\Algebra,Q}}} & } (0,0) \textnormal{ und }}
\mathe{\Aq{(0,h)} \Diagramm{   \ar@{<.>}[r]_>>>>{\MatheTiny{\Aq{\RRed_{\Algebra,Q}}}} & } \Aq{(0,0)}.
}
Nun ist $\Aq{\RRed_{\Algebra,Q}}$ konvergent, es gibt also Wege:
\mathe{\Diagramm{ 
\Aq{(0,h)} \ar@{->}[r] & \Aq{(v_1)} \ar@{->}[r] & \dots \Aq{(v_{n})} \ar@{=}[d]\\
\Aq{(0,0)} \ar@{->}[r] & \Aq{(u_1)} \ar@{->}[r] & \dots \Aq{(u_m)} 
}}
und wir erhalten folgendes Diagramm:
\mathe{\Diagramm{ 
(0,h)=:v'_0 
\ar@{<.>}[r]_>>>>{\MatheTiny{\Red_{\xi}}} & v''_0 \ar@{->}[r]_>{\MatheTiny{\RRed_{\Algebra,Q}}} & v'_1 \ar@{<.>}[r]_>>>>{\MatheTiny{\Red_{\xi}}} & v''_1 \ar@{->}[r]_>{\MatheTiny{\RRed_{\Algebra,Q}}} & v'_2 
\ar@{<.>}[r]_>>>>{\MatheTiny{\Red_{\xi}}} & v''_2 \ar@{->}[r]_>{\MatheTiny{\RRed_{\Algebra,Q}}} & \dots 
																		  	\ar@{->}[r]_>{\MatheTiny{\RRed_{\Algebra,Q}}} & v'_n \ar@{<.>}[d]_>>>>{\MatheTiny{\Red_{\xi}}} \\
(0,0)=:u'_0  
\ar@{<.>}[r]_>>>>{\MatheTiny{\Red_{\xi}}} & u''_0 \ar@{->}[r]_>{\MatheTiny{\RRed_{\Algebra,Q}}} & u'_1 \ar@{<.>}[r]_>>>>{\MatheTiny{\Red_{\xi}}} & u''_1 \ar@{->}[r]_>{\MatheTiny{\RRed_{\Algebra,Q}}} & u'_2 
\ar@{<.>}[r]_>>>>{\MatheTiny{\Red_{\xi}}} & u''_2 \ar@{->}[r]_>{\MatheTiny{\RRed_{\Algebra,Q}}} & \dots 
 																				  	\ar@{->}[r]_>{\MatheTiny{\RRed_{\Algebra,Q}}} & u'_n  \\
},}
wobei $v''_i,v'_i \in \Aq{v_i}$ und $u''_i,u'_i \in \Aq{u_i}$.

Nach Lemma \ref{LemmaIV} und Lemma \ref{LemmaII} liegen alle $v''_i,v'_i,u''_i,u'_i\in \FreieA(\Urbild)$.
Nach Lemma \ref{LemmaIII} gilt $\ProjA(v'_i )= \ProjA(v'_{i+1})$ und  $\ProjA(u'_i )= \ProjA(u'_{i+1})$.
Zusammen mit Lemma \ref{LemmaL} folgt:
\mathe{ \ProjA(v'_i )- \ProjA(v''_{i}), \ProjA(u'_i )- \ProjA(u''_{i}) \in L.
}

}%Beweis
}%Satz

Um ein Erzeugersystem für den Kern zu erhalten, reicht es also ein konvergentes $\Aq{\Red_{\Algebra,Q}}$ zu kennen. Wenn $\Red_{\Algebra,Q}$ konvergent ist, dann ist es auch $\Aq{\Red_{\Algebra,Q}}$. Um die Konvergenz von ${\Red_{\Algebra,Q}}$ nachzuweisen kann man zeigen, dass alle $\P$-minimalen Überschneidungen behebbar sind (vgl. \ref{noetherschModul}). Wir wollen nun zeigen, dass wir nicht alle diese Überschneidungen überprüfen müssen.

Zunächst werden wir zeigen, dass die Regeln aus $\Red_{\xi}$ mit denen aus $\Red_{\Xi}$ kommutieren. Dazu betrachten wir  folgendes Lemma.

\Lemma{$\Red_{\Xi}$ kommutiert mit $\Red_{\xi}$}{
\label{LemmaVI}
\AufzaehlungP{
\item
Sei $(g_1,g_2)\in \Red_{\xi}$ und $(g_2,h_2)\in \Red_{\Xi}$, dann gibt es ein $h_1 \in \RigF(1)$, so dass  $(g_1,h_1)\in \Red_{\Xi}$ und $(h_1,h_2)\in \Red_{\xi}$, also 
\mathe{\Diagramm{
 g_1 \ar@{->}[r]_>>{\MatheTiny{\Red_{\xi}}} \ar@{->}@/_1pc/[dr]_>>{\MatheTiny{\Red_{\Xi}}}
 & g_2 \ar@{->}[r]_>>{\MatheTiny{\Red_{\Xi}}} & h_2 \\
 																				 & \bs{h_1} \ar@{->}@/_1pc/[ur]_>>>{\MatheTiny{\Red_{\xi}}}
}.}
\item 
Sei $(g_1,g_2)\in \Red_{\xi}$ und $(g_1,h_1)\in \Red_{\Xi}$, dann gibt es ein $h_2 \in \RigF(1)$, so dass $(h_1,h_2)\in \Red_{\xi}$ und $(g_2,h_2)\in \Red_{\Xi}$, also
\mathe{\Diagramm{
h_1 \ar@/_1pc/[dr]_>>{\MatheTiny{\Red_{\xi}}} &g_1 \ar@{->}[r]_>>{\MatheTiny{\Red_{\xi}}} \ar@{->}[l]^>>{\MatheTiny{\Red_{\Xi}}}& g_2 \ar@/^1pc/[dl]^>>{\MatheTiny{\Red_{\Xi}}}& \\
																	 &\bs{h_2} 
}.}
\item 
Sei $(g_1,h_1)\in \Red_{\Xi}$ und sei $\Diagramm{ g_2 \ar@{<.>}[r]_>>>>{\MatheTiny{\Red_{\xi}}} & g_1}$, dann gibt es ein $h_2 \in \RigF(1)$ mit $(g_2,h_2)\in \Red_{\Xi}$ und $\Diagramm{ h_2 \ar@{<.>}[r]_>>>>{\MatheTiny{\Red_{\xi}}} & h_1}$, also
\mathe{\Diagramm{
h_1 \ar@/_1pc/@{<.>}[dr]_>>>>{\MatheTiny{\Red_{\xi}}} &g_1 \ar@{<.>}[r]_>>>>{\MatheTiny{\Red_{\xi}}} \ar@{->}[l]^>>{\MatheTiny{\Red_{\Xi}}}& g_2 \ar@/^1pc/[dl]^>>{\MatheTiny{\Red_{\Xi}}}&\\
																	 &\bs{h_2} 
}.}
}
\Beweis{Lemma \ref{LemmaVI}}{
Ein Element $g\in \RigF(1)$ können wir zerlegen in einen Anteil in $g_{\xi}\in \RigF(\Urbild)$ und einen Anteil $g_{\Xi}$, der nur Monome mit jeweils einem Erzeuger $\bs{e}_i$ des Bildmoduls enthält:
\mathe{ g = (g_{\Xi},g_{\xi})}
Sei $((g_{\Xi},g_{\xi}) \nach \bar{g}) \in \Red_{\xi}$ eine Regel, die mit $g$ beginnt; sie verändert nur die zweite Stelle, also $\bar{g}= (g_{\Xi},\bar{g}_{\xi})$.\\
Da  $\RRed_{\Algebra,Q}$ durch $\RSam$ induziert ist, gibt es für jedes $v \in \RigF(1)$ die Regel $((v,g_{\xi}) \nach 
(v,\bar{g}_{\xi},) ) \in \Red_{\xi}$. \\ 
~\\
Sei $((g_{\Xi},g_{\xi}) \nach h) \in \Red_{\Xi}$ eine Regel, die mit $g$ beginnt; sie verändert $g_{\xi}$ nicht, also ist $h=(h_{\Xi}, h_{\xi} +' g_{\xi})$. \\
Da  $\RRed_{\Algebra,Q}$ durch $\RSam$ induziert ist, gibt es für jedes $v \in \RigF(1)$ die Regel $((g_{\Xi},v) \nach (h_{\Xi}, h_{\xi} +' v)) \in \Red_{\Xi}$.

Also "`kommutiert"' jede Regel aus $\Red_{\Xi}$ mit jeder aus $\Red_{\xi}$, weshalb die beiden ersten Punkte gelten. Den letzten Punkt zeigt man durch mehrfaches Anwenden der ersten beiden.
}
}

\Def{schwach vollständig}{
Sei $\RS_{\Algebra,Q}:=\RSa \cup \RS_{Q}$ und sei $\P_{Q}$ das Modulprädikat. Falls jede $\P_{Q}$-minimale Überschneidung $(w,\pi, \bar{\pi})$ mit $\pi, \bar{\pi} \in \red_{\Algebra}\cup \red_{e}$ mit Regeln aus $\RRed_{\Algebra,Q}$ behebbar ist, dann nennen wir $\RS_{\Algebra,Q}$ schwach vollständig.
}

\Satz{Hauptsatz}{
\label{HauptsatzKern}
Sei $\RS_{\Algebra,Q}$ schwach vollständig, dann ist 
\mathe{ \Aq{\RRed_{\Algebra,Q}} \subset \RigF(1) \times \RigF(1) \textnormal{ lokal konfluent.}
}
\Beweis{Satz \ref{HauptsatzKern}}{
Sei $\RS_{\Algebra,Q}$ schwach vollständig und seien $(\Aq{g},\Aq{h}),(\Aq{g},\Aq{\bar{h}})$ zwei Reduzierungsregeln, die mit $\Aq{g}$ beginnen, dann existieren $g_1,g_2 \in \Aq{g}$ und $h_1\in \Aq{h}$ sowie  $\bar{h}_2\in \Aq{\bar{h}}$ mit:
\mathe{
\Diagramm{&&g_1 \ar@{<.>}[r]_>>>>{\MatheTiny{\Red_{\xi}}} \ar@{->}[ld]^>{\MatheTiny{\RRed_{\Algebra,Q}}}& g_2 \ar@{->}[rd]_>{\MatheTiny{\RRed_{\Algebra,Q}}}\\
h \ar@{<.>}[r]_>>>>{\MatheTiny{\Red_{\xi}}} &h_1&&& \bar{h}_2 \ar@{<.>}[r]_>>>>{\MatheTiny{\Red_{\xi}}} &\bar{h}
.}}%mathe
Nun gilt $\Aq{g}\in  \RigF(1) $ und damit auch $g_1,g_2,h_1, h, \bar{h}_1, \bar{h}\in \RigF(1)$.
Da $(\Aq{g_i},\Aq{h_i})\in \Aq{\RRed_{\Algebra,Q}}$, ist $\Aq{g_i}\neq\Aq{h_i}$, also liegen $(g_i,h_i)$ in $\Red_{\Xi}$ . Mit Lemma \ref{LemmaVI} vereinfacht sich das Diagramm zu:
\mathe{
\Diagramm{&&&  \ar@{->}[ld]^>{\MatheTiny{\Red_{\Xi}}} g_2 \ar@{->}[rd]_>{\MatheTiny{\Red_{\Xi}}}\\
h \ar@{<.>}[r]_>>>>{\MatheTiny{\Red_{\xi}}} &{h_1}&\bs{h_2}  \ar@{<.>}[l]^>>>>{\MatheTiny{\Red_{\xi}}}&& \bar{h}_2 \ar@{<.>}[r]_>>>>{\MatheTiny{\Red_{\xi}}} &\bar{h}.
}}%mathe
Da $(g_2,h_2),(g_2,\bar{h}_2)\in \Red_{\Xi}$, sind sie von einer Regel aus $\red_{\Algebra} \cup \red_e$ induziert, nach Vor\-aussetzung und ähnlich wie Satz \ref{konfluentModul} können wir zeigen, dass es Wege $h_2 \Weg z$ und $\bar{h} \Weg z$ gibt.  Diese Wege sind auf den Äquivalenzklassen von folgender Form:

\mathe{
\Diagramm{&  \ar@{->}[ld]^>{\MatheTiny{}} \Aq{g_2} \ar@{->}[rd]\\
\Aq{h}=\Aq{h_1}=\Aq{h_2} \ar@{.>}[rdd]_>{\MatheTiny{}} && ~~~~~~~\Aq{\bar{h}_2}=\Aq{\bar{h}} ~~~~~~~\ar@{.>}[ldd]_>{\MatheTiny{}}\\
~\\
&\Aq{z}.
}}%mathe

}%Beweis
}

Wir erhalten also ein Erzeugersystem für den Kern, indem wir nachweisen, dass alle minimalen Überschneidungen in $\RSam$, die höchstens einen Erzeuger des Bildmoduls enthalten, bezüglich $\Red_{\Algebra,Q}$ behebbar sind.

%% file: Haupt/RechnungGBA.tex
In diesem Kapitel werden wir für die orthogonale freie Quantengruppe ein vollständiges Reduktionssystem angeben.  

\section{Definition $\fQG(n)$}
\label{sec:GröbnerbasisFürDieFreieQuantenGruppeAN}

\begin{definition}[$\fQG(n)$]
Sei $\fQG(n)$ die $\K$-Algebra, die von $a_{i,j}$, wobei $i,j \in \{1, \dots , n\}$, und den Relationen 

\mathe{ \sum_{p=1}^n a_{i,p}a_{j,p} = \delta_{i,j}}
\mathe{ \sum_{p=1}^n a_{p,i}a_{p,j} = \delta_{i,j}}
erzeugt wird. Es gilt also $\fQG(n) = \FreieA_{n^2} \Bruchstrich I_{\fQG(n)}$, wobei \mathe{I_{\fQG(n)}:=\left< \sum_{p=1}^n a_{i,p}a_{j,p} - \delta_{i,j}, \sum_{p=1}^n a_{p,i}a_{p,j} = \delta_{i,j} ~|~ i,j = 1 \dots n \right>.}\\
Wir nennen $\fQG(n)$ die orthogonale freie Quantengruppe mit $n^2$-Erzeuger.
\end{definition}

\Frage{Literatur Thom fragen}

\begin{bemerkung}
Schreibt man die Erzeuger als Matrix:
\mathe{ \A:=\left( a_{i,j} \right)_{i,j},}
dann lassen sich die Relationen schreiben als:
\mathe{\A \A^t=\id,}
\mathe{\A^t \A=\id,}
wobei $\id$ die Identitätsmatrix und $\A^t$ die transponierte Matrix von $\A$ ist.
\end{bemerkung}

%
%\begin{beispiel}[$\fQG(1)$]
%Im Fall $\fQG(1)$ gibt es nur einen Erzeuger $a_{1,1}$ und nur eine Relation $a_{1,1}^2=1$. Also ist $\fQG(1)$ die Projektive Algebra mit nur einem Erzeuger. Sie ist halbeinfach, daher wird der Fall $n=1$ nicht  weiter betrachtet.
%\end{beispiel}
%

\Satz{vollständiges Reduktionssystem $\red_{\fQG(n)}$}{
\label{GBA}

Sei $\red_{\fQG(n)}:=\Menge{\widetilde{Z}_{p,q},\widetilde{S}_{p,q}, \widetilde{Z}_{p,q,r},\widetilde{S}_{p,q,r}}{p,q,r = 1, \dots , n}$, wobei
\mathe{\begin{array}{rl}
	\widetilde{Z}_{p,q}:=&\left(a_{p,1} a_{q,1} ~,~ \Zi{p}{q}\right) %&\textnormal{ mit } &Z_{p,q} :=& 
	,\\
	\widetilde{S}_{p,q}:=&\left(a_{1,p} a_{1,q}~,~ \Si{p}{q}\right) %&\textnormal{ mit } &S_{p,q} :=& \Si{p}{q}
	,\\
	\widetilde{Z}_{p,q,r}:=&\left(a_{p,1} a_{q,2} a_{r,2}~,~ \Zii{p}{q}{r} -\left(  \Zi{p}{q} \right)a_{r,1}\right) % &\textnormal{ mit } &Z_{p,q,r} :=& \Zii{p}{q}{r} - \Zi{p}{q} a_{r,1}
	,\\
	\widetilde{S}_{p,q,r}:=&\left(a_{1,p} a_{2,q} a_{2,r}~,~ \Sii{p}{q}{r} - \left( \Si{p}{q} \right)a_{1,r}\right)% &\textnormal{ mit } &S_{p,q,r} :=& \Sii{p}{q}{r} - \Si{p}{q} a_{1,r}
	,
\end{array}
}
dann ist $\red_{\fQG(n)}$ mit der kanonischen Ordnung ein vollständiges Wortersetzungssystem für $\fQG(n)$.
}%Satz

\Bem{Spezialfälle $n=1$ und $n=2$}{
Falls $n=1$, gibt es nur einen Erzeuger $a_{1,1}$ und nur eine Relation $(a_{1,1}^2~,~ 1)$.\\
Falls $n=2$, fällt in $\widetilde{Z}_{p,q,r}$ und $\widetilde{S}_{p,q,r}$ auf der rechten Seite jeweils die erste Summe weg.
}

\Ver{
Um die Lesbarkeit zu erhöhen, setzen wir:
\mathe{\begin{array}{rl}
Z_{p,q} :=& \Zi{p}{q},\\
S_{p,q} :=& \Si{p}{q},\\
Z_{p,q,r} :=& \Zii{p}{q}{r} ,\\
S_{p,q,r} :=& \Sii{p}{q}{r} .
\end{array}
}%mathe
Die Regeln in $\tilde{r}_{A(n)}$ verkürzen sich dann zu:
\mathe{\begin{array}{rl}
	a_{p,1} a_{q,1} & \Pfeil{\widetilde{Z}_{p,q}} Z_{p,q},\\
	a_{1,p} a_{1,q} & \Pfeil{\widetilde{S}_{p,q}} S_{p,q},\\
	a_{p,1} a_{q,2} a_{r,2} &\Pfeil{\widetilde{Z}_{p,q,r}} Z_{p,q,r} - Z_{p,q} a_{r,1},\\
	a_{1,p} a_{2,q} a_{2,r} &\Pfeil{\widetilde{S}_{p,q,r}} S_{p,q,r} - S_{p,q} a_{1,r}.
\end{array}
}%mathe
 }%Ver

Bevor wir den Satz \ref{GBA} beweisen, bemerken wir noch, dass wegen Satz \ref{vollstaendig} und der Definition \ref{DefOrdnungRigF}, das  durch $\red_{\fQG(n)}$ induzierte Reduktionssystem konvergent ist.

\Beweis{Satz \ref{GBA} }{
Es sind zwei Dinge zu zeigen:
\AufzaehlungZ{
\item Es gilt: $\left< x-y ~|~ (x,y) \in \red_{\fQG(n)} \right> = I_{\fQG(n)}$.
\item Das Ersetzungssystem $\red_{\fQG(n)}$ ist vollständig.
}%Aufzaehlung

Zu 1.: \\
Es gilt: 
Nach Definition ist \mathe{
I_{\fQG(n)}:=\left< - (\widetilde{Z}_{p,q} - a_{p,1}a_{q,1}), -(\widetilde{S}_{p,q} - a_{p,1}a_{q,1}) ~|~ p,q=1,\dots,n \right>
.}
Also gilt: $I_{\fQG(n)} \subset \tilde{r}_{\fQG(n)}$.
\\
Andererseits sehen wir, dass  $Z_{p,q,r},S_{p,q,r}\subset I_{\fQG(n)}$gilt durch folgende Gleichungen:
\Gleichung{ Z_{p,q}a_{r_1}&= a_{p,1}a_{q,1}a_{r_1} =& a_{p,1}Z_{q,r},\\
						S_{p,q}a_{1,r}&= a_{1,p}a_{1,q}a_{1,r} =& a_{1,p}S_{q,r}.}

~\\
%\begin{landscape}
Zu 2.:\\
Um Vollständigkeit zu zeigen müssen alle minimalen Überschneidungen behebbar sein. In der folgenden Tabelle sind alle %\Frage{warum gibt es keine weiteren} 
minimalen Überschneidungen aufgeführt. Ein Eintrag in der Tabelle besteht aus drei Teilen: der Länge der Überschneidung, dem Monom, für das es zwei Reduzierungsregeln ($\textnormal{Zeile} \cdot s$ und $p \cdot \textnormal{Spalte}$) gibt, und einer Seitenangabe, auf der die Überschneidung behoben wird.
\mathe{\begin{array}{|l|l|l|}\hline
 	& \LMatrix{~\\ \widetilde{Z}_{x,y}= a_{x,1}a_{y,1}}& \widetilde{S}_{x,y}=a_{1,x}a_{1,y} %&\widetilde{Z}_{x,y,z}=a_{x,1}a_{y,2}a_{z,2}& \widetilde{S}_{x,y,z}=a_{1,x}a_{2,y}a_{2,z}
\\ \hline %entspricht \cline{1-5}
 \LMatrix{~\\ \widetilde{Z}_{p,q} \\ ~ \\ ~ }
 		&\LMatrix{1: a_{p,1}a_{q,1}a_{y,1} &\Seite{ZiZil} \\ 2: \textnormal{beide gleich} } 
 		&\LMatrix{1: a_{p,1}a_{1,1}a_{1,y} &\Seite{ZiSil} \\ 2: a_{1,1}a_{1,1} &\Seite{ZiSill} } 
%		&\LMatrix{1: a_{p,1}a_{q,1}a_{y,2}a_{z,2} &\Seite{ZiZiil} \\ 2: \textnormal{nicht möglich}} 
%		&\LMatrix{1: a_{p,1}a_{1,1}a_{2,y}a_{2,z} &\Seite{ZiSiil} \\ 2: a_{1,1}a_{2,1}a_{1,z} &\Seite{ZiSiill}} 
  \\ \hline
 \LMatrix{~\\ \widetilde{S}_{p,q} \\ ~ \\ ~ }
  	&\LMatrix{1: a_{1,p}a_{1,1}a_{y,1} &\Seite{SiZil} \\ 2: a_{1,1}a_{1,1} &\Seite{SiZill} } 
		&\LMatrix{1: a_{1,p}a_{1,q}a_{1,y} &\Seite{SiSil} \\ 2: \textnormal{beide gleich} } 
%		&\LMatrix{1: a_{1,p}a_{1,1}a_{y,2}a_{z,2} &\Seite{SiZiil} \\ 2: a_{1,1}a_{1,2}a_{z,2} &\Seite{SiZiill} } 
%		&\LMatrix{1: a_{1,p}a_{1,1}a_{2,y}a_{2,z} &\Seite{SiSiil} \\ 2: \textnormal{nicht möglich} } 
 \\ \hline
 \LMatrix{~\\ \widetilde{Z}_{p,q,r}  \\ ~ \\ ~ }	
 		&\LMatrix{1: \textnormal{nicht möglich}  \\ 2: \textnormal{nicht möglich}} 
 		&\LMatrix{1: a_{p,1}a_{q,2}a_{1,2}a_{1,y} &\Seite{ZiiSil} \\ 2: a_{p,1}a_{1,2}a_{1,2} &\Seite{ZiiSill} } 
%		&\LMatrix{1: \textnormal{nicht möglich} \\ 2: \textnormal{nicht möglich} \\ 3: \textnormal{beide gleich} }
%		&\LMatrix{1: a_{p,1}a_{q,2}a_{1,2}a_{2,y}a_{2,z} &\Seite{ZiiSiil} \\ 2: a_{p,1}a_{1,2}a_{2,2}a_{2,z} &\Seite{ZiiSiill} \\3: a_{1,1}a_{2,2}a_{2,2} \Seite{ZiiSiilll} } 
 \\ \hline
 \LMatrix{~\\ \widetilde{S}_{p,q,r}\\ ~ \\ ~ }
 		&\LMatrix{1: a_{1,p}a_{2,q}a_{2,1}a_{y,1} &\Seite{SiiZil} \\ 2: a_{1,p}a_{2,1}a_{2,1} &\Seite{SiiZill} } 
		&\LMatrix{1: \textnormal{nicht möglich} \\ 2: \textnormal{nicht möglich} } 
%		&\LMatrix{1: a_{1,p}a_{2,q}a_{2,1}a_{y,2}a_{z,2} &\Seite{SiiZiil} \\ 2: a_{1,p}a_{2,1}a_{2,2}a_{z,2} &\Seite{SiiZiill}\\3: a_{1,1}a_{2,2}a_{2,2} &\Seite{SiiZiilll} } 
%		&\LMatrix{1: \textnormal{nicht möglich} \\ 2: \textnormal{nicht möglich} \\3:\textnormal{beide gleich}} 
\\ \hline%
%
%\end{array}}
%\end{landscape}
&&\\
%%%Hier beginnt die zweite Tabelle
 &\widetilde{Z}_{x,y,z}=a_{x,1}a_{y,2}a_{z,2}& \widetilde{S}_{x,y,z}=a_{1,x}a_{2,y}a_{2,z}
 	\\ \hline %entspricht \cline{1-5}
 \LMatrix{~\\ \widetilde{Z}_{p,q} \\ ~ \\ ~ }
 		&\LMatrix{1: a_{p,1}a_{q,1}a_{y,2}a_{z,2} &\Seite{ZiZiil} \\ 2: \textnormal{nicht möglich}} 
		&\LMatrix{1: a_{p,1}a_{1,1}a_{2,y}a_{2,z} &\Seite{ZiSiil} \\ 2: a_{1,1}a_{2,1}a_{2,z} &\Seite{ZiSiill}} 
  \\ \hline
 \LMatrix{~\\ \widetilde{S}_{p,q} \\ ~ \\ ~ }
 		&\LMatrix{1: a_{1,p}a_{1,1}a_{y,2}a_{z,2} &\Seite{SiZiil} \\ 2: a_{1,1}a_{1,2}a_{z,2} &\Seite{SiZiill} } 
		&\LMatrix{1: a_{1,p}a_{1,1}a_{2,y}a_{2,z} &\Seite{SiSiil} \\ 2: \textnormal{nicht möglich} } 
 \\ \hline
 \LMatrix{~\\ \widetilde{Z}_{p,q,r}  \\ ~ \\ ~ }	
 		&\LMatrix{1: \textnormal{nicht möglich} \\ 2: \textnormal{nicht möglich} \\ 3: \textnormal{beide gleich} }
		&\LMatrix{1: a_{p,1}a_{q,2}a_{1,2}a_{2,y}a_{2,z} &\Seite{ZiiSiil} \\ 2: a_{p,1}a_{1,2}a_{2,2}a_{2,z} &\Seite{ZiiSiill} \\3: a_{1,1}a_{2,2}a_{2,2} &\Seite{ZiiSiilll} } 
 \\ \hline
 \LMatrix{~\\ \widetilde{S}_{p,q,r}\\ ~ \\ ~ }
 		&\LMatrix{1: a_{1,p}a_{2,q}a_{2,1}a_{y,2}a_{z,2} &\Seite{SiiZiil} \\ 2: a_{1,p}a_{2,1}a_{2,2}a_{z,2} &\Seite{SiiZiill}\\3: a_{1,1}a_{2,2}a_{2,2} &\Seite{SiiZiilll} } 
		&\LMatrix{1: \textnormal{nicht möglich} \\ 2: \textnormal{nicht möglich} \\3:\textnormal{beide gleich}} 
 \\ \hline
\end{array}}

%Die Lable etsrechen der obigen Tabelle. Sie betehen aus drei  Teilen   Zi steht für Z_{p,q}. Zii für Z_{p,q,r}. Si und Sii entsprechen S_{p,q} und S_{p,q,r}, die l, ll bzw. lll stehen für eine Überschneidung der länge 1, 2 bzw. 3. 

\input{Haupt/RechnungGBA/zusammen}
}%Beweis

%% file: Haupt/RechnungGBA/zusammen.tex
%Nur lokal für die Rechnungen brauch ich eine andere Summenumgebung 
\renewcommand{\Sum}[1]{\sum_{#1}^{n}}

\section{Behebung der minimalen Überschneidungen }
\label{sec:RechnungenGBA}
In den folgenden Rechnungen steht $\delta_{x,y}$ immer für das Kroneckersymbol, \mathe{\delta_{x,y} = \left\{ \Matrix{ 1& \textnormal{ für } x=y\\ 0& \textnormal{ sonst}}. \right. }\\
\Bem{Symmetrie}{
Durch Vertauschen der Indizes aus $Z_{p,q}$ erhält man $S_{p,q}$, und durch Vertauschen der Indizes aus $Z_{p,q,r}$ erhält man $S_{p,q,r}$, weshalb  einige der Rechnungen symmetrisch sind.}

\MitRechnungen{
\input{Haupt/RechnungGBA/Rechenregel}\newpage

\input{Haupt/RechnungGBA/ZiZil}\newpage
\input{Haupt/RechnungGBA/SiSil}\newpage

\input{Haupt/RechnungGBA/ZiSil}\newpage
\input{Haupt/RechnungGBA/SiZil}\newpage %~\newpage

\input{Haupt/RechnungGBA/ZiSill} \newpage
% Identische Rechnung wie SiSill 

\input{Haupt/RechnungGBA/ZiZiil}\newpage
\input{Haupt/RechnungGBA/SiSiil}\newpage

\input{Haupt/RechnungGBA/ZiSiil}\newpage
\input{Haupt/RechnungGBA/SiZiil}\newpage

\input{Haupt/RechnungGBA/ZiSiill} \newpage
\input{Haupt/RechnungGBA/SiZiill} \newpage

\input{Haupt/RechnungGBA/ZiiSil}\newpage
\input{Haupt/RechnungGBA/SiiZil}\newpage

\input{Haupt/RechnungGBA/ZiiSill}\newpage
\input{Haupt/RechnungGBA/SiiZill}\newpage

\input{Haupt/RechnungGBA/ZiiSiil}\newpage
\input{Haupt/RechnungGBA/SiiZiil}\newpage

\input{Haupt/RechnungGBA/ZiiSiill}\newpage
\input{Haupt/RechnungGBA/SiiZiill}\newpage

\input{Haupt/RechnungGBA/ZiiSiilll}

% Identische Rechnung wie SiSill 
}

%% file: Haupt/RechnungGBA/Rechenregel.tex
\subsection{Rechenregeln}
\label{sec:Rechenregeln}

\Satz{Rechenregel}{
\label{Rechenregel}
\label{RechenRegeln}
In der freien Algebra $\F_{n^2}$ gilt:
\AufzaehlungZ{
\item \mathe{ \Sum{j=2} Z_{p,j}a_{j,y} = \Sum{i=2} a_{p,i} S_{i,y} + \left( \delta_{1,y} -\delta_{p,1}\right) a_{p,y}}
\item \mathe{ \Sum{j=3} Z_{p,1,j} a_{j,z} =\Sum{i=3} a_{p,1} S_{i,i,z} - (1-\delta_{1,z}-\delta_{2,z})a_{p,1}a_{1,z} }
\item \mathe{ \Sum{j=2} Z_{p,j,j} = \Sum{i=3} a_{p,1} S_{i,i}  + a_{p,1}}
\item \mathe{ \Sum{j=3}Z_{p,1} a_{j,1}a_{j,z} = \Sum{i=2}a_{p,i} S_{i,1,z} - \delta_{1,z} \Sum{i=2}a_{p,i}a_{1,i} + \delta_{p,1}\Sum{j=3}a_{j,1}a_{j,z} }
\item \mathe{ Z_{p,q,1}\left(\Sum{j=3} a_{j,y}a_{j,z} - \delta_{y,z}\right)= a_{p,1}\Sum{i=3}a_{q,i} S_{i,y,z} + \delta_{q,1} a_{p,1}\left(\Sum{j=3} a_{j,y}a_{j,z} - \delta_{y,z}\right)}
\item\mathe{\Sum{j=2} Z_{p,q,j} a_{j,y} = \Sum{i=3} a_{p,1}a_{q,i} S_{i,y} + \delta_{y,1}a_{p,1}a_{q,1}+ \delta_{y,2}a_{p,1}a_{q,2} -\delta_{q,1}a_{p,1}a_{1,y}
}
}%AufzaehlungZ

\Beweis{Satz \ref{Rechenregel}}{
\MitRechnungen{
~\\ Zu 1.:
\Gleichung{ &\Sum{j=2} Z_{p,j}a_{j,y} 
\\	&=	- \Sum{j=2} \left( \mZi{p}{j} \right) a_{j,y}
\\&= - \Sum{i=2} \Sum{j=2} a_{p,i} a_{j,i} a_{j,y}
 		+ \left( \Sum{i=2} a_{p,i} \delta_{i,y} 
 					- \underbrace{\Sum{i=2} a_{p,i} \delta_{i,y}}_{=a_{p,y}\left(1-\delta_{1,y}\right)}
 	 	  \right)
 		+ \underbrace{\Sum{j=2} \delta_{p,j} a_{j,y}}_{= \left(1-\delta_{p,1}\right) a_{p,y}} 
\\&= \Sum{i=2} a_{p,i} S_{i,y} + \left( \delta_{1,y} -\delta_{p,1} \right) a_{p,y} 
}%Gleichung

~\\ Zu 2.:
\Gleichung{
&\Sum{j=3} Z_{p,1,j} a_{j,z} 
\\=& -\Sum{j=3} \left(\mZii{p}{1}{j} \right)a_{j,z}
\\=& - \Sum{i=3} \Sum{j=3} a_{p,1}a_{1,i} a_{j,i}a_{j,z} 
\\=& - \Sum{i=3} a_{p,1} \left( \mSjj{i}{i}{z} + a_{1,i}\delta_{i,z} \right)
\\=& \Sum{i=3} a_{p,1} S_{i,i,z} - \Sum{i=3}a_{p,1}a_{1,i}\delta_{i,z}
\\=& \Sum{i=3} a_{p,1} S_{i,i,z} - (1- \delta_{1,z}-\delta_{2,z})a_{p,1}a_{1,z}
}%Gleichung

~\\ Zu 3.:
\Gleichung{
&\Sum{j=2} Z_{p,j,j} 
\\=& - \Sum{j=2} \left( \mZii{p}{j}{j} \right)
\\=& -\Sum{i=3} a_{p,1} \left( \mSj{i}{i} \right) + a_{p,1}
\\=& \Sum{i=3} a_{p,1} S_{i,i}  + a_{p,1}
}

~\\Zu 4.:
\Gleichung{
&\Sum{j=3}Z_{p,1} a_{j,1}a_{j,z} 
\\=& -\Sum{j=3} \left(\mZi{p}{1} \right)a_{j,1}a_{j,z}
\\=&- \Sum{i=2}\Sum{j=3} a_{p,i}a_{1,i}a_{j,1}a_{j,z} + \Sum{j=3}\delta_{p,1}a_{j,1}a_{j,z}
\\=& - \Sum{i=2}a_{p,i}\left( \mSjj{i}{1}{z} + a_{1,i}\delta_{1,z}\right)  + \Sum{j=3}\delta_{p,1}a_{j,1}a_{j,z}
\\=&  \Sum{i=2}a_{p,i} S_{i,1,z} - \Sum{i=2}a_{p,i}a_{1,i}\delta_{1,z}  + \Sum{j=3}\delta_{p,1}a_{j,1}a_{j,z}
}

~\\Zu 5.:
\Gleichung{
&Z_{p,q,1}\left(\Sum{j=3} a_{j,y}a_{j,z} - \delta_{y,z}\right)
\\=&\left( \Zii{p}{q}{1} \right)  \left(\Sum{j=3} a_{j,y}a_{j,z} - \delta_{y,z}\right)
\\=& - a_{p,1}\Sum{i=3} a_{q,i} \left( a_{1,i}  \Sum{j=3} a_{j,y}a_{j,z} - a_{1,i} \delta_{y,z}\right) + \delta_{q,1} a_{p,1}\left(\Sum{j=3} a_{j,y}a_{j,z} - \delta_{y,z}\right)
\\=& a_{p,1}\Sum{i=3}a_{q,i} S_{i,y,z} + \delta_{q,1} a_{p,1}\left(\Sum{j=3} a_{j,y}a_{j,z} - \delta_{y,z}\right)
}
~\\Zu 6.:
\Gleichung{
&\Sum{j=2} Z_{p,q,j} a_{j,y} 
\\=& -\Sum{j=2} \left(\mZii{p}{q}{j}\right) a_{j,y}
\\&=-\Sum{i=3}a_{p,1}a_{q,i} \left( \mSj{i}{y} + \delta_{i,y}\right) + \Sum{j=2}a_{p,1}\delta_{q,j}a_{j,y}
\\&= -\Sum{i=3} a_{p,1}a_{q,i} \left( \mSj{i}{y} \right) - \Sum{i=3}a_{p,1}a_{q,i}\delta_{i,y}+ \Sum{j=2}a_{p,1}\delta_{q,j}a_{j,y}
\\&= \Sum{i=3} a_{p,1}a_{q,i} S_{i,y} - (1-\delta_{y,1}-\delta_{y,2})a_{p,1}a_{q,y}+ (1-\delta_{q,1})a_{p,1}a_{q,y}
\\&=  \Sum{i=3} a_{p,1}a_{q,i} S_{i,y} + \delta_{y,1}a_{p,1}a_{q,1}+ \delta_{y,2}a_{p,1}a_{q,2} -\delta_{q,1}a_{p,1}a_{1,y}
}
}%MitRechnung
}%Beweis

}%Satz

%% file: Haupt/RechnungGBA/ZiZil.tex
\subsection{Reduktionswege für $a_{p,1}a_{q,1}a_{y,1}$}
\label{ZiZil}

\subsubsection{Beginnend mit $\left( a_{p,1}a_{q,1}a_{y,1}~,~ Z_{p,q} a_{y,1}\right)$:}
\label{sec:ZPQAPy1}
\Gleichung{
\underline{a_{p,1}a_{q,1}}a_{y,1} ~ \Pfeil{\widetilde{Z}_{p,q}} &Z_{p,q} a_{y,1} 
}%Geleichung
\subsubsection{Beginnend mit $( a_{p,1}a_{q,1}a_{y,1}~,~a_{p,1} Z_{q,y})$:}
\label{sec:aP1ZXY}
\Gleichung{
a_{p,1}\underline{a_{q,1}a_{y,1}} ~
 \Pfeil{\widetilde{Z}_{q,y}} 
 			&a_{p,1} Z_{q,y} = a_{p,1} \left( \Zi{q}{y} \right)\\
	 		&= - \underline{a_{p,1}a_{q,2}a_{y,2}} +  Z_{p,q,y} \\
 \Pfeil{\widetilde{Z}_{p,q,y}}
 			& - Z_{p,q,y} + Z_{p,q}a_{y,1} +   Z_{p,q,y}\\
 			&= Z_{p,q}a_{y,1}
}
\Bem{Symmetrie}{Diese Rechnungen sind nach Vertauschen der Indizes identisch mit denen in \Seite{SiSil}.}

%% file: Haupt/RechnungGBA/SiSil.tex
\subsection{Reduktionswege für $a_{1,p}a_{1,q}a_{1,y}$}
\label{SiSil}
\subsubsection{Beginnend mit $\left( a_{1,p}a_{1,q}a_{1,y}~,~ S_{p,q} a_{1,y}\right)$:}
\label{sec:SPQAPy1}
\Gleichung{
\underline{a_{1,p}a_{1,q}}a_{1,y} ~ \Pfeil{\widetilde{S}_{p,q}} &S_{p,q} a_{1,y} 
}%Geleichung
\subsubsection{Beginnend mit $( a_{1,p}a_{1,q}a_{1,y}~,~a_{1,p} S_{q,y})$:}
\label{sec:aP1SXY}
\Gleichung{
a_{1,p}\underline{a_{1,q}a_{1,y}} ~
 \Pfeil{\widetilde{S}_{q,y}} 
 			&a_{1,p} S_{q,y} = a_{1,p} \left( \Si{q}{y} \right)\\
	 		&= - \underline{a_{1,p}a_{2,q}a_{2,y}} +  S_{p,q,y} \\
 \Pfeil{\widetilde{S}_{p,q,y}}
 			& - S_{p,q,y} + S_{p,q}a_{1,y} +   S_{p,q,y}\\
 			&= S_{p,q}a_{1,y}
}
\Bem{Symmetrie}{Diese Rechnungen sind nach Vertauschen der Indizes identisch mit denen in \Seite{ZiZil}.}
\newpage

%% file: Haupt/RechnungGBA/ZiSil.tex
\subsection{Reduktionswege für $ a_{p,1}a_{1,1}a_{1,y}$}
\label{ZiSil}

\subsubsection{Beginnend mit $(a_{p,1}a_{1,1}a_{1,y} ~,~ Z_{p,1}a_{1,y})$:}
\Gleichung{
\underline{a_{p,1}a_{1,1} }a_{1,y}  ~ 
	\Pfeil{\widetilde{Z}_{p,1} } & \left( \Zi{p}{1} \right) a_{1,y}
	\\=&- \Sum{i=2} a_{p,i}\underline{a_{1,i} a_{1,y}} + \delta_{p,1}a_{1,y}
	\\ \Pfeil{\widetilde{S}_{i,y}} & - \Sum{i=2} a_{p,i} S_{i,y} + \delta_{p,1}a_{1,y}
	\\=& \Sum{i=2} a_{p,i} \left( \mSj{i}{y} \right) + \delta_{p,1}a_{1,y}
%	\\=& \Sum{i=2} \Sum{j=2}  a_{p,i}  a_{j,i}  a_{j,y} - \Sum{i=2} a_{p,i} \delta_{i,y} + \delta_{p,1}a_{1,y}
	\\=& \Sum{i=2} \Sum{j=2}  a_{p,i}  a_{j,i}  a_{j,y} - (1-\delta_{1,y})a_{p,y} + \delta_{p,1}a_{1,y}
	\\=& \Sum{i=2} \Sum{j=2}  a_{p,i}  a_{j,i}  a_{j,y} - a_{p,y}+ \delta_{1,y}a_{p,1} + \delta_{p,1} a_{1,y}
}%Geleichung
\subsubsection{Beginnend mit $(a_{p,1}a_{1,1}a_{1,y}  ~,~ a_{p,1} S_{1,y})$:}
\Gleichung{
  a_{p,1}\underline{a_{1,1}a_{1,y} }  ~ 
  \Pfeil{\widetilde{S}_{1,y} } 
    & - a_{p,1}\left( \mSj{1}{y}\right)  
 	 	\\=& - \Sum{j=2} \underline{a_{p,1} a_{j,1}} a_{j,y} + a_{p,1} \delta_{1,y}
 \\ \Pfeil{\widetilde{Z}_{p,j}}
  	& - \Sum{j=2} Z_{p,j} a_{j,y} + a_{p,1} \delta_{1,y}
  	\\=&  \Sum{j=2} \left( \mZi{p}{j} \right) a_{j,y} + a_{p,1} \delta_{1,y}
  %	\\=&  \Sum{j=2} \Sum{i=2} a_{p,i}  a_{j,i}  a_{j,y} - \Sum{j=2}\delta_{p,j}a_{j,y}  + a_{p,1} \delta_{1,y}
  	\\=& \Sum{j=2} \Sum{i=2} a_{p,i}  a_{j,i}  a_{j,y} -(1-\delta_{p,1})a_{p,y}  + a_{p,1} \delta_{1,y}
  	\\=& \Sum{j=2} \Sum{i=2} a_{p,i}  a_{j,i}  a_{j,y} - a_{p,y} + \delta_{p,1}a_{1,y}  + \delta_{1,y} a_{p,1} 
}%Gleichung

\Bem{Symmetrie}{Diese Rechnungen sind nach Vertauschen der Indizes identisch mit denen in \Seite{SiZil}.}

%% file: Haupt/RechnungGBA/SiZil.tex
\subsection{Reduktionswege für $ a_{1,p}a_{1,1}a_{y,1}$}
\label{SiZil}

\subsubsection{Beginnend mit $(a_{1,p}a_{1,1}a_{y,1} ~,~ S_{p,1}a_{y,1})$:}
\Gleichung{
\underline{a_{1,p}a_{1,1} }a_{y,1}  ~ 
	\Pfeil{\widetilde{S}_{p,1} } & S_{p,1} a_{y,1} = \left( \Si{p}{1} \right) a_{y,1}\\
		&=- \Sum{i=2} a_{i,p}\underline{a_{i,1} a_{y,1}} + \delta_{1,p}a_{y,1}\\
	\Pfeil{\widetilde{Z}_{i,y}} & - \Sum{i=2} a_{i,p} Z_{i,y} + \delta_{1,p}a_{y,1}\\
		&= \Sum{i=2} a_{i,p} \left( \mZj{i}{y} \right) + \delta_{1,p}a_{y,1}\\
		&= \Sum{i=2} \Sum{j=2}  a_{i,p}  a_{i,j}  a_{y,j} - \underbrace{\Sum{i=2} a_{i,p} \delta_{i,y}}_{=(1-\delta_{y,1})a_{1,p}} + \delta_{1,p}a_{y,1}
}%Geleichung
\subsubsection{Beginnend mit $(a_{1,p}a_{1,1}a_{y,1}  ~,~ a_{1,p} Z_{1,y})$:}
\Gleichung{
  a_{1,p}\underline{a_{1,1}a_{y,1} }  ~ 
  \Pfeil{\widetilde{Z}_{1,y} } 
    &a_{1,p} Z_{1,y} = - a_{1,p}\left( \mZj{1}{y}\right)  \\
 	 	&= - \Sum{j=2} \underline{a_{1,p} a_{j,1}} a_{i,y} + a_{1,p} \delta_{y,1}\\
  \Pfeil{\widetilde{S}_{p,j}}
  	& - \Sum{j=2} S_{p,j} a_{y,j} + a_{1,p} \delta_{y,1}\\
  	&=  \Sum{j=2} \left( \mSi{p}{j} \right) a_{y,j} + a_{1,p} \delta_{y,1}\\
  	&=  \Sum{j=2} \Sum{i=2} a_{i,p}  a_{i,j}  a_{y,j} - \underbrace{\Sum{j=2}\delta_{p,j}a_{y,j}}_{=(1-\delta_{p,1})a_{y,1} }  + a_{1,p} \delta_{y,1}
}%Gleichung

\Bem{Symmetrie}{Diese Rechnungen sind nach Vertauschen der Indizes identisch mit denen in \Seite{ZiSil}.}

%% file: Haupt/RechnungGBA/ZiSill.tex
\subsection{Reduktionswege für $a_{1,1}a_{1,1} $}
\label{SiZill}
\label{ZiSill}%Die Rechungen sind identisch!

\subsubsection{Beginnend mit $(a_{1,1}a_{1,1} ~,~ Z_{1,1})$:}

\Gleichung{
\underline{a_{1,1}a_{1,1}}  ~ 
	\Pfeil{\widetilde{Z}_{1,1} } 
		& - \Sum{i=2} a_{1,i}a_{1,i} + \delta_{1,1}\\
	\Pfeil{\widetilde{S}_{i,i}} 
		&  \Sum{i=2} \left( \mSj{i}{i} \right) + 1
}%Geleichung
\subsubsection{Beginnend mit $(a_{1,1}a_{1,1}  ~,~ S_{1,1})$:}
\Gleichung{
\underline{a_{1,1}a_{1,1}}  ~ 
	\Pfeil{\widetilde{S}_{1,1} } 
		& - \Sum{j=2} a_{j,1}a_{j,1} + \delta_{1,1}\\
	\Pfeil{\widetilde{Z}_{j,j}} 
		&  \Sum{j=2} \left( \mZi{j}{j}\right)  + 1
}%Geleichung

%% file: Haupt/RechnungGBA/ZiZiil.tex
\subsection{Reduktionswege für $a_{p,1}a_{q,1}a_{y,2}a_{z,2} $}
\label{ZiZiil}

\subsubsection{Beginnend mit $(a_{p,1}a_{q,1}a_{y,2}a_{z,2} ~,~ Z_{p,q} a_{y,2}a_{z,2}) $:}
\Gleichung{
\underline{a_{p,1}a_{q,1} }a_{y,2}a_{z,2}  ~ 
	\Pfeil{\widetilde{Z}_{p,q}} 
		& Z_{p,q} a_{y,2}a_{z,2} 
}%Geleichung
\subsubsection{Beginnend mit $(a_{p,1}a_{q,1}a_{y,2}a_{z,2}  ~,~ a_{p,1} \left( Z_{q,y,z} - Z_{q,y}a_{z,1} \right))$:}
\Gleichung{
a_{p,1} \underline{a_{q,1}a_{y,2}a_{z,2} }  ~
  \Pfeil{\widetilde{Z}_{q,y,z} } 
  	& a_{p,1} \left( Z_{q,y,z} - Z_{q,y}a_{z,1} \right) \\
  	&= - a_{p,1} \left( \mZii{q}{y}{z} \right) 
  	\\ &+ a_{p,1}\left( \mZi{q}{y} \right) a_{z,1}\\
  	&= -\underline{a_{p,1}a_{q,1}} \left( \Sum{i=3} a_{y,i}a_{z,i} - \delta_{y,z} \right) 
  	\\ &
  	+ \left( \underline{a_{p,1}a_{q,2}a_{y,2}} - Z_{p,q,y} \right) a_{z,1}\\
\Pfeil{ \widetilde{Z}_{p,q}, \widetilde{Z}_{p,q,y}}
		& - Z_{p,q}  \left( \Sum{i=3} a_{y,i}a_{z,i} - \delta_{y,z} \right)  
 	\\ &
		 + \left( Z_{p,q,y} - Z_{p,q}a_{y,1} - Z_{p,q,y}\right) a_{z,1}\\
		&= -  Z_{p,q} \left(  \Sum{i=3} a_{y,i}a_{z,i} - \delta_{y,z}  + \underline{a_{y,1} a_{z,1}}\right)\\
\Pfeil{\widetilde{Z}_{y,z}} 
		&= -  Z_{p,q} \left(  \Sum{i=3} a_{y,i}a_{z,i} - \delta_{y,z} + Z_{y,z} \right)\\
		&= - Z_{p,q}\left( ~- a_{2,y} a_{2,z} \right)
}%Gleichung
\Bem{Symmetrie}{Diese Rechnungen sind nach Vertauschen der Indizes identisch mit denen in \Seite{SiSiil}.}

%% file: Haupt/RechnungGBA/SiSiil.tex
\subsection{Reduktionswege für $a_{1,p}a_{1,q}a_{2,y}a_{2,z} $}
\label{SiSiil}

\subsubsection{Beginnend mit $(a_{1,p}a_{1,q}a_{2,y}a_{2,z} ~,~ S_{p,q} a_{2,y}a_{2,z}) $:}
\Gleichung{
\underline{a_{1,p}a_{1,q} }a_{2,y}a_{2,z}  ~ 
	\Pfeil{\widetilde{S}_{p,q}} 
		& S_{p,q} a_{2,y}a_{2,z} 
}%Geleichung
\subsubsection{Beginnend mit $(a_{1,p}a_{1,q}a_{2,y}a_{2,z}  ~,~ a_{1,p} \left( S_{q,y,z} - S_{q,y}a_{1,z} \right))$:}
\Gleichung{
a_{1,p} \underline{a_{1,q}a_{2,y}a_{2,z} }  ~
  \Pfeil{\widetilde{S}_{q,y,z} } 
  	& a_{1,p} \left( S_{q,y,z} - S_{q,y}a_{1,z} \right) \\
  	=& - a_{1,p} \left( \mSii{q}{y}{z} \right) 
  	\\ &
  	+ a_{1,p}\left( \mSi{q}{y} \right) a_{1,z}\\
  	=& -\underline{a_{1,p}a_{1,q}} \left( \Sum{i=3} a_{i,y}a_{i,z} - \delta_{y,z} \right) 
  	\\ &
  	+ \left( \underline{a_{1,p}a_{q,2}a_{2,y}} - S_{p,q,y} \right) a_{1,z}\\
\Pfeil{ \widetilde{S}_{p,q}, \widetilde{S}_{p,q,y}}
		& - S_{p,q}  \left( \Sum{i=3} a_{i,y}a_{i,z} - \delta_{y,z} \right)  
 	\\ &
		 + \left( S_{p,q,y} - S_{p,q}a_{y,1} - S_{p,q,y}\right) a_{1,z}\\
		=& -  S_{p,q} \left(  \Sum{i=3} a_{i,y}a_{i,z} - \delta_{y,z}  + \underline{a_{y,1} a_{1,z}}\right)\\
\Pfeil{\widetilde{S}_{y,z}} 
		=& -  S_{p,q} \left(  \Sum{i=3} a_{i,y}a_{i,z} - \delta_{y,z} + S_{y,z} \right)\\
		=& - S_{p,q}\left( ~- a_{2,y} a_{2,z} \right)
}%Gleichung
\Bem{Symmetrie}{Diese Rechnungen sind nach Vertauschen der Indizes identisch mit denen in \Seite{ZiZiil}.}

%% file: Haupt/RechnungGBA/ZiSiil.tex
\subsection{Reduktionswege für $\Erzeuger{p}{1}\Erzeuger{1}{1}\Erzeuger{2}{y}\Erzeuger{2}{z} $}
\label{ZiSiil}
\Bem{Symmetrie}{Diese Rechnungen sind nach Vertauschen der Indizes identisch mit denen in \ref{SiZiil}.}

\input{Haupt/RechnungGBA/RechnungZiSiil}

%% file: Haupt/RechnungGBA/RechnungZiSiil.tex
\subsubsection{Beginnend mit $(\Erzeuger{p}{1}\Erzeuger{1}{1}\Erzeuger{2}{y}\Erzeuger{2}{z} ~,~\Zeile_{p,1}\Erzeuger{2}{y}\Erzeuger{2}{z} )$:}

\Gleichung{
\underline{\Erzeuger{p}{1}\Erzeuger{1}{1}} \Erzeuger{2}{y}\Erzeuger{2}{z}  ~ 
 \Pfeil{ \widetilde{\Zeile}_{p,1}} &  \left( \Zi{p}{1} \right) \Erzeuger{2}{y}\Erzeuger{2}{z}\\
		= &- \Sum{i=2} \Erzeuger{p}{i} \underline{\Erzeuger{1}{i} \Erzeuger{2}{y} \Erzeuger{2}{z}} + \DELTA_{p,1} \Erzeuger{2}{y}\Erzeuger{2}{z}\\
	\Pfeil{ \widetilde{\Spalte}_{i,y,z} }& - \Sum{i=2} \Erzeuger{p}{i} \left( \Spalte_{i,y,z} - \Spalte_{i,y}\Erzeuger{1}{z} \right)+ \DELTA_{p,1} \Erzeuger{2}{y}\Erzeuger{2}{z}\\
		= &  \Sum{i=2} \Erzeuger{p}{i} \left( \mSjj{i}{y}{z}  \right)  
		\\& + \Sum{i=2} \Erzeuger{p}{i} \Spalte_{i,y} \Erzeuger{1}{z}  + \DELTA_{p,1} \Erzeuger{2}{y}\Erzeuger{2}{z}
  \\= & \Sum{i=2}\Sum{j=3}\Erzeuger{p}{i} \Erzeuger{1}{i} \Erzeuger{j}{y}\Erzeuger{j}{z} 
  			- \Sum{i=2}\Erzeuger{p}{i}\Erzeuger{1}{i}\DELTA_{y,z}  
	\\&+ \Sum{i=2} \Erzeuger{p}{i} \Spalte_{i,y} \Erzeuger{1}{z}  + \DELTA_{p,1} \Erzeuger{2}{y}\Erzeuger{2}{z}
	\\= & \Sum{i=2}\Sum{j=3}\Erzeuger{p}{i} \Erzeuger{1}{i} \Erzeuger{j}{y}\Erzeuger{j}{z} 
				+ \left( \Zeile_{p,1} - \DELTA_{p,1} \right) \DELTA_{y,z}
	\\& + \Sum{i=2} \Erzeuger{p}{i} \Spalte_{i,y} \Erzeuger{1}{z}  + \DELTA_{p,1} \Erzeuger{2}{y}\Erzeuger{2}{z}
	\\= & \Sum{i=2}\Sum{j=3}\Erzeuger{p}{i} \Erzeuger{1}{i} \Erzeuger{j}{y}\Erzeuger{j}{z} +  \Zeile_{p,1}\DELTA_{y,z} 
	\\&+ \Sum{i=2} \Erzeuger{p}{i} \Spalte_{i,y} \Erzeuger{1}{z}  - \DELTA_{p,1}\left( \DELTA_{y,z} -  \Erzeuger{2}{y}\Erzeuger{2}{z}\right)
}%Geleichung
\newpage
\subsubsection{Beginnend mit $(\Erzeuger{p}{1}\Erzeuger{1}{1}\Erzeuger{2}{y}\Erzeuger{2}{z}  ~,~ \Erzeuger{p}{1} \left( \Spalte_{1,y,z} - \Spalte_{1,y} \Erzeuger{1}{z} \right))$:}
\Gleichung{
&\Erzeuger{p}{1}\underline{\Erzeuger{1}{1}\Erzeuger{2}{y}\Erzeuger{2}{z} }  ~ 
\\  \Pfeil{ \widetilde{\Spalte}_{1,y,z} }&
  % \Erzeuger{p}{1} \left( \Spalte_{1,y,z} - \Spalte_{1,y} \Erzeuger{1}{z} \right) \\	=&
	 - \Erzeuger{p}{1} \left( \mSjj{1}{y}{z} \right)
\\ & + \Erzeuger{p}{1} \left( \mSj{1}{y} \right) \Erzeuger{1}{z} 
\\
	=& - \underline{ \Erzeuger{p}{1}\Erzeuger{1}{1}} \left( \Sum{j=3} \Erzeuger{j}{y} \Erzeuger{j}{z} -  \DELTA_{y,z} \right) 
\\ & + \Sum{j=2}\underline{\Erzeuger{p}{1}\Erzeuger{j}{1}} \Erzeuger{j}{y} \Erzeuger{1}{z} - \Erzeuger{p}{1}\DELTA_{1,y}\Erzeuger{1}{z}
\\ 
  \Pfeil{ \widetilde{\Zeile}_{p,1}, \widetilde{\Zeile}_{p,j}} &
			 \left( \mZi{p}{1} \right) \left( \Sum{j=3} \Erzeuger{j}{y} \Erzeuger{j}{z} - \DELTA_{y,z} \right)
\\&		 + \Sum{j=2} \Zeile_{p,j} \Erzeuger{j}{y} \Erzeuger{1}{z} - \Erzeuger{p}{1}\DELTA_{1,y}\Erzeuger{1}{z}
\\ 
	\stackrel{\ref{Rechenregel}}{=}&  \Sum{i=2} \Sum{j=3} \Erzeuger{p}{i} \Erzeuger{1}{i} \Erzeuger{j}{y} \Erzeuger{j}{z} 
		- \DELTA_{p,1} \Sum{j=3} \Erzeuger{j}{y} \Erzeuger{j}{z}
		+\Zeile_{p,1} \DELTA_{y,z}
\\&+ \left( \Sum{i=2} \Erzeuger{p}{i} \Spalte_{i,y} + \left( \DELTA_{1,y} -\DELTA_{p,1}\right) \Erzeuger{p}{y} \right) \Erzeuger{1}{z}
			 - \Erzeuger{p}{1}\DELTA_{1,y}\Erzeuger{1}{z}
\\ = & \Sum{i=2} \Sum{j=3} \Erzeuger{p}{i} \Erzeuger{1}{i} \Erzeuger{j}{y} \Erzeuger{j}{z} 
			 +\Zeile_{p,1} \DELTA_{y,z}
			 +  \Sum{i=2} \Erzeuger{p}{i} \Spalte_{i,y}\Erzeuger{1}{z} 
	\\&		 -\DELTA_{p,1} \left(\Sum{j=3} \Erzeuger{j}{y} \Erzeuger{j}{z} + \underline{\Erzeuger{1}{y} \Erzeuger{1}{z}} \right)	
\\ \Pfeil{ \widetilde{\Spalte}_{y,z}} & 
			\Sum{i=2} \Sum{j=3} \Erzeuger{p}{i} \Erzeuger{1}{i} \Erzeuger{j}{y} \Erzeuger{j}{z} 
			 +\Zeile_{p,1} \DELTA_{y,z}
			 +  \Sum{i=2} \Erzeuger{p}{i} \Spalte_{i,y}\Erzeuger{1}{z} 
	\\&		 -\DELTA_{p,1} \left(\Sum{j=3} \Erzeuger{j}{y} \Erzeuger{j}{z} - \left( \mSj{y}{z} \right) \right)	
	\\= & \Sum{i=2}\Sum{j=3}\Erzeuger{p}{i} \Erzeuger{1}{i} \Erzeuger{j}{y}\Erzeuger{j}{z} +  \Zeile_{p,1}\DELTA_{y,z} + \Sum{i=2} \Erzeuger{p}{i} \Spalte_{i,y} \Erzeuger{1}{z}  
	\\&- \DELTA_{p,1}\left( \DELTA_{y,z} -  \Erzeuger{2}{y}\Erzeuger{2}{z}\right)
}%Gleichung

%% file: Haupt/RechnungGBA/SiZiil.tex
\subsection{Reduktionswege für $\Erzeuger{1}{p}\Erzeuger{1}{1}\Erzeuger{y}{2}\Erzeuger{z}{2} $}
\label{SiZiil}
\Bem{Symmetrie}{Diese Rechnungen sind nach Vertauschen der Indizes identisch mit denen in \ref{ZiSiil}.}

\input{StyleFiles/VertauschteIndizes}

\input{Haupt/RechnungGBA/RechnungZiSiil}
\input{StyleFiles/NormaleIndizes}

%% file: StyleFiles/VertauschteIndizes.tex
\renewcommand{\Erzeuger}[2]{a_{#2,#1}}
\renewcommand{\Zeile}{S}
\renewcommand{\Spalte}{Z}

%% file: StyleFiles/NormaleIndizes.tex
\renewcommand{\Erzeuger}[2]{a_{#1,#2}}
\renewcommand{\Zeile}{Z}
\renewcommand{\Spalte}{S}

%% file: Haupt/RechnungGBA/ZiSiill.tex
\subsection{Reduktionswege für $\Erzeuger{1}{1}\Erzeuger{2}{1}\Erzeuger{2}{z}$}
\label{ZiSiill}

\input{Haupt/RechnungGBA/RechnungZiSiill}

\Bem{Symmetrie}{Diese Rechnungen sind nach Vertauschen der Indizes identisch mit denen in \Seite{SiZiill}.}

%% file: Haupt/RechnungGBA/RechnungZiSiill.tex
\subsubsection{Beginnend mit $( \Erzeuger{1}{1}\Erzeuger{2}{1}\Erzeuger{2}{z}~,~ \Zeile_{1,2}\Erzeuger{2}{z} )$:}
\Gleichung{
\underline{\Erzeuger{1}{1}\Erzeuger{2}{1} } \Erzeuger{2}{z} ~ \Pfeil{ \widetilde{\Zeile}_{1,2}} =& \left( \Zi{1}{2} \right) \Erzeuger{2}{z}\\
	=& -\Sum{i=2}\underline{\Erzeuger{1}{i}\Erzeuger{2}{i}\Erzeuger{2}{z}} 
\\	\Pfeil{\widetilde{\Spalte}_{i,i,z}} &- \Sum{i=2} \left(\Spalte_{i,i,z} - \Spalte_{i,i} \Erzeuger{1}{z} \right)
}%Geleichung

\subsubsection{Beginnend mit $( \Erzeuger{1}{1}\Erzeuger{2}{1}\Erzeuger{1}{z} ~,~ \Spalte_{1,1,z} - \Spalte_{1,1}\Erzeuger{1}{z} )$:}
\Gleichung{
&\underline{\Erzeuger{1}{1}\Erzeuger{2}{1}\Erzeuger{1}{z}}  ~
\\ \Pfeil{\widetilde{\Spalte}_{1,1,z} } &  \left( \Sii{1}{1}{z} \right) + \left( \Sum{j=2} \underline{\Erzeuger{j}{1}\Erzeuger{j}{1}} - \delta_{1,1} \right)\Erzeuger{1}{z}
\\ \Pfeil{\widetilde{\Zeile}_{j,j}} & -\Sum{j=3}\underline{\Erzeuger{1}{1}\Erzeuger{j}{1}}\Erzeuger{j}{z} 
	+ \Erzeuger{1}{1}\delta_{1,z}  
\\&	- \Sum{i=2} \left(  \Sum{j=2}\Erzeuger{j}{i} \Erzeuger{j}{i} - \delta_{j,j} \right) \Erzeuger{1}{z} - \Erzeuger{1}{z}
	\\
	\Pfeil{\widetilde{\Zeile}_{1,j}} & 
	\Sum{j=3} \left( \mZi{1}{j} \right) \Erzeuger{j}{z}+ \Erzeuger{1}{1}\delta_{1,z}
  \\& - \Sum{i=2} \left( \Sum{j=2}\Erzeuger{j}{i} \Erzeuger{j}{i} - \delta_{i,i} \right) \Erzeuger{1}{z} 
  - \Erzeuger{1}{z}\\ 
= &\Sum{i=2}\Sum{j=3} \Erzeuger{1}{i}\Erzeuger{j}{i}\Erzeuger{j}{z} + \Erzeuger{1}{1}\delta_{1,z} 
\\&+ \Sum{i=2} \Spalte_{i,i}\Erzeuger{1}{z} - \left(  \Erzeuger{1}{1}\delta_{1,z}  + \Sum{i=2}  \Erzeuger{1}{i} \delta_{i,z}\right)\\
=  &- \Sum{i=2}\left( \Spalte_{i,i,z} - \Spalte_{i,i}\Erzeuger{1}{z} \right)
}%Gleichung

%% file: Haupt/RechnungGBA/SiZiill.tex
\subsection{Reduktionswege für $\Erzeuger{1}{1}\Erzeuger{1}{2}\Erzeuger{z}{2}$}
\label{SiZiill}

\input{StyleFiles/VertauschteIndizes}
\input{Haupt/RechnungGBA/RechnungZiSiill}
\Bem{Symmetrie}{Diese Rechnungen sind nach Vertauschen der Indizes identisch mit denen in \Seite{ZiSiill}.}

\input{StyleFiles/NormaleIndizes}

%% file: Haupt/RechnungGBA/ZiiSil.tex
\subsection{Reduktionswege für $\Erzeuger{p}{1}\Erzeuger{q}{2}\Erzeuger{1}{2}\Erzeuger{1}{y} $}
\label{ZiiSil}
\input{Haupt/RechnungGBA/RechnungZiiSil}

\Bem{Symmetrie}{Diese Rechnungen sind nach Vertauschen der Indizes identisch mit denen in \ref{SiiZil}.}

%% file: Haupt/RechnungGBA/RechnungZiiSil.tex
\subsubsection{Beginnend mit $(\Erzeuger{p}{1}\Erzeuger{q}{2}\Erzeuger{1}{2}\Erzeuger{1}{y} ~,~ \left( \Zeile_{p,q,1} - \Zeile_{p,q} \Erzeuger{1}{1} \right)\Erzeuger{1}{y})$:}
\Gleichung{
&\underline{\Erzeuger{p}{1}\Erzeuger{q}{2}\Erzeuger{1}{2} }\Erzeuger{1}{y}  ~ 
\\\Pfeil{ \widetilde{Z}_{p,q,1}} 
&   \left( \Zii{p}{q}{1} \right)\Erzeuger{1}{y} - \Zeile_{p,q} \underline{\Erzeuger{1}{1} \Erzeuger{1}{y}}
\\ \Pfeil{ \widetilde{S}_{1,y}}& - \Erzeuger{p}{1} \Sum{i=3}\Erzeuger{q}{i}\underline{\Erzeuger{1}{i}\Erzeuger{1}{y}} + \Erzeuger{p}{1}\DELTA_{q,1}\Erzeuger{1}{y} - \Zeile_{p,q} \Spalte_{1,y}
\\ \Pfeil{ \widetilde{S}_{i,y}}&  \Sum{i=3} \Erzeuger{p}{1} \Erzeuger{q}{i} \left( \mSj{i}{y} \right) 
		+ \Erzeuger{p}{1}\DELTA_{q,1}\Erzeuger{1}{y} - \Zeile_{p,q} \Spalte_{1,y}
\\=& \Sum{i=3}\Sum{i=2} \Erzeuger{p}{1} \Erzeuger{q}{i}\Erzeuger{j}{i}\Erzeuger{j}{y} - \Zeile_{p,q} \Spalte_{1,y} 
	\\&- \Sum{i=3}\Erzeuger{p}{1}\Erzeuger{q}{i} \DELTA_{i,y} + \Erzeuger{p}{1}\DELTA_{q,1}\Erzeuger{1}{y}
\\=& \Sum{i=3}\Sum{i=2} \Erzeuger{p}{1} \Erzeuger{q}{i}\Erzeuger{j}{i}\Erzeuger{j}{y} - \Zeile_{p,q} \Spalte_{1,y} 
		\\ &- (1- (\DELTA_{1,y}+\DELTA_{2,y})) \Erzeuger{p}{1}\Erzeuger{q}{y} + \Erzeuger{p}{1}\DELTA_{q,1}\Erzeuger{1}{y}
\\=& \Sum{i=3}\Sum{i=2} \Erzeuger{p}{1} \Erzeuger{q}{i}\Erzeuger{j}{i}\Erzeuger{j}{y} - \Zeile_{p,q} \Spalte_{1,y} 
		\\ &- \Erzeuger{p}{1}\Erzeuger{q}{y} + \DELTA_{1,y}\Erzeuger{p}{1}\Erzeuger{q}{1}+\DELTA_{2,y} \Erzeuger{p}{1}\Erzeuger{q}{2} + \Erzeuger{p}{1}\DELTA_{q,1}\Erzeuger{1}{y}
}%Geleichung
\newpage

\subsubsection{Beginnend mit $(\Erzeuger{p}{1}\Erzeuger{q}{2}\Erzeuger{1}{2}\Erzeuger{1}{y}  ~,~ \Erzeuger{p}{1}\Erzeuger{q}{2}\Spalte_{2,y})$:}
\Gleichung{
&\Erzeuger{p}{1}\Erzeuger{q}{2}\underline{\Erzeuger{1}{2}\Erzeuger{1}{y} }  ~
\\ \Pfeil{\widetilde{S}_{2,y} } & - \Erzeuger{p}{1}\Erzeuger{q}{2} \left( \mSj{2}{y} \right)
\\=& - \Sum{j=2} \underline{\Erzeuger{p}{1}\Erzeuger{q}{2}\Erzeuger{j}{2}}\Erzeuger{j}{y} + \Erzeuger{p}{1}\Erzeuger{q}{2}\DELTA_{2,y}
\\ \Pfeil{\widetilde{Z}_{p,q,j}}& - \Sum{j=2} \left(\Zeile_{p,q,j} - \Zeile_{p,q}\Erzeuger{j}{1}\right)\Erzeuger{j}{y} 
	+ \Erzeuger{p}{1}\Erzeuger{q}{2}\DELTA_{2,y}
\\=& \Sum{j=2} \left( \mZii{p}{q}{j} \right)\Erzeuger{j}{y} + \Zeile_{p,q} \Sum{j=2}  \Erzeuger{j}{1}\Erzeuger{j}{y}
\\&+ \Erzeuger{p}{1}\Erzeuger{q}{2}\DELTA_{2,y}
\\=& \Sum{j=2} \left( \mZii{p}{q}{j} \right)\Erzeuger{j}{y} -  \Zeile_{p,q}  \Spalte_{1,y} 
\\&+  \Zeile_{p,q} \DELTA_{1,y} + \Erzeuger{p}{1}\Erzeuger{q}{2}\DELTA_{2,y}
\\=& \Sum{i=3}\Sum{j=2} \Erzeuger{p}{1}\Erzeuger{q}{i}\Erzeuger{j}{i}\Erzeuger{j}{y} - \Zeile_{p,q}\Spalte_{1,y} 
	\\&- \underbrace{\Sum{j=2} \Erzeuger{p}{1} \DELTA_{q,j} \Erzeuger{j}{y}}_{(1-\DELTA{1,q}) \Erzeuger{p}{1}\Erzeuger{q}{y}} + \Zeile_{p,q}\DELTA_{1,y} + \Erzeuger{p}{1}\Erzeuger{q}{2}\DELTA_{2,y}
\\=& \Sum{i=3}\Sum{j=2} \Erzeuger{p}{1}\Erzeuger{q}{i}\Erzeuger{j}{i}\Erzeuger{j}{y} - \Zeile_{p,q}\Spalte_{1,y} 
	\\&- \Erzeuger{p}{1} \Erzeuger{q}{y} + \DELTA_{1,q} \Erzeuger{p}{1}\Erzeuger{1}{y} + \Zeile_{p,q}\DELTA_{1,y} + \Erzeuger{p}{1}\Erzeuger{q}{2}\DELTA_{2,y}
}%Gleichung

%% file: Haupt/RechnungGBA/SiiZil.tex
\subsection{Reduktionswege für $\Erzeuger{1}{p}\Erzeuger{2}{q}\Erzeuger{2}{1}\Erzeuger{y}{1} $}
\label{SiiZil}

\input{StyleFiles/VertauschteIndizes}
\input{Haupt/RechnungGBA/RechnungZiiSil}
\Bem{Symmetrie}{Diese Rechnungen sind nach Vertauschen der Indizes identisch mit denen in \ref{ZiiSil}.}

\input{StyleFiles/NormaleIndizes}

%% file: Haupt/RechnungGBA/ZiiSill.tex
\subsection{Reduktionswege für $\Erzeuger{p}{1}\Erzeuger{1}{2}\Erzeuger{1}{2}$}
\label{ZiiSill}
\input{Haupt/RechnungGBA/RechnungZiiSill}
\Bem{Symmetrie}{Diese Rechnungen sind nach Vertauschen der Indizes identisch mit denen in \ref{SiiZill}.}

%% file: Haupt/RechnungGBA/RechnungZiiSill.tex
\subsubsection{Beginnend mit $(\Erzeuger{p}{1}\Erzeuger{1}{2}\Erzeuger{1}{2} ~,~ \Zeile_{p,1,1} - \Zeile_{p,1}\Erzeuger{1}{1})$:}
\Gleichung{
&\underline{\Erzeuger{p}{1}\Erzeuger{1}{2}\Erzeuger{1}{2} }  ~ 
\\ \Pfeil{\widetilde{\Zeile}_{p,1,1}} & -\Erzeuger{p}{1} \Sum{i=3} \underline{\Erzeuger{1}{i}\Erzeuger{1}{i}} + \Erzeuger{p}{1}\delta_{1,1} + \Sum{i=2}\Erzeuger{p}{i}\underline{\Erzeuger{1}{i}\Erzeuger{1}{1}} - \delta_{p,1} \Erzeuger{1}{1}
\\ \Pfeil{\widetilde{\Spalte}_{i,i}, \widetilde{\Spalte}_{i,1} } & -\Erzeuger{p}{1} \Sum{i=3} \Spalte_{i,i} + \Erzeuger{p}{1} + \Sum{i=2}\Erzeuger{p}{i}\Spalte_{i,1}- \delta_{p,1} \Erzeuger{1}{1}
\\= & -\Erzeuger{p}{1} \Sum{i=3} \Spalte_{i,i} + \Sum{i=2}\Erzeuger{p}{i}\Spalte_{i,1} + (1- \delta_{p,1}) \Erzeuger{p}{1}
}%Geleichung
\subsubsection{Beginnend mit $(\Erzeuger{p}{1}\Erzeuger{1}{2}\Erzeuger{1}{2}  ~,~ \Erzeuger{p}{1} \Spalte_{2,2} )$:}
\Gleichung{
\Erzeuger{p}{1}\underline{\Erzeuger{1}{2}\Erzeuger{1}{2} }  ~ \Pfeil{ \widetilde{\Spalte}_{2,2}} & - \Sum{j=2}\underline{\Erzeuger{p}{1}\Erzeuger{j}{2}\Erzeuger{j}{2}} + \Erzeuger{p}{1}\delta_{2,2}
\\ \Pfeil{\widetilde{\Zeile}_{p,j,j}} & - \Sum{j=2} \left( \Zeile_{p,j,j} - \Zeile_{p,j}\Erzeuger{j}{1}\right) + \Erzeuger{p}{1}
\\=&   \Sum{j=2} \left( \mZii{p}{j}{j}\right)  
\\&- \Sum{j=2} \left( \mZi{p}{j}\right) \Erzeuger{j}{1} +  \Erzeuger{p}{1} 
\\=& \Erzeuger{p}{1} \Sum{i=3}  \left( \Sum{j=2} \Erzeuger{j}{i} \Erzeuger{j}{i} - \delta_{i,i}\right) - \Erzeuger{p}{1} 
\\&- \Sum{i=2} \Erzeuger{p}{i}\Sum{j=2}  \left(  \Erzeuger{j}{i}\Erzeuger{j}{1} - 0 \right) + \underbrace{\Sum{j=2} \delta_{p,j}\Erzeuger{j}{1}}_{= (1-\delta_{p,1})\Erzeuger{p}{1}} + \Erzeuger{p}{1}
\\=& - \Erzeuger{p}{1} \Sum{i=3}\Spalte_{i,i} + \Sum{i=2} \Erzeuger{p}{i}\Spalte_{i,1}  + (1-\delta_{p,1})\Erzeuger{p}{1}
}%Gleichung

%% file: Haupt/RechnungGBA/SiiZill.tex
\subsection{Reduktionswege für $\Erzeuger{1}{p}\Erzeuger{2}{1}\Erzeuger{2}{1}$}
\label{SiiZill}

\input{StyleFiles/VertauschteIndizes}
\input{Haupt/RechnungGBA/RechnungZiiSill}
\Bem{Symmetrie}{Diese Rechnungen sind nach Vertauschen der Indizes identisch mit denen in \ref{ZiiSill}.}

\input{StyleFiles/NormaleIndizes}

%% file: Haupt/RechnungGBA/ZiiSiil.tex
\subsection{Reduktionswege für $\Erzeuger{p}{1}\Erzeuger{q}{2}\Erzeuger{1}{2}\Erzeuger{2}{y}\Erzeuger{2}{z}$}
\label{ZiiSiil}
\input{Haupt/RechnungGBA/RechnungZiiSiil}

\Bem{Symmetrie}{Diese Rechnungen sind nach Vertauschen der Indizes identisch mit denen in \ref{SiiZiil}.}

%% file: Haupt/RechnungGBA/RechnungZiiSiil.tex
\subsubsection{Beginnend mit $(\Erzeuger{p}{1}\Erzeuger{q}{2}\Erzeuger{1}{2}\Erzeuger{2}{y}\Erzeuger{2}{z} ~,~ \left( \Zeile_{p,q,1} - \Zeile_{p,q}\Erzeuger{1}{1} \right) \Erzeuger{2}{y}\Erzeuger{2}{z})$:}
\Gleichung{
&\underline{\Erzeuger{p}{1}\Erzeuger{q}{2}\Erzeuger{1}{2} }\Erzeuger{2}{y}\Erzeuger{2}{z}  ~
\\ \Pfeil{\widetilde{Z}_{p,q,1} } &  - \left(\mZii{p}{q}{1}\right)\Erzeuger{2}{y}\Erzeuger{2}{z} - \Zeile_{p,q} \underline{\Erzeuger{1}{1} \Erzeuger{2}{y}\Erzeuger{2}{z}}
\\ \Pfeil{\widetilde{S}_{1,y,z} } & - \Sum{i=3} \Erzeuger{p}{1}\Erzeuger{q}{i}\underline{\Erzeuger{1}{i}\Erzeuger{2}{y}\Erzeuger{2}{z}} + \Erzeuger{p}{1}\DELTA_{q,1}\Erzeuger{2}{y}\Erzeuger{2}{z}  
\\&- \Zeile_{p,q} \left( \Spalte_{1,y,z}-\Spalte_{1,y}\Erzeuger{1}{z}\right)
\\ \Pfeil{\widetilde{S}_{i,y,z} } & - \Sum{i=3} \Erzeuger{p}{1}\Erzeuger{q}{i}\left( \Spalte_{i,y,z} - \Spalte_{i,y}\Erzeuger{1}{z} \right) 
+ \Erzeuger{p}{1}\DELTA_{q,1}\Erzeuger{2}{y}\Erzeuger{2}{z}  
\\&- \Zeile_{p,q} \left( \Spalte_{1,y,z}-\Spalte_{1,y}\Erzeuger{1}{z}\right)
}%Gleichung 
\subsubsection{Beginnend mit $(\Erzeuger{p}{1}\Erzeuger{q}{2}\Erzeuger{1}{2}\Erzeuger{2}{y}\Erzeuger{2}{z}  ~,~ \Erzeuger{p}{1}\Erzeuger{q}{2}\left( \Spalte_{2,y,z} - \Spalte_{2,y}\Erzeuger{1}{z}\right))$:}
\Gleichung{
&\Erzeuger{p}{1}\Erzeuger{q}{2}\underline{\Erzeuger{1}{2}\Erzeuger{2}{y}\Erzeuger{2}{z} }  ~
\\ \Pfeil{\widetilde{S}_{2,y,z} } &  - \Erzeuger{p}{1}\Erzeuger{q}{2}\left( \mSjj{2}{y}{z} \right) 
\\&+ \Erzeuger{p}{1}\Erzeuger{q}{2} \left( \mSj{2}{y} \right) \Erzeuger{1}{z}
\\=& -\underline{\Erzeuger{p}{1}\Erzeuger{q}{2}\Erzeuger{1}{2}} \left(\Sum{j=3} \Erzeuger{i}{y}\Erzeuger{i}{z} - \DELTA_{y,z}\right) 
\\&+ \Sum{j=2}\underline{\Erzeuger{p}{1}\Erzeuger{q}{2}\Erzeuger{j}{2}} \Erzeuger{j}{y}\Erzeuger{1}{z} - \Erzeuger{p}{1}\Erzeuger{q}{2}\DELTA_{2,y}\Erzeuger{1}{z}
\\
\Pfeil{\widetilde{Z}_{p,q,1}, \widetilde{Z}_{p,q,j}   } & - \left( \Zeile_{p,q,1}-\Zeile_{p,q}\Erzeuger{1}{1}\right) \left(\Sum{j=3} \Erzeuger{j}{y}\Erzeuger{j}{z} - \DELTA_{y,z}\right) 
\\&+\Sum{j=2}\left( \Zeile_{p,q,j}-\Zeile_{p,q}\Erzeuger{j}{1}\right) \Erzeuger{j}{y}\Erzeuger{1}{z} - \Erzeuger{p}{1}\Erzeuger{q}{2}\DELTA_{2,y}\Erzeuger{1}{z}
%%
%
%
%Seitenumbruch in der Gleichung
\displaybreak \\
\stackrel{\ref{Rechenregel}}{=}&-  \left(\Erzeuger{p}{1}\Sum{i=3}\Erzeuger{q}{i} \Spalte_{i,y,z} + \DELTA_{q,1} \Erzeuger{p}{1}\left(\Sum{j=3} \Erzeuger{j}{y}\Erzeuger{j}{z} - \DELTA_{y,z}\right)\right)
\\& - \Zeile_{p,q}\Spalte_{1,y,z} +\Biggl(\Sum{i=3} a_{p,1}a_{q,i} S_{i,y} + \delta_{y,1}a_{p,1}a_{q,1}.
\\& ~~~~~~~~~~~~~~~~~~~~~~~~~~~~~~~~~~~~+ \delta_{y,2}a_{p,1}a_{q,2} -\delta_{q,1}a_{p,1}a_{1,y}\Biggr)\Erzeuger{1}{z}
\\&- \Sum{j=2}\Zeile_{p,q}  \Erzeuger{j}{1} \Erzeuger{j}{y} \Erzeuger{1}{z} - \Erzeuger{p}{1}\Erzeuger{q}{2}\DELTA_{2,y}\Erzeuger{1}{z}
\\=&- \Sum{i=3}\Erzeuger{p}{1}\Erzeuger{q}{i}\left( \Spalte_{i,y,z}-  \Spalte_{i,y}\Erzeuger{1}{z}\right)
\\& - \DELTA_{q,1} \Erzeuger{p}{1}\left(\Sum{j=3} \Erzeuger{j}{y}\Erzeuger{j}{z} - \DELTA_{y,z}\right) - \DELTA_{q,1}\Erzeuger{p}{1}\underline{\Erzeuger{1}{y}\Erzeuger{1}{z}}
\\& +\left(\DELTA_{y,1}\underline{\Erzeuger{p}{1}\Erzeuger{q}{1}}- \Sum{j=2}\Zeile_{p,q}  \Erzeuger{j}{1} \Erzeuger{j}{y}  \right)\Erzeuger{1}{z} - \Zeile_{p,q}\Spalte_{1,y,z}
\\\Pfeil{\widetilde{Z}_{p,q},\widetilde{S}_{p,q}}&- \Sum{i=3}\Erzeuger{p}{1}\Erzeuger{q}{i}\left( \Spalte_{i,y,z}-  \Spalte_{i,y}\Erzeuger{1}{z}\right)
\\& - \DELTA_{q,1} \Erzeuger{p}{1}\left(\Sum{j=3} \Erzeuger{j}{y}\Erzeuger{j}{z} - \DELTA_{y,z}\right) - \DELTA_{q,1}\Erzeuger{p}{1}\Spalte_{p,q}
\\& +\left(\DELTA_{y,1}\Zeile_{p,q}- \Sum{j=2}\Zeile_{p,q}  \Erzeuger{j}{1} \Erzeuger{j}{y}  \right)\Erzeuger{1}{z} - \Zeile_{p,q}\Spalte_{1,y,z}
\\=&- \Sum{i=3} \Erzeuger{p}{1}\Erzeuger{q}{i}\left( \Spalte_{i,y,z} - \Spalte_{i,y}\Erzeuger{1}{z} \right) + \Erzeuger{p}{1}\DELTA_{q,1}\Erzeuger{2}{y}\Erzeuger{2}{z}  
\\&- \Zeile_{p,q} \left( \Spalte_{1,y,z}-\Spalte_{1,y}\Erzeuger{1}{z}\right)
}

%% file: Haupt/RechnungGBA/SiiZiil.tex
\subsection{Reduktionswege für $\Erzeuger{1}{p}\Erzeuger{2}{q}\Erzeuger{2}{1}\Erzeuger{y}{2}\Erzeuger{z}{2}$}
\label{SiiZiil}

\input{StyleFiles/VertauschteIndizes}
\input{Haupt/RechnungGBA/RechnungZiiSiil}
\Bem{Symmetrie}{Diese Rechnungen sind nach Vertauschen der Indizes identisch mit denen in \ref{ZiiSiil}.}

\input{StyleFiles/NormaleIndizes}

%% file: Haupt/RechnungGBA/ZiiSiill.tex
\subsection{Reduktionswege für $\Erzeuger{p}{1}\Erzeuger{1}{2}\Erzeuger{2}{2}\Erzeuger{2}{z}$}
\label{ZiiSiill}

\input{Haupt/RechnungGBA/RechnungZiiSiill}

\Bem{Symmetrie}{Diese Rechnungen sind nach Vertauschen der Indizes identisch mit denen in \ref{SiiZiill}.}

%% file: Haupt/RechnungGBA/RechnungZiiSiill.tex
\subsubsection{Beginnend mit $(\Erzeuger{p}{1}\Erzeuger{1}{2}\Erzeuger{2}{2}\Erzeuger{2}{z} ~,~  \left( \Zeile_{p,1,2}- \Zeile_{p,1}\Erzeuger{2}{1} \right) \Erzeuger{2}{z})$:}
\Gleichung{
\underline{\Erzeuger{p}{1}\Erzeuger{1}{2}\Erzeuger{2}{2} }\Erzeuger{2}{z}  ~ 
\Pfeil{\widetilde{Z}_{p,1,2} } &  -\left( \mZii{p}{1}{2}\right) \Erzeuger{2}{z} 
\\&+ \left( \mZi{p}{1} \right) \Erzeuger{2}{1} \Erzeuger{2}{z}
\\=& - \Sum{i=3} \Erzeuger{p}{1}\underline{\Erzeuger{1}{i}\Erzeuger{2}{i}\Erzeuger{2}{z}} 
\\&+ \Sum{i=2}\Erzeuger{p}{i}\underline{\Erzeuger{1}{i}\Erzeuger{2}{1}\Erzeuger{2}{z}} - \DELTA_{p,1}\Erzeuger{2}{1}\Erzeuger{2}{z}
\\ 
\Pfeil{\widetilde{S}_{i,i,z}, \widetilde{S}_{i,1,z} } 
	& - \Sum{i=3} \Erzeuger{p}{1}\left( \Spalte_{i,i,z}- \Spalte_{i,i}\Erzeuger{1}{z}\right)  
\\&+ \Sum{i=2}\Erzeuger{p}{i}\left( \Spalte_{i,1,z}- \Spalte_{i,1}\Erzeuger{1}{z}\right) - \DELTA_{p,1}\Erzeuger{2}{1}\Erzeuger{2}{z}
}%Geleichung
%\newpage
%
%
%
\subsubsection{Beginnend mit $(\Erzeuger{p}{1}\Erzeuger{1}{2}\Erzeuger{2}{2}\Erzeuger{2}{z}  ~,~ \Erzeuger{p}{1}\left( \Spalte_{2,2,z}-\Spalte_{2,2}\Erzeuger{1}{z} \right))$:}
\Gleichung{
&\Erzeuger{p}{1}\underline{\Erzeuger{1}{2}\Erzeuger{2}{2}\Erzeuger{2}{z} }  ~
\\ \Pfeil{ \widetilde{S}_{2,2,z}} &	-\Erzeuger{p}{1} \left( \mSjj{2}{2}{z} \right) 
\\&+ \Erzeuger{p}{1}\left( \mSj{2}{2} \right)\Erzeuger{1}{z}
\\=& -\Sum{j=3}\underline{ \Erzeuger{p}{1} \Erzeuger{1}{2}\Erzeuger{j}{2}} \Erzeuger{j}{z} + \Erzeuger{p}{1}\Erzeuger{1}{2} \DELTA_{2,z} 
\\&+ \Sum{j=2} \underline{\Erzeuger{p}{1} \Erzeuger{j}{2} \Erzeuger{j}{2}} \Erzeuger{1}{z} - \Erzeuger{p}{1} \Erzeuger{1}{z}
\displaybreak\\ \Pfeil{ \widetilde{Z}_{p,1,j},\widetilde{Z}_{p,j,j} } &
-\Sum{j=3}\left( \Zeile_{p,1,j}-\Zeile_{p,1}\Erzeuger{j}{1}\right) \Erzeuger{j}{z} -( 1- \DELTA_{2,z}) \Erzeuger{p}{1} \Erzeuger{1}{z}
\\&+ \Sum{j=2} \left( \Zeile_{p,j,j}-\Zeile_{p,j}\Erzeuger{j}{1}\right)\Erzeuger{1}{z} 
\\\stackrel{\ref{Rechenregel}}{=}& - \left( \Sum{i=3}\Erzeuger{p}{1}\Spalte_{i,i,z} - (1-\DELTA_{1,z}-\DELTA_{2,z})\Erzeuger{p}{1}\Erzeuger{1}{z} \right) 
\\& + \left( \Sum{i=2}\Erzeuger{p}{i} \Spalte_{i,1,z} - \Sum{i=2}\Erzeuger{p}{i}\Erzeuger{1}{i}\DELTA_{1,z}  + \Sum{j=3}\DELTA_{p,1}\Erzeuger{j}{1}\Erzeuger{j}{z} \right)
\\& -( 1- \DELTA_{2,z}) \Erzeuger{p}{1} \Erzeuger{1}{z}+  \left(\Sum{i=3} \Erzeuger{p}{1}\Spalte_{i,i} + \Erzeuger{p}{1}\right)\Erzeuger{1}{z}
\\& - \left( \Sum{i=2} \Erzeuger{p}{i} \Spalte_{i,1} + \left( \DELTA_{1,1} -\DELTA_{p,1}\right) \Erzeuger{p}{1}\right) \Erzeuger{1}{z}
\\=& - \Sum{i=3} \Erzeuger{p}{1}\left( \Spalte_{i,i,z}- \Spalte_{i,i}\Erzeuger{1}{z}\right)  + \Sum{i=2}\Erzeuger{p}{i}\left( \Spalte_{i,1,z}- \Spalte_{i,1}\Erzeuger{1}{z}\right) 
\\& -\DELTA_{p,1}\Erzeuger{2}{1}\Erzeuger{2}{z}
	-\DELTA_{1,z}%\underbrace{
			\left(\underline{\Sum{i=2}\Erzeuger{p}{i}\Erzeuger{1}{i}+\Erzeuger{p}{1}\Erzeuger{1}{1}}\right)
			%}_{\Pfeil{\widetilde{Z}_{p,1}}0} 
\\& + \DELTA_{p,1}%\underbrace{
			\left(\Sum{j=2}\Erzeuger{j}{1}\Erzeuger{j}{z} +\underline{\Erzeuger{1}{1} \Erzeuger{1}{z}}
 		\right)
 		%}_{\Pfeil{\widetilde{Z}_{1,z}}0} 
%
\\\Pfeil{\widetilde{\Zeile}_{p,1},\widetilde{\Zeile}_{1,z}}
& - \Sum{i=3} \Erzeuger{p}{1}\left( \Spalte_{i,i,z}- \Spalte_{i,i}\Erzeuger{1}{z}\right)  + \Sum{i=2}\Erzeuger{p}{i}\left( \Spalte_{i,1,z}- \Spalte_{i,1}\Erzeuger{1}{z}\right) 
\\& -\DELTA_{p,1}\Erzeuger{2}{1}\Erzeuger{2}{z}
	-\DELTA_{1,z}%\underbrace{
			\left(\Sum{i=2}\Erzeuger{p}{i}\Erzeuger{1}{i}+\Zeile_{p,1}\right)
			%}_{\Pfeil{\widetilde{Z}_{p,1}}0} 
\\& + \DELTA_{p,1}%\underbrace{
			\left(\Sum{j=2}\Erzeuger{j}{1}\Erzeuger{j}{z} +\Zeile_{1,z}
 		\right)
 		%}_{\Pfeil{\widetilde{Z}_{1,z}}0} 
\\=	& - \Sum{i=3} \Erzeuger{p}{1}\left( \Spalte_{i,i,z}- \Spalte_{i,i}\Erzeuger{1}{z}\right)  
\\&+ \Sum{i=2}\Erzeuger{p}{i}\left( \Spalte_{i,1,z}- \Spalte_{i,1}\Erzeuger{1}{z}\right) - \DELTA_{p,1}\Erzeuger{2}{1}\Erzeuger{2}{z}
}%Gleichung

%% file: Haupt/RechnungGBA/SiiZiill.tex
\subsection{Reduktionswege für $\Erzeuger{1}{p}\Erzeuger{2}{1}\Erzeuger{2}{2}\Erzeuger{z}{2}$}
\label{SiiZiill}

\input{StyleFiles/VertauschteIndizes}
\input{Haupt/RechnungGBA/RechnungZiiSiill}
\Bem{Symmetrie}{Diese Rechnungen sind nach Vertauschen der Indizes identisch mit denen in \ref{ZiiSiill}.}

\input{StyleFiles/NormaleIndizes}

%% file: Haupt/RechnungGBA/ZiiSiilll.tex
\subsection{Reduktionswege für $\Erzeuger{1}{1}\Erzeuger{2}{2}\Erzeuger{2}{2} $}
\label{ZiiSiilll}
\label{SiiZiilll} %Die beiden Rechungen sind identisch.

\subsubsection{Beginnend mit $( \Erzeuger{1}{1}\Erzeuger{2}{2}\Erzeuger{2}{2} ~,~ \Zeile_{1,2,2}-\Zeile_{1,2}\Erzeuger{2}{1} )$:}
\Gleichung{
&\underline{ \Erzeuger{1}{1}\Erzeuger{2}{2}\Erzeuger{2}{2} }  ~
\\\Pfeil{\widetilde{Z}_{1,2,2} } & \Zeile_{1,2,2}-\Zeile_{1,2}\Erzeuger{2}{1}
\\=&  -\left(\mZii{1}{2}{2}\right)
\\&+\left(\mZi{1}{2}\right)\Erzeuger{2}{1}
\\=&  -\Sum{i=3}\underline{\Erzeuger{1}{1}\Erzeuger{2}{i}\Erzeuger{2}{i}} + \Erzeuger{1}{1} - \Sum{i=2}\underline{\Erzeuger{1}{i}\Erzeuger{2}{i}\Erzeuger{2}{1}}
\\ \Pfeil{\widetilde{S}_{1,i,i}, \widetilde{S}_{i,i,1} } & -\Sum{i=3} \left( \Spalte_{1,i,i} - \Spalte_{1,i}\Erzeuger{1}{i} \right)+ \Erzeuger{1}{1} - \Sum{i=2}\left( \Spalte_{i,i,1} - \Spalte_{i,i}\Erzeuger{1}{1}\right)
\\=& -\Sum{i=3} \left( \mSjj{1}{i}{i}\right) 
\\& -\Sum{i=3} \left(  \mSj{1}{i}\right)\Erzeuger{1}{i} + \Erzeuger{1}{1} 
\\& - \Sum{i=2}\left( \mSjj{i}{i}{1}\right)
\\& + \Sum{i=2}\left(\mSj{i}{i}\right)\Erzeuger{1}{1}
\\=&-\Sum{i=3} \left( \mSjj{1}{i}{i}\right) -\Sum{i=3} \Sum{j=2}\Erzeuger{j}{1}\Erzeuger{j}{i}\Erzeuger{1}{i}
\\& + \Erzeuger{1}{1} - \Sum{i=2}\Sum{j=3}\Erzeuger{1}{i}\Erzeuger{j}{i} \Erzeuger{j}{1}
\\&+ \Sum{i=2}\left(\mSj{i}{i}\right)\Erzeuger{1}{1}
}%Geleichung
\newpage

\subsubsection{Beginnend mit $(\Erzeuger{1}{1}\Erzeuger{2}{2}\Erzeuger{2}{2}   ~,~ \Spalte_{1,2,2}-\Spalte_{1,2}\Erzeuger{1}{2} )$:}
\Gleichung{
&\underline{ \Erzeuger{1}{1}\Erzeuger{2}{2}\Erzeuger{2}{2} }  ~ 
\\ \Pfeil{\widetilde{S}_{1,2,2} } & \Spalte_{1,2,2}-\Spalte_{1,2}\Erzeuger{1}{2} 
\\=&  -\left(\mSjj{1}{2}{2}\right)
\\&+ \left(\mSj{1}{2}\right)\Erzeuger{1}{2} 
\\=& -\Sum{j=3}\underline{\Erzeuger{1}{1}\Erzeuger{j}{2}\Erzeuger{j}{2}} + \Erzeuger{1}{1} + \Sum{j=2} \underline{\Erzeuger{j}{1}\Erzeuger{j}{2}\Erzeuger{1}{2} }
\\ \Pfeil{\widetilde{Z}_{1,j,j}, \widetilde{Z}_{j,j,1} }& -\Sum{j=3} \left(\Zeile_{1,j,j} - \Zeile_{1,j}\Erzeuger{j}{1}\right) + \Erzeuger{1}{1} 
+ \Sum{j=2} \left(\Zeile_{j,j,1} - \Zeile_{j,j}\Erzeuger{1}{1}\right)
\\ = &  -\Sum{j=3} \Zeile_{1,j,j} - \Sum{j=3}\left(\mZi{1}{j}\right)\Erzeuger{j}{1}
\\& + \Erzeuger{1}{1} -\Sum{j=2} \left( \mZii{j}{j}{1}\right)  - \Sum{j=2} \Zeile_{j,j} \Erzeuger{1}{1}
\\=&-\Sum{i=3} \left( \mSjj{1}{i}{i}\right) -\Sum{i=3} \Sum{j=2}\Erzeuger{j}{1}\Erzeuger{j}{i}\Erzeuger{1}{i}
\\& + \Erzeuger{1}{1} - \Sum{i=2}\Sum{j=3}\Erzeuger{1}{i}\Erzeuger{j}{i} \Erzeuger{j}{1}
\\&+ \Sum{i=2}\left(\mSj{i}{i}\right)\Erzeuger{1}{1}
}%Gleichung

%% file: Haupt/EDA-Beweis.tex
\section{Lineare Basis für $\fQG(2)$}
\label{sec:LineareBasisFürN2}

In diesem Abschnitt werden wir für den Fall $\fQG(2)$ einen Automaten angeben, der prüft, ob ein Wort ein Basiselement ist.
Nach Bemerkung \ref{BemLineareBasis} bilden die Wörter, die sich nicht reduzieren lassen, eine Basis. 

Dazu betrachten wir zunächst die folgende Tabelle. Sie enthält die linke Seite aller Reduzierungsregeln aus $\red_{\fQG(n)}$ aus Satz \ref{GBA} für $n=2$, also die Teilwörter, die reduziert werden können.

\mathe{
\begin{array}{|l|c|c|c|c|} \hline
\Matrix{~\\~\\~}&\Matrix{\widetilde{Z}_{p,q}\\a_{p,1}a_{q,1}} & \Matrix{\widetilde{S}_{p,q}\\a_{1,p}a_{1,q}} & \Matrix{\widetilde{Z}_{p,q,r}\\a_{p,1}a_{q,2}a_{r,2}} & \Matrix{\widetilde{S}_{p,q,r}\\a_{1,p}a_{2,q}a_{2,r}}\\ 
\hline
\Matrix{p=1,q=1, r=1} &\Matrix{~\\~}a_{1,1}a_{1,1}& {a_{1,1}a_{1,1}} & a_{1,1}a_{1,2}a_{1,2}& a_{1,1}a_{2,1}a_{2,1}\\
\hline
\Matrix{p=1,q=2, r=1} &\Matrix{~\\~} a_{1,1}a_{2,1}& a_{1,1}a_{1,2} & a_{1,1}a_{2,2}a_{1,2}& a_{1,1}a_{2,2}a_{2,1}\\
\hline
\Matrix{p=2,q=1,r=1} & \Matrix{~\\~}a_{2,1}a_{1,1}& a_{1,2}a_{1,1} & a_{2,1}a_{1,2}a_{1,2}& a_{1,2}a_{2,1}a_{2,1}\\
\hline
\Matrix{p=2,q=2, r=1} &\Matrix{~\\~} a_{2,1}a_{2,1}& a_{1,2}a_{1,2} & a_{2,1}a_{2,2}a_{1,2}& a_{1,2}a_{2,2}a_{2,1}\\
\hline
\Matrix{p=1,q=1, r=2} &\Matrix{~\\~}&															 & a_{1,1}a_{1,2}a_{2,2}& a_{1,1}a_{2,1}a_{2,2}\\
\hline
\Matrix{p=1,q=2, r=2} & \Matrix{~\\~}&															 & a_{1,1}a_{2,2}a_{2,2}& a_{1,1}a_{2,2}a_{2,2}\\
\hline
\Matrix{p=2,q=1, r=2} & \Matrix{~\\~}&															 & a_{2,1}a_{1,2}a_{2,2}& a_{1,2}a_{2,1}a_{2,2}\\
\hline
\Matrix{p=2,q=2, r=2} & \Matrix{~\\~}&															 & a_{2,1}a_{2,2}a_{2,2}& a_{1,2}a_{2,2}a_{2,2}\\
\hline
\end{array}
}

Um zu prüfen ob ein Wort $w$ unreduzierbar ist, kann man wie folgt vorgehen. Zuerst überprüft man, ob es eine Regel in der obigen Tabelle gibt, die mit dem ersten Buchstaben von $w$ beginnt. Da die Regeln höchstens Wörter aus drei Buchstaben reduzieren, reicht es, nur die nächsten zwei Buchstaben zu betrachten. Falls es keine entsprechende Regel gibt, wiederholt man das Verfahren vom nächsten Buchstaben aus, solange bis man am Wortende angelangt ist.

 %In einem Unreduzierbaren Wort, darf beispielsweise auf ein $a_{1,1}$ nur ein $a_{2,2}a_{1,1}$ folgen, da alle anderen Teilworte nach obiger Tabelle reduziert werden können. 

\goodbreak

Man kann auch andersherum vorgehen. Dazu schreiben wir alle unreduzierbaren Teilwörter, die aus drei Buchstaben bestehen, in eine neue Tabelle. Sie erhält man aus obiger Tabelle indem man in jede Spalte die Wörter schreibt, die mit dem entsprechenden Buchstaben beginnen, und nicht in der obigen Tabelle vorkommen.

\subsubsection{
Tabelle der unreduzierbaren Wörter}
\mathe{
\begin{array}{|l|c|c|c|c|} \hline
				& a_{1,1} & a_{1,2} & a_{2,1} & a_{2,2}\\ 
\hline
\Matrix{~\\~\\~\\~} 
		a_{1,1}	& 		&					& 				& \Matrix{a_{2,2}a_{1,1}\bs{a_{2,2}}}
\\ 
\hline
\Matrix{~\\~\\~\\~} 	
	a_{1,2}	&					& 				& 
a_{2,1}a_{1,2}\bs{a_{2,1}} & \Matrix{a_{2,2}a_{1,2}\bs{a_{2,1}}\\a_{2,2}a_{1,2}\bs{a_{2,2}}}
\\ 
\hline
\Matrix{~\\~\\~\\~} 	
	a_{2,1}&	& 
a_{1,2}a_{2,1}\bs{a_{1,2}}	& & \Matrix{a_{2,2}a_{2,1}\bs{a_{1,2}}\\a_{2,2}a_{2,1}\bs{a_{2,2}}}
\\ 
\hline
%\Matrix{~\\~\\~\\~\\~} 	
	a_{2,2}&  a_{1,1}a_{2,2}\bs{a_{1,1}}		& \Matrix{a_{1,2}a_{2,2}\bs{a_{1,1}}\\a_{1,2}a_{2,2}\bs{a_{1,2}}} &  \Matrix{a_{2,1}a_{2,2}\bs{a_{1,1}}\\a_{2,1}a_{2,2}\bs{a_{2,1}}} & \Matrix{a_{2,2}a_{2,2}\bs{a_{1,1}}\\a_{2,2}a_{2,2}\bs{a_{1,2}}\\a_{2,2}a_{2,2}\bs{a_{2,1}}\\a_{2,2}a_{2,2}\bs{a_{2,2}}}
\\ 
\hline
\end{array}
}

~\\
~\\
~\\
Um festzustellen, ob ein Wort unreduzierbar ist, können wir nun vom ersten Buchstaben aus die nächsten zwei Buchstaben betrachten und überprüfen, ob das entsprechende Teilwort in der neuen Tabelle vorkommt. Falls es vorkommt wiederholen wir das Verfahren vom nächsten Buchstaben aus, solange bis wir am Wortende angelangt sind.
~\\
~\\
\Bem{$a_{2,2}$}{
\label{BemAii}
Es gibt in der ersten Tabelle kein Wort, das mit $a_{2,2}$ beginnt. Es können also am Beginn eines unreduzierbaren Wortes beliebig viele $a_{2,2}$ stehen. Daher können wir die letzte Spalte in der neuen Tabelle ignorieren.
}
~\\
~\\

Wir wollen nun diese Verfahren durch einen endlichen Automaten abbilden. 

Für ein Wort beginnt man in einem Startknoten und wandert jeweils entsprechend des nächsten Buchstaben über Pfeile zu einem neuen Knoten. Falls dieser Knoten nicht der Knoten "`reduzierbar"' ist, so ist das Wort unreduzierbar.

%Bevor wir den Automaten konstruieren bemerken wir, dass es  in der ersten Tabelle kein Wort gibt, dass mit einem $a_{2,2}$ beginnt. Es können also am Beginn eines unreduzierbaren Wortes beliebig viele $a_{2,2}$ stehen. Wir können also die gesamte letzte Spalte der Tabelle der unreduzierbaren 

%Im folgenden  können also führende $a_{2,2}$ stets weglassen und einen Pfeil "`$a_{2,2}$"' vom Startknoten zum %Startknoten hinzufügen.
\newpage
\subsubsection{Konstruktion}
\label{sec:Konstruktion}

Wir starten mit einem Startknoten $S$ und einem Knoten "`reduzierbar"', sowie mit einem Knoten für jedes nicht leere Feld der Tabelle, benannt nach der Spalte und der Zeile. Außerdem fügen wir Knoten für jeden Buchstaben hinzu, sowie entsprechende Pfeile vom Startknoten ausgehend.

Wegen Bemerkung \ref{BemAii} können wir jeweils folgende Knoten zusammenfassen:
\AufzaehlungP{
\item "`$a_{2,2}$"', "`$a_{2,2}a_{2,2}$"' und "`$S$"',
\item "`$a_{1,1}$"' und "`$a_{2,2}a_{1,1}$"', 
\item "`$a_{1,2}$"' und "`$a_{2,2}a_{1,2}$"',
\item "`$a_{2,1}$"' und "`$a_{2,2}a_{2,1}$"'.
}
Für jeden Knoten, der einem Feld in der Tabelle entspricht, fügen wir für jeden Eintrag  einen Pfeil, benannt nach dem fett gedrucken Buchstaben, zum Knoten der letzten zwei Buchstaben des Eintrages ein. 
Als nächstes fügen wir von jedem Knoten für alle fehlenden Buchstaben Pfeile zum Knoten "`unreduzierbar"' hinzu.

Sei beispielsweise der erste Buchstabe ein $a_{1,1}$, dann gibt es gemäß der Tabelle der unreduzierbaren Wörter nur zwei weiterführende Knoten "`$a_{1,1}$"', falls ein $a_{2,2}$ folgt und "`reduzierbar"' sonst.
\begin{center}$
\entrymodifiers={+++[O][F-]}
\SelectTips{cm}{}
\xymatrix  {
 ~~~~~S~~~~~% \ar@(ur,ul) _{a_{2,2}}
		\ar[dd] _{a_{1,1}}
\\
*\txt{}
\\
 a_{1,1} \ar  `dr[dd] [dd]  ^>>>>>>>>>>{a_{2,2}}
 \ar@{->}[rrr]^>>>>>>>>>>>>>>{a_{1,1},a_{1,2},a_{2,1}}
& *\txt{}
	&*\txt{}
	&\txt{reduzierbar}
\\
*\txt{}
\\
 a_{1,1}a_{2,2}%\ar `ul[uu]^>{a_{1,1}} [uu]  
}
$
\end{center}

~\\
~\\

Alle Pfeile aus dem Knoten "`reduzierbar"' führen wieder in diesen Knoten zurück. Da wir uns nur für die unreduzierbaren Wörter interessieren, lassen wir den Knoten "`reduzierbar"' und die zu ihm führenden Pfeile im Folgenden weg und interpretieren fehlende Pfeile als Pfeile zu diesem Knoten.

~\\
~\\
So erhalten wir aus der Tabelle der unreduzierbaren Wörter den auf folgender Seite stehenden endlichen Automaten.

\begin{landscape}

\begin{center}$
\entrymodifiers={+++[O][F-]}
\SelectTips{cm}{}
\xymatrix{
*\txt{} 
& *\txt{}
& *\txt{}
& \txt{\tiny{$a_{1,2}a_{2,1}$} }
	\ar `d^r[rr] `^u[rr] _{a_{1,2}} [rr] 
& *\txt{}
& \txt{\tiny{$a_{2,1}a_{1,2}$} }
	\ar `u^l[ll] `^d[ll] _{a_{2,1}} [ll] 
\\
*\txt{}
\\
*\txt{}
\\
*\txt{}
\\ \txt{\tiny{$a_{1,2}a_{2,2}$}}
%\ar@/_/(dr,dl)  [r]
\ar `d^r[r] `^u[r] _{a_{1,2}} [r] 
\ar `dl^r[ddr]  [ddrrrr] ^{a_{1,1}}
& \txt{\tiny{$~~~~~~a_{1,2}~~~~~~$}} %\\ \txt{\tiny{$a_{2,2}a_{1,2}$}} }
	\ar `u^l[l] `^d[l] _{a_{2,2}} [l] 
	\ar[uuuurr] _{a_{2,1}}
& *\txt{}
& *\txt{}
& S \ar@(ur,ul) _{a_{2,2}}
		\ar[lll] _{a_{1,2}}
		\ar[rrr] ^{a_{2,1}}
		\ar[dd] _{a_{1,1}}
& *\txt{}
& *\txt{}
& \txt{\tiny{$~~~~~~a_{2,1}~~~~~~$}}
	\ar `d^r[r] `^u[r] _{a_{2,2}} [r] 
	\ar[uuuull] ^{a_{1,2}}
	& 
	\txt{\tiny{$a_{2,1}a_{2,2}$}}
	\ar `u^l[l] `^d[l] _{a_{2,1}} [l] 
	\ar `dr_l[ddl]  [ddllll] _{a_{1,1}}
	&*\txt{}
\\
*\txt{}
\\
*\txt{}
&*\txt{}
&*\txt{}
&*\txt{}
& \txt{\tiny{$~~~~~~a_{1,1}~~~~~~$}} \ar  `dr[dd] [dd]  ^>>>>>>>>>>{a_{2,2}}
&*\txt{}
& *\txt{}
	&*\txt{}
		&*\txt{}
			&*\txt{}
				&*\txt{}	&*\txt{}	&*\txt{}
\\
*\txt{}
\\
*\txt{}
&*\txt{}
&*\txt{}
&*\txt{}
& \txt{\tiny{$a_{1,1}a_{2,2}$}}\ar `ul[uu]^>{a_{1,1}} [uu]  
}
$\\~\\ Endlicher Automat zum Verifizieren von Basiselementen
\end{center}
\end{landscape}

%% file: Haupt/BiAufloesung.tex
\chapter{Projektive Auflösungen}
\label{sec:AuflösungVonFQG}
\newpage
\section{Auflösung von $\fQG(n)$ als Bimodul}
\label{sec:AuflösungVonFQGN}

In diesem Kapitel konstruieren wir für $n\geq 3$ %\footnote{Für $n=1$ siehe Bemerkung \ref{sec:SpezialfallN1}}
%, für $n=2$ müssten ähnliche Rechnungen  
eine projektive Auflösung für die  orthogonale freie Quantengruppe $\fQG(n)$ als $\fQG(n)\otimes \fQG^{op}(n)$-Modul.

Dazu sei $\Aev(n) := \fQG(n) \otimes_{\K} \fQG^{op}(n)$, wobei $\fQG^{op}(n)$ die Algebra mit vertauschter Multiplikation ist.
Die Elemente in $\Aev(n)$ sind erzeugt von Elementen der Form $a \otimes b$, wobei $a,b \in \fQG(n)$. Um zu verdeutlichen, dass wir $\Aev(n)$ als $\fQG(n)$-Bimodul betrachten, schreiben wir anstelle von $a \otimes b$ oft auch $a\boldsymbol{e}b$, wobei $\boldsymbol{e}$ der Erzeuger des Moduls sein soll. Im Folgenden schreiben wir die Erzeuger der Moduln stets fett gedruckt.

~\\
Wir wollen nun zeigen, dass die folgende Sequenz exakt ist:

\xymatrix{				 
	\ar[d] 0  
\\
 	\ar@{_{(}->}[d]_{\Phi_3} \Aev(n) 		 				&  \ar@{|->}[d] \bs{f}								&			 															
\\
	\ar[d]_{\Phi_2} \left(\Aev(n)\right)^{n^2} 	&  -\sum\limits_{i,j,k=1}^n a_{j,i}\bs{f}_{j,k} a_{k,i}
																								+ \sum\limits_{i=1}^n \bs{f}_{i,i}	& \ar@{|->}[d] \bs{f}_{p,q}
\\
	\ar[d]_{\Phi_1} \left(\Aev(n)\right)^{n^2}  &  \ar@{|->}[d] \bs{e}_{p,q}	&\sum\limits_{i=1}^n\left(a_{p,i}\bs{e}_{q,i}
																																								+\bs{e}_{p,i}a_{q,i}\right)
\\
	\ar@{->>}[d]_{\Phi_0} \Aev(n) 							& a_{p,q}\bs{e}- \bs{e}a_{p,q}& \ar@{|->}[d] a \otimes b	
\\
	\ar[d] \fQG(n)															&															& ab
 \\
	0 				 
		%					%
							%	
	}
Die Abbildung $\Phi_0$ ist die übliche Multiplikationsabbildung. Mit $\Phi_0(a\otimes 1)=a$ sehen wir, dass sie surjektiv ist. In \cite{MR674652} 
wird gezeigt, dass  für jede Algebra mit $1$ der Kern der Multiplikationsabbildung von $a\otimes 1 - 1 \otimes a$ aufgespannt wird, wobei $a$ über die gesamte Algebra läuft. In Abschnitt 10.2. in \cite{MR674652} wird sogar gezeigt, dass 
\mathe{\Kern{\Phi_0}= \sum_{a \in \textnormal{Erzeuger}} \left(a \otimes 1 - 1 \otimes a\right) \times \Aev(n).} Somit ist auch die Exaktheit an der zweiten Stelle bewiesen. Im Folgenden betrachten wir nur die anderen Stellen.

\newpage
Wenn wir die Erzeuger der Moduln und der Algebra $\fQG(n)$ als generische Matrizen  $\bs{E}:=(\bs{e}_{p,q})_{p,q=1\dots n}$, $\bs{F}:=(\bs{f}_{p,q})_{p,q=1\dots n}$ und $A:=(a_{p,q})_{p,q=1\dots n}$ schreiben, können wir die Abbildungen durch folgende Matrizen beschreiben:

\xymatrix{				 
	\ar[d] 0  
\\
 	\ar@{_{(}->}[d]_{\Phi_3} \Aev(n) 		 				&   \ar@{|->}[d]\bs{f}			&
\\
	\ar[d]_{\Phi_2} \left(\Aev(n)\right)^{n^2}  &   -\tr(A^t\bs{F}A)+\tr(\bs{F})    	&\ar@{|->}[d]\bs{F}
\\
	\ar[d]_{\Phi_1} \left(\Aev(n)\right)^{n^2}  &  \ar@{|->}[d] \bs{E}			&A\bs{E}^t+\bs{E}A^t		
\\
	\ar@{->>}[d]_{\Phi_0} \Aev(n) 							& A\bs{e}-\bs{e}A 							&
\\
	\ar[d] \fQG(n)															& 													&	
 \\
	0 				 
		%					%
							%	
	}
	
%Sequenz horizontal	
%\xymatrix{ 0 	& \ar[l] \fQG(n) 
%							& \ar@{->>}[l]_{\Phi_0} \Aev(n) 
%							& \ar[l]_{\Phi_1} \left(\Aev(n)\right)^{n^2}
%							& \ar[l]_{\Phi_2} \left(\Aev(n)\right)^{n^2}
%							& \ar@{_{(}->}[l]_-{\Phi_3} \Aev 
%							& \ar[l] 0   \\
%							%	
%							&	
%							& Ae-eA
%							& \ar@{|->}[l] E
%							& A^tFA+\tr(F)
%							&  \ar@{|->}[l]f
%							\\
%							%
%							&
%							&
%							&
%							AE^t+EA^t
%							& \ar@{|->}[l]F
%	}

Wir werden  für jede Abbildung mit der in \ref{sec:KernEinesAoAHomomorphismus} vorgestellten Methode für den Kern ein Erzeugersystem $S_{\Kern{\Phi_i}}$ berechnen, das mit dem Bild der Erzeuger unter der Abbildung $\Phi_{i+1}$ übereinstimmt. Dazu wählen wir auf $\Aev(n)$ die Silbenordnung zu den Trennungsbuchstaben 
\Gleichung{
&\bs{e}_{n,n} > \dots >\bs{e}_{n,1}> \bs{e}_{n-1,n}> \dots >\bs{e}_{1,1} \\
\textnormal{bzw. }~~ & \bs{f}_{n,n} > \dots >\bs{f}_{n,1}> \bs{f}_{n-1,n}> \dots  >\bs{f}_{1,1}.
}

Sei $\RSa:=\red_{\fQG(n)}$ wie in Abschnitt \ref{sec:GröbnerbasisFürDieFreieQuantenGruppeAN} definiert.
Um ein Erzeugersystem für den $\Kern{\Phi}$ zu erhalten müssen wir ein Erzeugersystem $S=S_{\Phi} \cup S_{es} \cup S_{\Kern{\Phi}} \cup S_{fs}$ für den Graphen $\Gamma(\Phi)$ angeben, so dass das induzierte Wortersetzungssystem $\red_{\Algebra,\Gamma(\Phi)}$ schwach vollständig ist (vgl. \ref{sec:KernEinesAoAHomomorphismus}).

Für jede Abbildung betrachten wir zunächst eine Übersichtstabelle, aus der hervorgeht, welche einzelnen Beweisschritte gemacht werden müssen. Die Tabelle ist wie folgt aufgebaut:

\AufzaehlungZ{
\item{Zeile:} Die Abbildung in Matrixschreibweise.
\item{Zeile:} Die Abbildung in Komponentenschreibweise.
\item{Zeile:} Die Regeln $\RS_{\Phi}$. Sie erhält man, indem man die Abbildung so umstellt, dass das größte Monom auf der linken Seite steht. Schreibt man sie als Gleichung, so erhält man $S_{\Phi}$.
\item{Zeile:} Die Regeln $\RS_{es}$. Falls $\RS_{es}$ nicht benötigt wird, lassen wir die Zeile leer. Schreibt man sie als Gleichung, so erhält man $S_{es}$. Zusammen mit $\RS_{\Phi}$ sind es die einzigen Regeln, die ein Monom mit einem Erzeuger des Bildmoduls reduzieren, also ist $\red_e:=\RS_{{\Phi}} \cup \RS_{es}$.
 Auf der dahinter angegebenen Seite findet sich ein Beweis, dass $S_{es}\subset <S_{\Phi}>$. 
\item{Zeile:} Die Regeln $\RS_{\Kern{\Phi}}$. 
Schreibt man sie als Gleichung, so erhält man $S_{\Kern{\Phi}}$. Auf der dahinter angegebenen Seite findet sich ein Beweis, dass $S_{\Kern{\Phi}}\subset <S_{\Phi} \cup S_{es}>$. 
\item{Zeile:} Die Regeln $\RS_{fs}$. Falls $\RS_{es}$ nicht benötigt wird, lassen wir die Zeile leer. Zusammen mit $\RS_{\Kern{\Phi}}$ sind es die einzigen Regeln die Monome mit einem Erzeuger des Urbildmoduls reduzieren, also ist $\red_f:=\RS_{\Kern{\Phi}} \cup \RS_{fs}$. Nach Satz \ref{SatzKern} und \ref{HauptsatzKern} ist $\left< \ProjA(x-y) ~|~ (x,y) \in \red_f \right> = \Kern{\Phi}$, falls $\RSa \cup \red_f \cup \red_e$ schwach vollständig ist.
Schreibt man sie als Gleichung, so erhält man $S_{fs}$. Auf der dahinter angegebenen Seite findet sich ein Beweis, dass $S_{fs}\subset  <S_{\Phi} \cup S_{es}>$ und ein Beweis, dass gilt:
\mathe{
\left< \ProjA(x-y) ~|~ (x,y) \in \red_{fs} \right> \subset 
\left< \ProjA(x-y) ~|~ (x,y) \in \red_{\Kern{\Phi}} \right>.
}
Also ist $S_{\Kern{\Phi}}$ ein Erzeugersystem für den Kern.
\item{Zeile:} Eine Tabelle mit allen minimalen Überschneidungen $\RSa\cup \RS_{e} \cup \RS_{f}$, in denen genau ein Erzeuger des Bildmoduls vorkommt. Um Satz \ref{HauptsatzKern} nutzen zu können müssen alle minimalen Überschneidungen, die höchstens einen Erzeuger des Bildmoduls enthalten, betrachtet werden. Die Überschneidungen, in denen kein Erzeuger vorkommt, wurden jedoch schon in Abschnitt \ref{sec:RechnungenGBA} überprüft. Ein Eintrag $w$ in der Tabelle  ist als Wort einer Überschneidung $(w,\T{Zeile},\T{Spalte})$ zu verstehen.
}

\Bem{Computerunterstützung}{
Die Regeln $\RS_{es},\RS_{\Kern{\Phi}}$ und $\RS_{fs}$ wurden mit Computerunterstützung geraten.
Benutzt wurden neben viel selbst geschriebener Software die Programme: Plural \cite{GLS03}, GAP \cite{GAP4} und Magma \cite{Magma1997}.}
\newpage

\input{Haupt/RechnungAufloesung/BiPhi1}

\input{Haupt/RechnungAufloesung/BiPhi2}

\input{Haupt/RechnungAufloesung/BiPhi3}

%	

%% file: Haupt/RechnungAufloesung/BiPhi1.tex
\section{$\Phi_1$}
\label{sec:Phi1}

\subsection{Übersichtstabelle für $\Phi_1$}
\label{sec:ÜbersichtstabelleFürPhi1}

\mathe{\begin{array}{|l|l|l|}
 \hline
						& \LMatrix{~\\ \Phi_1 : \Aev(n)^{n^2} \rightarrow \Aev(n)\\ ~ } 	
\\ \hline
 \LMatrix{~\\ \textnormal{Matrix} \\ ~ } 			& ~~~E \mapsto Ae - eA 
\\ \hline
 \LMatrix{~\\ \textnormal{Komponenten} \\ ~ }	& \bs{e}_{p,q} \mapsto a_{p,q}\bs{e} - \bs{e} a_{p,q}
\\ \hline
 \LMatrix{~\\ \RS_{\Phi_1} \\ ~ } & \LMatrix{
 \widetilde{V}_{p,q}= \left( a_{p,q}\bs{e} ~,~ \bs{e}a_{p,q} + \bs{e}_{p,q} \right)
 }
\\ \hline
&\\ \hline
 \LMatrix{~\\ \RS_{\Kern{\Phi_1}} \\ ~ } & \LMatrix{  \widetilde{K}_{p,q}= \left( a_{p,n}\bs{e}_{q,n} ~,~ -\Sum{i=1}{n-1} a_{p,i} \bs{e}_{q,i}- \Sumn{i=1}  \bs{e}_{p,i}a_{q,i} \right) & \Seite{sec:RegelnRSkernPhi1}}\\
 \hline
 \LMatrix{~\\ \RS_{fs} \\ ~ } & \LMatrix{
  \widetilde{K}'_{p,q}= \left( a_{n,p}\bs{e}_{n,q} ~,~ -\Sum{i=1}{n-1} a_{i,p} \bs{e}_{i,q}- \Sumn{i=1}  \bs{e}_{i,p}a_{i,q} \right) &\Seite{sec:RegelnRSfsPhi1} }
\\ \hline
  \textnormal{Konflikte}  & \LMatrix{~\\ \begin{array}{|l|ll|}\hline

						&\LMatrix{~\\ \widetilde{V}_{p,q}\\ ~ } &
						\\ \hline
						\LMatrix{~\\\widetilde{Z}_{p,q}\\ ~ } & a_{p,1}a_{q,1}\bs{e} &\ref{ZpqV1}
						\\ \hline
						\LMatrix{~\\\widetilde{S}_{p,q}\\ ~ } & a_{1,p}a_{1,q}\bs{e} &\ref{SpqV1}
						\\ \hline
						\LMatrix{~\\\widetilde{Z}_{p,q,r}\\ ~ }& a_{p,1}a_{q,2}a_{r,2}\bs{e} &\ref{ZpqrV1}
						\\ \hline
						\LMatrix{~\\\widetilde{S}_{p,q,r}\\ ~ }	& a_{1,p}a_{2,q}a_{2,r}\bs{e} &\ref{SpqrV1}
						 \\\hline
					  \end{array}\\ ~ }
\\ \hline
\end{array}
}

\newpage
\subsection{$\RS_{\Kern{\Phi_1}}$}
\label{sec:RegelnRSkernPhi1}

Wir werden zeigen, dass $S_{\Kern{\Phi_1}} \subset \left< S_{{\Phi_1}} \right>$. Dazu  schreiben wir $S_{\Phi_1}$ in Matrizenschreibweise:
\mathe{ \bs{E}- \left( A\bs{e} - \bs{e}A\right).}
Wenn wir von rechts mit $A^t$ multiplizieren, erhalten wir wegen$AA^t=\id$:
\mathe{  \bs{E}A^t - A \bs{e}A^t  + \bs{e}.}

Wenn wir zunächst transponieren und dann von links mit $A$ multiplizieren, erhalten wir wegen $A^tA=\id$:
\mathe{ A \bs{E}^t - \bs{e} + A \bs{e} A^t.}

Durch Addieren erhalten wir $\bs{E}A^t - A \bs{E}^t \subset \left< S_{{\Phi_1}} \right>$. 

Da $\bs{E}A^t - A \bs{E}^t$ genau $S_{\Kern{\Phi_1}}$ in Matrixschreibweise ist, folgt die Behauptung.

\subsection{$\RS_{fs}$}
\label{sec:RegelnRSfsPhi1}
Wir werden zeigen,  dass $S_{fs} \subset \left< S_{\Kern{\Phi_1}} \right>$. Zusammen mit dem vorhergehenden Abschnitt folgt dann  auch $S_{fs} \subset \left< S_{{\Phi_1}} \right>$.

Dazu schreiben wir $S_{\Kern{\Phi_1}}$ in Matrizenschreibweise:
\mathe{\bs{E}A^t - A \bs{E}^t.}
Wenn wir von links mit $A^t$ und von rechts mit $A$ multiplizieren, erhalten wir:
\mathe{ A^t \bs{E} -  \bs{E}^t A,}
was  $S_{fs}$ in Matrizenschreibweise entspricht; also ist 
$S_{fs} \subset \left< S_{\Kern{\Phi_1}} \right>$.

 \newpage
\MitRechnungen{
\input{Haupt/RechnungAufloesung/Phi1/Rechenregeln}\newpage

\input{Haupt/RechnungAufloesung/Phi1/Zpq}\newpage

\input{Haupt/RechnungAufloesung/Phi1/Spq}\newpage

\input{Haupt/RechnungAufloesung/Phi1/Zpqr}\newpage

\input{Haupt/RechnungAufloesung/Phi1/Spqr}\newpage
}%\MitRechnungen

%% file: Haupt/RechnungAufloesung/Phi1/Rechenregeln.tex
%Nur lokal für die Rechnungen brauch ich eine andere Summenumgebung 

\subsection{Rechenregeln}
\label{Rechenregeln1}

\Satz{Rechenregeln}{
Es gilt:
\Gleichung{
1. ~ Z_{p,q}\bs{e} \Pfeil{}\dots\Pfeil{} & \bs{e}Z_{p,q}+ \bs{e}_{p,1}a_{q,1}+ a_{p,1}e_{q,1}\\
2. ~ S_{p,q}\bs{e} \Pfeil{}\dots\Pfeil{} & \bs{e}S_{p,q}+ \bs{e}_{1,p}a_{1,q}+ a_{1,p}e_{1,q}\\
3. ~ Z_{p,q,r}\bs{e} \Pfeil{}\dots\Pfeil{} &  \bs{e} Z_{p,q,r } +  \bs{e}_{p,1}\left(Z_{q,r}+ a_{q,2}a_{r,2} \right)+ a_{p,1}\bs{e}_{q,1} a_{r,1} \\
									&+a_{p,1}\bs{e}_{q,2} a_{r,2} +Z_{p,q} \bs{e}_{r,1} + a_{p,1}a_{q,2} \bs{e}_{r,2} \\
4. ~ S_{p,q,r}\bs{e} \Pfeil{}\dots\Pfeil{} & \bs{e} S_{p,q,r } +  \bs{e}_{1,p}\left(S_{q,r}+ a_{2,q}a_{2,r} \right)+ a_{1,p}\bs{e}_{1,q} a_{1,r}\\
									&+a_{1,p}\bs{e}_{2,q} a_{2,r}+S_{p,q} \bs{e}_{1,r} + a_{1,p}a_{2,q} \bs{e}_{2,r} \\
}

\Beweis{}{
Zu 1.:\\
\Gleichung{
Z_{p,q}\bs{e} =& -\Sumn{i=2} a_{p,i} \underline{ a_{q,i}\bs{e} } + \delta_{p,q}\bs{e}\\
	\Pfeil{\widetilde{V}_{q,i} } &-\Sumn{i=2} \left( \underline{a_{p,i} \bs{e}} a_{q,i} + a_{p,i}\bs{e}_{q,i}\right) + \delta_{p,q}\bs{e}\\
	\Pfeil{\widetilde{V}_{p,i} } &-\Sumn{i=2} \left( \bs{e}a_{p,i}a_{q,i} + \bs{e}_{p,i}a_{q,i}
	  + a_{p,i}\bs{e}_{q,i}\right) + \delta_{p,q}\bs{e}\\
	  =~ & \bs{e}Z_{p,q} - \Sumn{i=2}\bs{e}_{p,i}a_{q,i}  - \sum_{i=2}^{n-1}a_{p,i}\bs{e}_{q,i}- \underline{a_{p,n}\bs{e}_{q,n}} \\
	\Pfeil{\widetilde{K}_{p,q} } & \bs{e}Z_{p,q} - \Sumn{i=2}\bs{e}_{p,i}a_{q,i}  - \sum_{i=2}^{n-1}a_{p,i}\bs{e}_{q,i} + \sum_{i=1}^{n-1}a_{p,i}\bs{e}_{q,i} +\Sumn{i=1}\bs{e}_{p,i}a_{q,i}  \\
	=~& \bs{e}Z_{p,q} + \bs{e}_{p,1} a_{q,1} + a_{p,1} \bs{e}_{q,1}
}

\newpage
Zu 2.: Diese Rechnung ist nach Vertauschen der Indizes identisch zu der vorhergehenden:
\Gleichung{
S_{p,q}\bs{e} =&  -\Sumn{i=2} a_{i,p} \underline{ a_{i,q}\bs{e} } + \delta_{p,q}\bs{e}\\
	\Pfeil{\widetilde{V}_{i,q} } &-\Sumn{i=2} \left( \underline{a_{i,p} \bs{e}} a_{i,q} + a_{i,p}\bs{e}_{i,q}\right) + \delta_{p,q}\bs{e}\\
	\Pfeil{\widetilde{V}_{i,p} } &-\Sumn{i=2} \left( \bs{e}a_{i,p}a_{i,q} + \bs{e}_{i,p}a_{i,q}
	  + a_{i,p}\bs{e}_{i,q}\right) + \delta_{p,q}\bs{e}\\
	  =&\bs{e}S_{p,q} - \Sumn{i=2}\bs{e}_{i,p}a_{i,q}  - \sum_{i=2}^{n-1}a_{i,p}\bs{e}_{i,q}- \underline{a_{n,p}\bs{e}_{n,q}} \\
	\Pfeil{\widetilde{K}'_{p,q} } & \bs{e}S_{p,q} - \Sumn{i=2}\bs{e}_{i,p}a_{i,q}  - \sum_{i=2}^{n-1}a_{i,p}\bs{e}_{i,q} + \sum_{i=1}^{n-1}a_{i,p}\bs{e}_{i,q} +\Sumn{i=1}\bs{e}_{i,p}a_{i,q}  \\
	=& \bs{e}S_{p,q} + \bs{e}_{1,p} a_{1,q} + a_{1,p} \bs{e}_{1,q}
}
\newpage
Zu 3.:
\Gleichung{
Z_{p,q,r}\bs{e}=& - \Sumn{i=3} a_{p,1} a_{q,i}\underline{a_{r,i}\bs{e}} + \delta_{q,r}\underline{a_{p,1}\bs{e}} \\
	\Pfeil{\widetilde{V}_{r,i},\widetilde{V}_{p,1} }& 
	- \Sumn{i=3} a_{p,1} \underline{a_{q,i}\bs{e}}a_{r,i}  - \Sumn{i=3} a_{p,1} a_{q,i} \bs{e}_{r,i} 
	+ \delta_{q,r}\left(\bs{e}a_{p,1} +\bs{e}_{p,1} \right) \\
 \Pfeil{\widetilde{V}_{q,i} }& - \Sumn{i=3} \underline{a_{p,1} \bs{e}}a_{q,i}a_{r,i} - \Sumn{i=3} a_{p,1} \bs{e}_{q,i} a_{r,i}  \\ &
  - \sum_{i=3}^{n-1} a_{p,1} a_{q,i} \bs{e}_{r,i} -a_{p,1} \underline{a_{q,n} \bs{e}_{r,n}}
	+ \delta_{q,r}\bs{e}a_{p,1}   +\delta_{q,r}\bs{e}_{p,1}\\ 
	\Pfeil{\widetilde{V}_{p,1}, \widetilde{K}_{q,r} }& 
	- \Sumn{i=3} \bs{e}a_{p,1}a_{q,i}a_{r,i} - \Sumn{i=3} \bs{e}_{p,1} a_{q,i}a_{r,i} - a_{p,1}\Sumn{i=3}  \bs{e}_{q,i} a_{r,i}  \\ &- \sum_{i=3}^{n-1} a_{p,1} a_{q,i} \bs{e}_{r,i}   + a_{p,1}\left(\Sum{i=1}{n-1} a_{q,i} \bs{e}_{r,i}+ \Sumn{i=1}  \bs{e}_{q,i}a_{r,i}\right)\\&
	+ \bs{e}a_{p,1}\delta_{q,r}   +\bs{e}_{p,1}\delta_{q,r}\frac{}{}\\ 
	= ~& \bs{e} Z_{p,q,r } +  \bs{e}_{p,1}\left(Z_{q,r}+ a_{q,2}a_{r,2} \right)+ a_{p,1}\left(\bs{e}_{q,1} a_{r,1}+\bs{e}_{q,2} a_{r,2}\right)\frac{}{} \\&  +\underline{a_{p,1} a_{q,1}} \bs{e}_{r,1}  + a_{p,1}a_{q,2} \bs{e}_{r,2}\frac{}{}\\
	\Pfeil{ \widetilde{Z}_{p,q}  }&  \bs{e} Z_{p,q,r } +  \bs{e}_{p,1}\left(Z_{q,r}+ a_{q,2}a_{r,2} \right)+ a_{p,1}\left(\bs{e}_{q,1} a_{r,1}+\bs{e}_{q,2} a_{r,2}\right)\frac{}{} \\&+Z_{p,q} \bs{e}_{r,1} + a_{p,1}a_{q,2} \bs{e}_{r,2} \frac{}{}
}%Gleichung

\newpage

Zu 4.:\\Diese Rechnung ist nach Vertauschen der Indizes identisch zu der vorhergehenden:
\Gleichung{
S_{p,q,r}\bs{e}=& - \Sumn{i=3} a_{1,p} a_{i,q}\underline{a_{i,r}\bs{e}} + \delta_{q,r}\underline{a_{1,p}\bs{e}} \\
	\Pfeil{\widetilde{V}_{i,r},\widetilde{V}_{1,p} }& 
	- \Sumn{i=3} a_{1,p} \underline{a_{i,q}\bs{e}}a_{i,r}  - \Sumn{i=3} a_{1,p} a_{i,q} \bs{e}_{i,r} 
	+ \delta_{q,r}\left(\bs{e}a_{1,p} +\bs{e}_{1,p} \right) \\
 \Pfeil{\widetilde{V}_{i,q} }& - \Sumn{i=3} \underline{a_{1,p} \bs{e}}a_{i,q}a_{i,r} - \Sumn{i=3} a_{1,p} \bs{e}_{i,q} a_{i,r}  \\ &
  - \sum_{i=3}^{n-1} a_{1,p} a_{i,q} \bs{e}_{i,r} -a_{1,p} \underline{a_{n,q} \bs{e}_{n,r}}
	+ \delta_{q,r}\bs{e}a_{1,p}   +\delta_{q,r}\bs{e}_{1,p}\\ 
	\Pfeil{\widetilde{V}_{1,p}, \widetilde{K}'_{q,r} }& 
	- \Sumn{i=3} \bs{e}a_{1,p}a_{i,q}a_{i,r} - \Sumn{i=3} \bs{e}_{1,p} a_{i,q}a_{i,r} - a_{1,p}\Sumn{i=3}  \bs{e}_{i,q} a_{i,r}  \\ &- \sum_{i=3}^{n-1} a_{1,p} a_{i,q} \bs{e}_{i,r}   + a_{1,p}\left(\sum\limits_{i=1}^{n-1} a_{i,q} \bs{e}_{i,r}+ \sum\limits_{i=1}^n  \bs{e}_{i,q}a_{i,r}\right)\\&
	+ \bs{e}a_{1,p}\delta_{q,r}   +\bs{e}_{1,p}\delta_{q,r}\frac{}{}\\ 
	= ~& \bs{e} S_{p,q,r } +  \bs{e}_{1,p}\left(S_{q,r}+ a_{2,q}a_{2,r} \right)+ a_{1,p}\left(\bs{e}_{1,q} a_{1,r}+\bs{e}_{2,q} a_{2,r}\right)\frac{}{} \\&  +\underline{a_{1,p} a_{1,q}} \bs{e}_{1,r}  + a_{1,p}a_{2,q} \bs{e}_{2,r}\frac{}{}\\
	\Pfeil{ \widetilde{Z}_{p,q}  }&  \bs{e} S_{p,q,r } +  \bs{e}_{1,p}\left(S_{q,r}+ a_{2,q}a_{2,r} \right)+ a_{1,p}\left(\bs{e}_{1,q} a_{1,r}+\bs{e}_{2,q} a_{2,r}\right)\frac{}{} \\&+S_{p,q} \bs{e}_{1,r} + a_{1,p}a_{2,q} \bs{e}_{2,r} \frac{}{}
	}
}%Beweis
}%Satz

%% file: Haupt/RechnungAufloesung/Phi1/Zpq.tex
\subsection{Reduktionswege für $a_{p,1}a_{q,1}\bs{e}$}
\label{ZpqV1}

\subsubsection{Beginnend mit $(a_{p,1}a_{q,1}\bs{e} ~,~ Z_{p,q}\bs{e})$:}

\Gleichung{
\underline{a_{p,1}a_{q,1} }\bs{e} ~ 
	\Pfeil{\widetilde{Z}_{p,q} } & \underline{Z_{p,q}\bs{e}}  \\
	\Pfeil{\ref{Rechenregeln1}} ~& \bs{e}Z_{p,q} + \bs{e}_{p,1} a_{q,1} + a_{p,1} \bs{e}_{q,1}
}%Gleichung

\subsubsection{Beginnend mit $(a_{p,1}a_{q,1}\bs{e} ~,~ a_{p,1}{V}_{q,1})$:}
\Gleichung{
a_{p,1}\underline{a_{q,1}\bs{e} } ~ 
	\Pfeil{\widetilde{V}_{q,1} } & \underline{a_{p,1} \bs{e}} a_{q,1} + a_{p,1} \bs{e}_{q,1}\\
	\Pfeil{\widetilde{V}_{p,1} } & \bs{e} \underline{a_{p,1}a_{q,1}} + \bs{e}_{p,1} a_{q,1} + a_{p,1} \bs{e}_{q,1}\\
	\Pfeil{\widetilde{Z}_{p,q} }&  \bs{e}Z_{p,q} + \bs{e}_{p,1} a_{q,1} + a_{p,1} \bs{e}_{q,1}
}%Geleichung
\Bem{Symmetrie}{Diese Rechnungen sind nach Vertauschen der Indizes identisch mit denen in \ref{SpqV1}.}

%% file: Haupt/RechnungAufloesung/Phi1/Spq.tex
\subsection{Reduktionswege für $a_{1,p}a_{1,q}\bs{e}$}
\label{SpqV1}

\subsubsection{Beginnend mit $(a_{1,p}a_{1,q}\bs{e} ~,~ S_{p,q}\bs{e})$:}

\Gleichung{
\underline{a_{1,p}a_{1,q} }\bs{e} ~ 
	\Pfeil{\widetilde{S}_{p,q} } & \underline{S_{p,q}\bs{e} } \\
	\Pfeil{\ref{Rechenregeln1} } & \bs{e}S_{p,q} + \bs{e}_{1,p} a_{1,q} + a_{1,p} \bs{e}_{1,q}
}%Gleichung

\subsubsection{Beginnend mit $(a_{1,p}a_{1,q}\bs{e} ~,~ a_{1,p}V_{1,q})$:}
\Gleichung{
a_{1,p}\underline{a_{1,q}\bs{e} } ~ 
	\Pfeil{\widetilde{V}_{1,q} } & \underline{a_{1,p} \bs{e}} a_{1,q} + a_{1,p} \bs{e}_{1,q}\\
	\Pfeil{\widetilde{V}_{1,p} } & \bs{e} \underline{a_{1,p}a_{1,q}} + \bs{e}_{1,p} a_{1,q} + a_{1,p} \bs{e}_{1,q}\\
	\Pfeil{\widetilde{S}_{p,q} }&  \bs{e}S_{p,q} + \bs{e}_{1,p} a_{1,q} + a_{1,p} \bs{e}_{1,q}
}%Geleichung

\Bem{Symmetrie}{Diese Rechnungen sind nach Vertauschen der Indizes identisch mit denen in \ref{ZpqV1}.}

%% file: Haupt/RechnungAufloesung/Phi1/Zpqr.tex
\subsection{Reduktionswege für $a_{p,1}a_{q,2}a_{r,2}\bs{e}$}
\label{ZpqrV1}

\subsubsection{Beginnend mit $(a_{p,1}a_{q,2}a_{r,2}\bs{e} ~,~ \left(Z_{p,q,r}-Z_{p,q}a_{r,1}\right)\bs{e})$:}
%
%V
%\bs{e}a_{p,q} + \bs{e}_{p,q}
%\bs{e}a_{p,1} + \bs{e}_{p,1}
%
%
\Gleichung{
&\underline{a_{p,1}a_{q,2}a_{r,2} }\bs{e} ~ \\
\Pfeil{\widetilde{Z}_{p,q,r} } &\underline{Z_{p,q,r}\bs{e}}  - Z_{p,q}  \underline{a_{r,1}\bs{e}}  \\
	\Pfeil{\ref{Rechenregeln1} ,\widetilde{V}_{r,1} }&  \bs{e} Z_{p,q,r } +  \bs{e}_{p,1}\left(Z_{q,r}+ a_{q,2}a_{r,2} \right)+ a_{p,1}\left(\bs{e}_{q,1} a_{r,1}+\bs{e}_{q,2} a_{r,2}\right)
	\\&+Z_{p,q} \bs{e}_{r,1} + a_{p,1}a_{q,2} \bs{e}_{r,2} \frac{}{} - \underline{Z_{p,q} \bs{e}}a_{r,1}  - Z_{p,q} \bs{e}_{r,1}
 \\
 \Pfeil{ \ref{Rechenregeln1} }& 
	 \bs{e} Z_{p,q,r } +  \bs{e}_{p,1}\left(Z_{q,r}+ a_{q,2}a_{r,2} \right)+ a_{p,1}\left(\bs{e}_{q,1} a_{r,1}+\bs{e}_{q,2} a_{r,2}\right) \\&+ a_{p,1}a_{q,2} \bs{e}_{r,2}  - \left(\bs{e}Z_{p,q}+ \bs{e}_{p,1}a_{q,1}+ a_{p,1}e_{q,1} \right) a_{r,1} \frac{}{} \\
=~\frac{\frac{}{}}{}&	\frac{}{} \bs{e} Z_{p,q,r } +  \bs{e}_{p,1}\left(Z_{q,r} - \underline{a_{q,1}a_{r,1}} + a_{q,2}a_{r,2} \right) + a_{p,1}\bs{e}_{q,2} a_{r,2} \\
	 & + a_{p,1}a_{q,2} \bs{e}_{r,2} \frac{}{} - \bs{e}Z_{p,q} a_{r,1} \frac{}{}\\ 
\Pfeil{\widetilde{Z}_{q,r} }&	 \bs{e} Z_{p,q,r } +  \bs{e}_{p,1}\left(Z_{q,r} - Z_{q,r} + a_{q,2}a_{r,2} \right) + a_{p,1}\bs{e}_{q,2} a_{r,2}\frac{}{}\\
	 &   + a_{p,1}a_{q,2} \bs{e}_{r,2} \frac{}{} - \bs{e}Z_{p,q} a_{r,1} 
	}%Gleichung
\subsubsection{Beginnend mit $(a_{p,1}a_{q,2}a_{r,2}\bs{e} ~,~ a_{p,1}a_{q,2}{V}_{r,2})$:}
\Gleichung{
&a_{p,1}a_{q,2}\underline{a_{r,2}\bs{e} } ~ \\
	\Pfeil{\widetilde{V}_{r,2} } & a_{p,1}\underline{ a_{q,2}\bs{e}} a_{r,2} + a_{p,1} a_{q,2}\bs{e}_{r,2}\\
	\Pfeil{\widetilde{V}_{q,2} } & \underline{ a_{p,1}\bs{e}}a_{q,2} a_{r,2} + a_{p,1}\bs{e}_{q,2} a_{r,2}+ a_{p,1} a_{q,2}\bs{e}_{r,2}\\
	\Pfeil{\widetilde{V}_{p,1} } & \bs{e}\underline{a_{p,1}a_{q,2} a_{r,2}}  + \bs{e}_{p,1}a_{q,2} a_{r,2}+ a_{p,1}\bs{e}_{q,2} a_{r,2}+ a_{p,1} a_{q,2}\bs{e}_{r,2}\\
		\Pfeil{\widetilde{Z}_{p,q,r} } & \bs{e}Z_{p,q,r}- \bs{e}Z_{p,q}a_{r,1}  + \bs{e}_{p,1}a_{q,2} a_{r,2}+ a_{p,1}\bs{e}_{q,2} a_{r,2}+ a_{p,1} a_{q,2}\bs{e}_{r,2}
}%Geleichung
\Bem{Symmetrie}{Diese Rechnungen sind nach Vertauschen der Indizes identisch mit denen in \ref{SpqrV1}.}

%% file: Haupt/RechnungAufloesung/Phi1/Spqr.tex
\subsection{Reduktionswege für $a_{1,p}a_{2,q}a_{2,r}\bs{e}$}
\label{SpqrV1}

\subsubsection{Beginnend mit $(a_{1,p}a_{2,q}a_{2,r}\bs{e} ~,~ \left(S_{p,q,r}-S_{p,q}a_{1,r}\right)\bs{e})$:}

%V
%\bs{e}a_{p,q} + \bs{e}_{p,q}
%\bs{e}a_{1,p} + \bs{e}_{1,p}

\Gleichung{
&\underline{a_{1,p}a_{2,q}a_{2,r} }\bs{e} ~ \\
\Pfeil{\widetilde{S}_{p,q,r} } &\underline{S_{p,q,r}\bs{e}}  - S_{p,q}  \underline{a_{1,r}\bs{e}}  \\
	\Pfeil{\ref{Rechenregeln1} ,\widetilde{V}_{1,r} }&  \bs{e} S_{p,q,r } +  \bs{e}_{1,p}\left(S_{q,r}+ a_{2,q}a_{2,r} \right)+ a_{1,p}\left(\bs{e}_{1,q} a_{1,r}+\bs{e}_{2,q} a_{2,r}\right)
	\\&+S_{p,q} \bs{e}_{1,r} + a_{1,p}a_{2,q} \bs{e}_{2,r} \frac{}{} - \underline{S_{p,q} \bs{e}}a_{1,r}  - S_{p,q} \bs{e}_{1,r}
 \\
 \Pfeil{ \ref{Rechenregeln1} }& 
	 \bs{e} S_{p,q,r } +  \bs{e}_{1,p}\left(S_{q,r}+ a_{2,q}a_{2,r} \right)+ a_{1,p}\left(\bs{e}_{1,q} a_{1,r}+\bs{e}_{2,q} a_{2,r}\right) \\& + a_{1,p}a_{2,q} \bs{e}_{2,r}  - \left(\bs{e}S_{p,q}+ \bs{e}_{1,p}a_{1,q}+ a_{1,p}e_{1,q} \right) a_{1,r} \\
=~\frac{\frac{}{}}{}&	 \bs{e} S_{p,q,r } +  \bs{e}_{1,p}\left(S_{q,r} - \underline{a_{1,q}a_{1,r}} + a_{2,q}a_{2,r} \right) + a_{1,p}\bs{e}_{2,q} a_{2,r}\frac{}{}  \\
	 & + a_{1,p}a_{2,q} \bs{e}_{2,r} \frac{}{} - \bs{e}S_{p,q} a_{1,r} \\ 
\Pfeil{\widetilde{S}_{q,r} }&	 \bs{e} S_{p,q,r } +  \bs{e}_{1,p}\left(S_{q,r} - S_{q,r} + a_{2,q}a_{2,r} \right) + a_{1,p}\bs{e}_{2,q} a_{2,r}\frac{}{}\\
	 &   + a_{1,p}a_{2,q} \bs{e}_{2,r} \frac{}{} - \bs{e}S_{p,q} a_{1,r} 
	}%Gleichung

\subsubsection{Beginnend mit $(a_{1,p}a_{2,q}a_{2,r}\bs{e} ~,~ a_{1,p}a_{2,q}{V}_{2,r})$:}
\Gleichung{
&a_{1,p}a_{2,q}\underline{a_{2,r}\bs{e} } ~ \\
	\Pfeil{\widetilde{V}_{2,r} } & a_{1,p}\underline{ a_{2,q}\bs{e}} a_{2,r} + a_{1,p} a_{2,q}\bs{e}_{2,r}\\
	\Pfeil{\widetilde{V}_{2,q} } & \underline{ a_{1,p}\bs{e}}a_{2,q} a_{2,r} + a_{1,p}\bs{e}_{2,q} a_{2,r}+ a_{1,p} a_{2,q}\bs{e}_{2,r}\\
	\Pfeil{\widetilde{V}_{1,p} } & \bs{e}\underline{a_{1,p}a_{2,q} a_{2,r}}  + \bs{e}_{1,p}a_{2,q} a_{2,r}+ a_{1,p}\bs{e}_{2,q} a_{2,r}+ a_{1,p} a_{2,q}\bs{e}_{2,r}\\
		\Pfeil{\widetilde{S}_{p,q,r} } & \bs{e}S_{p,q,r}- \bs{e}S_{p,q}a_{1,r}  + \bs{e}_{1,p}a_{2,q} a_{2,r}+ a_{1,p}\bs{e}_{2,q} a_{2,r}+ a_{1,p} a_{2,q}\bs{e}_{2,r}
}%Geleichung
\Bem{Symmetrie}{Diese Rechnungen sind nach Vertauschen der Indizes identisch mit denen in \ref{ZpqrV1}.}

%% file: Haupt/RechnungAufloesung/BiPhi2.tex
\section{$\Phi_2$}
\label{sec:Phi2}

\subsection{Übersichtstabelle für $\Phi_2$}
\label{sec:ÜbersichtstabelleFürPhi2}

\mathe{\begin{array}{|l|l|l|}
 \hline
						& \LMatrix{~\\ \Phi_2 : \Aev(n)^{n^2} \rightarrow \Aev(n)^{n^2}\\ ~ } 	
\\ \hline
 \LMatrix{~\\ \textnormal{Matrix} \\ ~ } 			& ~~~F \mapsto AE^t+ EA^t 
\\ \hline
 \LMatrix{ \textnormal{Kompo-}\\ \textnormal{nenten}  }	& \bs{f}_{p,q} \mapsto \sum\limits_{i=1}^{n} \left(a_{p,i}\bs{e}_{q,i} + \bs{e}_{p,i} a_{q,i}\right)\LMatrix{~\\~\\}
\\ \hline
 \LMatrix{~\\ \RS_{\Phi_2} \\ ~ } &\LMatrix{
 \widetilde{V}_{p,q}= \left( a_{p,n}\bs{e}_{q,n} ~,~ -  \sum\limits_{i=1}^{n-1}a_{p,i}\bs{e}_{q,i} - \sum\limits_{i=1}^{n}\bs{e}_{p,i} a_{q,i} + \bs{f}_{p,q}\right)
}
\\ \hline
\LMatrix{~\\ \RS_{es} \\ ~ } &\LMatrix{
  \widetilde{W}_{p,q}= \left(a_{n,p}\bs{e}_{n,q}~,~ 
  \LMatrix{~\\~} \right. 
  & -\Sum{j=1}{n-1}a_{j,p}\bs{e}_{j,q} -\Sumn{j=1}\bs{e}_{j,p}a_{j,q}  \\
  &\left. + \Sumn{k,j=1}a_{j,p}\bs{f}_{j,k}a_{k,q}\right) &\ref{sec:RegelnRSesPhi2}
  }%LMatrix
\\ \hline
 \LMatrix{~\\ \RS_{\Kern{\Phi_2}} \\ ~ } & \LMatrix{
 \widetilde{K}= \left( a_{n,1}\bs{f}_{n,n}a_{n,1} ~,~
  \LMatrix{~\\~}\right.& 
 -\Sumn{i,j,k=1\\(i,j,k)\neq (1,n,n)}a_{j,i} \bs{f}_{j,k}a_{k,i} 
  \\
  &\left. + \sum\limits_{i=1}^{n}\bs{f}_{i,i} \right) &\ref{sec:RegelnRSkernPhi2}
 }
\\ \hline
&\\ \hline %r_fs
  \textnormal{Konflikte}  & \LMatrix{~\\ \begin{array}{|l|l|l|}\hline

						&\LMatrix{~\\ \widetilde{V}_{p,q}: a_{p,n}\bs{e}_{q,n}\\ ~ } &\widetilde{W}_{p,q}: a_{n,p}\bs{e}_{n,q}
						\\ \hline
						\LMatrix{~\\\widetilde{Z}_{p,q}:\\ a_{p,1} a_{q,1}   \\ ~ } &	 \textnormal{keine Konflikte}  &\LMatrix{a_{p,1} a_{n,1}\bs{e}_{n,q}	&\ref{ZpqW2}		}
						\\ \hline
						\LMatrix{~\\\widetilde{S}_{p,q}:\\ a_{1,p} a_{1,q} \\ ~ } &\LMatrix{ a_{1,p}a_{1,n}\bs{e}_{q,n} &\ref{SpqV2}}&\textnormal{keine Konflikte }				%
						\\ \hline
						\LMatrix{~\\\widetilde{Z}_{p,q,r}:\\a_{p,1} a_{q,2} a_{r,2}  \\~ }&\LMatrix{ \textnormal{keine Konflikte}  %für}&n\neq2
						}&\LMatrix{a_{p,1} a_{q,2} a_{n,2}\bs{e}_{n,r}&\ref{ZpqrW2} }
						\\ \hline
						\LMatrix{~\\\widetilde{S}_{p,q,r}:\\a_{1,p} a_{2,q} a_{2,r} \\ ~ }	& \LMatrix{a_{1,p}a_{2,q}a_{2,n}\bs{e}_{q,n} &\ref{SpqrV2}} &\LMatrix{\textnormal{keine Konflikte}% für}& n\neq2
						}
						 \\\hline
						 \LMatrix{~\\\widetilde{W}_{p,q}:\\a_{n,p} \bs{e}_{n,q} \\ ~ }	& \LMatrix{a_{n,n}\bs{e}_{n,n} &\ref{WV2}} &
						 \\\hline
					  \end{array}\\ ~ 
					 }%LMatrix
\\ \hline
\end{array}
}%mathe

\newpage
\subsection{$\RS_{es}$}
\label{sec:RegelnRSesPhi2}

Wir werden zeigen, dass $S_{es} \subset \left< S_{{\Phi_2}} \right>$. Dazu  schreiben wir $S_{\Phi_2}$ in Matrizenschreibweise:
\mathe{ \bs{F}-\left(A\bs{E}^t + \bs{E}A^t\right).}
Wenn wir von links mit $A^t$ und von rechts mit $A$ multiplizieren, erhalten wir, da $A^tA=\id$:
\mathe{  A^t \bs{F} A - \bs{E}^t A - A^t \bs{E}.}
Also ist  $A^t \bs{F} A - \bs{E}^t A - A^t \bs{E} \subset \left< S_{{\Phi_2}} \right>$. 

Da $A^t \bs{F} A - \bs{E}^t A - A^t \bs{E}$ gerade $S_{es}$ in Matrixschreibweise ist, folgt die Behauptung.

\subsection{$\RS_{\Kern{\Phi_2}}$}
\label{sec:RegelnRSkernPhi2}

Wir werden zeigen, dass $S_{\Kern{\Phi_2}} \subset \left< S_{{\Phi_2}} \right>$. Dazu  schreiben wir $S_{\Phi_2}$ in Matrizenschreibweise:
\mathe{ \bs{F}-\left( A\bs{E}^t + \bs{E}A^t\right).}

Wenn wir die Spur bilden, erhalten wir:
\mathe{ \tr(\bs{F})-\tr(A\bs{E}^t) - \tr(\bs{E}A^t).}
Also ist  $\tr(\bs{F})-\tr(A\bs{E}^t) - \tr(\bs{E}A^t) \subset \left< S_{{\Phi_2}} \right>$. 

Da $S_{\Kern{\Phi_2}}$ in Matizenschreibweise gerade
\mathe{\tr(\bs{F})-\tr\left(A^tFA\right)} ist, folgt die Behauptung mit folgender Gleichung:
\Gleichung{
\tr\left(AE^t +EA^t\right) 
&= \Sum{i=1}{n}\left(\left(AE^t\right)_{i,i} + \left(EA^t\right)_{i,i} \right) \\
&=  \Sum{i=1}{n} \Sum{j=1}{n} \left(a_{i,j}e_{i,j} + e_{i,j}a_{i,j}\right)\\
&=  \Sum{i=j}{n}\left(\left(E^tA\right)_{j,j} + \left(A^tE\right)_{j,j} \right) \\
&= \tr\left(A^tAE^tA+A^tEA^tA\right)
&=\tr\left(A^tFA\right)
.}
 
\newpage

\MitRechnungen{

\input{Haupt/RechnungAufloesung/Phi2/Rechenregeln}\newpage
\input{Haupt/RechnungAufloesung/Phi2/SpqV2}\newpage
\input{Haupt/RechnungAufloesung/Phi2/SpqrV2}\newpage
\input{Haupt/RechnungAufloesung/Phi2/ZpqW2}\newpage
\input{Haupt/RechnungAufloesung/Phi2/ZpqrW2}\newpage
\input{Haupt/RechnungAufloesung/Phi2/WV2}\newpage
}%MitRechnungen

%% file: Haupt/RechnungAufloesung/Phi2/Rechenregeln.tex
\subsection{Rechenregeln}
\label{Rechenregeln2}

\Satz{Rechenregeln}{
Sei $ m \in \left\{1 \dots n\right\} $, dann gilt:
1.:
\Gleichung{ \Sumn{i=1}a_{p,i}a_{q,i} \mapsto \delta_{p,q} 
}%Gleichung
und 2.:
\Gleichung{& \Sumn{j=m}a_{j,p}\underline{a_{j,n} \bs{e}_{q,n}} \mapsto \dots \mapsto	\\&	-\Sumn{j=m} \sum\limits_{i=1}^{n-1}a_{j,p}a_{j,i}\bs{e}_{q,i} 	
		 + \Sum{j=1}{m-1}\Sumn{i=1}a_{j,p}\bs{e}_{j,i} a_{q,i} 
			+\bs{e}_{q,p} \\&
					- \Sum{j=1}{m-1}a_{j,p}\bs{f}_{j,q}
}%Gleichung

und 3.:
\Gleichung{&\Sumn{i=m}a_{p,i}a_{n,i}\bs{e}_{n,q} \mapsto \dots \mapsto	\\&     -\Sumn{i=m}\Sum{j=1}{n-1}a_{p,i}a_{j,i}\bs{e}_{j,q} 
+\Sum{i=1}{m-1}\Sumn{j=1}a_{p,i}\bs{e}_{j,i} a_{j,q}
  +	\bs{e}_{p,q}\\&
  -	\Sumn{j=1}\bs{f}_{p,j}a_{j,q}
  + \Sumn{i=m}\Sumn{k,j=1}a_{p,i} a_{j,i}\bs{f}_{j,k}a_{k,q}
}%Gleichung

								.
}%Satz

\beweis{
Zu 1.:
\Gleichung{\Sumn{i=1}a_{p,i}a_{q,i} =& \underline{a_{p,1}a_{q,1}} +\Sumn{i=2}a_{p,i}a_{q,i}\\
\Pfeil{\widetilde{Z}_{p,q}}&-\Sumn{i=2}a_{p,i}a_{q,i} + \delta_{p,q} +\Sumn{i=2}a_{p,i}a_{q,i} \\=& \delta_{p,q}
}%Gleichung
\newpage
Zu 2.:
\Gleichung{
&  \Sumn{j=m}a_{j,p}\underline{a_{j,n} \bs{e}_{q,n}} \\
	\Pfeil{\widetilde{V}_{j,q} } & 
	 \Sumn{j=m}a_{j,p}\left( -\sum\limits_{i=1}^{n-1}a_{j,i}\bs{e}_{q,i} - \sum\limits_{i=1}^{n}\bs{e}_{j,i} a_{q,i} + \bs{f}_{j,q}\right) \\
	=& -\Sumn{j=m} \sum\limits_{i=1}^{n-1}a_{j,p}a_{j,i}\bs{e}_{q,i} 
			+\Sumn{j=m} a_{j,p}\bs{f}_{j,q}  \\&
	-  \Sum{i,j=1}{n-1} a_{j,p}\bs{e}_{j,i} a_{q,i}	
	-  \Sum{j=1}{n-1} a_{j,p}\bs{e}_{j,n} a_{q,n}	\\&
	-  \Sum{i=1}{n-1} \underline{a_{n,p}\bs{e}_{n,i}} a_{q,i}	
		+ \Sum{j=1}{m-1}\Sumn{i=1}a_{j,p}\bs{e}_{j,i} a_{q,i} - \underline{a_{n,p}\bs{e}_{n,n}} a_{q,n} \\	
	\Pfeil{\widetilde{W2}_{p,i} ,\widetilde{W2}_{p,n} } &	-\Sumn{j=m} \sum\limits_{i=1}^{n-1}a_{j,p}a_{j,i}\bs{e}_{q,i}
	+\Sumn{j=m} a_{j,p}\bs{f}_{j,q} \\&
	-  \Sum{i,j=1}{n-1} a_{j,p}\bs{e}_{j,i} a_{q,i}	
	-  \Sum{j=1}{n-1} a_{j,p}\bs{e}_{j,n} a_{q,n}	\\&
	-  \Sum{i=1}{n-1} \left( -\Sum{j=1}{n-1}a_{j,p}\bs{e}_{j,i} -\Sumn{j=1}\bs{e}_{j,p}a_{j,i} + \Sumn{k,j=1}a_{j,p}\bs{f}_{j,k}a_{k,i}\right)a_{q,i}		 \\&
	 + \Sum{j=1}{m-1}\Sumn{i=1}a_{j,p}\bs{e}_{j,i} a_{q,i} 
	\\&
	 -\left( -\Sum{j=1}{n-1}a_{j,p}\bs{e}_{j,n} -\Sumn{j=1}\bs{e}_{j,p}a_{j,n} + \Sumn{k,j=1}a_{j,p}\bs{f}_{j,k}a_{k,n}
	 \right)a_{q,n}
	\\
	=& -\Sumn{j=m} \sum\limits_{i=1}^{n-1}a_{j,p}a_{j,i}\bs{e}_{q,i} 	
		 + \Sum{j=1}{m-1}\Sumn{i=1}a_{j,p}\bs{e}_{j,i} a_{q,i} \\&
			+\Sumn{i,j=1}\bs{e}_{j,p}\underline{a_{j,i} a_{q,i}} 
	+\Sumn{j=m} a_{j,p}\bs{f}_{j,q} 	
					- \Sumn{i,k,j=1}a_{j,p}\bs{f}_{j,k}\underline{a_{k,i}a_{q,i}}
					\NeueSeite
	\Pfeil{\ref{Rechenregeln2} }&
			-\Sumn{j=m} \sum\limits_{i=1}^{n-1}a_{j,p}a_{j,i}\bs{e}_{q,i} 	
		 + \Sum{j=1}{m-1}\Sumn{i=1}a_{j,p}\bs{e}_{j,i} a_{q,i} 
			+\Sumn{j=1}\bs{e}_{j,p}\delta_{j,q} \\&
	+\Sumn{j=m} a_{j,p}\bs{f}_{j,q} 	
					- \Sumn{k,j=1}a_{j,p}\bs{f}_{j,k}\delta_{k,q}
					\\
					=
					&
					-\Sumn{j=m} \sum\limits_{i=1}^{n-1}a_{j,p}a_{j,i}\bs{e}_{q,i} 	
		 + \Sum{j=1}{m-1}\Sumn{i=1}a_{j,p}\bs{e}_{j,i} a_{q,i} \\&
			+\bs{e}_{q,p} 
					- \Sum{j=1}{m-1}a_{j,p}\bs{f}_{j,q}
}%Gleichung
zu 3.:
\Gleichung{
& \Sumn{i=m}a_{p,i}\underline{a_{n,i}\bs{e}_{n,q}}
\\
\Pfeil{\widetilde{W_{i,q}}}&  \Sumn{i=m}a_{p,i}\left(
-\Sum{j=1}{n-1}a_{j,i}\bs{e}_{j,q} -\Sumn{j=1}\bs{e}_{j,i}a_{j,q} + \Sumn{k,j=1}a_{j,i}\bs{f}_{j,k}a_{k,q}
\right)
\\
=
&
 	-\Sumn{i=m}\Sum{j=1}{n-1}a_{p,i}a_{j,i}\bs{e}_{j,q} 
  -	\Sum{i=m}{n-1}\Sumn{j=1}a_{p,i}\bs{e}_{j,i}a_{j,q}
   \\
  &
  - \Sumn{j=1}\underline{a_{p,n}\bs{e}_{j,n}}a_{j,q} 
  + \Sumn{i=m}\Sumn{k,j=1}a_{p,i} a_{j,i}\bs{f}_{j,k}a_{k,q}
\\
\Pfeil{\widetilde{V_{p,j}}}&
 -\Sumn{i=m}\Sum{j=1}{n-1}a_{p,i}a_{j,i}\bs{e}_{j,q} 
  -	\Sum{i=m}{n-1}\Sumn{j=1}a_{p,i}\bs{e}_{j,i}a_{j,q} 
  \\&
  -	\Sumn{j=1}\left(-  \sum\limits_{i=1}^{n-1}a_{p,i}\bs{e}_{j,i} - \sum\limits_{i=1}^{n}\bs{e}_{p,i} a_{j,i} + \bs{f}_{p,j}\right)a_{j,q}
   \\&
  + \Sumn{i=m}\Sumn{k,j=1}a_{p,i} a_{j,i}\bs{f}_{j,k}a_{k,q}
  \NeueSeite
  =
  &
   -\Sumn{i=m}\Sum{j=1}{n-1}a_{p,i}a_{j,i}\bs{e}_{j,q} 
  +\Sum{i=1}{m-1}\Sumn{j=1}a_{p,i}\bs{e}_{j,i} a_{j,q}
  +	\sum\limits_{i=1}^{n}\bs{e}_{p,i}\underline{ \Sumn{j=1} a_{j,i}  a_{j,q}}\\&
  -	\Sumn{j=1}\bs{f}_{p,j}a_{j,q}
  + \Sumn{i=m}\Sumn{k,j=1}a_{p,i} a_{j,i}\bs{f}_{j,k}a_{k,q}
  \\
  \Pfeil{\ref{Rechenregeln2}}
    &
   -\Sumn{i=m}\Sum{j=1}{n-1}a_{p,i}a_{j,i}\bs{e}_{j,q} 
  +\Sum{i=1}{m-1}\Sumn{j=1}a_{p,i}\bs{e}_{j,i} a_{j,q}
  +	 \sum\limits_{i=1}^{n}\bs{e}_{p,i}\delta_{i,q}\\&
  -	\Sumn{j=1}\bs{f}_{p,j}a_{j,q}
  + \Sumn{i=m}\Sumn{k,j=1}a_{p,i} a_{j,i}\bs{f}_{j,k}a_{k,q}
  \\
  =
  &
     -\Sumn{i=m}\Sum{j=1}{n-1}a_{p,i}a_{j,i}\bs{e}_{j,q} 
+\Sum{i=1}{m-1}\Sumn{j=1}a_{p,i}\bs{e}_{j,i} a_{j,q}
  +	\bs{e}_{p,q}\\&
  -	\Sumn{j=1}\bs{f}_{p,j}a_{j,q}
  + \Sumn{i=m}\Sumn{k,j=1}a_{p,i} a_{j,i}\bs{f}_{j,k}a_{k,q}
}%Gleichung

%&\Sumn{i=m}a_{p,i}\underline{a_{n,i}\bs{e}_{n,n}a_{q,n}}\\
%\Pfeil{\widetilde{W}_{i,q} } & \Sumn{i=m}a_{p,i}\left(
%-\Sumn{k,j=1\\(k,j)\neq (n,n)}a_{j,i}\bs{e}_{j,k}a_{q,k} - \bs{e}_{q,i}+\Sumn{j=1}  a_{j,i}\bs{f}_{j,q}\right)
%\\
%&
%= 				-\Sumn{i=m}\Sumn{k,j=1\\(k,j)\neq (n,n)}a_{p,i}a_{j,i}\bs{e}_{j,k}a_{q,k} 
%					- \Sum{i=m}{n-1}a_{p,i}\bs{e}_{q,i} -\underline{a_{p,n}\bs{e}_{q,n}}\\&
%					+\Sumn{i=m}\Sumn{j=1} a_{p,i} a_{j,i}\bs{f}_{j,q}
%\\
%\Pfeil{\widetilde{V}_{p,q} }	
%				&	-\Sumn{i=m}\Sumn{k,j=1\\(k,j)\neq (n,n)}a_{p,i}a_{j,i}\bs{e}_{j,k}a_{q,k} 
%					- \Sum{i=m}{n-1}a_{p,i}\bs{e}_{q,i}\\
%					&- \left(
%								- \Sum{i=1}{n-1}a_{p,i}\bs{e}_{q,i} 
%								- \Sum{i=1}{n}\bs{e}_{p,i} a_{q,i}
%								+f_{p,q}
%						 \right)\\
%					&			+\Sumn{i=m}\Sumn{j=1}  a_{p,i}a_{j,i}\bs{f}_{j,q}
%\\
%=	&	-\Sumn{i=m}\Sumn{k,j=1\\(k,j)\neq (n,n)}a_{p,i}a_{j,i}\bs{e}_{j,k}a_{q,k} 
%					+ \Sum{i=1}{m-1}a_{p,i}\bs{e}_{q,i}\\
%					&+ \Sum{i=1}{n}\bs{e}_{p,i} a_{q,i}
%								-f_{p,q}			+\Sumn{i=m}\Sumn{j=1}  a_{p,i}a_{j,i}\bs{f}_{j,q}
%}

}%Beweis

%% file: Haupt/RechnungAufloesung/Phi2/SpqV2.tex
%Nur lokal für die Rechnungen brauch ich eine andere Summenumgebung 

\subsection{Reduktionswege für $a_{1,p}a_{1,n}\bs{e}_{q,n}$}
\label{SpqV2}

\subsubsection{Beginnend mit $(a_{1,p}a_{1,n}\bs{e}_{q,n} ~,~ S_{p,n}\bs{e}_{q,n})$:}

\Gleichung{
&\underline{a_{1,p}a_{1,n} }\bs{e}_{q,n} ~ \\
	\Pfeil{\widetilde{S}_{p,n} } &-  \underline{\Sumn{j=2}a_{j,p}a_{j,n} \bs{e}_{q,n}} + \delta_{p,n} \bs{e}_{q,n}\\
	\Pfeil{\ref{Rechenregeln2}} &	 
				\Sumn{j=2} \sum\limits_{i=1}^{n-1}a_{j,p}a_{j,i}\bs{e}_{q,i} 	
		 -\Sumn{i=1}a_{1,p}\bs{e}_{1,i} a_{q,i} 
			-\bs{e}_{q,p} \\&
					+ a_{1,p}\bs{f}_{1,q}
	  + \delta_{p,n} \bs{e}_{q,n}
}%Gleichung

% \LMatrix{~\\ \RSkern \\ ~ } &
% \widetilde{K}_{p,q}= 
% a_{n,n}\bs{f}_{n,n}a_{n,n} ~,~ -\Sumn{i,j,k=1\\(i,j,k)\neq (n,n,n)}{n-1} a_{i,k} \bs{f}_{i,j}a_{j,k} + \sum\limits_{i=1}^{n}\bs{f}_{i,i} 

% \widetilde{V}_{p,q}= \left( 
%a_{p,n}\bs{e}_{q,n} ~,~ - \sum\limits_{i=1}^{n-1}a_{p,i}\bs{e}_{q,i} - \sum\limits_{i=1}^{n}\bs{e}_{p,i} a_{q,i}+f_{p,q}
%
%W_{p,q} : a_{n,p}e_{n,n}a_{q,n} \mapsto 
%-\Sumn{i,j=1\\(i,j)\neq (n,n)}a_{j,p}\bs{e}_{j,i}a_{q,i} - \bs{e}_{q,p}+\Sumn{j=1}  a_{j,p}f_{j,q}
%
\subsubsection{Beginnend mit $(a_{1,p}a_{1,n}\bs{e}_{q,n} ~,~ a_{1,p}V_{1,q})$:}
\Gleichung{
&a_{1,p}\underline{a_{1,n}\bs{e}_{q,n} } ~ \\
	\Pfeil{\widetilde{V}_{1,q} } &   - \sum\limits_{i=1}^{n-1} \underline{a_{1,p}a_{1,i}}\bs{e}_{q,i} - \sum\limits_{i=1}^{n}a_{1,p}\bs{e}_{1,i} a_{q,i} + a_{1,p}\bs{f}_{1,q}\\
 \Pfeil{\widetilde{S}_{p,i} } &  - \sum\limits_{i=1}^{n-1} \left(-\Sumn{j=2}a_{j,p}a_{j,i} +\delta_{p,i}  \right)\bs{e}_{q,i}- \sum\limits_{i=1}^{n}a_{1,p}\bs{e}_{1,i} a_{q,i} + a_{1,p}\bs{f}_{1,q}\\
 =& \sum\limits_{i=1}^{n-1} \Sumn{j=2}a_{j,p}a_{j,i}\bs{e}_{q,i} -(1-\delta_{p,n}) \bs{e}_{q,p}- \sum\limits_{i=1}^{n}a_{1,p}\bs{e}_{1,i} a_{q,i} + a_{1,p}\bs{f}_{1,q}
	}%Geleichung

%\Bem{Symmetrie}{Diese Rechnungen sind nach Vertauschen der Indizes identisch zu denen in \ref{ZpqV1}.}

%% file: Haupt/RechnungAufloesung/Phi2/SpqrV2.tex
%Nur lokal für die Rechnungen brauch ich eine andere Summenumgebung 

\subsection{Reduktionswege für $a_{1,p}a_{2,q}a_{2,n}\bs{e}_{r,n}$}
\label{SpqrV2}

\subsubsection{Beginnend mit $(a_{1,p}a_{2,q}a_{2,n}\bs{e}_{r,n} ~,~ \left(S_{p,q,n}-S_{p,q}a_{1,n}\right)\bs{e}_{r,n})$:}

%-	\Sumn{j=k} \sum\limits_{i=1}^{n-1}a_{j,p}a_{j,i}\bs{e}_{r,i} 			 +  \Sum{j=1}{k-1}\Sumn{i=1}a_{j,p}\bs{e}_{j,i} a_{q,i} % + \bs{e}_{r,p} \\ & + \Sum{j=1}{k-1}a_{j,p}f_{j,q}  %

%-\Sumn{i=m}\Sumn{k,j=1 \\(k,j)\neq (n,n)}a_{p,i}a_{j,i}\bs{e}_{j,k}a_{q,k} 					+ \Sum{i=1}{m-1}a_{p,i}\bs{e}_{q,i}					+ \Sum{i=1}{n}\bs{e}_{p,i} a_{q,i} -f_{p,q}			+\Sumn{i=m}\Sumn{j=1}  a_{p,i}a_{j,i}f_{j,q}

\Gleichung{
&\underline{a_{1,p}a_{2,q}a_{2,n} }\bs{e}_{r,n} ~ \\
	\Pfeil{\widetilde{S}_{p,q,n} } &  
		-a_{1,p}\underline{\Sumn{j=3} a_{j,q}a_{j,n}\bs{e}_{r,n}}		
		+a_{1,p}\delta_{q,n}\bs{e}_{r,n}
		+\left(\Sumn{j=2}a_{j,p}a_{j,q}- \delta_{p,q}\right)\underline{a_{1,n}\bs{e}_{r,n} } 
	\\
	\Pfeil{\ref{Rechenregeln2},\widetilde{V}_{1,r} } &  -a_{1,p}\left( 
  	-\Sumn{j=3} \sum\limits_{i=1}^{n-1}a_{j,q}a_{j,i}\bs{e}_{r,i} 	
		 + \Sum{j=1}{2}\Sumn{i=1}a_{j,q}\bs{e}_{j,i} a_{r,i} 
			+\bs{e}_{r,q} 
					- \Sum{j=1}{2}a_{j,q}\bs{f}_{j,r}
	\right)    \\&
		+a_{1,p}\delta_{q,n}\bs{e}_{r,n}\\&
		+\left(\Sumn{j=2}a_{j,p}a_{j,q}- \delta_{p,q}\right)\left(
		- \sum\limits_{i=1}^{n-1}a_{1,i}\bs{e}_{r,i} - \sum\limits_{i=1}^{n}\bs{e}_{1,i} a_{r,i}+f_{1,r}\right)
		\\
		=&
		\sum\limits_{i=1}^{n-1}\left(
		\Sumn{j=3} a_{1,p}a_{j,q}a_{j,i}
		- \Sumn{j=2}a_{j,p}a_{j,q}a_{1,i} 
		+ \delta_{p,q}a_{1,i}
		\right)\bs{e}_{r,i}\\
		&
			 	-\Sumn{i=1}a_{1,p} a_{2,q}\bs{e}_{2,i} a_{r,i} 
				-a_{1,p}\bs{e}_{r,q} 
%
				%+a_{1,p} a_{1,q}\bs{f}_{1,r}
				+a_{1,p} a_{2,q}\bs{f}_{2,r}
				+\delta_{q,n}a_{1,p}\bs{e}_{r,n}\\&
		- \Sumn{i,j=1}\underline{a_{j,p}a_{j,q}}\bs{e}_{1,i} a_{r,i}
		+\Sumn{j=1}\underline{a_{j,p}a_{j,q}}f_{1,r}
	+ \delta_{p,q} \sum\limits_{i=1}^{n}\bs{e}_{1,i} a_{r,i}
	- \delta_{p,q}f_{1,r}
\\
&
	\Pfeil{\ref{Rechenregeln2}}
				\sum\limits_{i=1}^{n-1}\left(
		\Sumn{j=3} a_{1,p}a_{j,q}a_{j,i}
		- \Sumn{j=2}a_{j,p}a_{j,q}a_{1,i} 
		+ \delta_{p,q}a_{1,i}
		\right)\bs{e}_{r,i}\\
		&
			 	-\Sumn{i=1}a_{1,p} a_{2,q}\bs{e}_{2,i} a_{r,i} 
				-a_{1,p}\bs{e}_{r,q} 
				+a_{1,p} a_{2,q}\bs{f}_{2,r}
				+\delta_{q,n}a_{1,p}\bs{e}_{r,n}\\&
		- \Sumn{i=1}\delta_{p,q}\bs{e}_{1,i} a_{r,i}
		+\delta_{p,q}f_{1,r}
	+ \delta_{p,q} \sum\limits_{i=1}^{n}\bs{e}_{1,i} a_{r,i}
	- \delta_{p,q}f_{1,r}
		%		
	%	%=&	\Sumn{j=3} \sum\limits_{i=1}^{n-1}a_{1,p}a_{j,q}a_{j,i}\bs{e}_{r,i} 
%		- a_{1,p}\bs{e}_{r,q}	+a_{1,p}\delta_{q,n}\bs{e}_{r,n}\\&
%		- \sum\limits_{i=1}^{n-1}\Sumn{j=2}a_{j,p}a_{j,q}a_{1,i}\bs{e}_{r,i}
%		   +\delta_{p,q}\sum\limits_{i=1}^{n-1}a_{1,i}\bs{e}_{r,i}\\&
%		   -\Sumn{i=1}a_{1,p}a_{2,q}\bs{e}_{2,i} a_{r,i}		 -a_{1,p}a_{2,q}f_{2,r}
%		 \\&
%		+ \left(-\underline{a_{1,p}a_{1,q}} - \Sumn{j=2}a_{j,p}a_{j,q}+\delta_{p,q} \right)\Sumn{i=1}\bs{e}_{1,i} a_{r,i}
%		\\
%		& 
%		+\left(- \underline{a_{1,p}a_{1,q}}-\Sumn{j=2}a_{j,p}a_{j,q}+\delta_{p,q}\right)f_{1,r}
%		\\
%		%
%		\Pfeil{\widetilde{S}_{p,q} } 
%		&	\Sumn{j=3} \sum\limits_{i=1}^{n-1}a_{1,p}a_{j,q}a_{j,i}\bs{e}_{r,i} 
%		- a_{1,p}\bs{e}_{r,q}	+a_{1,p}\delta_{q,n}\bs{e}_{r,n}\\&
%		- \sum\limits_{i=1}^{n-1}\Sumn{j=2}a_{j,p}a_{j,q}a_{1,i}\bs{e}_{r,i}
%		   +\delta_{p,q}\sum\limits_{i=1}^{n-1}a_{1,i}\bs{e}_{r,i}\\&
%		   -\Sumn{i=1}a_{1,p}a_{2,q}\bs{e}_{2,i} a_{r,i}		 -a_{1,p}a_{2,q}f_{2,r}
%		 +0+0
}%Gleichung

%\newpage

\subsubsection{Beginnend mit $(a_{1,p}a_{2,q}a_{2,n}\bs{e}_{r,n} ~,~ a_{1,p}a_{2,q}V_{n,q})$:}
\Gleichung{
&a_{1,p}a_{2,q}\underline{a_{2,n}\bs{e}_{r,n} } ~ \\
	\Pfeil{\widetilde{V}_{2,q} } & a_{1,p}a_{2,q}\left(	- \sum\limits_{i=1}^{n-1}a_{2,i}\bs{e}_{r,i} 
					- \sum\limits_{i=1}^{n}\bs{e}_{2,i} a_{q,i}+f_{n,q}\right)\\
	=& 	- \sum\limits_{i=1}^{n-1} \underline{a_{1,p}a_{2,q}a_{2,i}}\bs{e}_{r,i} 
					- \sum\limits_{i=1}^{n} a_{1,p}a_{2,q}\bs{e}_{2,i} a_{q,i}+ a_{1,p}a_{2,q}f_{2,q}\\
 \Pfeil{\widetilde{S}_{p,q,i} } &	
  - \sum\limits_{i=1}^{n-1} \left(-\Sumn{j=3}a_{1,p}a_{j,q}a_{j,i}+a_{1,p}\delta_{q,i}+\Sumn{j=2}a_{j,p}a_{j,q}a_{1,i}-
  \delta_{p,q}a_{1,i}\right)\bs{e}_{r,i} 
  \\
  &
					- \sum\limits_{i=1}^{n} a_{1,p}a_{2,q}\bs{e}_{2,i} a_{r,i}+ a_{1,p}a_{2,q}f_{2,r}
}%Geleichung
\newpage

%\Bem{Symmetrie}{Diese Rechnungen sind nach Vertauschen der Indizes identisch zu denen in \ref{ZpqV1}.}

%% file: Haupt/RechnungAufloesung/Phi2/ZpqW2.tex
%Nur lokal für die Rechnungen brauch ich eine andere Summenumgebung 

\subsection{Reduktionswege für $a_{p,1}a_{n,1}\bs{e}_{n,q}$}
\label{ZpqW2}

\subsubsection{Beginnend mit $(a_{p,1}a_{n,1}\bs{e}_{n,q} ~,~ Z_{p,n}\bs{e}_{n,q})$:}

\Gleichung{
&\underline{a_{p,1}a_{n,1}}\bs{e}_{n,q}  ~ \\
\Pfeil{\widetilde{Z}_{p,n} } & - \underline{\Sumn{i=2}a_{p,i}a_{n,i}\bs{e}_{n,q}}+ \delta_{p,n}\bs{e}_{n,q}
\\
\Pfeil{\widetilde{W_{i,q}}}& 
	+\Sumn{i=2}\Sum{j=1}{n-1}a_{p,i}a_{j,i}\bs{e}_{j,q} 
  -\Sumn{j=1}a_{p,1}\bs{e}_{j,1} a_{j,q}
  -	\bs{e}_{p,q}\\&
  +	\Sumn{j=1}\bs{f}_{p,j}a_{j,q}
  - \Sumn{i=2}\Sumn{k,j=1}a_{p,i} a_{j,i}\bs{f}_{j,k}a_{k,q}
  + \delta_{p,n}\bs{e}_{n,q}
}%Gleichung

\subsubsection{Beginnend mit $(a_{p,1}a_{n,1}\bs{e}_{n,q} ~,~ a_{p,1}W_{1,q})$:}
\Gleichung{
&a_{p,1}\underline{a_{n,1}\bs{e}_{n,q} } ~ \\
	\Pfeil{\widetilde{W}_{1,q} } &  -\Sum{j=1}{n-1} \underline{a_{p,1}a_{j,1}}\bs{e}_{j,q} -\Sumn{j=1} a_{p,1}\bs{e}_{j,1}a_{j,q} 
	\\
	& +\Sumn{k,j=1}\underline{a_{p,1}a_{j,1}}\bs{f}_{j,k}a_{k,q}
\\
	\Pfeil{\widetilde{Z}_{p,j}, \widetilde{Z}_{p,j}} &
	\Sum{j=1}{n-1} \Sumn{i=2}a_{p,i}a_{j,i}\bs{e}_{j,q} -\Sum{j=1}{n-1} \delta_{p,j}\bs{e}_{j,q}
	-\Sumn{j=1} a_{p,1}\bs{e}_{j,1}a_{j,q} 
	\\
	&
	- \Sumn{k,j=1}\Sumn{i=2}a_{p,i}a_{j,i}\bs{f}_{j,k}a_{k,q}+ \Sumn{k,j=1}\delta_{p,j}\bs{f}_{j,k}a_{k,q}
		\\
		=
		&
			\Sum{j=1}{n-1} \Sumn{i=2}a_{p,i}a_{j,i}\bs{e}_{j,q} 
			-(1-\delta_{p,n})\bs{e}_{p,q}
			-\Sumn{j=1} a_{p,1}\bs{e}_{j,1}a_{j,q} 
	\\
	&
			- \Sumn{k,j=1}\Sumn{i=2}a_{p,i}a_{j,i}\bs{f}_{j,k}a_{k,q}
			+ \Sumn{k=1}\bs{f}_{p,k}a_{k,q}
	}%Geleichung

%\Bem{Symmetrie}{Diese Rechnungen sind nach Vertauschen der Indizes identisch zu denen in \ref{ZpqV1}.}

%% file: Haupt/RechnungAufloesung/Phi2/ZpqrW2.tex
%Nur lokal für die Rechnungen brauch ich eine andere Summenumgebung 

\subsection{Reduktionswege für $a_{p,1}a_{q,2}a_{n,2}\bs{e}_{n,r}$}
\label{ZpqrW2}

\subsubsection{Beginnend mit $(a_{n,n}\bs{e}_{n,n} ~,~ Z_{p,q,n}\bs{e}_{n,r}-Z_{p,q}a_{n,1}\bs{e}_{n,r})$:}

\Gleichung{
&\underline{a_{p,1}a_{q,2}a_{n,2}}\bs{e}_{n,r} ~ \\
\Pfeil{\widetilde{Z}_{p,q,n} } & 
-a_{p,1}\underline{\Sumn{i=3}a_{q,i}a_{n,i}\bs{e}_{n,r}}
+a_{p,1}\delta_{q,n}\bs{e}_{n,r}
+\left(\Sumn{i=2}a_{p,i}a_{q,i}-\delta_{p,q}\right)\underline{a_{n,1}\bs{e}_{n,r}}
\\
\Pfeil{\ref{Rechenregeln2},\widetilde{W_{1,r}}}& 
-a_{p,1}\left(
	-\Sumn{i=3}\Sum{j=1}{n-1}a_{q,i}a_{j,i}\bs{e}_{j,r} 
	+\Sumn{j=1}a_{q,1}\bs{e}_{j,1} a_{j,r} \right.\\& \left.
	+\Sumn{j=1}a_{q,2}\bs{e}_{j,2} a_{j,r}
  +	\bs{e}_{q,r}
  -	\Sumn{j=1}\bs{f}_{q,j}a_{j,r}
  + \Sumn{i=3}\Sumn{k,j=1}a_{q,i} a_{j,i}\bs{f}_{j,k}a_{k,r}
\right)\\&
+a_{p,1}\delta_{q,n}\bs{e}_{n,r}
+\left(\Sumn{i=2}a_{p,i}a_{q,i}-\delta_{p,q}\right)\left(
	-\Sum{j=1}{n-1}a_{j,1}\bs{e}_{j,r} \right.\\& \left.
	-\Sumn{j=1}\bs{e}_{j,1}a_{j,r} 
	+ \Sumn{k,j=1}a_{j,1}\bs{f}_{j,k}a_{k,r}
\right)
\\
=&
	\Sumn{i=3}\Sum{j=1}{n-1}a_{p,1}a_{q,i}a_{j,i}\bs{e}_{j,r} 
	-\Sumn{j=1}\underline{a_{p,1}a_{q,1}}\bs{e}_{j,1} a_{j,r}\\&
	-\Sumn{j=1}a_{p,1}a_{q,2}\bs{e}_{j,2} a_{j,r}
  -	a_{p,1}\bs{e}_{q,r}
  +	\Sumn{j=1}a_{p,1}\bs{f}_{q,j}a_{j,r}\\&
  - \Sumn{i=3}\Sumn{k,j=1}a_{p,1}a_{q,i} a_{j,i}\bs{f}_{j,k}a_{k,r}
+a_{p,1}\delta_{q,n}\bs{e}_{n,r}\\&
-\Sumn{i=2}\Sum{j=1}{n-1}a_{p,i}a_{q,i}a_{j,1}\bs{e}_{j,r} 
-\Sumn{i=2}\Sumn{j=1}a_{p,i}a_{q,i}\bs{e}_{j,1}a_{j,r} \\&
+ \Sumn{i=2}\Sumn{k,j=1}a_{p,i}a_{q,i}a_{j,1}\bs{f}_{j,k}a_{k,r}\\&
+\delta_{p,q}\left(\Sum{j=1}{n-1}a_{j,1}\bs{e}_{j,r} 
+\Sumn{j=1}\bs{e}_{j,1}a_{j,r} 
- \Sumn{k,j=1}a_{j,1}\bs{f}_{j,k}a_{k,r}\right)
\NeueSeite
\Pfeil{\widetilde{Z}_{p,q}}&
	\Sumn{i=3}\Sum{j=1}{n-1}a_{p,1}a_{q,i}a_{j,i}\bs{e}_{j,r} \\&
	+\Sumn{j=1}\Sumn{i=2}a_{p,i}a_{q,i}\bs{e}_{j,1} a_{j,r}
		-\Sumn{j=1}\delta_{p,q}\bs{e}_{j,1} a_{j,r}\\&
	-\Sumn{j=1}a_{p,1}a_{q,2}\bs{e}_{j,2} a_{j,r}
  -	a_{p,1}\bs{e}_{q,r}
  +	\Sumn{j=1}a_{p,1}\bs{f}_{q,j}a_{j,r}\\&
  - \Sumn{i=3}\Sumn{k,j=1}a_{p,1}a_{q,i} a_{j,i}\bs{f}_{j,k}a_{k,r}
+a_{p,1}\delta_{q,n}\bs{e}_{n,r}\\&
-\Sumn{i=2}\Sum{j=1}{n-1}a_{p,i}a_{q,i}a_{j,1}\bs{e}_{j,r} 
-\Sumn{i=2}\Sumn{j=1}a_{p,i}a_{q,i}\bs{e}_{j,1}a_{j,r} \\&
+ \Sumn{i=2}\Sumn{k,j=1}a_{p,i}a_{q,i}a_{j,1}\bs{f}_{j,k}a_{k,r}\\&
+\delta_{p,q}\left(\Sum{j=1}{n-1}a_{j,1}\bs{e}_{j,r} 
+\Sumn{j=1}\bs{e}_{j,1}a_{j,r} 
- \Sumn{k,j=1}a_{j,1}\bs{f}_{j,k}a_{k,r}\right)\\
=
&
	\Sumn{i=3}\Sum{j=1}{n-1}a_{p,1}a_{q,i}a_{j,i}\bs{e}_{j,r} 	
  -	a_{p,1}(1-\delta_{q,n})\bs{e}_{q,r}\\&
	-\Sumn{i=2}\Sum{j=1}{n-1}\left(a_{p,i}a_{q,i} -\delta_{p,q}\right) a_{j,1}\bs{e}_{j,r}
	 -\Sumn{j=1}a_{p,1}a_{q,2}\bs{e}_{j,2} a_{j,r}
	 \\&
  +	\Sumn{j=1}a_{p,1}\bs{f}_{q,j}a_{j,r}\\&
  - \Sumn{i=3}\Sumn{k,j=1}a_{p,1}a_{q,i} a_{j,i}\bs{f}_{j,k}a_{k,r} 
\\&
+ \Sumn{i=2}\Sumn{k,j=1}\left( a_{p,i}a_{q,i}-\delta_{p,q} \right) a_{j,1}\bs{f}_{j,k}a_{k,r}
\\
=
&
	-\Sum{j=1}{n-1}Z_{p,q,j}\bs{e}_{j,r} 
	+\Sum{j=1}{n-1}Z_{p,q} a_{j,1}\bs{e}_{j,r} \\&
	-\Sumn{j=1}a_{p,1}a_{q,2}\bs{e}_{j,2} a_{j,r}
	 \\&
  +\Sumn{k,j=1}Z_{p,q,j}\bs{f}_{j,k}a_{k,r} 
- \Sumn{k,j=1}Z_{p,q} a_{j,1}\bs{f}_{j,k}a_{k,r}
}%Gleichung
\newpage
\subsubsection{Beginnend mit $(a_{p,1}a_{q,2}a_{n,2}\bs{e}_{n,r} ~,~ a_{p,1}a_{q,2}W_{2,r})$:}
\Gleichung{
&a_{p,1}a_{q,2}\underline{a_{n,2}\bs{e}_{n,r}} \\
\Pfeil{\widetilde{W}_{2,r}} & 
-\Sum{j=1}{n-1}\underline{a_{p,1}a_{q,2}a_{j,2}}\bs{e}_{j,r} -\Sumn{j=1}a_{p,1}a_{q,2}\bs{e}_{j,2}a_{j,r} + \Sumn{k,j=1}\underline{a_{p,1}a_{q,2}a_{j,2}}\bs{f}_{j,k}a_{k,r}\\
=
&
-\Sum{j=1}{n-1}Z_{p,q,j}\bs{e}_{j,r}
+\Sum{j=1}{n-1}Z_{p,q}a_{j,1}\bs{e}_{j,r}  \\&
-\Sumn{j=1}a_{p,1}a_{q,2}\bs{e}_{j,2}a_{j,r} \\&
+ \Sumn{k,j=1}Z_{p,q,j}\bs{f}_{j,k}a_{k,r}
- \Sumn{k,j=1}Z_{p,q}a_{j,1}\bs{f}_{j,k}a_{k,r}
	}%Geleichung

%\Bem{Symmetrie}{Diese Rechnungen sind nach Vertauschen der Indizes identisch zu denen in \ref{ZpqV1}.}

%% file: Haupt/RechnungAufloesung/Phi2/WV2.tex
%Nur lokal für die Rechnungen brauch ich eine andere Summenumgebung 

\subsection{Reduktionswege für $a_{n,n}\bs{e}_{n,n}$}
\label{WV2}

\subsubsection{Beginnend mit $(a_{n,n}\bs{e}_{n,n} ~,~ W_{n,n})$:}

\Gleichung{
&\underline{a_{n,n}\bs{e}_{n,n}} ~ \\
\Pfeil{\widetilde{W}_{n,n} } & 
-\Sum{j=1}{n-1}\underline{a_{j,n}\bs{e}_{j,n}} -\Sumn{j=1}\bs{e}_{j,n}a_{j,n} + \Sumn{k,j=1}a_{j,n}\bs{f}_{j,k}a_{k,n}\\
\Pfeil{\widetilde{V}_{j,j} } & 
-\Sum{j=1}{n-1}\left(-  \sum\limits_{i=1}^{n-1}a_{j,i}\bs{e}_{j,i} - \sum\limits_{i=1}^{n}\bs{e}_{j,i} a_{j,i} + \bs{f}_{j,j}
\right)
\\&-\Sumn{j=1}\bs{e}_{j,n}a_{j,n} + \Sumn{k,j=1}a_{j,n}\bs{f}_{j,k}a_{k,n}
\\
=
&
\Sum{i,j=1}{n-1}a_{j,i}\bs{e}_{j,i} 
+\Sum{i,j=1}{n-1} \bs{e}_{j,i} a_{j,i} 
-\bs{e}_{n,n}a_{n,n} 
-\Sum{j=1}{n-1}\bs{f}_{j,j}
\\&
+ \Sumn{k,j=1}a_{j,n}\bs{f}_{j,k}a_{k,n}\\
}%Gleichung

\newpage

\subsubsection{Beginnend mit $(a_{n,n}\bs{e}_{n,n} ~,~ V_{n,n})$:}
\Gleichung{
&\underline{a_{n,n}\bs{e}_{n,n}} ~ \\
\Pfeil{\widetilde{V}_{n,n} } & 
-  \sum\limits_{i=1}^{n-1}\underline{a_{n,i}\bs{e}_{n,i}} - \sum\limits_{i=1}^{n}\bs{e}_{n,i} a_{n,i} + \bs{f}_{n,n}
\\
\Pfeil{\widetilde{W}_{i,i} } & -  \sum\limits_{i=1}^{n-1}\left(-\Sum{j=1}{n-1}a_{j,i}\bs{e}_{j,i} -\Sumn{j=1}\bs{e}_{j,i}a_{j,i} + \Sumn{k,j=1}a_{j,i}\bs{f}_{j,k}a_{k,i}
\right)
\\& - \sum\limits_{i=1}^{n}\bs{e}_{n,i} a_{n,i} + \bs{f}_{n,n}
\\
=
&
 \Sum{i,j=1}{n-1}a_{j,i}\bs{e}_{j,i} 
	+\Sum{i,j=1}{n-1}\bs{e}_{j,i}a_{j,i} -\bs{e}_{n,n} a_{n,n} + \bs{f}_{n,n}\\&
	-\Sum{i=2}{n-1}\Sumn{k,j=1}a_{j,i}\bs{f}_{j,k}a_{k,i}
	- \Sumn{k,j=1\\(k,j)\neq(n,n)}a_{j,1}\bs{f}_{j,k}a_{k,1}  \\&
	- \underline{a_{n,1}\bs{f}_{n,n}a_{n,1}}
\\
\Pfeil{\widetilde{K}}&
 \Sum{i,j=1}{n-1}a_{j,i}\bs{e}_{j,i} 
	+\Sum{i,j=1}{n-1}\bs{e}_{j,i}a_{j,i} -\bs{e}_{n,n} a_{n,n} + \bs{f}_{n,n}\\&	-\Sum{i=2}{n-1}\Sumn{k,j=1}a_{j,i}\bs{f}_{j,k}a_{k,i}
	- \Sumn{k,j=1\\(k,j)\neq(n,n)}a_{j,1}\bs{f}_{j,k}a_{k,1} 	\\&
+\Sumn{i,j,k=1\\(i,j,k)\neq (1,n,n)}a_{j,i} \bs{f}_{j,k}a_{k,i} 
- \sum\limits_{i=1}^{n}\bs{f}_{i,i} 
\\
=
&
 \Sum{i,j=1}{n-1}a_{j,i}\bs{e}_{j,i} 
	+\Sum{i,j=1}{n-1}\bs{e}_{j,i}a_{j,i} -\bs{e}_{n,n} a_{n,n} - \Sum{i=1}{n-1}\bs{f}_{i,i}\\&
%		-\Sum{i=2}{n-1}\Sumn{k,j=1}a_{j,i}\bs{f}_{j,k}a_{k,i}
%	- \Sumn{k,j=1\\(k,j)\neq(n,n)}a_{j,1}\bs{f}_{j,k}a_{k,1} 	\\&
%+\Sumn{i,j,k=1\\(i,j,k)\neq (1,n,n)}a_{j,i} \bs{f}_{j,k}a_{k,i} 
+\Sumn{k,j=1}a_{j,n} \bs{f}_{j,k}a_{k,n}
	}%Geleichung

%\Bem{Symmetrie}{Diese Rechnungen sind nach Vertauschen der Indizes identisch zu denen in \ref{ZpqV1}.}

%% file: Haupt/RechnungAufloesung/BiPhi3.tex
\section{$\Phi_3$}
\label{sec:Phi3}

\subsection{Übersichtstabelle für $\Phi_3$}
\label{sec:ÜbersichtstabelleFürPhi3}

Um die Lesbarkeit zu erhöhen setzen wir:
\mathe{ Q := \Sumn{i=3}\Sumn{j,k=1} a_{j,i}\bs{f}_{j,k}a_{k,i} - \Sumn{i=1}\bs{f}_{i,i}-\bs{f}.}
%\kleinerText
{

\mathe{\begin{array}{|l|l|l|}
 \hline
						& \LMatrix{~\\ \Phi_3 : \Aev(n) \rightarrow \Aev(n)^{n^2}\\ ~ } 	
\\ \hline
 \LMatrix{~\\ \textnormal{Matrix} \\ ~ } 			& ~~~\bs{f} \mapsto -\tr(A^t\bs{F}A)+\tr(\bs{F}) 
\\ \hline
 \LMatrix{~\\ \textnormal{Kompo-} \\ \textnormal{nenten} \\ ~ }	& ~~~\bs{f} \mapsto -\sum\limits_{i,j,k=1}^{n} a_{j,i}\bs{f}_{j,k}a_{k,i} + \sum\limits_{i=1}^{n}\bs{f}_{i,i} 
\\ \hline
 \LMatrix{~\\ \RS_{\Phi_3} \\ ~ } & \LMatrix{ 
 	\widetilde{V}= \left( a_{n,1}\bs{f}_{n,n}a_{n,1} ~,~ 
 	%\LMatrix{~\\~} \right.&
 	 -Q - \Sumn{j,k=1} a_{j,2} \bs{f}_{j,k}a_{k,2}  - \Sumn{j,k=1\\(j,k)\neq(n,n)} a_{j,1} \bs{f}_{j,k}a_{k,1}  \right)
 }
 \\ \hline
 \LMatrix{~\\ \RS_{es} \\ ~ } &
 \LMatrix{
 \widetilde{V_p}= \left( \LMatrix{~\\a_{p,1}a_{n,2}\bs{f}_{n,n}a_{n,2} ~,~ \\~}	  \right.
  &
  -\Sumn{j,k=1} Z_{p,j}\bs{f}_{j,k}a_{k,1}  \\
  &
  \left. 	                           
  -a_{p,1}\left( Q+ \Sumn{j,k=1\\(j,k)\neq (n,n)}a_{j,2} \bs{f}_{j,k}a_{k,2}  \right) \LMatrix{~\\~\\~}\right)
  \\~\\
  \widetilde{W_p}= \left(	\LMatrix{~\\ a_{n,1}\bs{f}_{n,n}a_{n,2}a_{p,2}~,~ \\~}
   \right.& 
  - a_{n,1}\bs{f}_{n,n}\left(\Sumn{l=3}a_{n,l}a_{p,l}  -\delta_{n,p} \right)
  \\& 
  +	\Sumn{j,k=1\\(j,k)\neq (n,n)}a_{j,1} \bs{f}_{j,k} Z_{k,p}
	\\& 
	\left.
		+\left(Q+\Sumn{j,k=1}a_{j,2} \bs{f}_{j,k}a_{k,2}\right)a_{p,1} \LMatrix{~\\~\\~\\} \right)
 	}%LMatrix
\\ \hline
 \LMatrix{ \RSkern  } &  \LMatrix{~\\~  }
\\ \hline
& \\ \hline
\textnormal{Konflikte}& \LMatrix{~\\ \textnormal{siehe Tabelle } \ref{sec:Konflikte_Phi3_n}\\~}
\\ \hline
\end{array}\\
 }%Mathe
 }% kleinerText

 % \subsection{Allgemeine Bemerkungen über die Konflikte in $\Phi_3$}
%\label{nichtImModul}
%
%Die Konflikte die zwischen den Regeln die durch $\Phi_3$ erzeugt werden, haben mehrere Besonderheiten.
%\begin{itemize}
%	\item{Führende Terme:} Das größte Monom der Regeln $\widetilde{V}$, $\widetilde{V_p}$ und $\widetilde{W}_p$ hat jeweils einen Anteil in $\fQG^{op}$. Deshalb können auch Überschneidungen zwischen diesen Regeln entstehen, z.B. zwischen $\widetilde{V}$ und $\widetilde{V_n}$:
%	$a_{n,1}\bs{f}_{n,n}a_{n,1}a_{n,2}\bs{f}_{n,n}a_{n,2}.$\\
%	In diesen Überschneidungen kommen aber jeweils zwei Erzeuger des Moduls vor. Da alle Regeln die Anzahl der Modulerzeuger konstant lassen und Monome mit zwei Modulerzeuger Elemente des Moduls sind, brauchen wir diese Überschneidungen nicht betrachten.
%	\item{Spezialfall $n=2$:} Hier gibt es deutlich mehr Konflikte als für alle anderen Fälle. Wir betrachten den Fall $n=2$ daher gesondert.
%\end{itemize}
%~\\
%\Frage{Achtung: In Reduktionssysteme für Modul muß erklärt werden, warum 	\mathe{a_{n,1}\bs{f}_{n,n}a_{n,1}a_{n,2}\bs{f}_{n,n}a_{n,2}.
%	}
%nicht betrachtet werden muß}
% 

%\newpage

\newpage
\subsection{$\RS_{es}$}
\label{sec:RegelnRSesPhi3}

Wir werden zeigen, dass $S_{es} \subset \left< S_{{\Phi_3}} \right>$. Dazu  betrachten wir $S_{\Phi_3}$:
\mathe{
a_{n,1}\bs{f}_{n,n}a_{n,1} -\left( 
 	 -Q - \Sumn{j,k=1} a_{j,2} \bs{f}_{j,k}a_{k,2}  - \Sumn{j,k=1\\(j,k)\neq(n,n)} a_{j,1} \bs{f}_{j,k}a_{k,1} \right).
}
Wenn wir von links mit $a_{p,1}$ multiplizieren, erhalten wir:
\kleinerText{
\Gleichung{%& ~~~ 
& \underline{a_{p,1}a_{n,1}}\bs{f}_{n,n}a_{n,1} + {a_{p,1}}Q + {a_{p,1}}\Sumn{j,k=1} a_{j,2} \bs{f}_{j,k}a_{k,2}  + \Sumn{j,k=1\\(j,k)\neq(n,n)}\underline{a_{p,1}a_{j,1}} \bs{f}_{j,k}a_{k,1} \\
= &  %& ~~~ 
a_{p,1}a_{n,2}\bs{f}_{n,n}a_{n,2} +	                           
  a_{p,1} Q+  a_{p,1}\Sumn{j,k=1\\(j,k)\neq (n,n)}a_{j,2} \bs{f}_{j,k}a_{k,2}   +\Sumn{j,k=1} Z_{p,j}\bs{f}_{j,k}a_{k,1}  .
}
}
Das entspricht genau dem Teil von $S_{es}$, der von $\widetilde{V}_p$ induziert ist.
\\
Wenn wir von rechts mit $a_{p,1}$ multiplizieren, erhalten wir:
\kleinerText{
\Gleichung{%& ~~~
 &a_{n,1}\bs{f}_{n,n}\underline{a_{n,1}a_{p,1}} + Q a_{p,1} + \Sumn{j,k=1} a_{j,2} \bs{f}_{j,k}a_{k,2}a_{p,1}  + \Sumn{j,k=1\\(j,k)\neq(n,n)}a_{j,1} \bs{f}_{j,k}\underline{a_{k,1}a_{p,1}}\\
=&
%& ~~~
 a_{n,1}\bs{f}_{n,n}a_{n,2}a_{p,2}+ a_{n,1}\bs{f}_{n,n}\left(\Sumn{l=3}a_{n,l}a_{p,l}  +\delta_{n,p} \right)
	\\
	%	~~~~~~~~~~~~~~~~~~~~~~~~~~~~~~~~~~~~~~~~~~~
		&-\left(Q+\Sumn{j,k=1}a_{j,2} \bs{f}_{j,k}a_{k,2}\right)a_{p,1}-
	\Sumn{j,k=1\\(j,k)\neq (n,n)}a_{j,1} \bs{f}_{j,k} Z_{k,p}.
}
}
Das entspricht genau dem Teil von $S_{es}$, der von $\widetilde{W}_p$ induziert ist.\\

Zusammen folgt: $S_{es}\subset  \left< S_{{\Phi_3}} \right>$.%
%Es kommen in $\RSkern$ keine Regeln vor, die Abbildung $\Phi_3$ ist also injektiv.
% 
% 
% 
% 
% 
% 
% 
% \newpage
% 
% 
% 
% 
% 
% 
% 
% 
% 
% 
% 
% 
\subsection{Konflikte:}
\label{sec:Konflikte_Phi3_n}

Wir wollen nun eine Tabelle aller möglichen Konflikte betrachten. Dazu überlegen wir uns zunächst, dass es keine Konflikte innerhalb von $\red_{e}$ gibt. Dazu müsste es eine Überschneidung geben, so dass der Erzeuger des Moduls $\bs{f}_{n,n}$ in dem Wort  nur einmal vorkommt. Dies ist aber durch die festgesetzten Indizes nicht möglich.
\\
Zwischen den Regeln aus $\RSa$ und denen aus $\red_{e}$ kann es Überschneidungen geben. Anders als in den vorhergehenden Fällen können die Regeln aus  $\RSa$ sogar rechts von denen aus $\red_{e}$ vorkommen. Es kommen Überschneidungen mit einem oder mit zwei Buchstaben vor, in den meisten Fällen ist jedoch durch die festgesetzten Indizes keine Überschneidung möglich. In der folgenden Tabelle stehen alle minimalen Überschneidungen, wobei oberhalb des Striches die Regeln aus $\RSa$ von links und unterhalb von rechts angewendet werden. Wenn es Überschneidungen mit mehr als einem Buchstaben gibt, stehen beide Fälle übereinander.

\newpage
\begin{landscape}
Tabelle der minimalen Überschneidungen in $\red_{\Phi}\cup\red_{es}$:
\input{Haupt/RechnungAufloesung/Phi3/TabelleAllgemeinerFall}
\end{landscape}
\newpage

%
%\subsection{Konflikte für $n=2$:}
%\label{sec:Konflikte_Phi3_2}
% Für diesen Spezialfall müssen vereinfachen sich die Regeln zu:
% \kleinerText{
% \mathe{\LMatrix{
% 	\widetilde{V}= ( a_{2,1}\bs{f}_{2,2}a_{2,1} ~,&~  \Sum{i=1}{2}\bs{f}_{i,i}+\bs{f}-\Sum{j,k=1}{2} a_{j,2} \bs{f}_{j,k}a_{k,2}  - \Sum{j,k=1\\(j,k)\neq(2,2)}{2}a_{j,1} \bs{f}_{j,k}a_{k,1}  )\\
%  \widetilde{V_p}= ( 	a_{p,1}a_{2,2}\bs{f}_{2,2}a_{2,2} ~,&~
%  \Sum{j,k=1}{2} \left( a_{p,2}a_{j,2}-\delta_{p,j}\right)\bs{f}_{j,k}a_{k,1}  +a_{p,1} \Sum{i=1}{2}\bs{f}_{i,i}+a_{p,1}\bs{f}	                           
%  -a_{p,1} \Sum{j,k=1\\(j,k)\neq (2,2)}{2}a_{j,2} \bs{f}_{j,k}a_{k,2}  
%  )\\
%  \widetilde{W_p}= (	a_{2,1}\bs{f}_{2,2}a_{2,2}a_{p,2}~,&~ 
%  -a_{n,1}f_{n,n}\delta_{n,p} -	\Sum{j,k=1\\(j,k)\neq (2,2)}{2}a_{j,1} \bs{f}_{j,k} \left(a_{k,2}a_{p,2}-\delta_{k,p}\right)
%		-\Sum{i=1}{2}\bs{f}_{i,i}a_{p,1}-\bs{f}a_{p,1}	 +\Sum{j,k=1}{2}a_{j,2} \bs{f}_{j,k}a_{k,2}a_{p,1} )
%}}
%}
% Wir müssen zusätzlich zu den Konflikten in \ref{sec:Konflikte_Phi3_n} noch die Konflikte:
% 
%\begin{itemize}
%	\item{$\widetilde{S}_{p,q,1}$ und $\widetilde{V}$:} $a_{1,p}a_{2,q}a_{2,1}\bs{f}_{2,2}a_{2,1}$ ~~~~~~\ref{Spq1V_phi3}
% \item{$\widetilde{S}_{p,q,1}$ und $\widetilde{W}_p$:} $a_{1,p}a_{2,q}a_{2,1}\bs{f}_{2,2}a_{2,2}a_{p,2}$ ~~~~~~\ref{Spq1Wp_phi3}
%\end{itemize}
%betrachten.
%
%%\kleinerText{
%%\input{Haupt/RechnungAufloesung/Phi3/TabelleSpezialfalln=2}
%%}

\newpage

\MitRechnungen{
\input{Haupt/RechnungAufloesung/Phi3/ZpnV_phi3}\newpage
\input{Haupt/RechnungAufloesung/Phi3/VZnq_phi3}\newpage
\input{Haupt/RechnungAufloesung/Phi3/VZnqr_phi3}\newpage
\input{Haupt/RechnungAufloesung/Phi3/ZpnWq_phi3}\newpage
\input{Haupt/RechnungAufloesung/Phi3/ZpnVq_phi3}\newpage
\input{Haupt/RechnungAufloesung/Phi3/Sp1V1_phi3}\newpage
\input{Haupt/RechnungAufloesung/Phi3/Spq1V2_phi3}\newpage
\input{Haupt/RechnungAufloesung/Phi3/W1S2q_phi3}\newpage
\input{Haupt/RechnungAufloesung/Phi3/W1S2qr_phi3}\newpage
}

%% file: Haupt/RechnungAufloesung/Phi3/TabelleAllgemeinerFall.tex
  \mathe{
\LMatrix{~\\ \begin{array}{|l|l|l|l|}\hline

						&\LMatrix{~\\ \widetilde{V}: a_{n,1}\bs{f}_{n,n} a_{n,1} \\ ~ } & \widetilde{V_p}:a_{p,1}a_{n,2}\bs{f}_{n,n}a_{n,2}&
\widetilde{W_p}:	a_{n,1}\bs{f}_{n,n}a_{n,2}a_{p,2}
						\\ \hline
\LMatrix{~\\  \widetilde{Z}_{p,q}:\\ ~a_{p,1} a_{q,1} 
\\ ~ \\~} 
						&  \LMatrix{ ~\\a_{p,1}a_{n,1}\bs{f}_{n,n} a_{n,1} ~~~\ref{ZpnV3} 
												\\ \hline 
														a_{n,1}\bs{f}_{n,n}a_{n,1}a_{p,1} ~~~ \ref{VZnq3} \\~} 
						&
							 \LMatrix{ a_{p,1}a_{q,1}a_{n,2}\bs{f}_{n,n}a_{n,2} ~~~ \ref{ZpqVq3} \\ \textnormal{keine Konflikte}
							 					\\ \hline	
							 							\textnormal{keine Konflikte}\\~}
						&
							 \LMatrix{ ~\\ a_{p,1}a_{n,1}\bs{f}_{n,n}a_{n,2}a_{q,2} ~~~ \ref{ZpnWq3}	
							 					\\ \hline	
							 							\textnormal{keine Konflikte}			 \\ \textnormal{keine Konflikte}}
						\\ \hline
\LMatrix{~\\\widetilde{S}_{p,q}:\\ ~a_{1,p} a_{1,q} 
\\~\\~ } 
						&  \LMatrix{ ~\\ \textnormal{keine Konflikte} 
												\\ \hline 
												 \textnormal{keine Konflikte}\\~ }
						& \LMatrix{ a_{1,p} a_{1,1}a_{n,2}\bs{f}_{n,n}a_{n,2} ~~~ \ref{Sp1V1_phi3} \\ \textnormal{keine Konflikte}
												\\ \hline 
													\textnormal{keine Konflikte}\\~}
						& \LMatrix{~\\ \textnormal{keine Konflikte}
												\\ \hline	
							 						a_{n,1}\bs{f}_{n,n}a_{n,2}a_{1,2}a_{1,q}~~~ \ref{W1S2q_phi3}			 \\ \textnormal{keine Konflikte}}
						\\ \hline
\LMatrix{~\\\widetilde{Z}_{p,q,r}:\\~a_{p,1} a_{q,2} a_{r,2}
\\ ~\\~ }
						& \LMatrix{~\\ \textnormal{keine Konflikte} 
												\\ \hline  
												 a_{n,1}\bs{f}_{n,n}a_{n,1}a_{q,2}a_{r,2} ~~~\ref{VZnqr3}\\~}
						& \LMatrix{ \textnormal{keine Konflikte} \\ \textnormal{keine Konflikte}
												\\ \hline	
							 						\textnormal{keine Konflikte}	\\~		}
						& \LMatrix{~\\ \textnormal{keine Konflikte} 
												\\ \hline	
							 						\textnormal{keine Konflikte}		\\ \textnormal{keine Konflikte}}	
						\\ \hline
\LMatrix{~\\ \widetilde{S}_{p,q,r}:\\~a_{1,p} a_{2,q} a_{2,r}
\\ ~\\~}	
						& \LMatrix{ ~\\ \textnormal{keine Konflikte} 
												\\ \hline 
												 \textnormal{keine Konflikte}\\~ }
					 	& \LMatrix{	a_{1,p} a_{2,q}  a_{2,1}a_{n,2}\bs{f}_{n,n}a_{n,2} ~~~ \ref{Spq1V2_phi3}	\\ \textnormal{keine Konflikte}
												\\ \hline	
							 						\textnormal{keine Konflikte}	\\~		}	
						& \LMatrix{ ~\\ \textnormal{keine Konflikte} 
												\\ \hline 
												 a_{n,1}\bs{f}_{n,n}a_{n,2}a_{1,2} a_{2,q} a_{2,r} ~~~ \ref{W1S2qr_phi3}\\ \textnormal{keine Konflikte} }
						\\ \hline 
%&&&\\ \hline
						%
						%
%\LMatrix{~\\\widetilde{V}:\\ \\~ }
%						& \LMatrix{ \textnormal{beide gleich}
%											}
%						&\LMatrix{ 
%						 					 \textnormal{keine Konflikte}\\
%										 }
%						&\LMatrix{ 
%											 \textnormal{keine Konflikte}\\
%										 }
%						 \\\hline
%						 %
%						 %
%\LMatrix{~\\\widetilde{V_p}:\\ \\~ } 
%						&\LMatrix{ 
%											 \textnormal{keine Konflikte}\\
%										 }
%						&\LMatrix{ 
%											 \textnormal{beide gleich}
%										 }
%						&\LMatrix{ 
%											 \textnormal{keine Konflikte}
%										 }
%						 \\\hline
%						 %
%						 %
%\LMatrix{~\\\widetilde{W_p}:\\ \\~ }	
%						&\LMatrix{ 
%											 \textnormal{keine Konflikte}
%										 }
%						&\LMatrix{ 
%											 \textnormal{keine Konflikte}
%										 }
%						&\LMatrix{ 
%											 \textnormal{beide gleich}
%										 }
%						 \\\hline
%						 %
%						 %
\end{array}\\ ~ 
}%LMatrix			 
}%mathe

%% file: Haupt/RechnungAufloesung/Phi3/ZpnV_phi3.tex
%phi3 V & Z
%V
%  a_{n,1}\bs{f}_{n,n}a_{n,1}
% -\Sumn{i,j,k=1\\(i,j,k)\neq (1,n,n)}a_{j,i} \bs{f}_{j,k}a_{k,i} + \Sumn{i=1}\bs{f}_{i,i} 

%Z
% a_{p,1} a_{q,1}
% \Zi{p}{q}

\subsection{Reduktionswege für $a_{p,1}a_{n,1}\bs{f}_{n,n}a_{n,1} $}
\label{ZpnV3}

\subsubsection{Beginnend mit $(a_{p,1}a_{n,1}\bs{f}_{n,n}a_{n,1} ~,~ a_{p,1}V)$:}
\Gleichung{a_{p,1}
&\underline{a_{n,1}\bs{f}_{n,n}a_{n,1}} \\ ~ 
	\Pfeil{\widetilde{V}_{} } & 
	-a_{p,1}Q - \Sumn{j,k=1} a_{p,1}a_{j,2} \bs{f}_{j,k}a_{k,2}  - \Sumn{j,k=1\\(j,k)\neq(n,n)}\underline{a_{p,1}a_{j,1}} \bs{f}_{j,k}a_{k,1}  \\
\Pfeil{\widetilde{Z}_{p,j} } & 
		-a_{p,1}Q - \Sumn{j,k=1} a_{p,1}a_{j,2} \bs{f}_{j,k}a_{k,2}  - \Sumn{j,k=1\\(j,k)\neq(n,n)}Z_{p,j} \bs{f}_{j,k}a_{k,1}  \\
		=&	-a_{p,1}Q - \Sumn{j,k=1\\(j,k)\neq(n,n)} a_{p,1}a_{j,2} \bs{f}_{j,k}a_{k,2}  - \underline{a_{p,1}a_{n,2} \bs{f}_{n,n}a_{n,2}}\\& - \Sumn{j,k=1\\(j,k)\neq(n,n)}Z_{p,j} \bs{f}_{j,k}a_{k,1}  \\
\Pfeil{\widetilde{V}_{p} } & 
	-a_{p,1}Q - \Sumn{j,k=1\\(j,k)\neq(n,n)} a_{p,1}a_{j,2} \bs{f}_{j,k}a_{k,2}   
	\\& -\left(-\Sumn{j,k=1} Z_{p,j}\bs{f}_{j,k}a_{k,1}  -a_{p,1}\left( Q+ \Sumn{j,k=1\\(j,k)\neq (n,n)}a_{j,2} \bs{f}_{j,k}a_{k,2}  \right) \right) 
	\\& - \Sumn{j,k=1\\(j,k)\neq(n,n)}Z_{p,j} \bs{f}_{j,k}a_{k,1}
	\\
	=&~ Z_{p,n}\bs{f}_{n,n}a_{n,1}
}%Geleichung
\subsubsection{Beginnend mit $(a_{p,1}a_{n,1}\bs{f}_{n,n}a_{n,1} ~,~ Z_{p,n}\bs{f}_{n,n}a_{n,1} )$:}
\Gleichung{
\underline{a_{p,1}a_{n,1}}\bs{f}_{n,n}a_{n,1}  ~ 
	\Pfeil{\widetilde{Z}_{p,n} } & Z_{p,n}\bs{f}_{n,n}a_{n,1} 
}

%% file: Haupt/RechnungAufloesung/Phi3/VZnq_phi3.tex
%phi3 V & Z
%V
%  a_{n,1}\bs{f}_{n,n}a_{n,1}
% -\Sumn{i,j,k=1\\(i,j,k)\neq (1,n,n)}a_{j,i} \bs{f}_{j,k}a_{k,i} + \Sumn{i=1}\bs{f}_{i,i} 

%Z
% a_{p,1} a_{q,1}
% \Zi{p}{q}

\subsection{Reduktionswege für $a_{n,1}\bs{f}_{n,n}a_{n,1}a_{q,1} $}
\label{VZnq3}

\subsubsection{Beginnend mit $(a_{n,1}\bs{f}_{n,n}a_{n,1}a_{q,1} ~,~ Va_{q,1})$:}
\Gleichung{
&\underline{a_{n,1}\bs{f}_{n,n}a_{n,1}}a_{q,1} \\ ~ 
	\Pfeil{\widetilde{V}_{} } & -Qa_{q,1} - \Sumn{j,k=1} a_{j,2} \bs{f}_{j,k}a_{k,2}a_{q,1}  - \Sumn{j,k=1\\(j,k)\neq(n,n)}a_{j,1} \bs{f}_{j,k}\underline{a_{k,1} a_{q,1}} \\
	\Pfeil{\widetilde{Z}_{k,q} } &
	-Qa_{q,1} - \Sumn{j,k=1} a_{j,2} \bs{f}_{j,k}a_{k,2}a_{q,1}  - \Sumn{j,k=1\\(j,k)\neq(n,n)}a_{j,1} \bs{f}_{j,k}Z_{k,q} 
	}

\subsubsection{Beginnend mit $(a_{n,1}\bs{f}_{n,n}a_{n,1}a_{q,1} ~,~ a_{n,1}\bs{f}_{n,n} \left(Z_{n,q}\right))$:}
\Gleichung{&a_{n,1}\bs{f}_{n,n}
\underline{a_{n,1}a_{q,1} } ~ \\
	\Pfeil{\widetilde{Z}_{n,q} } & 
	a_{n,1}\bs{f}_{n,n}Z_{n,q} \\
		=&
	-\underline{a_{n,1}\bs{f}_{n,n}a_{n,2}a_{q,2}} 
	- \Sumn{i=3}a_{n,1}\bs{f}_{n,n}a_{n,i}a_{q,i} + \delta_{n,q}a_{n,1}\bs{f}_{n,n}\\
		\Pfeil{\widetilde{W}_{q} } & 
		 a_{n,1}\bs{f}_{n,n}\left(\Sumn{l=3}a_{n,l}a_{q,l}  -\delta_{n,q} \right)
		-	\Sumn{j,k=1\\(j,k)\neq (n,n)}a_{j,1} \bs{f}_{j,k} Z_{k,q}\\&
		-\left(Q+\Sumn{j,k=1}a_{j,2} \bs{f}_{j,k}a_{k,2}\right)a_{q,1} \\
		&- \Sumn{i=3}a_{n,1}\bs{f}_{n,n}a_{n,i}a_{q,i} + \delta_{n,q}a_{n,1}\bs{f}_{n,n}\\
		=&-	\Sumn{j,k=1\\(j,k)\neq (n,n)}a_{j,1} \bs{f}_{j,k} Z_{k,q}	-Qa_{q,1} -\Sumn{j,k=1}a_{j,2} \bs{f}_{j,k}a_{k,2}a_{q,1} 
}%Geleichung

%% file: Haupt/RechnungAufloesung/Phi3/VZnqr_phi3.tex
%phi3 V & Zii

\subsection{Reduktionswege für $a_{n,1}\bs{f}_{n,n}a_{n,1}a_{q,2} a_{r,2} $}
\label{VZnqr3}

\subsubsection{Beginnend mit $(a_{n,1}\bs{f}_{n,n}a_{n,1}a_{q,2} a_{r,2} ~,~ Va_{q,2} a_{r,2})$:}
\Gleichung{
&\underline{a_{n,1}\bs{f}_{n,n}a_{n,1}}a_{q,2} a_{r,2} \\ ~ 
	\Pfeil{\widetilde{V}_{} } & 
-Qa_{q,2} a_{r,2}
- \Sumn{j,k=1} a_{j,2} \bs{f}_{j,k}a_{k,2}a_{q,2} a_{r,2}  \\&
- \Sumn{j,k=1\\(j,k)\neq(n,n)}a_{j,1} \bs{f}_{j,k}\underline{a_{k,1}a_{q,2} a_{r,2}}  \\
%
%%%
\Pfeil{\widetilde{Z}_{k,q,r} } & 
	-Qa_{q,2} a_{r,2}
- \Sumn{j,k=1} a_{j,2} \bs{f}_{j,k}a_{k,2}a_{q,2} a_{r,2}  \\&
- \Sumn{j,k=1\\(j,k)\neq(n,n)} a_{j,1} \bs{f}_{j,k}\left( Z_{k,q,r}-Z_{k,q}a_{r,1} \right)  \\
}%Geleichung
\subsubsection{Beginnend mit $(a_{n,1}\bs{f}_{n,n}a_{n,1}a_{q,2} a_{r,2} ~,~ a_{n,1}\bs{f}_{n,n}\left(Z_{n,q,r}-Z_{n,q}a_{r,1}\right) )$:}
\Gleichung{a_{n,1}\bs{f}_{n,n}
&\underline{a_{n,1}a_{q,2} a_{r,2}} \\ ~ 
	\Pfeil{\widetilde{Z}_{n,q,r} } &
	- \underline{a_{n,1}\bs{f}_{n,n}a_{n,1}} \left(\Sumn{l=3}a_{q,i}a_{r,i} -\delta_{q,r}\right) 
	\\
	&+ a_{n,1}\bs{f}_{n,n}	\left(\Sumn{i=2} a_{n,i}a_{q,i}-\delta_{n,q} \right) 	a_{r,1}
	\\
		\Pfeil{\widetilde{V}_{} } & 
	 \left(Q + \Sumn{j,k=1} a_{j,2} \bs{f}_{j,k}a_{k,2} 
	 				\right.\\&~~~~~~~~~~~~~~\left. + \Sumn{j,k=1\\(j,k)\neq(n,n)}a_{j,1} \bs{f}_{j,k}a_{k,1} \right)\left(\Sumn{l=3}a_{q,l}a_{r,l} -\delta_{q,r}\right) \\
	&+ \underline{a_{n,1}\bs{f}_{n,n}a_{n,2}a_{q,2}}a_{r,1}+ 
	a_{n,1}\bs{f}_{n,n}\left(\Sumn{i=3} a_{n,i}a_{q,i}-\delta_{n,q} \right) a_{r,1}
	\NeueSeite
	\Pfeil{\widetilde{W}_{q} } &
		 \left(Q + \Sumn{j,k=1} a_{j,2} \bs{f}_{j,k}a_{k,2} 
	 				\right.\\&~~~~~~~~~~~~~~\left. + \Sumn{j,k=1\\(j,k)\neq(n,n)}a_{j,1} \bs{f}_{j,k}a_{k,1} \right)\left(\Sumn{l=3}a_{q,l}a_{r,l} -\delta_{q,r}\right) \\
	& - a_{n,1}\bs{f}_{n,n}\left(\Sumn{l=3}a_{n,l}a_{q,l}  -\delta_{n,q} \right)a_{r,1}
	\\&
	+	\Sumn{j,k=1\\(j,k)\neq (n,n)}a_{j,1} \bs{f}_{j,k} Z_{k,q}a_{r,1}
	\\&
		+\left(Q+\Sumn{j,k=1}a_{j,2} \bs{f}_{j,k}a_{k,2}\right)\underline{a_{q,1} 		a_{r,1}}
	\\&+	a_{n,1}\bs{f}_{n,n}\left(\Sumn{i=3} a_{n,i}a_{q,i}-\delta_{n,q} \right) a_{r,1}\\
\Pfeil{\widetilde{Z}_{q,r}}&
	 \left(Q + \Sumn{j,k=1} a_{j,2} \bs{f}_{j,k}a_{k,2} \right)\left(\Sumn{l=3}a_{q,l}a_{r,l} -\delta_{q,r}\right)\\&
	 				 + \Sumn{j,k=1\\(j,k)\neq(n,n)}a_{j,1} \bs{f}_{j,k}a_{k,1} \left(\Sumn{l=3}a_{q,l}a_{r,l} -\delta_{q,r}\right) \\
	&+	\Sumn{j,k=1\\(j,k)\neq (n,n)}a_{j,1} \bs{f}_{j,k} Z_{k,q}a_{r,1}
	\\&
		+\left(Q+\Sumn{j,k=1}a_{j,2} \bs{f}_{j,k}a_{k,2}\right)Z_{q,r}\\
=& 	 -\left(Q + \Sumn{j,k=1} a_{j,2} \bs{f}_{j,k}a_{k,2} \right)a_{q,2}a_{r,2} \\&
	 		- \Sumn{j,k=1\\(j,k)\neq(n,n)}a_{j,1} \bs{f}_{j,k}\left(Z_{k,q,r}-Z_{k,q}a_{r,1}\right)
}%Gleichung

%% file: Haupt/RechnungAufloesung/Phi3/ZpnWq_phi3.tex
%phi3 Z & W 

\subsection{Reduktionswege für $a_{p,1}a_{n,1}\bs{f}_{n,n}a_{n,2}a_{q,2} $}
\label{ZpnWq3}

\subsubsection{Beginnend mit $(a_{p,1}a_{n,1}\bs{f}_{n,n}a_{n,2}a_{q,2}~,~a_{p,1}W_{q} )$:}
\Gleichung{a_{p,1}
&\underline{a_{n,1}\bs{f}_{n,n}a_{n,2}a_{q,2}} \\ ~ 
	\Pfeil{\widetilde{W}_{q} } & 
	- \underline{a_{p,1}a_{n,1}}\bs{f}_{n,n}\left(\Sumn{l=3}a_{n,l}a_{q,l}  -\delta_{n,q} \right)
	\\
	&+	\Sumn{j,k=1\\(j,k)\neq (n,n)}\underline{a_{p,1}a_{j,1}} \bs{f}_{j,k} Z_{k,q}
	\\&
		+a_{p,1}\left(Q+\Sumn{j,k=1}a_{j,2} \bs{f}_{j,k}a_{k,2}\right)a_{q,1} \\
	\Pfeil{\widetilde{Z}_{p,n}, \widetilde{Z}_{p,j}} & 	
		- Z_{p,n}\bs{f}_{n,n}\left(\Sumn{l=3}a_{n,l}a_{q,l}  -\delta_{n,q} \right)
		\\&+	\Sumn{j,k=1\\(j,k)\neq (n,n)}Z_{p,j} \bs{f}_{j,k} Z_{k,q}
	\\&
		+a_{p,1}\left(Q+\Sumn{j,k=1\\(j,k)\neq(n,n)}a_{j,2} \bs{f}_{j,k}a_{k,2}\right)a_{q,1} \\&
		+ \underline{a_{p,1}a_{n,2} \bs{f}_{n,n}a_{n,2}}a_{q,1} \\
	\Pfeil{\widetilde{V}_{p} } & 
	%- Z_{p,n}\bs{f}_{n,n}\left(\Sumn{l=3}a_{n,l}a_{q,l}  -\delta_{n,q} \right)\\&
	+	\Sumn{j,k=1}Z_{p,j} \bs{f}_{j,k} Z_{k,q} + Z_{p,n}\bs{f}_{n,n}a_{n,2}a_{q,2}
	\\&
		+a_{p,1}\left(Q+\Sumn{j,k=1\\(j,k)\neq(n,n)}a_{j,2} \bs{f}_{j,k}a_{k,2}\right)a_{q,1} \\
		&
  -\Sumn{j,k=1} Z_{p,j}\bs{f}_{j,k}\underline{a_{k,1}a_{q,1}}\\&
  -a_{p,1}\left( Q+ \Sumn{j,k=1\\(j,k)\neq (n,n)}a_{j,2} \bs{f}_{j,k}a_{k,2}  \right)a_{q,1} 
	\NeueSeite
	\Pfeil{\widetilde{Z}_{k,q}} &
		+	\Sumn{j,k=1}Z_{p,j} \bs{f}_{j,k} Z_{k,q} 
		+ Z_{p,n}\bs{f}_{n,n}a_{n,2}a_{q,2}
	\\&
		%+a_{p,1}\left(Q+\Sumn{j,k=1\\(j,k)\neq(n,n)}a_{j,2} \bs{f}_{j,k}a_{k,2}\right)a_{q,1} \\
		%
  -\Sumn{j,k=1} Z_{p,j}\bs{f}_{j,k}Z_{k,q}\\
  %-a_{p,1}\left( Q+ \Sumn{j,k=1\\(j,k)\neq (n,n)}a_{j,2} \bs{f}_{j,k}a_{k,2}  \right)a_{q,1} 
		%%
= &~ Z_{p,n}\bs{f}_{n,n}a_{n,2}a_{q,2}
}%Geleichung
\subsubsection{Beginnend mit $(a_{p,1}a_{n,1}\bs{f}_{n,n}a_{n,2}a_{q,2}~,~Z_{p,n}\bs{f}_{n,n}a_{n,2}a_{q,2} )$:}
\Gleichung{
\underline{a_{p,1}a_{n,1}}\bs{f}_{n,n}a_{n,2}a_{q,2} ~ 
	\Pfeil{\widetilde{Z}_{p,n} } & Z_{p,n}\bs{f}_{n,n}a_{n,2}a_{q,2}
}

%% file: Haupt/RechnungAufloesung/Phi3/ZpnVq_phi3.tex
%phi3 Z & V_p
%V
%  a_{n,1}\bs{f}_{n,n}a_{n,1}
% -\Sumn{i,j,k=1\\(i,j,k)\neq (1,n,n)}a_{j,i} \bs{f}_{j,k}a_{k,i} + \Sumn{i=1}\bs{f}_{i,i} 

%Z
% a_{p,1} a_{q,1}
% \Zi{p}{q}

%W1
%	a_{n,1}\bs{f}_{n,n}a_{n,2}a_{q,2}
% - a_{n,1}\bs{f}_{n,n}\left(\Sumn{l=3}a_{n,l}a_{q,l}  +\delta_{n,q} \right)+ 	\Sumn{j,k=1\\(j,k)\neq (n,n)}a_{j,1} \bs{f}_{j,k} Z_{k,q}		+\left(\Sumn{i=2}\Sumn{j,k=1}a_{j,i} \bs{f}_{j,k}a_{k,i}	- \Sumn{i=1}\bs{f}_{i,i}-  \bs{f}\right)a_{q,1} 

%W2
%	a_{p,1}a_{n,2}\bs{f}_{n,n}a_{n,2}
% -\Sumn{j,k=1} Z_{p,j}\bs{f}_{j,k}a_{k,1}  	                                                                           -a_{p,1}\left( \Sumn{i=3}\Sumn{j,k=1} a_{j,i} \bs{f}_{j,k}a_{k,i} 		+ \Sumn{j,k=1\\(j,k)\neq (n,n)}a_{j,2} \bs{f}_{j,k}a_{k,2} - \Sumn{i=1}\bs{f}_{i,i} -  \bs{f} \right)

\subsection{Reduktionswege für $a_{p,1}a_{q,1}a_{n,2}\bs{f}_{n,n}a_{n,2}$}
\label{ZpqVq3}

\subsubsection{Beginnend mit $(a_{p,1}a_{q,1}a_{n,2}\bs{f}_{n,n}a_{n,2}~,~a_{p,1}V_{q} )$:}
\Gleichung{a_{p,1}
&\underline{a_{q,1}a_{n,2}\bs{f}_{n,n}a_{n,2}} \\ ~ 
	\Pfeil{\widetilde{V}_{q} } &  
  -\Sumn{j,k=1} a_{p,1}Z_{q,j}\bs{f}_{j,k}a_{k,1}\\&  	                           
  -\underline{a_{p,1}a_{q,1}}\left( Q+ \Sumn{j,k=1\\(j,k)\neq (n,n)}a_{j,2} \bs{f}_{j,k}a_{k,2}  \right) 
			\\
	\Pfeil{\widetilde{Z}_{p,q}} & 	  
	-\Sumn{j,k=1} a_{p,1}Z_{q,j}\bs{f}_{j,k}a_{k,1}  \\&	                           
  -Z_{p,q}\left( Q+ \Sumn{j,k=1\\(j,k)\neq (n,n)}a_{j,2} \bs{f}_{j,k}a_{k,2}  \right) 
	\\
	=&
	 \Sumn{j,k=1}\left(\underline{a_{p,1}a_{q,2}a_{j,2}}-Z_{p,q,j} \right)\bs{f}_{j,k}a_{k,1} 
	\\
	&	-Z_{p,q}\left( Q+ \Sumn{j,k=1\\(j,k)\neq (n,n)}a_{j,2} \bs{f}_{j,k}a_{k,2}  \right) 
			\\
	\Pfeil{\widetilde{Z}_{p,q,j}} &  	 
	\Sumn{j,k=1}\left( Z_{p,q,j}-Z_{p,q}a_{j,1}-Z_{p,q,j} \right) \bs{f}_{j,k}a_{k,1} 
	\\
	&	-Z_{p,q}\left( Q+ \Sumn{j,k=1\\(j,k)\neq (n,n)}a_{j,2} \bs{f}_{j,k}a_{k,2}  \right) \\
	&
	=	Z_{p,q}\left(\underline{a_{n,1}\bs{f}_{n,n}a_{n,1}} 
	+\Sumn{j,k=1\\(j,k)\neq(n,n)}a_{j,1}  \bs{f}_{j,k}a_{k,1} 
	\right.\\& \left.
		+ Q+ \Sumn{j,k=1\\(j,k)\neq (n,n)}a_{j,2} \bs{f}_{j,k}a_{k,2}  \right) 
		\NeueSeite%\\
	\Pfeil{\widetilde{V}}& 
	Z_{p,q}\left(-Q - \Sumn{j,k=1} a_{j,2} \bs{f}_{j,k}a_{k,2}  - \Sumn{j,k=1\\(j,k)\neq(n,n)}a_{j,1} \bs{f}_{j,k}a_{k,1}
	\right.\\& \left.
	+\Sumn{j,k=1\\(j,k)\neq(n,n)}a_{j,1}  \bs{f}_{j,k}a_{k,1} 
			+ Q+ \Sumn{j,k=1\\(j,k)\neq (n,n)}a_{j,2} \bs{f}_{j,k}a_{k,2}  \right) 
			\\
	=& Z_{p,q}  a_{n,2} \bs{f}_{n,n}a_{n,2}
}%Gleichung

\subsubsection{Beginnend mit $(a_{p,1}a_{q,1}a_{n,2}\bs{f}_{n,n}a_{n,2}~,~Z_{p,q}a_{n,2}\bs{f}_{n,n}a_{n,2} )$:}
\Gleichung{
&\underline{a_{p,1}a_{q,1}}a_{n,2}\bs{f}_{n,n}a_{n,2} \\ ~ 
	\Pfeil{\widetilde{Z}_{p,q} } & Z_{p,q}a_{n,2}\bs{f}_{n,n}a_{n,2}
}%Gleichung

%% file: Haupt/RechnungAufloesung/Phi3/Sp1V1_phi3.tex
%Sp1 & V1 phi3

\subsection{Reduktionswege für $a_{1,p} a_{1,1}a_{n,2}\bs{f}_{n,n}a_{n,2}
 $}
\label{Sp1V1_phi3}

\subsubsection{Beginnend mit $(a_{1,p} a_{1,1}a_{n,2}\bs{f}_{n,n}a_{n,2} ~ , ~ a_{1,p}V_1)$:}

\Gleichung{a_{1,p} 
&\underline{a_{1,1}a_{n,2}\bs{f}_{n,n}a_{n,2} } \\ ~
	\Pfeil{\widetilde{V}_{1} } &  
	-\Sumn{j,k=1} a_{1,p}Z_{1,j}\bs{f}_{j,k}a_{k,1}  	\\&                           
	-\underline{a_{1,p}a_{1,1}}\left( Q+ \Sumn{j,k=1\\(j,k)\neq (n,n)}a_{j,2} \bs{f}_{j,k}a_{k,2}  \right) 
	\\
	\Pfeil{\widetilde{S}_{p,1} } &  
	\Sumn{j,k=1}\Sumn{l=2} \underline{a_{1,p}a_{1,l}}a_{j,l}\bs{f}_{j,k}a_{k,1}  -\Sumn{j,k=1} a_{1,p}\delta_{1,j}\bs{f}_{j,k}a_{k,1}
	\\&                        
	-S_{p,1}\left( Q+ \Sumn{j,k=1\\(j,k)\neq (n,n)}a_{j,2} \bs{f}_{j,k}a_{k,2}  \right) 
		\\
\Pfeil{\widetilde{S}_{p,l} } &
	\Sumn{j,k=1}\underline{\Sumn{l=2} S_{p,1}a_{j,l}}\bs{f}_{j,k}a_{k,1}  -\Sumn{j,k=1} a_{1,p}\delta_{1,j}\bs{f}_{j,k}a_{k,1}
	\\&                        
	-S_{p,1}\left( Q+ \Sumn{j,k=1\\(j,k)\neq (n,n)}a_{j,2} \bs{f}_{j,k}a_{k,2}  \right)
\\
\stackrel{\ref{sec:Rechenregeln}}{=}&	
	\Sumn{j,k=1}\left( \Sumn{l=2}a_{l,p}Z_{l,j} + \left(\delta_{j,1}-\delta_{1,p}\right)a_{j,p} \right)\bs{f}_{j,k}a_{k,1} \\&
		  - \Sumn{k=1} a_{1,p} \bs{f}_{1,k}a_{k,1} 	-S_{p,1}\left( Q+ \Sumn{j,k=1\\(j,k)\neq (n,n)}a_{j,2} \bs{f}_{j,k}a_{k,2}  \right)
\\
=&
\Sumn{j,k=1} \Sumn{l=2}a_{l,p}Z_{l,j}\bs{f}_{j,k}a_{k,1}  -\Sumn{j,k=1\\(j,k)\neq(n,n)}\delta_{1,p}a_{j,1}\bs{f}_{j,k}a_{k,1} \\& -\delta_{1,p}\underline{a_{n,1}\bs{f}_{n,n}a_{n,1} }
		 \\
    & 	-S_{p,1}\left( Q+ \Sumn{j,k=1\\(j,k)\neq (n,n)}a_{j,2} \bs{f}_{j,k}a_{k,2}  \right)
\NeueSeite
\Pfeil{\widetilde{V} } &
\Sumn{j,k=1} \Sumn{l=2}a_{l,p}Z_{l,j}\bs{f}_{j,k}a_{k,1}  -\Sumn{j,k=1\\(j,k)\neq(n,n)}\delta_{1,p}a_{j,1}\bs{f}_{j,k}a_{k,1} \\& 
+\delta_{1,p}\left(Q + \Sumn{j,k=1} a_{j,2} \bs{f}_{j,k}a_{k,2}  + \Sumn{j,k=1\\(j,k)\neq(n,n)}a_{j,1} \bs{f}_{j,k}a_{k,1} \right)
		 \\
    & 	-S_{p,1}\left( Q+ \Sumn{j,k=1\\(j,k)\neq (n,n)}a_{j,2} \bs{f}_{j,k}a_{k,2}  \right)
    \\
    &= 
    \Sumn{j,k=1} \Sumn{l=2}a_{l,p}Z_{l,j}\bs{f}_{j,k}a_{k,1}  
+\delta_{1,p}\left(Q + \Sumn{j,k=1} a_{j,2} \bs{f}_{j,k}a_{k,2}  \right)
		 \\
    & 	-S_{p,1}\left( Q+ \Sumn{j,k=1\\(j,k)\neq (n,n)}a_{j,2} \bs{f}_{j,k}a_{k,2}  \right)
}

\subsubsection{Beginnend mit $(a_{1,p} a_{1,1}a_{n,2}\bs{f}_{n,n}a_{n,2} ~ , ~ S_{p,1}a_{n,2}\bs{f}_{n,n}a_{n,2})$:}

\Gleichung{
&\underline{a_{1,p} a_{1,1} } a_{n,2}\bs{f}_{n,n}a_{n,2}\\ ~
		\Pfeil{\widetilde{S}_{p,1} } &  - \Sumn{l=2}a_{l,p}\underline{a_{l,1} a_{n,2}\bs{f}_{n,n}a_{n,2}} + \delta_{p,1}a_{n,2}\bs{f}_{n,n}a_{n,2}
		\\
		\Pfeil{\widetilde{V}_{l} } & 
		+ \Sumn{l=2}a_{l,p}\Sumn{j,k=1} Z_{l,j}\bs{f}_{j,k}a_{k,1}  	  \\&                         
		+ \Sumn{l=2}a_{l,p} a_{l,1}\left( Q+ \Sumn{j,k=1\\(j,k)\neq (n,n)}a_{j,2} \bs{f}_{j,k}a_{k,2}  \right) 
		 \\&+ \delta_{p,1}a_{n,2}\bs{f}_{n,n}a_{n,2}
		\\
		=&
		+ \Sumn{l=2}a_{l,p}\Sumn{j,k=1} Z_{l,j}\bs{f}_{j,k}a_{k,1}  	  \\&                         
		-\left(S_{p,1}-\delta_{p,1}\right)\left( Q+ \Sumn{j,k=1\\(j,k)\neq (n,n)}a_{j,2} \bs{f}_{j,k}a_{k,2}  \right) 
		 \\&+ \delta_{p,1}a_{n,2}\bs{f}_{n,n}a_{n,2}
}

%% file: Haupt/RechnungAufloesung/Phi3/Spq1V2_phi3.tex
%Spq1 & V2 phi3

\subsection{Reduktionswege für $a_{1,p} a_{2,q}a_{2,1}a_{n,2}\bs{f}_{n,n}a_{n,2}
 $}
\label{Spq1V2_phi3}

\subsubsection{Beginnend mit $(a_{1,p} a_{2,q}a_{2,1}a_{n,2}\bs{f}_{n,n}a_{n,2} ~ ,a_{1,p} a_{2,q}V_2 ~ )$:}

\Gleichung{a_{1,p} a_{2,q}
&\underline{a_{2,1}a_{n,2}\bs{f}_{n,n}a_{n,2} } \\ ~
\Pfeil{\widetilde{V}_{2} } & 
-\Sumn{j,k=1} a_{1,p} a_{2,q}Z_{2,j}\bs{f}_{j,k}a_{k,1}  \\&	                           
  -\underline{a_{1,p} a_{2,q}a_{2,1}}\left( Q+ \Sumn{j,k=1\\(j,k)\neq (n,n)}a_{j,2} \bs{f}_{j,k}a_{k,2}  \right) 
  \\
  \Pfeil{\widetilde{S}_{p,q,1} } & 
  \Sumn{j,k=1} \Sumn{i=2}\underline{a_{1,p} a_{2,q}a_{2,i}}a_{j,i}\bs{f}_{j,k}a_{k,1}  \\&	 
  -\Sumn{j,k=1} a_{1,p} a_{2,q}\delta_{2,j}\bs{f}_{j,k}a_{k,1}  \\&                          
  -\left(S_{p,q,1}-S_{p,q}a_{1,1}\right)\left( Q+ \Sumn{j,k=1\\(j,k)\neq (n,n)}a_{j,2} \bs{f}_{j,k}a_{k,2}  \right) \\
   \Pfeil{\widetilde{S}_{p,q,i} } &
   \Sumn{j,k=1} %\underline{
   		\Sumn{i=2}S_{p,q,i}a_{j,i}
   %}%underline
   \bs{f}_{j,k}a_{k,1}  \\&	 
     -\Sumn{j,k=1} \Sumn{i=2}S_{p,q}a_{1,i}a_{j,i}\bs{f}_{j,k}a_{k,1}  \\&
  -\Sumn{k=1} a_{1,p} a_{2,q}\bs{f}_{2,k}a_{k,1}  \\&                          
  -\left(S_{p,q,1}-S_{p,q}a_{1,1}\right)\left( Q+ \Sumn{j,k=1\\(j,k)\neq (n,n)}a_{j,2} \bs{f}_{j,k}a_{k,2}  \right) 
   \NeueSeite
   =&
     \Sumn{j,k=1} \underline{
     		\Sumn{i=2}S_{p,q,i}a_{j,i}
     }%underline
     \bs{f}_{j,k}a_{k,1}  \\&	 
     +\Sumn{j,k=1} S_{p,q}Z_{1,j}\bs{f}_{j,k}a_{k,1} -\Sumn{j,k=1} S_{p,q}\delta_{1,j}\bs{f}_{j,k}a_{k,1} \\&
  -\Sumn{k=1} a_{1,p} a_{2,q}\bs{f}_{2,k}a_{k,1}  \\&                          
  -\left(S_{p,q,1}-S_{p,q}a_{1,1}\right)\left( Q+ \Sumn{j,k=1\\(j,k)\neq (n,n)}a_{j,2} \bs{f}_{j,k}a_{k,2}  \right) \\
  %\NeueSeite
  %
  \stackrel{\ref{RechenRegeln}}{=}&
     \Sumn{j,k=1}
     \left(% \underline{\Sumn{i=2}S_{p,q,i}a_{j,i}}
     \Sumn{i=3}a_{1,p}a_{i,q}Z_{i,j}+\delta_{j,1}\underline{a_{1,p}a_{1,q}}\right.\\&\left.+ \delta_{j,2}a_{1,p}a_{2,q}-\delta_{q,1}a_{1,p}a_{j,1}
     \Matrix{~\\~}\right)
     \bs{f}_{j,k}a_{k,1}  \\&	
       +\Sumn{j,k=1} S_{p,q}Z_{1,j}\bs{f}_{j,k}a_{k,1} -\Sumn{k=1} S_{p,q}\bs{f}_{1,k}a_{k,1} \\&
  -\Sumn{k=1} a_{1,p} a_{2,q}\bs{f}_{2,k}a_{k,1}  \\&                          
  -\left(S_{p,q,1}-S_{p,q}a_{1,1}\right)\left( Q+ \Sumn{j,k=1\\(j,k)\neq (n,n)}a_{j,2} \bs{f}_{j,k}a_{k,2}  \right) \\ 
\Pfeil{\widetilde{S}_{p,q} } &
     \Sumn{j,k=1}\Sumn{i=3}a_{1,p}a_{i,q}Z_{i,j} \bs{f}_{j,k}a_{k,1} 
     +     \Sumn{k=1}S_{p,q} \bs{f}_{1,k}a_{k,1} 
     \\&
   %  + \Sumn{k=1}a_{1,p}a_{2,q} \bs{f}_{2,k}a_{k,1} 
     -\delta_{q,1}a_{1,p}\left(\Sumn{j,k=1\\(j,k)\neq(n,n)}a_{j,1} \bs{f}_{j,k}a_{k,1} 
     +\underline{a_{n,1} \bs{f}_{n,n}a_{n,1}}\right) 
     \\&	
     +\Sumn{j,k=1} S_{p,q}Z_{1,j}\bs{f}_{j,k}a_{k,1} -\Sumn{k=1} S_{p,q}\bs{f}_{1,k}a_{k,1} \\&
  %-\Sumn{k=1} a_{1,p} a_{2,q}\bs{f}_{2,k}a_{k,1}  \\&                          
  -\left(S_{p,q,1}-S_{p,q}a_{1,1}\right)\left( Q+ \Sumn{j,k=1\\(j,k)\neq (n,n)}a_{j,2} \bs{f}_{j,k}a_{k,2}  \right)  
\NeueSeite
\Pfeil{\widetilde{V}} &
     \Sumn{j,k=1}\Sumn{i=3}a_{1,p}a_{i,q}Z_{i,j} \bs{f}_{j,k}a_{k,1} 
     %+     \Sumn{k=1}S_{p,q} \bs{f}_{1,k}a_{k,1} 
     \\&
     -\delta_{q,1}a_{1,p}\left(\Sumn{j,k=1\\(j,k)\neq(n,n)}a_{j,1} \bs{f}_{j,k}a_{k,1} 
     -Q - \Sumn{j,k=1} a_{j,2} \bs{f}_{j,k}a_{k,2}  \right.\\&\left.-\Sumn{j,k=1\\(j,k)\neq(n,n)}a_{j,1} \bs{f}_{j,k}a_{k,1} \right)
     \\&	
     +\Sumn{j,k=1} S_{p,q}Z_{1,j}\bs{f}_{j,k}a_{k,1} %-\Sumn{k=1} S_{p,q}\bs{f}_{1,k}a_{k,1} 
     \\&
  -\left(S_{p,q,1}-S_{p,q}a_{1,1}\right)\left( Q+ \Sumn{j,k=1\\(j,k)\neq (n,n)}a_{j,2} \bs{f}_{j,k}a_{k,2}  \right) \\ 
=&
    \Sumn{j,k=1}\Sumn{i=3}a_{1,p}a_{i,q}Z_{i,j} \bs{f}_{j,k}a_{k,1}     +\Sumn{j,k=1} S_{p,q}Z_{1,j}\bs{f}_{j,k}a_{k,1} 
     \\&
     +\delta_{q,1}a_{1,p}\left(  Q +\Sumn{j,k=1} a_{j,2} \bs{f}_{j,k}a_{k,2}  \right)      
     \\&
  -\left(S_{p,q,1}-S_{p,q}a_{1,1}\right)\left( Q+ \Sumn{j,k=1\\(j,k)\neq (n,n)}a_{j,2} \bs{f}_{j,k}a_{k,2}  \right)  
}

\subsubsection{Beginnend mit $(a_{1,p} a_{2,q}a_{2,1}a_{n,2}\bs{f}_{n,n}a_{n,2} ~ , ~ (S_{p,q,1}-S_{p,q}a_{1,1})a_{n,2}\bs{f}_{n,n}a_{n,2} )$:}

\Gleichung{
&\underline{ a_{1,p} a_{2,q}a_{2,1}}a_{n,2}\bs{f}_{n,n}a_{n,2} \\ ~
  \Pfeil{\widetilde{S}_{p,q,1} } & 
  -\Sumn{i=3}a_{1,p}a_{i,q}\underline{a_{i,1}a_{n,2}\bs{f}_{n,n}a_{n,2}}
  +\delta_{q,1}a_{1,p}a_{n,2}\bs{f}_{n,n}a_{n,2}\\&
  -S_{p,q}\underline{a_{1,1}a_{n,2}\bs{f}_{n,n}a_{n,2}}\\
  \Pfeil{\widetilde{V}_{i},\widetilde{V}_{1} } &   
     \Sumn{j,k=1}\Sumn{i=3}a_{1,p}a_{i,q} Z_{i,j}\bs{f}_{j,k}a_{k,1}  	\\ &                          
  +\Sumn{i=3}a_{1,p}a_{i,q}a_{i,1}\left( Q+ \Sumn{j,k=1\\(j,k)\neq (n,n)}a_{j,2} \bs{f}_{j,k}a_{k,2}  \right) 
  \\&+\delta_{q,1}a_{1,p}a_{n,2}\bs{f}_{n,n}a_{n,2}\\&
  +\Sumn{j,k=1} S_{p,q}Z_{1,j}\bs{f}_{j,k}a_{k,1}  	                           
  \\&+S_{p,q}a_{1,1}\left( Q+ \Sumn{j,k=1\\(j,k)\neq (n,n)}a_{j,2} \bs{f}_{j,k}a_{k,2}  \right) 
\\=&
     \Sumn{j,k=1}\Sumn{i=3}a_{1,p}a_{i,q} Z_{i,j}\bs{f}_{j,k}a_{k,1}   +\Sumn{j,k=1} S_{p,q}Z_{1,j}\bs{f}_{j,k}a_{k,1}  	 	\\ &                          
  \\&+\delta_{q,1}a_{1,p}\left(Q+ \Sumn{j,k=1}a_{j,2} \bs{f}_{j,k}a_{k,2}  \right) \\&
  +\left(S_{p,q}a_{1,1}-S_{p,q,1}\right)\left( Q+ \Sumn{j,k=1\\(j,k)\neq (n,n)}a_{j,2} \bs{f}_{j,k}a_{k,2}  \right)
}

%% file: Haupt/RechnungAufloesung/Phi3/W1S2q_phi3.tex
%W1 & S2q phi3

\subsection{Reduktionswege für $a_{n,1}\bs{f}_{n,n}a_{n,2}a_{1,2} a_{1,q} $}
\label{W1S2q_phi3}

\subsubsection{Beginnend mit $(a_{n,1}\bs{f}_{n,n}a_{n,2}a_{1,2} a_{1,q}  ~ , ~ a_{n,1}\bs{f}_{n,n}a_{n,2}\left(S_{2,q}\right))$:}
\Gleichung{
&a_{n,1}\bs{f}_{n,n}a_{n,2}\underline{a_{1,2} a_{1,q}}\\ ~
\Pfeil{\widetilde{S}_{2,q} } & -a_{n,1}\bs{f}_{n,n}a_{n,2}\left(\Sumn{i=2}a_{i,2}a_{i,q} -\delta_{2,q}\right)\\
=&-\Sumn{i=2}\underline{a_{n,1}\bs{f}_{n,n}a_{n,2}a_{i,2}}a_{i,q}
+a_{n,1}\bs{f}_{n,n}a_{n,2}\delta_{2,q}\\
\Pfeil{\widetilde{W}_{i} } &
-\Sumn{i=2}\left(
- a_{n,1}\bs{f}_{n,n}\left(\Sumn{l=3}a_{n,l}a_{i,l}  -\delta_{n,i} \right)
\right.\\&\left.
+ 	\Sumn{j,k=1\\(j,k)\neq (n,n)}a_{j,1} \bs{f}_{j,k} Z_{k,i}
		+\left(Q+\Sumn{j,k=1}a_{j,2} \bs{f}_{j,k}a_{k,2}\right)a_{i,1}
\Matrix{~\\~\\~\\~\\}\right)a_{i,q}\\&
+a_{n,1}\bs{f}_{n,n}a_{n,2}\delta_{2,q}\\
=&
\Sumn{i=2} a_{n,1}\bs{f}_{n,n}\left(\Sumn{l=3}a_{n,l}a_{i,l} \right)a_{i,q}
-\Sumn{i=2} a_{n,1}\bs{f}_{n,n}\left(\delta_{n,i} \right)a_{i,q}\\&
- \Sumn{j,k=1\\(j,k)\neq (n,n)}a_{j,1} \bs{f}_{j,k} \underline{\Sumn{i=2}Z_{k,i}a_{i,q}}\\&
-\left(Q+\Sumn{j,k=1}a_{j,2} \bs{f}_{j,k}a_{k,2}\right)\Sumn{i=2}a_{i,1}a_{i,q}
+a_{n,1}\bs{f}_{n,n}a_{n,2}\delta_{2,q}\\
\stackrel{\ref{Rechenregel}}{=}&
 \Sumn{l=3}a_{n,1}\bs{f}_{n,n}a_{n,l}\Sumn{i=2}a_{i,l}a_{i,q} - a_{n,1}\bs{f}_{n,n}a_{n,q}\\&
- \Sumn{j,k=1\\(j,k)\neq (n,n)}a_{j,1} \bs{f}_{j,k} \left( \Sumn{i=2}a_{k,i}S_{i,q}+\left(\delta_{1,q}-\delta_{k,1}\right)a_{k,q}\right)\\&
+\left(Q+\Sumn{j,k=1}a_{j,2} \bs{f}_{j,k}a_{k,2}\right)\left(S_{1,q}-\delta_{1,q}\right)
+a_{n,1}\bs{f}_{n,n}a_{n,2}\delta_{2,q}
\NeueSeite
 &-\Sumn{l=3}a_{n,1}\bs{f}_{n,n}a_{n,l}\left(S_{l,q}-\delta_{l,q}\right) - a_{n,1}\bs{f}_{n,n}a_{n,q}\\&
- \Sumn{j,k=1\\(j,k)\neq (n,n)}a_{j,1} \bs{f}_{j,k} \left( \Sumn{i=2}a_{k,i}S_{i,q}+\left(\delta_{1,q}-\delta_{k,1}\right)a_{k,q}\right)\\&
+\left(Q+\Sumn{j,k=1}a_{j,2} \bs{f}_{j,k}a_{k,2}\right)\left(S_{1,q}-\delta_{1,q}\right)
+a_{n,1}\bs{f}_{n,n}a_{n,2}\delta_{2,q}
\\
=
&-\Sumn{l=3}a_{n,1}\bs{f}_{n,n}a_{n,l}\left(S_{l,q}\right)\\&
- \Sumn{j,k=1\\(j,k)\neq (n,n)}a_{j,1} \bs{f}_{j,k} \left( \Sumn{i=2}a_{k,i}S_{i,q}-\delta_{k,1}a_{k,q}\right)\\&
+\left(Q+\Sumn{j,k=1}a_{j,2} \bs{f}_{j,k}a_{k,2}\right)\left(S_{1,q}\right)\\&
+\Sumn{l=3}a_{n,1}\bs{f}_{n,n}a_{n,l}\left(\delta_{l,q}\right) - a_{n,1}\bs{f}_{n,n}a_{n,q}
\\&
- \Sumn{j,k=1\\(j,k)\neq (n,n)}a_{j,1} \bs{f}_{j,k} \left(\delta_{1,q}a_{k,q}\right)
\\&
+\left(Q+\Sumn{j,k=1}a_{j,2} \bs{f}_{j,k}a_{k,2}\right)\left(-\delta_{1,q}\right)
+a_{n,1}\bs{f}_{n,n}a_{n,2}\delta_{2,q}
\\
=
&-\Sumn{l=3}a_{n,1}\bs{f}_{n,n}a_{n,l}S_{l,q}\\&
- \Sumn{j,k=1\\(j,k)\neq (n,n)}a_{j,1} \bs{f}_{j,k} \left( \Sumn{i=2}a_{k,i}S_{i,q}-\delta_{k,1}a_{k,q}\right)\\&
+\left(Q+\Sumn{j,k=1}a_{j,2} \bs{f}_{j,k}a_{k,2}\right)\left(S_{1,q}\right)\\&
-\delta_{1,q}\left(\Matrix{~\\~\\~\\~}\underline{a_{n,1}\bs{f}_{n,n}a_{n,1}}
\right.\\&\left.+\Sumn{j,k=1\\(j,k)\neq (n,n)}a_{j,1} \bs{f}_{j,k} a_{k,1} +\left(Q+\Sumn{j,k=1}a_{j,2} \bs{f}_{j,k}a_{k,2}\right)\right)
\\&
\NeueSeite
\Pfeil{\widetilde{V}}
&-\Sumn{l=3}a_{n,1}\bs{f}_{n,n}a_{n,l}S_{l,q}\\&
- \Sumn{j,k=1\\(j,k)\neq (n,n)}a_{j,1} \bs{f}_{j,k} \left( \Sumn{i=2}a_{k,i}S_{i,q}-\delta_{k,1}a_{k,q}\right)\\&
+\left(Q+\Sumn{j,k=1}a_{j,2} \bs{f}_{j,k}a_{k,2}\right)\left(S_{1,q}\right)\\&
-\delta_{1,q}\left(
 -\left(Q + \Sumn{j,k=1} a_{j,2} \bs{f}_{j,k}a_{k,2}  + \Sumn{j,k=1\\(j,k)\neq(n,n)}a_{j,1} \bs{f}_{j,k}a_{k,1}\right)
\right.\\&\left.
+\Sumn{j,k=1\\(j,k)\neq (n,n)}a_{j,1} \bs{f}_{j,k} a_{k,1} +\left(Q+\Sumn{j,k=1}a_{j,2} \bs{f}_{j,k}a_{k,2}\right)
\right)
\\
=
&-\Sumn{l=3}a_{n,1}\bs{f}_{n,n}a_{n,l}S_{l,q}\\&
- \Sumn{j,k=1\\(j,k)\neq (n,n)}a_{j,1} \bs{f}_{j,k} \left( \Sumn{i=2}a_{k,i}S_{i,q}-\delta_{k,1}a_{k,q}\right)\\&
+\left(Q+\Sumn{j,k=1}a_{j,2} \bs{f}_{j,k}a_{k,2}\right)\left(S_{1,q}\right)
}

\subsubsection{Beginnend mit $(a_{n,1}\bs{f}_{n,n}a_{n,2}a_{1,2} a_{1,q} ~ , ~ W_1a_{1,q})$:}

\Gleichung{ 
&\underline{a_{n,1}\bs{f}_{n,n}a_{n,2}a_{1,2} } a_{1,q} \\ ~
\Pfeil{\widetilde{W}_{1} } & 
 - a_{n,1}\bs{f}_{n,n}\left(\Sumn{l=3}a_{n,l}a_{1,l}  -\delta_{n,1} \right)a_{1,q}
 \\&
 +	\Sumn{j,k=1\\(j,k)\neq (n,n)}a_{j,1} \bs{f}_{j,k}Z_{k,1}a_{1,q} 
	\\&
		+\left(Q+\Sumn{j,k=1}a_{j,2} \bs{f}_{j,k}a_{k,2}\right)a_{1,1} a_{1,q} 
\\
=&
- a_{n,1}\bs{f}_{n,n}\Sumn{l=3}a_{n,l}\underline{a_{1,l}  a_{1,q}}
 \\&
 -	\Sumn{j,k=1\\(j,k)\neq (n,n)}a_{j,1} \bs{f}_{j,k}\left(\Sumn{i=2}a_{k,i}\underline{a_{1,i}a_{1,q}}
 -\delta_{k,1}a_{1,q} \right)
	\\&
		+\left(Q+\Sumn{j,k=1}a_{j,2} \bs{f}_{j,k}a_{k,2}\right)\underline{a_{1,1} a_{1,q}}
\\
\Pfeil{\widetilde{S}_{l,q},\widetilde{S}_{i,q}, \widetilde{S}_{1,q} } & 
- a_{n,1}\bs{f}_{n,n}\Sumn{l=3}a_{n,l}S_{l,q}
 \\&
 -	\Sumn{j,k=1\\(j,k)\neq (n,n)}a_{j,1} \bs{f}_{j,k}\left(\Sumn{i=2}a_{k,i}S_{i,q}
 -\delta_{k,1}a_{1,q} \right)
	\\&
		+\left(Q+\Sumn{j,k=1}a_{j,2} \bs{f}_{j,k}a_{k,2}\right)S_{1,q}
}

%% file: Haupt/RechnungAufloesung/Phi3/W1S2qr_phi3.tex
%W1 & S2qr phi3

\subsection{Reduktionswege für $a_{n,1}\bs{f}_{n,n}a_{n,2}a_{1,2} a_{2,q} a_{2,r}$}
\label{W1S2qr_phi3}

Um die Lesbarkeit zu erhöhen setzen wir:
\Gleichung{\bs{X}:=&-\Sumn{l=3}a_{n,1}\bs{f}_{n,n}a_{n,l}\left(S_{l,q,r}-S_{1,q}a_{1,r}\right)\\&
-	\Sumn{j,k=1\\(j,k)\neq (n,n)}a_{j,1} \bs{f}_{j,k} Z_{k,1}\left(\Sumn{i=3}a_{i,q}a_{i,r}-\delta_{q,r}\right)\\&
	+\left(Q+\Sumn{j,k=1}a_{j,2}\bs{f}_{j,k}a_{k,2}\right)\left(S_{1,q,r}-S_{1,q}a_{1,r}\right).}

\subsubsection{Beginnend mit $(a_{n,1}\bs{f}_{n,n}a_{n,2}a_{1,2} a_{2,q} a_{2,r} ~ , ~ a_{n,1}\bs{f}_{n,n}a_{n,2}\left(S_{2,q,r}-S_{2,q}a_{1,r}\right))$:}

\Gleichung{
&a_{n,1}\bs{f}_{n,n}a_{n,2}\underline{a_{1,2} a_{2,q} a_{2,r}} \\ 
\Pfeil{\widetilde{S}_{2,q,r} } & 
-\underline{a_{n,1}\bs{f}_{n,n}a_{n,2}a_{1,2}}\left(\Sumn{i=3}a_{i,q}a_{i,r}
-\delta_{q,r}\right)
 \\&
+\Sumn{i=2}\underline{a_{n,1}\bs{f}_{n,n}a_{n,2}a_{i,2}}a_{i,q}a_{1,r}
\\&
-a_{n,1}\bs{f}_{n,n}a_{n,2}\delta_{2,q}a_{1,r}\\
\Pfeil{\widetilde{W}_{1}, \widetilde{W}_{i}  } & 
 \left(\Matrix{~\\~\\~\\~\\~}a_{n,1}\bs{f}_{n,n}\left(\Sumn{l=3}a_{n,l}a_{1,l}  -\delta_{n,1} \right)
 \right.\\&\left.
 -	\Sumn{j,k=1\\(j,k)\neq (n,n)}a_{j,1} \bs{f}_{j,k} Z_{k,1} \right.\\&\left.
	-\left(Q+\Sumn{j,k=1}a_{j,2} \bs{f}_{j,k}a_{k,2}\right)a_{1,1}\Matrix{~\\~\\~\\~\\~} \right)\left(\Sumn{i=3}a_{i,q}a_{i,r}
-\delta_{q,r}\right)
 \NeueSeite
 &- \Sumn{i=2}a_{n,1}\bs{f}_{n,n}\left(\Sumn{l=3}a_{n,l}a_{i,l}  -\delta_{n,i} \right)a_{i,q}a_{1,r}
 \\&
 +\Sumn{i=2}	\Sumn{j,k=1\\(j,k)\neq (n,n)}a_{j,1} \bs{f}_{j,k} Z_{k,i}a_{i,q}a_{1,r}
	\\&
		+\Sumn{i=2}\left(Q+\Sumn{j,k=1}a_{j,2} \bs{f}_{j,k}a_{k,2}\right)a_{i,1} 
a_{i,q}a_{1,r}
\\&
-a_{n,1}\bs{f}_{n,n}a_{n,2}\delta_{2,q}a_{1,r}
\\
 =&
 -\Sumn{l=3}a_{n,1}\bs{f}_{n,n}a_{n,l}S_{l,q,r}\\&
-	\Sumn{j,k=1\\(j,k)\neq (n,n)}a_{j,1} \bs{f}_{j,k} Z_{k,1}\left(\Sumn{i=3}a_{i,q}a_{i,r}-\delta_{q,r}\right)\\&
	+\left(Q+\Sumn{j,k=1}a_{j,2}\bs{f}_{j,k}a_{k,2}\right)S_{1,q,r}
 \\&
 - a_{n,1}\bs{f}_{n,n}\Sumn{l=3}a_{n,l}\left(\Sumn{i=2}a_{i,l}  a_{i,q}\right)a_{1,r}\\&
 + \Sumn{i=2}a_{n,1}\bs{f}_{n,n}\left(\delta_{n,i} \right)a_{i,q}a_{1,r}
 \\&
 +\Sumn{i=2}	\Sumn{j,k=1\\(j,k)\neq (n,n)}a_{j,1} \bs{f}_{j,k} Z_{k,i}a_{i,q}a_{1,r}
	\\&
		+\left(Q+\Sumn{j,k=1}a_{j,2} \bs{f}_{j,k}a_{k,2}\right)\left(\Sumn{i=2}a_{i,1} 
a_{i,q}\right)a_{1,r}
\\&
-a_{n,1}\bs{f}_{n,n}a_{n,2}\delta_{2,q}a_{1,r}\\
=&
\bs{X}	%
 -\Sumn{l=3}a_{n,1}\bs{f}_{n,n}a_{n,l}\delta_{l,q}a_{1,r} -a_{n,1}\bs{f}_{n,n}a_{n,2}\delta_{2,q}a_{1,r}\\&
  + \Sumn{i=2}a_{n,1}\bs{f}_{n,n}\left( \delta_{n,i} \right)a_{i,q}a_{1,r}
 \\&
 +	\Sumn{j,k=1\\(j,k)\neq (n,n)}a_{j,1} \bs{f}_{j,k} \underline{\Sumn{i=2}Z_{k,i}a_{i,q}}a_{1,r}
	\\&
	+\left(Q+\Sumn{j,k=1}a_{j,2} \bs{f}_{j,k}a_{k,2}\right)\left(\delta_{1,q}\right)a_{1,r}
\NeueSeite
\stackrel{\ref{Rechenregel}}{=}&
\bs{X}	+\delta_{q,1}\underline{a_{n,1}\bs{f}_{n,n}a_{n,1}}a_{1,r}\\
&
 +	\Sumn{j,k=1\\(j,k)\neq (n,n)}a_{j,1} \bs{f}_{j,k} \left(\Sumn{i=2}a_{k,i}S_{i,q}+(\delta_{1,q}-\delta_{k,1})a_{k,q}\right)a_{1,r}
	\\&
+\delta_{1,q}	\left(Q+\Sumn{j,k=1}a_{j,2} \bs{f}_{j,k}a_{k,2}\right)a_{1,r}
\\
\Pfeil{\widetilde{V}}&
\bs{X}	-\delta_{q,1}    \left(Q + \Sumn{j,k=1} a_{j,2} \bs{f}_{j,k}a_{k,2}  
	\right.\\&\left. ~~~~~~~~~~~~~~~~
+ \Sumn{j,k=1\\(j,k)\neq(n,n)}a_{j,1} \bs{f}_{j,k}a_{k,1} \right) a_{1,r}\\&
 +	\Sumn{j,k=1\\(j,k)\neq (n,n)}a_{j,1} \bs{f}_{j,k} \left(\Sumn{i=2}a_{k,i}S_{i,q}a_{1,r}-\delta_{k,1} \underline{a_{1,q}a_{1,r}}\right)\\&
	+\delta_{1,q}\left(	\Sumn{j,k=1\\(j,k)\neq (n,n)}a_{j,1} \bs{f}_{j,k} a_{k,1}
	\right.\\&\left. ~~~~~~~~~~~~~~~~
	+Q+\Sumn{j,k=1}a_{j,2} \bs{f}_{j,k}a_{k,2}\right)a_{1,r}
\\
\Pfeil{\widetilde{Z}_{q,r}}&
\bs{X}	
 +	\Sumn{j,k=1\\(j,k)\neq (n,n)}a_{j,1} \bs{f}_{j,k} \left(\Sumn{i=2}a_{k,i}S_{i,q}a_{1,r}-\delta_{k,1} Z_{q,r}\right)
}

\subsubsection{Beginnend mit $(a_{n,1}\bs{f}_{n,n}a_{n,2}a_{1,2} a_{2,q} a_{2,r} ~ , ~ )$:}

\Gleichung{ 
&\underline{a_{n,1}\bs{f}_{n,n}a_{n,2}a_{1,2} } a_{2,q} a_{2,r}\\ ~
\Pfeil{\widetilde{W}_{1} } & 
 - a_{n,1}\bs{f}_{n,n}\left(\Sumn{l=3}a_{n,l}a_{1,l}  -\delta_{n,1} \right)a_{2,q} a_{2,r}
 \\&
 +	\Sumn{j,k=1\\(j,k)\neq (n,n)}a_{j,1} \bs{f}_{j,k} Z_{k,1}a_{2,q} a_{2,r}
	\\&
		+\left(Q+\Sumn{j,k=1}a_{j,2} \bs{f}_{j,k}a_{k,2}\right)\underline{a_{1,1} a_{2,q} a_{2,r}}
\\
\Pfeil{\widetilde{S}_{1,q,r} } & 
 - a_{n,1}\bs{f}_{n,n}\Sumn{l=3}a_{n,l}\underline{a_{1,l}  a_{2,q} a_{2,r}}
 %+\delta_{n,1} a_{n,1}\bs{f}_{n,n}a_{2,q} a_{2,r}
 \\&
 +	\Sumn{j,k=1\\(j,k)\neq (n,n)}a_{j,1} \bs{f}_{j,k} Z_{k,1}a_{2,q} a_{2,r}
	\\&
		+\left(Q+\Sumn{j,k=1}a_{j,2} \bs{f}_{j,k}a_{k,2}\right)\left(S_{1,q,r}-S_{1,q}a_{1,r}\right)
\\
\Pfeil{\widetilde{S}_{l,q,r} } & 
 - a_{n,1}\bs{f}_{n,n}\Sumn{l=3}a_{n,l}\left(S_{l,q,r}-S_{l,q}a_{1,r}\right)
 %+\delta_{n,1} a_{n,1}\bs{f}_{n,n}a_{2,q} a_{2,r}
 \\&
 +	\Sumn{j,k=1\\(j,k)\neq (n,n)}a_{j,1} \bs{f}_{j,k} \underline{Z_{k,1}a_{2,q} a_{2,r}}
	\\&
		+\left(Q+\Sumn{j,k=1}a_{j,2} \bs{f}_{j,k}a_{k,2}\right)\left(S_{1,q,r}-S_{1,q}a_{1,r}\right)\\
\stackrel{\ref{ZiSiil}}{=}&
 - a_{n,1}\bs{f}_{n,n}\Sumn{l=3}a_{n,l}\left(S_{l,q,r}-S_{l,q}a_{1,r}\right)
 \\&
 +	\Sumn{j,k=1\\(j,k)\neq (n,n)}a_{j,1} \bs{f}_{j,k}
\left(\Sumn{i=2}\Sumn{l=3}a_{k,i}a_{1,i}a_{l,q}a_{l,r}+Z_{k,1}\delta_{q,r}
\right.\\&\left.+\Sumn{i=2}a_{k,i}S_{i,q}a_{1,r}-\delta_{k,1}\left(\delta_{q,r}-a_{2,q}a_{2,r}\right) \right)
	\\&
		+\left(Q+\Sumn{j,k=1}a_{j,2} \bs{f}_{j,k}a_{k,2}\right)\left(S_{1,q,r}-S_{1,q}a_{1,r}\right)
\NeueSeite
=&
 - a_{n,1}\bs{f}_{n,n}\Sumn{l=3}a_{n,l}\left(S_{l,q,r}-S_{l,q}a_{1,r}\right)
 \\&
 +	\Sumn{j,k=1\\(j,k)\neq (n,n)}a_{j,1} \bs{f}_{j,k}
\left(\Sumn{i=2}\Sumn{l=3}a_{k,i}a_{1,i}a_{l,q}a_{l,r}+Z_{k,1}\delta_{q,r}
\right.\\&\left.+\Sumn{i=2}a_{k,i}S_{i,q}a_{1,r}-\delta_{k,1}\left(\delta_{q,r}-a_{2,q}a_{2,r}\right) \right)
	\\&
		+\left(Q+\Sumn{j,k=1}a_{j,2} \bs{f}_{j,k}a_{k,2}\right)\left(S_{1,q,r}-S_{1,q}a_{1,r}\right)\\
		=&
		 - a_{n,1}\bs{f}_{n,n}\Sumn{l=3}a_{n,l}\left(S_{l,q,r}-S_{l,q}a_{1,r}\right)
 \\&
 		+\left(Q+\Sumn{j,k=1}a_{j,2} \bs{f}_{j,k}a_{k,2}\right)\left(S_{1,q,r}-S_{1,q}a_{1,r}\right)
\\&
 -	\Sumn{j,k=1\\(j,k)\neq (n,n)}a_{j,1} \bs{f}_{j,k}
\left(\Sumn{l=3}\left(Z_{k,1}-\delta_{k,1}\right)a_{l,q}a_{l,r}-Z_{k,1}\delta_{q,r}\right)\\&
+\Sumn{j,k=1\\(j,k)\neq (n,n)}a_{j,1} \bs{f}_{j,k}\Sumn{i=2}a_{k,i}S_{i,q}a_{1,r}\\&
-\Sumn{j,k=1\\(j,k)\neq (n,n)}a_{j,1} \bs{f}_{j,k}\delta_{k,1}\left(\delta_{q,r}-a_{2,q}a_{2,r}\right)
	\\
	=&
	\bs{X}
 +	\Sumn{j,k=1\\(j,k)\neq (n,n)}a_{j,1} \bs{f}_{j,k}
 \delta_{k,1}\left(\Sumn{l=3}a_{l,q}a_{l,r}-\left(\delta_{q,r}-a_{2,q}a_{2,r}\right)\right)\\&
+\Sumn{j,k=1\\(j,k)\neq (n,n)}a_{j,1} \bs{f}_{j,k}\Sumn{i=2}a_{k,i}S_{i,q}a_{1,r}\\
	=&
	\bs{X}
 -	\Sumn{j,k=1\\(j,k)\neq (n,n)}a_{j,1} \bs{f}_{j,k}
 \left(\delta_{k,1}Z_{q,r}
-\Sumn{i=2}a_{k,i}S_{i,q}a_{1,r}\right)
	}

%% file: Haupt/HochschildAufloesung/HochschildAufloesung.tex
%\section{Übersichtstabellen für $\K_{\Lambda}$}
%\label{Tabellen_K_Lambda}
%
%
%
%hier fehlen noch die Übersichtstabellen

\newpage

\chapter{(Ko)-Homologie}
\label{sec:HochschildKoHomologie}

Wir wollen im Folgenden 1-dimensionale Moduln über der orthogonalen freien Quantengruppe $\fQG$ betrachten.\\

Sei $\fQG(n)$ die in Abschnitt \ref{sec:GröbnerbasisFürDieFreieQuantenGruppeAN} definierte freie Quantengruppe. Zur Erinnerung: Wir können die Erzeuger von $\fQG$  schreiben als $n\times n$-Matrix $A:=(a_{p,q})_{p,q=1,\dots, n}$ und die Relationen als $AA^t=A^tA=\id$.

\Def{1-dimensionaler $\fQG$-Modul $\K_{\Lambda}$}{
 Ein $\fQG$-Modul, zu dem es eine 1-dimensionale Darstellung gibt, die gegeben ist durch einen Algebra-Homomorphismus :
\MatheAbbildung{ & \fQG}{&a_{j,i}}{\K}{\lambda_{j,i},}
heißt 1-dimensionaler Modul. Die Modulstruktur ist gegeben durch:
\MatheAbbildung{ \fQG \times \K}{a_{j,i},x}{\K}{ \lambda_{j,i} \cdot x.}
Diesen Modul bezeichnen wir mit $\K_{\Lambda}$. Da die Multiplikation in $\K$ kommutiert, ist jeder 1-dimensionale Modul immer auch ein Bimodul.
}

Die Bilder der Erzeuger schreiben wir wieder als Matrix: $\Lambda:=(\lambda_{p,q})_{p,q=1,\dots, n}$. 
Für diese Matrix gelten%, induziert durch den Algebren-Homomorphismus,
 ähnliche Relationen wie für die Erzeuger-Matrix $A$: 
\mathe{\Lambda\Lambda^t=\Lambda^t \Lambda = \id;}
sie ist also orthogonal.\\
Umgekehrt entspricht jede orthogonale Matrix $M$ einem 1-dimensionalen $\fQG$-Modul $\K_{M}$.\\
~\\
\goodbreak

\newpage
\input{Haupt/HochschildAufloesung/Speziallfall_n=1}
\newpage

\subsection{Auflösung für $n\geq 3$}
\label{sec:FallN1}

Als Reduktionssystem $\red_{\Algebra}$ für die Algebra $\fQG(n)$ nehmen wir das in Satz \ref{GBA} definierte:

\mathe{\begin{array}{rl}
	\widetilde{Z}_{p,q}:=&\left(a_{p,1} a_{q,1} ~,~ \Zi{p}{q}\right) %&\textnormal{ mit } &Z_{p,q} :=& 
	,\\
	\widetilde{S}_{p,q}:=&\left(a_{1,p} a_{1,q}~,~ \Si{p}{q}\right) %&\textnormal{ mit } &S_{p,q} :=& \Si{p}{q}
	,\\
	\widetilde{Z}_{p,q,r}:=&\left(a_{p,1} a_{q,2} a_{r,2}~,~ \Zii{p}{q}{r} -\left(  \Zi{p}{q} \right)a_{r,1}\right) % &\textnormal{ mit } &Z_{p,q,r} :=& \Zii{p}{q}{r} - \Zi{p}{q} a_{r,1}
	,\\
	\widetilde{S}_{p,q,r}:=&\left(a_{1,p} a_{2,q} a_{2,r}~,~ \Sii{p}{q}{r} - \left( \Si{p}{q} \right)a_{1,r}\right)% &\textnormal{ mit } &S_{p,q,r} :=& \Sii{p}{q}{r} - \Si{p}{q} a_{1,r}
	.
\end{array}
}

Alle Moduln in der Bimodulauflösung in Abschnitt \ref{sec:AuflösungVonFQGN} sind freie $\fQG$-Moduln. Wenn wir also 
 die Bimodulauflösung  von rechts mit $\tensor_{\fQG} \K_{\lambda}$ multiplizieren, ist die resultierende Sequenz wieder exakt. Mit $\fQG \tensor_{\fQG} \K_{\lambda} \iso \K_{\lambda}$ erhalten wir die folgende Sequenz, bzw. die Auflösung von $\K_{\lambda}$:

\label{Aufloesung_K_Lambda}
\xymatrix{				 
	\ar[d] 0  
\\
 	\ar@{_{(}->}[d]_{\phi_3} \fQG 		 				&  \ar@{|->}[d] \bs{f}								&			 															
\\
	\ar[d]_{\phi_2} \left(\fQG\right)^{n^2} 	&  -\sum\limits_{i,j,k=1}^n a_{j,i}\bs{f}_{j,k} \lambda_{k,i}
																								+ \sum\limits_{i=1}^n \bs{f}_{i,i}	& \ar@{|->}[d] \bs{f}_{p,q}
\\
	\ar[d]_{\phi_1} \left(\fQG\right)^{n^2}  &  \ar@{|->}[d] \bs{e}_{p,q}	&\sum\limits_{i=1}^n\left(a_{p,i}\bs{e}_{q,i}																																								+\bs{e}_{p,i}\lambda_{q,i}\right)
\\
	\ar@{->>}[d]_{\phi_0} \fQG 							& a_{p,q}\bs{e}- \bs{e}\lambda_{p,q}& \ar@{|->}[d] a_{p,q}	
\\
	\ar[d] \K_{\Lambda}															&															& \lambda_{p,q}
 \\
	0 				 
		%					%
							%	
	}
	
%Die Exaktheit beweisen wir mit den gleichen Rechnungen, wie für $\fQG \otimes \fQG^{op}$-Moduln \ref{sec:AuflösungVonFQGN}. Die folgenden Tabellen geben die Möglichen Konflikte an. Um sie zu beheben brauchen wir nur in den alten Rechnungen die $a_{i,j}$ rechts vom Tesorzeichen durch $\lambda_{i,j}$ ersetzen.

\MitText{
Die fett geschriebenen Buchstaben $\bs{e},\bs{e}_{p,q}, \bs{f}_{p,q}$ und $\bs{f}$ stehen für die Erzeuger der Moduln. 
Die Elemente aus $\K_{\Lambda}$ kommutieren mit denen aus $\fQG$, daher können wir sie vorziehen.

Zur besseren Übersicht schreiben wir die Auflösung auch in Matrizen-Schreibweise.
Mit den Buchstaben $\bs{E}$ und $\bs{F}$ bezeichnen wir die Matrizen der entsprechenden Erzeuger, also $\bs{E}:=(\bs{e}_{p,q})_{p,q=1, \dots , n}$ und $\bs{F}:=(\bs{f}_{p,q})_{p,q=1, \dots , n}$.
}%mitText
\mathe{\begin{array}{|l|l|}
\hline
 \Matrix{~\\~}\textnormal{Koeffizienten vorne} & \textnormal{Koeffizienten hinten} 
\\ \hline
 \Matrix{~\\~}\phi_3: \bs{f} \mapsto -\tr(\Lambda A^t\bs{F})+\tr(\bs{F})& \phi_3: \bs{f} \mapsto -\tr(A^t\bs{F}\Lambda)+\tr(\bs{F})
\\ \hline
 \Matrix{~\\~}	\phi_2: \bs{F} \mapsto A\bs{E}^t+\left(\Lambda\bs{E}^t\right)^t& \phi_2: \bs{F} \mapsto A\bs{E}^t+\bs{E}\Lambda^t
\\ \hline
 \Matrix{~\\~}	\phi_1: \bs{E} \mapsto \left(A-\Lambda\right)\bs{e}		& \phi_1: \bs{E} \mapsto A\bs{e}-\bs{e}\Lambda
\\ \hline
 \Matrix{~\\~}	 \phi_0: a_{i,j}\bs{e} \mapsto \lambda_{i,j}				& \phi_0: a_{i,j}\bs{e} \mapsto \lambda_{i,j}
\\ \hline
\end{array}
}
\MitText{
Hierbei beachten wir, dass das Transponieren nicht wie üblich mit der Multiplikation verträglich ist:
\mathe{ \left( A B \right)^t_{p,q}=\sum_{i=1}^{n} a_{q,i}b_{i,p} \textnormal{ aber } 
 \left( B^t A^t \right)_{p,q} = \sum_{i=1}^n b_{i,p}a_{q,i}. }

Im Spezialfall $\Lambda=\id$ erhalten wir, dass die in \cite{ThomCollins} vorgestellte Sequenz exakt ist:

\xymatrix{ 0 	& \ar[l]  \K_{\id}
							& \ar@{->>}[l]_{} \fQG
							& \ar[l]_{} \fQG^{n^2}
							& \ar[l]_{} \fQG^{n^2}
							& \ar@{_{(}->}[l]_-{} \fQG 
							& \ar[l] 0   
							\\
							& \delta_{p,q}
							& \ar@{|->}[l] a_{p,q}
							& \sum\limits_{i=1}^n a_{p,i} \bs{e}_{q,i} + \bs{e}_{p,q}
							& \ar@{|->}[l]\bs{f}
							\\
							&	
							& a_{p,q}-\delta_{p,q}
							& \ar@{|->}[l] \bs{e}_{p,q}
							& -\sum\limits_{j,k=1}^n a_{j,k} \bs{f}_{j,k} + \sum\limits_{i=1}^n \bs{f}_{i,i}
							&  \ar@{|->}[l] \bs{f}
	}

Um die Dimension der Homologie bzw. Kohomologie zu berechnen betrachten wir die Auflösung von $\K_{\Lambda}$, ohne die erste nicht triviale Stelle unter dem Funktor $\left( \K_{\Omega} \otimes_{\fQG} \bullet \right)$ bzw. dem Funktor $\Hom_{\fQG}\left( \bullet, \K_{\Omega} \right)$. Mit anderen Worten, wir berechnen $\Tor{\K_{\Omega}, \K_{\Lambda}}$ bzw. $\Ext{\K_{\Omega},\K_{\Lambda}}$.

\newpage
}%mitText
 \input{Haupt/HochschildAufloesung/Tor}

%% file: Haupt/HochschildAufloesung/Speziallfall_n=1.tex
\subsection{Spezialfall $n=1$}
\label{sec:SpezialfallN1}

Falls $n=1$, dann hat die Algebra $\fQG(1)$ nur einen Erzeuger $a_{1,1}$ und nur eine Relation $a_{1,1}a_{1,1}=1$. Diese Algebra ist isomorph zu der Algebra $\K[z]/(z^2-1)$.

Wir wollen uns zunächst überlegen, warum diese Algebra  halbeinfach ist, falls $char(\K) \neq 2$.

Seien $e_1 =\frac{z+1}{2}$ und $e_2=\frac{z-1}{2}$, dann sind $e_1$ und $e_2$ orthogonale Idempotenten, so dass $e_1 + e_2 =1$.

Insbesondere
\mathe{
\fQG(1)\cong \K<e_1> \oplus \K<e_2>.
}

Wir bemerken weiterhin, dass es zu $\fQG(1)$ nur zwei eindimensionale Moduln gibt:
Einen Modul $\K_1$ zu $e_1$ und einen Modul $\K_2$ zu $e_2$. Diese Moduln sind beide projektiv als direkte Summanden des freien  Moduls $\fQG(1)$.
%Also ist die Algebra $\fQG(1)$ halbeinfach und eine projektive Auflösung dieser eindimensionalen Moduln ist von folgender Form:
%\Gleichung{
%0 \leftarrow \K_1 \leftarrow \fQG(1) \leftarrow \K_2 \leftarrow 0\\
%0 \leftarrow \K_2 \leftarrow \fQG(1) \leftarrow \K_1 \leftarrow 0
%}

Da sie projektiv sind, verschwinden alle höheren Ext- und Torgruppen zwischen ihnen. % Die Ko- und Homologie verschwindet also auch.
 
~\\

Falls $char(\K)= 2$, dann gibt es nur einen eindimensionalen Modul $\K$, der durch folgende eindimensionale  Darstellung gegeben ist:

\mathe{\LMatrix{
\phi\colon &\K[z]/(z^2-1) &\to \K\\
& z &\mapsto 1\\
& 1 &\mapsto 1. }
}
 
Dieser Modul hat eine unendliche freie Auflösung
\mathe{
0 \leftarrow \K \leftarrow \fQG(1) \leftarrow \cdots \leftarrow \fQG(1) \leftarrow \cdots,
}
wobei jedes Differential durch Multiplikation mit $z+1$
gegeben ist.

Für die Tor- und Extgruppen gilt also:
\Gleichung{
 Ext^k(\K,\K)\cong \K &\textnormal{ für alle } k\geq0 , \\
 Tor_k(\K,\K) \cong \K &\textnormal{ für alle } k\geq0
 .}

%% file: Haupt/HochschildAufloesung/Tor.tex
\section{Tor}
\label{sec:Tor}
In diesem Kapitel wollen wir die Dimensionen der Homologie über die Tor-Gruppen berechnen.
\Satz{Tor}{
Sei  $\HH_i$ die Homologie des eindimensionalen Moduls $\K_{\Lambda}$, wobei $i=1,2,3$. Dann gelten für die Dimensionen folgende Werte:
\mathe{\begin{array}{|l|c|c|c|}
\hline
\Matrix{~\\ ~\\~} & \Omega= \Lambda & \Omega = -\Lambda & \textnormal{sonst} \\ \hline
\Matrix{~\\ \Dim{\HH_0}\\~}	& 1 & 0 & 0 \\ \hline
\Matrix{~\\ \Dim{\HH_1}\\~} & \frac{n^2-n}{2} & \frac{n^2+n-2}{2} & k_{-1}+k_{\bigwedge}-1 \\ \hline
\Matrix{~\\ \Dim{\HH_2}\\~} & \frac{n^2-n}{2} & \frac{n^2+n-2}{2} & k_{-1}+k_{\bigwedge}-1\\ \hline
\Matrix{~\\ \Dim{\HH_3}\\~} & 1& 0&0 \\ \hline
\end{array}
,
}%mathe
mit $k_{-1}:=\# \Menge{i}{r_i=-1}$ und $k_{\bigwedge}:=\# \Menge{(i,j)}{i>j, r_ir_j=1}$, wobei $r_i$ die Vielfachheit der Nullstellen des Charakteristischenpolynoms von $\Omega\Lambda^t$ sind.
}%Satz

Den Komplex erhalten wir durch Anwenden des Funktors $\left( \K_{\Omega} \otimes_{\fQG} \bullet \right)$ auf  die in \ref{Aufloesung_K_Lambda} gegebene Auflösung. Sie sieht wegen $\K_{\Omega} \otimes_{\fQG} \fQG \iso \K$ wie folgt aus:
\mathe{\begin{array}{|l|}
\hline  \Matrix{~\\~~\\~} 
0 \leftarrow \K
\stackrel{\phi_{1*}}{\longleftarrow} \K^{n^2} 
\stackrel{\phi_{2*}}{\longleftarrow} \K^{n^2} 
\stackrel{\phi_{3*}}{\longleftarrow} \K 
 \leftarrow 0
\\ \hline
 \Matrix{~\\~~\\~} \phi_{1*}: \bs{E} \mapsto \left(\Omega-\Lambda\right)\bs{e}
\\ \hline
 \Matrix{~\\~~\\~}  \phi_{2*}: \bs{F} \mapsto \Omega\bs{E}^t+\left(\Lambda\bs{E}^t\right)^t
\\ \hline
 \Matrix{~\\~~\\~}  \phi_{3*}: \bs{f} \mapsto -\tr(\Lambda \Omega^t\bs{F})+\tr(\bs{F})
\\ \hline
\end{array}
}%mathe

Die Homologie bekommen wir aus dem Komplex durch:
\mathe{ \HH_i:=  \Kern{\Phi_{i*}}/\Bild{\Phi_{i+1*}}.}
Wir berechnen die Dimensionen über:
\mathe{
\Dim{\Kern{\Phi_{i*}}/ \Bild{\Phi_{i+1*}}}  = \Dim{\Kern{\Phi_{i*}}}-\Rang{i+1} }

\newpage
 \input{Haupt/HochschildAufloesung/Phi1.tex}\newpage
\input{Haupt/HochschildAufloesung/Phi2.tex}\newpage
\input{Haupt/HochschildAufloesung/Phi3.tex}

%% file: Haupt/HochschildAufloesung/Phi1.tex
\subsection{$\Phi_{1*}: \bs{E} \mapsto \Omega\bs{e}-\Lambda\bs{e}$}
\label{sec:Phi1Stern}

\mathe{\begin{array}{|c|c|}
\hline  
\Dim{\Kern{\Phi_{1*}}} & \Rang{1} \Matrix{~\\~}\\
\hline
\begin{array}{c|c} 
 \Omega=\Lambda & \textnormal{sonst} \\
\hline 
n^2 &n^2-1
\end{array}
& \begin{array}{c|c} 
\Omega=\Lambda & \textnormal{sonst} \\
\hline 
0 &1
\end{array}\Matrix{~\\~\\~}\\
 \hline
\end{array}
}%mathe

Sei zunächst $\Omega=\Lambda$, dann ist $\Bild{\Phi_{1*}}=0$ und der $\Kern{\Phi_{1*}}$ ist das gesamte Urbild, 
also ist $\Rang{1}=0$ und $\Dim{\Kern{\Phi_{1*}}}=n^2$.

Falls $\Omega\neq\Lambda$, dann existiert wenigstens ein $(p,q)$ mit $\omega_{p,q}\neq\lambda_{p,q}$. Dann ist  $\Phi_{1*}(e_{p,q})\neq 0$, also gilt $\Rang{1}>0$ . Da der Bildraum nur eindimensional ist, gilt $\Rang{1}=1$ und mit $\Dim{\Kern{\Phi_{1*}}}= \Dim{\textnormal{Urbild}}-\Dim{\textnormal{Bild}}=n^2-1$.

%% file: Haupt/HochschildAufloesung/Phi2.tex
\subsection{$\Phi_{2*}:  \bs{F} \mapsto \Omega\bs{E}^t+\left(\Lambda\bs{E}^t\right)^t$}
\label{sec:Phi2Stern}

\mathe{\begin{array}{|c|c|}
\hline  
\Dim{\Kern{\Phi_{2*}}} & \Rang{2}\Matrix{~\\~}\\
\hline
\begin{array}{c|c|c} 
\Omega=\Lambda &\Omega=-\Lambda & \textnormal{sonst} \\
\hline 
\frac{n^2-n}{2}&\frac{n^2+n}{2}&k_{-1}+k_{\bigwedge}
\end{array}
& \begin{array}{c|c|c} 
\Omega=\Lambda & \Omega=-\Lambda &\textnormal{sonst} \\
\hline 
\frac{n^2+n}{2}&\frac{n^2-n}{2}&n^2-(k_{-1}+k_{\bigwedge})
\end{array}\Matrix{~\\~\\~}\\
 \hline
\end{array}
}%mathe

Mit $k_{-1}:=\# \Menge{i}{r_i=-1}$ und $k_{\bigwedge}:=\# \Menge{(i,j)}{i>j, r_ir_j=1}$, die Vielfachheit der Nullstellen des Charakteristischenpolynoms von $\Omega\Lambda^t$ sind.

Wir werden nun $\Phi_{2*}$%:\Omega E + (\Lambda^t E)^t$ 
als $n^2 \times n^2$-Matrix schreiben. Da wir von links multiplizieren, steht in der Zeile $(p-1)n+q$ das Bild von $f_{p,q}$ bezüglich der Basis des Bildraumes:
\mathe{ \Phi_{2*}(f_{p,q}) = \Sumn{i=1} \left( \omega_{p,i} e_{q,i}+ \lambda_{q,i} e _{p,i} \right) .}

Um uns Schreibarbeit zu sparen schreiben wir $\lambda_{p,\cdot}$ bzw. $\omega_{p,\cdot}$ für den Zeilenvektor $\left(\lambda_{p,1}, \lambda_{p,2},\dots ,\lambda_{p,n}\right)$ bzw. $\left(\omega_{p,1}, \omega_{p,2},\dots ,\omega_{p,n}\right)$. Den Nullvektor schreiben wir als $0_n:=(\underbrace{0,\dots,0}_{n})$.\\
Die Matrix $\Lambda$ können wir dann schreiben als:
\mathe{ \Lambda= \left(\Matrix{\lambda_{1,\cdot}\\ \lambda_{2,\cdot} \\ \vdots}\right).}

Die Matrix $\Phi_{2*}$ ist dann:

\mathe{\Matrix{ 
\Matrix{
f_{1,\cdot} & \Matrix{f_{1,1}\\ f_{1,2} \\ \vdots \\ f_{1,n}}
\\ \hline
f_{2,\cdot} & \Matrix{f_{2,1}\\ f_{2,2} \\ \vdots \\ f_{2,n}}
\\ \hline
						& \vdots \\ & \vdots	
\\\hline
f_{n,\cdot} & \Matrix{f_{n,1} \\ f_{n,2}\\ \vdots \\ f_{n,n}}
}
&
\left(\Matrix{
\BlockIV{ 
\BlockII{\omega_{1,\cdot}+&\lambda_{1,\cdot}} 	& 		 																	& 				&	 \\
\BlockII{~~~~~~~~~~&\lambda_{2,\cdot}}					& \BlockII{\omega_{1,\cdot}&~~~~~~~~~~}	&					& 	\\
	~~~~~~ ~\vdots 			 													&			 																	& \ddots	&		\\
\BlockII{~~~~~~~~~~&\lambda_{n,\cdot}}			 		&																				&					&\BlockII{\omega_{1,\cdot}&~~~~~~~~~~}
}
\\
\BlockIV{ 
\BlockII{\omega_{2,\cdot}&~~~~~~~~~~} &\BlockII{~~~~~~~~~~&\lambda_{1,\cdot}}	 			& 				&	 \\
																			&\BlockII{\omega_{2,\cdot}+&\lambda_{2,\cdot}}&					& 	\\
&			 													~~~~~~ ~\vdots 			 							& \ddots	&		\\
																								&\BlockII{~~~~~~~~~~&\lambda_{n,\cdot}}&				&\BlockII{\omega_{2,\cdot}&~~~~~~~~~~}
}
\\
\vdots \\ \vdots
\\
\BlockIV{ 
\BlockII{\omega_{n,\cdot}&~~~~~~~~~~}	 & 		 												& 				&	\BlockII{~~~~~~~~~~&\lambda_{1,\cdot}} \\
					& \BlockII{\omega_{n,\cdot}&~~~~~~~~~~}	&					& 	\BlockII{~~~~~~~~~~&\lambda_{2,\cdot}}\\
			 							&			 												& \ddots	&\vdots		\\
			 		&															&					&\BlockII{\omega_{n,\cdot}+&\lambda_{n,\cdot}}
}
}\right)\\
& \Matrix{\BlockIoffen{e_{1,1} \dots e_{1,n}} &\BlockIoffen{e_{2,1} \dots e_{2,n}}&\dots&\BlockIoffen{e_{1,n} \dots e_{n,n}}}
}
,}%mathe
wobei die leeren Stellen jeweils mit Nullen gefüllt sein sollen.
Alternativ lässt sich die Matrix schreiben als:
\mathe{
\Phi_{2*}= 
\underbrace{
\left(\Matrix{
\BlockI{\Lambda} & & &\\
	& \BlockI{\Lambda}& &\\
	& & \ddots &\\
	&	&	& \BlockI{\Lambda}
}\right)
}_{=: L_{\Lambda} }
 + 
\underbrace{
\left(\Matrix{
\BlockI{\omega_{1,\cdot}\id_{n \times n^2} }\\
\BlockI{\omega_{2,\cdot}\id_{n \times n^2} } \\
\vdots \\
\BlockI{\omega_{n,\cdot}\id_{n \times n^2} }
}\right)}_{=:D_{\Omega}}
}

Da $\Lambda$ mit $\Lambda^t$ invertierbar, ist die Matrix $L_{\Lambda}$ mit $L_{\Lambda^t}$ invertierbar. Das Multiplizieren mit einer invertierbaren Matrix ändert den Rang einer Matrix nicht, daher können wir anstelle von $\Phi_{2*}$ auch folgende Matrix betrachten:
\mathe{\Phi_{2*}L_{\Lambda^t} = \id + D_{\Omega}L_{\Lambda^t}=\id +D_{\Omega\Lambda^t},}

da für das Produkt aus den Matrizen $D_{\Omega}$ und $L_{\Lambda^t}$ gilt:
\mathe{D_{\Omega}L_{\Lambda^t}= D_{\Omega\Lambda^t}.}

\input{Haupt/HochschildAufloesung/Phi2_Fall_1.tex}
\input{Haupt/HochschildAufloesung/Phi2_Fall_2.tex}\newpage
\input{Haupt/HochschildAufloesung/Phi2_Fall_3.tex}\newpage
\input{Haupt/HochschildAufloesung/Phi2_Fall_sonst.tex}

%% file: Haupt/HochschildAufloesung/Phi2_Fall_1.tex
\subsection{1. Fall: $\K=\R$ und $n=2$}
\label{sec:1FallKRUndN2}

Dann sind $\Lambda$ und $\Omega$ Drehmatrizen oder Spiegelungsmatrizen. Also gibt es ein $\psi$ mit:
\mathe{\Omega\Lambda^t = \left( \Matrix{
\cos{\psi} & -\sin{\psi}\\
\sin{\psi} & \cos{\psi}
} \right) \textnormal{ oder }\Omega\Lambda^t = \left( \Matrix{
\cos{\psi} & \sin{\psi}\\
\sin{\psi} & -\cos{\psi}
} \right) 
.}
Die Matrix $\id + D_{\Omega\Lambda^t}$ ist dann:
\mathe{\left( \Matrix{
\cos{\psi}+1 & -\sin{\psi} & 0 &0\\
0 & 1& \cos{\psi} & -\sin{\psi}\\
\sin{\psi} & \cos{\psi}&1 & 0\\
0&0&\sin{\psi} & \cos{\psi}+1
}\right) 
\textnormal{ oder }
\left( \Matrix{
\cos{\psi}+1 & \sin{\psi} & 0 &0\\
0 & 1& \cos{\psi} & \sin{\psi}\\
\sin{\psi} & -\cos{\psi}&1 & 0\\
0&0&\sin{\psi} & -\cos{\psi}+1
}\right) 
.}
Für den Spezialfall $\cos{\psi}=-1$ vereinfacht sich die Matrix zu:
\mathe{\left( \Matrix{
0 & 0 & 0 &0\\
0 & 1& -1 & 0\\
0 & -1&1 & 0\\
0&0&0 & 0
}\right)
\textnormal{ oder }
\left( \Matrix{
0 & 0 & 0 &0\\
0 & 1& -1 & 0\\
0 & 1&1 & 0\\
0&0&0 & 2
}\right)
.}
Der Rang ist also $1$ oder $3$. Dieser Spezialfall (für $\Omega\Lambda^t$ ist eine Drehmatrix) tritt genau dann ein, wenn $\Lambda=-\Omega$.

Falls $\cos{\psi}\neq-1$, kann man die Matrizen mittels des Gaussalgorithmus in folgende Form bringen:
\mathe{
\left( \Matrix{
\sin{\psi} & -\frac{\sin^2{\psi}}{\cos{\psi}+1 } & 0 &0\\
0 & 1& \cos{\psi} & -\sin{\psi}\\
0 & 0&0 & 0\\
0&0&\frac{\sin^2{\psi}}{\cos{\psi}+1} & \sin{\psi}
}\right)
\textnormal{ oder }
\left( \Matrix{
\sin{\psi} & \frac{\sin^2{\psi}}{\cos{\psi}+1} & 0 &0\\
0 & 1& \cos{\psi} & \sin{\psi}\\
0 & 0&sin{\psi} & \frac{\sin^2{\psi}}{\cos{\psi}+1}\\
0&0&0 & 0
}\right) 
}

Der Rang ist also $3$.
\newpage
\subsubsection{Fall $\Omega\Lambda^t$ ist eine Drehmatrix:}
\label{sec:FallDrechmatrix}

\mathe{\left( \Matrix{
\sin{\psi} & -\frac{\sin^2{\psi}}{\cos{\psi}+1 } & 0 &0\\
0 & 1& \cos{\psi} & -\sin{\psi}\\
\sin{\psi} & \cos{\psi}&1 & 0\\
0&0&\sin{\psi} & \cos{\psi}+1
}\right)
\rightarrow 
\left( \Matrix{
\sin{\psi} & -\frac{\sin^2{\psi}}{\cos{\psi}+1 } & 0 &0\\
0 & 1& \cos{\psi} & -\sin{\psi}\\
0 & \cos{\psi}+\frac{\sin^2{\psi}}{\cos{\psi}+1 }&1 & 0\\
0&0&\sin{\psi} & \cos{\psi}+1
}\right)
}

mit $\cos{\psi}+\frac{\sin^2{\psi}}{\cos{\psi}+1 }=\frac{\cos^2{\psi}+\cos{\psi}+\sin^2{\psi}}{\cos{\psi}+1 }=1$

\mathe{
\left( \Matrix{
\sin{\psi} & -\frac{\sin^2{\psi}}{\cos{\psi}+1 } & 0 &0\\
0 & 1& \cos{\psi} & -\sin{\psi}\\
0 & 1&1 & 0\\
0&0&\frac{\sin{\psi}}{\cos{\psi}+1} & 1
}\right)
\rightarrow
\left( \Matrix{
\sin{\psi} & -\frac{\sin^2{\psi}}{\cos{\psi}+1 } & 0 &0\\
0 & 1& \cos{\psi} & -\sin{\psi}\\
0 & 0&1-\cos{\psi} & \sin{\psi}\\
0&0&\frac{\sin^2{\psi}}{\cos{\psi}+1} & \sin{\psi}
}\right)
}
mit $\frac{\sin^2{\psi}}{\cos{\psi}+1}-(1-\cos{\psi})=\frac{\sin^2{\psi}-(1-\cos^2{\psi})}{\cos{\psi}+1}=0$
\mathe{\rightarrow
\left( \Matrix{
\sin{\psi} & -\frac{\sin^2{\psi}}{\cos{\psi}+1 } & 0 &0\\
0 & 1& \cos{\psi} & -\sin{\psi}\\
0 & 0&0 & 0\\
0&0&\frac{\sin^2{\psi}}{\cos{\psi}+1} & \sin{\psi}
}\right)
}
\subsubsection{Fall $\Omega\Lambda^t$ ist eine Spiegelungsmatrix:}
\label{sec:FallSpiegelungsmatrix}

\mathe{
\left( \Matrix{
\cos{\psi}+1 & \sin{\psi} & 0 &0\\
0 & 1& \cos{\psi} & \sin{\psi}\\
\sin{\psi} & -\cos{\psi}&1 & 0\\
0&0&\sin{\psi} & -\cos{\psi}+1
}\right) 
\rightarrow
\left( \Matrix{
\sin{\psi} & \frac{\sin^2{\psi}}{\cos{\psi}+1} & 0 &0\\
0 & 1& \cos{\psi} & \sin{\psi}\\
\sin{\psi} & -\cos{\psi}&1 & 0\\
0&0&\sin{\psi} & -\cos{\psi}+1
}\right) 
}
\mathe{
\rightarrow
\left( \Matrix{
\sin{\psi} & \frac{\sin^2{\psi}}{\cos{\psi}+1} & 0 &0\\
0 & 1& \cos{\psi} & \sin{\psi}\\
0 & -\cos{\psi}-\frac{\sin^2{\psi}}{\cos{\psi}+1}&1 & 0\\
0&0&\sin{\psi} & -\cos{\psi}+1
}\right) }
mit $-\cos{\psi}-\frac{\sin^2{\psi}}{\cos{\psi}+1 }=-\frac{\cos^2{\psi}+\cos{\psi}+\sin^2{\psi}}{\cos{\psi}+1 }=-1$

\mathe{
\rightarrow
\left( \Matrix{
\sin{\psi} & \frac{\sin^2{\psi}}{\cos{\psi}+1} & 0 &0\\
0 & 1& \cos{\psi} & \sin{\psi}\\
0 & -1&1 & 0\\
0&0&\sin{\psi} & -\cos{\psi}+1
}\right) 
\rightarrow
\left( \Matrix{
\sin{\psi} & \frac{\sin^2{\psi}}{\cos{\psi}+1} & 0 &0\\
0 & 1& \cos{\psi} & \sin{\psi}\\
0 &0 &1+ \cos{\psi} & \sin{\psi}\\
0&0&\sin{\psi} & -\cos{\psi}+1
}\right) 
}
mit $\frac{\sin^2{\psi}}{\cos{\psi}+1}-(1-\cos{\psi})=\frac{\sin^2{\psi}-(1-\cos^2{\psi})}{\cos{\psi}+1}=0$
\mathe{
\rightarrow
\left( \Matrix{
\sin{\psi} & \frac{\sin^2{\psi}}{\cos{\psi}+1} & 0 &0\\
0 & 1& \cos{\psi} & \sin{\psi}\\
0 & 0&sin{\psi} & \frac{\sin^2{\psi}}{\cos{\psi}+1}\\
0&0&0 & 0
}\right) }

%% file: Haupt/HochschildAufloesung/Phi2_Fall_2.tex
\subsection{2. Fall: $\Omega=\Lambda$, also $\Omega\Lambda^t=\id$}
\label{sec:2FallOmegaLambdaAlsoOmegaLambdaTId}

Die Matrix $\id + D_{\id}$ ist in diesem Fall:

\mathe{%\Matrix{ 
%\Matrix{
%f_{1,\cdot} & \Matrix{f_{1,1}\\ f_{1,2} \\ \vdots \\ f_{1,n}}
%\\ \hline
%f_{2,\cdot} & \Matrix{f_{2,1}\\ f_{2,2} \\ \vdots \\ f_{2,n}}
%\\ \hline
%						& \vdots \\ & \vdots	
%\\\hline
%f_{n,\cdot} & \Matrix{f_{n,1} \\ \vdots \\ f_{n,n-1}\\ f_{n,n}}
%}
%&
\left(\Matrix{
\BlockIV{ 
\LMatrix{
\BlockI{2, 0, \dots, 0} & & &\\
~~~~~~1												& \BlockI{1,0,\dots, 0} & &\\
~~~~~~~~~~~\ddots							&												&\ddots&\\
~~~~~~~~~~~~~~~~~1				&												&			 & \BlockI{1,0,\dots, 0}
}}
\\
\BlockIV{ 
\LMatrix{
\BlockI{0, 1, \dots, 0} & ~~1& &\\
												& \BlockI{0,2,\dots, 0} & &\\
						&		~~~~~~~~~~~\ddots											&\ddots&\\
				&				~~~~~~~~~~~~~~~~~1								&			 & \BlockI{0,1,\dots, 0}
}}\\
\vdots \\ \vdots
\\
\BlockIV{ 
\LMatrix{
\BlockI{0,  \dots,0, 1} & & &~~1\\
												& \BlockI{0,\dots,0, 1} & &~~~~~~1	\\
						&												&\ddots&~~~~~~~~~~~\ddots	\\
				&												&			 & \BlockI{0,\dots, 0,2}
}}
}\right)
%& \Matrix{\BlockIoffen{e_{1,1} \dots e_{1,n}} &\BlockIoffen{e_{2,1} \dots e_{2,n}}&\dots&\BlockIoffen{e_{1,n} \dots e_{n,n}}}
%}
.}%mathe

Die Zeilen zu den Basisvektoren $f_{p,q}$ und $f_{q,p}$ stimmen für jedes Paar $(p,q)$ überein. Alle anderen Zeilen sind linear unabhängig. Der Rang ist daher $\Sumn{i=1}i=\frac{n^2+n}{2}$.

%% file: Haupt/HochschildAufloesung/Phi2_Fall_3.tex
\subsection{3. Fall: $\Omega=-\Lambda$, also $\Omega\Lambda^t=-\id$}
\label{sec:3FallOmegaLambdaAlsoOmegaLambdaTId}

Die Matrix $\id + D_{-\id}$ ist in diesem Fall:

\mathe{%\Matrix{ 
%\Matrix{
%f_{1,\cdot} & \Matrix{f_{1,1}\\ f_{1,2} \\ \vdots \\ f_{1,n}}
%\\ \hline
%f_{2,\cdot} & \Matrix{f_{2,1}\\ f_{2,2} \\ \vdots \\ f_{2,n}}
%\\ \hline
%						& \vdots \\ & \vdots	
%\\\hline
%f_{n,\cdot} & \Matrix{f_{n,1} \\ \vdots \\ f_{n,n-1}\\ f_{n,n}}
%}
%&
\left(\Matrix{
\BlockIV{ 
\LMatrix{
\BlockI{0, 0, \dots, 0} & & &\\
~~~~~~1												& \BlockI{-1,0,\dots, 0} & &\\
~~~~~~~~~~~\ddots							&												&\ddots&\\
~~~~~~~~~~~~~~~~~1				&												&			 & \BlockI{-1,0,\dots, 0}
}}
\\
\BlockIV{ 
\LMatrix{
\BlockI{0, -1, \dots, 0} & ~~1& &\\
												& \BlockI{0,0,\dots, 0} & &\\
						&		~~~~~~~~~~~\ddots											&\ddots&\\
				&				~~~~~~~~~~~~~~~~~1								&			 & \BlockI{0,-1,\dots, 0}
}}\\
\vdots \\ \vdots
\\
\BlockIV{ 
\LMatrix{
\BlockI{0,  \dots,0, -1} & & &~~1\\
												& \BlockI{0,\dots,0, -1} & &~~~~~~1	\\
						&												&\ddots&~~~~~~~~~~~\ddots	\\
				&												&			 & \BlockI{0,\dots, 0,0}
}}
}\right)
%& \Matrix{\BlockIoffen{e_{1,1} \dots e_{1,n}} &\BlockIoffen{e_{2,1} \dots e_{2,n}}&\dots&\BlockIoffen{e_{1,n} \dots e_{n,n}}}
%}
.}%mathe

Die Zeilen zu den Basisvektoren $f_{i,i}$ bestehen für jedes $i$ nur aus Nullen.
Die Zeilen zu den Basisvektoren $f_{p,q}$ und $f_{q,p}$ stimmen für jedes Paar $(p,q)$ bis auf das Vorzeichen überein. Alle anderen Zeilen sind linear unabhängig. Der Rang ist daher $\Sumn{i=1}i-n=\frac{n^2-n}{2}$.

%% file: Haupt/HochschildAufloesung/Phi2_Fall_sonst.tex
\newpage
\subsection{Allgemeiner Fall}
\label{sec:AllgemeinerFall}

\newcommand{\NS}{\textnormal{NS}}
\newcommand{\e}[1]{e_{\textnormal{\tiny{$#1$} }}}

In diesem Abschnitt wollen wir den Rang von $\id+D_{\Omega\Lambda^t}$ für eine beliebige orthogonale Matrix $\Omega\Lambda^t$ berechnen.

Zuerst stellen wir fest, dass wir den Rang von $\id+D_{\Omega \Lambda^t}$ über die Vielfachheit der Nullstellen des Charakteristischen Polynoms an der Stelle $-1$ berechnen können:
\mathe{ \Rang{\id+D_{\Omega\Lambda^t} }= n^2-\NS_{-1}(\Det{x\id-D_{\Omega\Lambda^t}})
.} 

\Frage{Das geht nur wenn bekannt ist, dass $D_{\Omega\Lambda^t}$ diagonalisierbar ist. 
In $\R$ ist das klar, weil: Dann ist $\Omega\Lambda^t$ diagonalisierbar, also gibt es ein diagonales $N$. Aus Beweis von Satz(Charakteristische Polynom) folgt das für Diagonales $N$ auch $D_{N}$ diagonalisierbar. Aus Lemma(D ist verträglich mit konjugation) folgt das dann auch $D_{M}$ diagonalisierbar. }

Wir werden nun das Charakteristische Polynom von $D_{\Psi}$ für eine beliebige orthogonale Matrix $\Psi$ berechnen.

Dazu wollen wir uns überlegen, warum wir annehmen können, dass die Matrix $\Psi$ in Jordanform vorliegt.

Falls $\K$ nicht algebraisch abgeschlossen ist,  verändert sich das Charakteristische Polynom von $D_{\Psi}$ über einem Abschluss von $\K$ nicht. Wir können also ohne Beschränkung der Allgemeinheit annehmen, dass $\K$ algebraisch abgeschlossen ist.

Für eine Matrix $\Psi$ gibt es in einem algebraisch abgeschlossenen Körper $\K$ eine Jordanform. Es gibt also eine Matrix $N$ und eine Jordan-Matrix $J_{\Psi}$, so dass gilt: $\Psi N = N J_{\Psi}.$

\Lemma{Konjugieren ist verträglich mit $D$}{
\label{LemmaKonjugierenD}
Seien $M,N$ und $X$ quadratische Matrizen, die $M X =  X N$ erfüllen. Dann gilt:
\mathe{
D_{M}(X\otimes X) = (X\otimes X)D_{N}.
}

\Beweis{Lemma \ref{LemmaKonjugierenD} }{

Um Schreibarbeit zu sparen bezeichnen wir die $p$-te Zeile der Matrix $M$ bzw. $N$ mit $m_{p,\cdot}$ bzw. $n_{p,\cdot}$ und die $p$-te Spalte entsprechend $m_{\cdot,p}$ bzw. $n_{\cdot,p}$.
Die Voraussetzung $M X = X N$ lässt sich dann schreiben als:
\mathe{ \Sumn{i=1}m_{p,i}x_{i,q}= \Sumn{i=1}x_{p,i}n_{i,q}~~~\forall (p,q)}
oder alternativ:
\mathe{\left< m_{p,\cdot}~,~x_{\cdot,q}\right>=\left< x_{p,\cdot}~,~n_{\cdot,q}\right>~~~\forall (p,q).}

Die Matrix $(X\otimes X)$ ist definiert durch:
\mathe{(X\otimes X):= \left(
\Matrix{
x_{1,1}\BlockI{X}& \dots &x_{1,n}\BlockI{X}\\
\vdots&& \vdots \\
x_{n,1}\BlockI{X}& \dots &x_{n,n}\BlockI{X}\\
}\right)
}
Für eine $n^2\times n^2$-Matrix $Y$ bezeichnen wir die Spalten und  Zeilen wie zuvor nach den zugehörigen Basisvektoren. An der Stelle $((p,q),(i,j))$, also im $p$-ten Zeilenblock im $i$-ten Spaltenblock an der $(q,j)$-Stelle, steht der Anteil des $(i,j)$-ten Basisvektor des Bildes zum Basisvektor $e_{p,q}$.
\mathe{\Matrix{
		&\Matrix{
						&1 & \dots &i\textnormal{-ter Spaltenblock} &\dots &n\\
						& & &\BlockIV{~~~~~~~~\dots &j& \dots &~~~~~~~~~~~~~} &\\
						}\\\Matrix{1 & \\ \vdots &  \\p \textnormal{-ter Block} & \BlockI{ \vdots \\ q \\ \vdots \\ } \\ \vdots \\n
			 } &
			 	\left(\Matrix{
			 			&  && \\
			 			&	\vdots & & \\
			 				\dots & \BlockIV{
			 												&\vdots&&\\
			 												 \dots &~~~~~~y_{(p,q),(i,j)} &	\dots&\\
			 												 &\vdots&&\\
			 												 &&&~
			 												} & \dots\\
			 										~~~~~~~~~~~~~~~~~~~~~~~&\vdots & & ~~~~~~~~~~~~\\~\\		 	
			 			 } \right)		
			 	 }
}

 Da in jeder Zeile von $D_{M}$ nur ein $n$-Tupel steht, lässt sich das Produkt $D_{M}(X\otimes X)$ leicht berechnen:
\mathe{ (D_{M}(X\otimes X))_{(p,q),(i,j)}:= \left<m_{p,\cdot}~,~x_{i,q} x_{\cdot,j}\right>~~~\forall ((p,q),(i,j)).}

In jeder Zeile von $D_{N}$ stehen wiederum nur wenige Einträge, nämlich ein Spaltenvektor auf verschiedene Blöcke verteilt. Daher gilt:
\mathe{ ((X\otimes X)D_{N})_{(p,q),(i,j)}:= \left<x_{p,\cdot}x_{i,q}~,~m_{\cdot,j}  \right>~~~\forall ((p,q),(i,j)).}
Mit der Voraussetzung $\left< m_{p,\cdot}~,~x_{\cdot,q}\right>=\left< x_{p,\cdot}~,~n_{\cdot,q}\right>~~~\forall (p,q)$ erhalten wir für jedes $(p,q),(i,j)$:
\Gleichung{
((X\otimes X)D_{N})_{(p,q),(i,j)}:=& \left<x_{p,\cdot}x_{i,q}~,~n_{\cdot,j}  \right>\\
=&x_{i,q}\left<x_{p,\cdot}~,~n_{\cdot,j}  \right>\\
=&x_{i,q}\left<m_{p,\cdot}~,~x_{\cdot,j}  \right>\\
=&\left<m_{p,\cdot}~,~x_{i,q} x_{\cdot,j}\right> \\
=:&(D_{M}(X\otimes X))_{(p,q),(i,j)}.
}

}%Beweis

}%Lemma

Als Folgerung daraus erhalten wir, dass das Charakteristische Polynom  von $D_{\Psi}$ mit dem Charakteristischen Polynom von $D_{J_{\Psi}}$ übereinstimmt.

\Satz{Charakteristisches Polynom von $D_{J}$}{
\label{SatzCharakteristischePolynom}
Sei $J$ eine Jordan-Matrix, dann gilt für das Charakteristische Polynom:
\mathe{
\chi_{D_{J}}(x)= \chi_{J}(x) \cdot  \chi_{\bigwedge^2 J}(x^2).
}
\Beweis{Satz \ref{SatzCharakteristischePolynom}}{

Sei 
\mathe{J:=\left( \Matrix{
r_1&\epsilon_1&&&\\
&\ddots&\ddots&\\
&&r_{n-1}&\epsilon_{n-1}\\
&&&r_{n}
}\right)
,}
wobei die leeren Stellen mit Nullen gefüllt sein sollen.
Dann ist die Matrix $D_{J}$  ein Endomorphismus auf einem $n^2$-dimensionalen Vektorraum. Sie ist von folgender Form:
\mathe{\Matrix{ 
\Matrix{
e_{1,\cdot} & \Matrix{e_{1,1}\\ e_{1,2} \\ \vdots \\ e_{1,n}}
\\ \hline
e_{2,\cdot} & \Matrix{e_{2,1}\\ e_{2,2} \\ \vdots \\ e_{2,n}}
\\ \hline
						& \vdots \\ & \vdots	
\\\hline
e_{n,\cdot} & \Matrix{e_{n,1} \\ \vdots \\ e_{n,n-1}\\ e_{n,n}}
}
&
\left(\Matrix{
\BlockIV{ 
\BlockII{r_1, \epsilon_1, \dots &~~~~~~~} 	& 		 																	& 				&	 \\
															& \BlockII{r_1, \epsilon_1, \dots &~~~~~~~}	&					& 	\\
				 											&			 												& \ddots	&		\\
			 		&																				&					&\BlockII{r_1, \epsilon_1, \dots &~~~~~~~}
}
\\
\BlockIV{ 
\BlockII{0, r_2, \epsilon_2,\dots &~~~} 	& 		 																	& 				&	 \\
															& \BlockII{0, r_2, \epsilon_2,\dots &~~~}&					& 	\\
				 											&			 												& \ddots	&		\\
			 		&																				&					&\BlockII{0, r_2, \epsilon_2,\dots &~~~}
}
\\
\vdots \\ \vdots
\\
\BlockIV{ 
\BlockII{0,~\dots~ ~~0,&r_n} 	& 		 																	& 				&	 \\
															& \BlockII{0,~\dots~ ~~0,&r_n}	&					& 	\\
				 											&			 												& \ddots	&		\\
			 		&																				&					&\BlockII{0,~\dots~ ~~0,&r_n}
}
}\right)\\
& \Matrix{\BlockIoffen{e_{1,1}, e_{1,2}, \dots~~~~} &\BlockIoffen{e_{2,1}, e_{2,2},  \dots~~~~}&\dots&\BlockIoffen{e_{n,1}, e_{n,2},  \dots ~~~~}}
}
.}%mathe

In den Beispielen \ref{sec:2FallOmegaLambdaAlsoOmegaLambdaTId} und \ref{sec:3FallOmegaLambdaAlsoOmegaLambdaTId} haben wir gesehen, dass die Zeilen zu den Basisvektoren $e_{p,p}$, $e_{p,q}$ und $e_{q,p}$ eine besondere Rolle spielen. Daher sortieren wir die Basisvektoren $e_{p,q}$ nach der Summe aus $p$ und $q$ um. Falls die Summe gleich ist, sortieren wir nach dem zweiten Index. Auf diese Weise folgen die Paare $e_{p,q}$ und $e_{q,p}$ aufeinander, gefolgt von $e_{p,p}$.

\mathe{\Matrix{ 
\Matrix{
\Matrix{e_{1,1}}
\\ \hline
\Matrix{e_{2,1}\\ e_{1,2}}
\\ \hline
 \Matrix{e_{3,1}\\ e_{1,3}}
\\ \hline
\Matrix{e_{2,2}}
\\ \hline
 \vdots
\\ \hline
 \Matrix{e_{n,n}}
}
&
\left(\Matrix{ 
\BlockI{~r_1~} &\Matrix{~&~~~~\epsilon_1}	
\\ \hline
&~\BlockII{&r_{2}~~\\~~r_1&}~&&\Matrix{\\ \epsilon_1} & \Matrix{\dots \\~}\\
\hline
&&~\BlockII{&r_3~~\\~~r_1&}~&&\Matrix{\dots \\ \dots}\\
\hline
&&&~\BlockI{~r_{2}~}&\dots\\
\hline
&&&&\ddots&\dots\\
&&&&&\BlockI{~r_n~}
}\right)\\
& \Matrix{
|\e{1,1}|&|\e{2,1}&\e{1,2}|&|\e{3,1}& \e{1,3}|&|\e{2,2}|&\dots &|\e{n,n}|
%\BlockIoffen{\e{1,1}} &\BlockIIoffen{\e{2,1}&\e{1,2}}&\BlockIIoffen{\e{3,1}& \e{1,3}} &\BlockIoffen{\e{2,2}} &\BlockIoffen{\e{3,2},\e{2,3}}&\BlockIoffen{\e{3,3}}
}
}
.}%mathe
In der Zeile zum Basisvektor $\e{p,q}$ kommen nur $r_p$ und $\epsilon_p$ vor, wobei $\epsilon_p$ weiter links steht.
Die Matrix $D_{J}$ sieht fast wie eine obere Dreiecksmatrix aus. Auf der Diagonalen stehen $1 \times 1$-und $2 \times 2$-Blöcke und in der oberen rechten Ecke steht an einigen Stellen ein $\epsilon$.
\\
Das Charakteristische Polynom berechnen wir über $\Det{x \id - D_{J}}$, also:
\mathe{
\Det{
\Matrix{ 
\BlockI{x-r_1~} &\Matrix{~&~~~~~~~-\epsilon_1}	
\\ \hline
&~\BlockII{x&-r_{2}~\\~-r_1&x}~&&\Matrix{\\ -\epsilon_1} & \Matrix{\dots \\~}\\
\hline
&&~\BlockII{x&-r_3~\\~-r_1&x}~&&\Matrix{\dots \\ \dots}\\
\hline
&&&~\BlockI{x-r_{2}~}&\dots\\
\hline
&&&&\ddots&\dots\\
&&&&&\BlockI{x-r_n~}
}%Matrix
}
.}
Für die Determinante von Matrizen, die sich in Blockmatrizen zerlegen lassen, wobei der untere Block leer ist, ist die Determinante gleich dem Produkt der Determinanten auf der Diagonalen:
\mathe{\Det{x \id - D_{J}}= (x-r_1)\Det{\Matrix{x&-r_{2}\\-r_1&x}}\Det{\Matrix{x&-r_{3}\\-r_1&x}}(x-r_2)\cdots (x-r_n).}
Da $ \Det{\Matrix{x&-r_{p}\\-r_q&x}}=x^2-r_pr_q$,  können wir das Charakteristische Polynom von $D_{J}$ schreiben als:
\mathe{\chi_{D_{J}}(x)=\left(\prod_{i=1}(x-r_i)\right)\left(\prod_{i>j}(x^2-r_i r_j) \right).}
Dies entspricht gerade der Behauptung: $\chi_{D_{J}}(x)= \chi_{J}(x) \cdot  \chi_{\bigwedge^2 J}(x^2)$.
}%Beweis
}%Satz

Zusammen mit den Vorüberlegungen erhalten wir, dass $\Rang{\Phi_{2*}}$ nur von Eigenwerten $r_i$ der Matrix $\Omega\Lambda^t$ abhängt. Sei dazu 
\Gleichung{
k_{-1}&:=\# \Menge{i}{r_i=-1}\\
k_{\bigwedge}&:=\# \Menge{(i,j)}{i>j, r_ir_j=1}
.}

Der Rang ist dann:
\mathe{\Rang{\Phi_{2*}}=n^2-(k_{-1}+k_{\bigwedge} )}
und der Kern:
\mathe{\Kern{\Phi_{2*}}=(k_{-1}+k_{\bigwedge} ).}

%% file: Haupt/HochschildAufloesung/Phi3.tex
\subsection{$\Phi_{3*}: \bs{f} \mapsto -\tr(\Lambda \Omega^t\bs{F})+\tr(\bs{F})$}
\label{sec:Phi3Stern}

\mathe{\begin{array}{|c|c|}
\hline  
\Dim{\Kern{\Phi_{3*}}} & \Rang{3}\Matrix{~\\~}\\
\hline 
\begin{array}{c|c} 
\Omega=\Lambda & \textnormal{sonst} \\
\hline 
1&0
\end{array}
& 
\begin{array}{c|c} 
\Omega=\Lambda & \textnormal{sonst}
\\
\hline 
0&1
\end{array}\Matrix{~\\~\\~} \\
 \hline
\end{array}
}%mathe

Sei zunächst $\Omega=\Lambda$, dann ist $\Lambda^t\Omega=\id$ und der $\Kern{\Phi_{3*}}$ das gesamte Urbild. Da der Urbildraum eindimensional ist, gilt $\Rang{3}=0$ und $\Dim{\Kern{\Phi_{3*}}}=1$.

Falls $\Omega\neq\Lambda$, dann existiert wenigstens ein $(p,q)$ mit $\Sumn{i=1}\lambda_{q,i}\omega_{p,i}+\delta_{p,q}\neq0$; also gilt $\Rang{3}>0$ . Da der Urbildraum nur eindimensional ist, gilt $\Rang{3}=1$ und mit $\Dim{\Kern{\Phi_{1*}}}= \Dim{\textnormal{Urbild}}-\Dim{\textnormal{Bild}}=1-1=0$.

%% file: Haupt/Ext.tex
\section{Ext}
\label{sec:Ext}
Die Extgruppen sind isomorph zu den Torgruppen. Hier erhält man die selben Matrizen, nur transponiert.

Wir betrachten, ähnlich wie in Kapitel \ref{sec:Tor}, die induzierte, an der ersten Stelle nicht exakte, Sequenz:
\mathe{
0 \leftarrow \Algebra  \stackrel{\phi_1}{\leftarrow} \Algebra^{n^2} \stackrel{\phi_2}{\leftarrow} \Algebra^{n^2}  \stackrel{\phi_3}{\leftarrow} \Algebra  \leftarrow 0
.}%mathe
Wir untersuchen zunächst den Vektorraum Isomorphismus $\Hom_{\Algebra}\left(\Algebra^n, \K_{\Omega}\right) \iso \K^n$.
 Sei $\bs{e}_1, \dots, \bs{e}_n$ eine Basis des Moduls $\Algebra^n$, dann ist  
 $\bs{\bar{e}}_1, \dots, \bs{\bar{e}}_1$ eine Basis von $\Hom_{\Algebra}\left(\Algebra^n, \K_{\Omega}\right)$ mit $\bs{\bar{e}}_j (\bs{e}_i)= \delta_{i,j}$. Lineares Fortsetzen liefert für $a_i \in \Algebra$:
\mathe{
\bs{\bar{e}}_j\left(\Sumn{i=1} a_{i} \bs{e}_i \right)= a_j
.}
So erhalten wir aus der induzierten Sequenz unter dem kontravarianten Funktor $\Hom_{\Algebra}\left(\cdot, \K_{\Omega}\right) $:
\mathe{ 0 \stackrel{}{\rightarrow}  \K_{\Omega}  
					\stackrel{\phi_1^*}{\rightarrow}  \K^{n^2}  
					\stackrel{\phi_2^*}{\rightarrow}  \K^{n^2}  
					\stackrel{\phi_3^*}{\rightarrow}  \K
					\stackrel{}{\rightarrow}  0.
}

Sei $\bar{\bs{e}}$ der Basisvektor von $\K_{\Omega} $, dann sind die Abbildungen folgendermaßen definiert: 
\AufzaehlungP{
\item{$\phi_1^*: \bar{\bs{e}} \mapsto \left(\Omega-\Lambda\right)\bar{\bs{E}}$}, wobei $\bar{\bs{E}}:=\left(\bar{\bs{e}}_{i,j}\right)$ die Matrix der Basisvektoren des Bildraumes ist. Die Abbildung erhalten wir durch:
\mathe{
	\bs{\bar{e}} \mapsto \left(
	\bs{e}_{p,q} 
	\stackrel{\Phi_1}{\mapsto} \left( a_{p,q}-\lambda_{p,q} \right) \bs{e} 		
	\stackrel{\bs{\bar{e}}}{\mapsto}\omega_{p,q}-\lambda_{p,q} \right),
	}
	Diese Abbildung stimmt mit $\phi_{1*}$ überein.

	\item{$\phi_2^*: \bar{\bs{E}} \mapsto \left(\Omega^t\bar{\bs{F}}\right)^t +  \left(\Lambda^t\bar{\bs{F}^t}\right)^t$} , wobei $\bar{\bs{F}}:=\left(\bar{\bs{f}}_{i,j}\right)$ die Matrix der Basisvektoren des Bildraumes ist. Die Abbildung erhalten wir durch:
	\mathe{\bs{\bar{e}}_{p,q} \mapsto \left(
	\bs{f}_{x,y} 
	\stackrel{\Phi_2}{\mapsto} \Sumn{i=1}\left(a_{x,i}\bs{e}_{y,i} + \lambda_{y,i}\bs{e}_{x,i}\right)
  \stackrel{\bs{\bar{e}}_{p,q}}{\mapsto} \underbrace{\Sumn{i=1} \left( \omega_{x,i}\delta_{p,y}\delta_{q,i} + \lambda_{y,i}\delta_{p,x}\delta_{q,i}\right)}_{=\omega_{x,q}\delta_{p,y} + \lambda_{y,q}\delta_{p,x}} \right),
  }
  also:
  \mathe{
   \bs{\bar{e}}_{p,q} \mapsto \Sumn{x,y=1} \left(
   				\omega_{x,q}\delta_{p,y} + \lambda_{y,q}\delta_{p,x} \right) \bar{\bs{f}}_{x,y} 
   				= \Sumn{x=1} \omega_{x,q}\bar{\bs{f}}_{x,p} +\Sumn{y=1} \lambda_{y,q}\bar{\bs{f}}_{p,y}
  .}
  \item{$\phi_3^*: \bar{\bs{F}} \mapsto \left(\left(\Lambda\Omega^t\right)^t+\id\right)\bar{\bs{f}} $}, wobei $\bar{\bs{f}}$ der Basisvektor des Bildraumes ist. Die Abbildung erhalten wir durch:
  \mathe{\bs{\bar{f}}_{p,q} \mapsto \left(
  \bs{f}
  \stackrel{\Phi_3}{\mapsto} \Sumn{i,j,k=1} \lambda_{k,i} a_{j,i}\bs{f}_{j,k} + \Sumn{i=1}\bs{f}_{i,i} 
  \stackrel{\bar{\bs{f}}_{p,q}}{\mapsto}
  \underbrace{\Sumn{i,j,k=1} \lambda_{k,i} \omega_{j,i}\delta_{p,j}\delta_{q,k} + \Sumn{i=1}\delta_{p,i}\delta_{q,i} }_{=\Sumn{i=1} \lambda_{q,i} \omega_{p,i}  + \delta_{p,q}}
  \right)
  }
   also:
  \mathe{
   \bs{\bar{f}}_{p,q} \mapsto \left( \Sumn{i=1} \lambda_{q,i} \omega_{p,i} + \delta_{p,q}\right)\bar{\bs{f}}
   .}
 
}%Aufzaehlung

  	Betrachten wir als Nächstes den Rang dieser Abbildungen.
\subsection{$\phi_1^*: \bar{\bs{e}} \mapsto \left(\Omega-\Lambda\right)\bar{\bs{E}}$}
\label{sec:Phi1Stern:Ext}
Die Abbildung $\phi_1^*$ hat höchstens Rang $1$, da der Urbildraum eindimensional ist. Sie hat Rang $0$, wenn das Bild nur die $0$ ist. Dies ist der Fall wenn:
\mathe{\Omega-\Lambda = 0,}
also gilt: 
\mathe{\RanG{1}=\left\{\LMatrix{ 0 \textnormal{ für } \Omega=\Lambda \\ 1 \textnormal{ sonst}} \right.
.}

\subsection{$\phi_2^*: \bar{\bs{E}} \mapsto \left(\Omega^t\bar{\bs{F}}\right)^t +  \left(\Lambda^t\bar{\bs{F}^t}\right)^t$}
\label{sec:Phi2Stern:Ext}
Hier beobachten wir, dass $\phi_2^*=\phi_{2*}^t$:
\mathe{\phi_2^*: \bar{\bs{E}} \mapsto \left(\Omega^t\bar{\bs{F}}\right)^t +  \left(\Lambda^t\bar{\bs{F}^t}\right)^t = \left(\Omega^t\bar{\bs{F}}\right)^t +  \left(\Lambda\bar{\bs{F}^t}\right)^t = \phi_{2*}^t
}
Wir betrachten die Abbildung in Matrizenschreibweise, wozu wir das Bild des Basisvektors $\bar{\bs{e}}_{p,q}$ in der Basis $\bar{\bs{f}}_{x,y}$ untersuchen:
\Gleichung{
	\bar{\bs{e}}_{(p,q)}{(x,y)} &= \delta_{p,y}\omega_{x,q} + \delta_{p,x}\lambda_{y,q}\\
		& = \Sumn{p=1} \Sumn{q=1} \left(\delta_{p,y}\omega_{x,q} +
					 \delta_{p,x}\lambda_{y,q}\right)\bs{e}_{p,q}\\
		& = \Sumn{q=1}\omega_{x,y}\bs{e}_{y,q} + \Sumn{p=1}\lambda_{y,q}\bs{e}_{x,q}\\
		&= \phi_{2*}.
		}

\Frage{Koeffizienten Rechnung}

Also gilt :
\mathe{ \Rang{2}=\RanG{2}
}
\Frage{Achtung hier ist noch ein Fehler $\phi_2^*=\phi_{2*}^t$}

\subsection{$\phi_3^*: \bar{\bs{F}} \mapsto \left(\left(\Lambda\Omega^t\right)^t+\id\right)\bar{\bs{f}} $}
\label{sec:Phi3Stern:Ext}
Die Abbildung $\phi_3^*$ hat höchstens Rang $1$, da der Bildraum eindimensional ist. Sie hat Rang $0$, wenn das Bild nur die $0$ ist. Dies ist nur dann der Fall, wenn:
\mathe{\left(\Lambda\Omega^t\right)^t=-\id,}
also gilt:
\mathe{\RanG{1}=\left\{\LMatrix{ 0 \textnormal{ für } \Omega=\Lambda \\ 1 \textnormal{ sonst}} \right.
.}

%% file: Bezeichner.tex
\newcommand{\Eintrag}[4]{ \textnormal{#1} & \Matrix{ #2 }     &  {\backslash\textnormal{#3}} & \ref{#4} & \textnormal{#4} \\ \hline }

\mathe{\begin{array}{|l|c|l|c|l|}
\hline
			\textnormal{Name}			& \textnormal{Zeichen}	& \textnormal{command} 	& \textnormal{Seite}	& \backslash \{ \textnormal{label} \}
\\\hline
%
%
%Reduktionssysteme
\Eintrag{Reduktionssystem}{\Red }		 						{Red}			{DefReduktionssystem}
\T{induziertes Reduktionssystem}	&{\Red}	&\backslash \T{Red}& \ref{DefInduziertesRedAlgebra}& \T{DefInduziertesRedAlgebra}\\
																	&{\RRed}&\backslash \T{RRed}& \ref{DefInduziertesRedRig}	&\T{DefInduziertesRedRig}\\
																	\hline
\Eintrag{$\red$-Reduktionsregeln auf $\Rig$}{\Rel}{Rel}{DefInduziertesRedRig}
\Eintrag{$\FreieA$-Reduktionsregeln auf $\Rig$}{\RelF}{RelF}{DefInduziertesRedRig}
%
%
%Pfeile
\T{Reduktionsregel}	& (x,y) 				& & \ref{DefReduktionsregel} & \T{DefReduktionsregel}\\
										& x \Regel y		&\backslash\T{Regel} &&\\
										& x\redRegel y  &\backslash\T{redRegel}&&\\ 
										& x\RedRegel y  &\backslash\T{RedRegel}&&\\ 
										& x\RRedRegel y  &\backslash\T{RRedRegel}&&\\ 
										& x\RedpRegel y  &\backslash\T{RedpRegel}&&\\ 
						  \hline
\T{Reduktionsweg}		& (x,y)_{i=1\dots n} 				& & \ref{DefReduktionsweg} & \T{DefReduktionsweg}\\
										& x \Weg y		&\backslash\T{Weg} &&\\ 	
										& x\RedWeg y  &\backslash\T{RedWeg}&&\\ 									
 										&	x\RRedWeg y& \backslash\T{RRedWeg} &&\\	
 										&	x\AqRedWeg y& \backslash\T{AqRedWeg} &&\\						
 		          \hline
\T{unreduzierbares Elment} &z		&&			\ref{DefUnreduzierbar} & \T{DefUnreduzierbar}\\ \hline
\Eintrag{beliebige Menge} { \MengeE}  						{MengeE}	{DefMenge}
\Eintrag{Normalform}{\Nf(x)}										{Nf}			{DefNormalform}
\Eintrag{Menge der Unreduzierbaren Elemente}{\Nf}{Nf}{NgleichTN}
\Eintrag{Church-Rosser-äquivalent}{x \RedAq \bar{x}}{RedAq}{DefChurchRosserAq}
\Eintrag{$\Red$-Prädikat}{\P}{P}{DefRPraedikat}
\T{Teilmenge $\wahr$} &{\Rp}					& \backslash \T{Rp}&\ref{DefRp}&\T{DefRp}\\
											&{\MengeE_{\P}} & 		&\ref{DefRp}&\T{DefRp}\\ 
								\hline
\T{Äquivalenzklasse}& {\Aq{ ~\cdot~ }}&\backslash \T{Aq}&\ref{DefAqR}&\T{DefAqR}\\
										& {\Aq{ ~\cdot~ }}_{\P} &&\\
								\hline
\Eintrag{$\Rp$-äquivalent}{\RpAq}{RpAq}{DefRpAq}
\Eintrag{Körper}{\lambda \in \K}{K}{secDefAlgebra}
\Eintrag{Alphabeth}{a\in\Alphabet}{Alphabet}{secDefAlgebra}
\Eintrag{Terme}{\lambda w \in \Terme }{Terme}{DefTerm}
\Eintrag{Monome}{w \in \Monome}{Monome}{DefMonom}
\T{freie Algebra}	& {f \in \FreieA} &\backslash\T{FreieA} &\ref{DefFreieAlgebra} &\T{DefFreieAlgebra} \\
									& \FreieA_n &\backslash\T{FreieA}  &&\\ 
									& \FreieAA & \backslash\T{FreieAA }&&\\
									& \FreieAM  & \backslash\T{FreieAM }&&\\
						\hline
\Eintrag{Algebra}{\Algebra }{Algebra}{DefAlgebra}
\Eintrag{Wortersetzungssystem}{\red}{red}{DefWortersetzungssystem}
\Eintrag{vollständiges $\red$}{\vRS}{vRS}{}
\Eintrag{}{\vRSa}{vRSa}{}
\Eintrag{}{\RSm}{vRSm}{}
\Eintrag{Rig}{\RigF}{RigF}{DefRigF}
\Eintrag{Modul}{\Modul}{Modul}{}
\Eintrag{Unter Modul}{\UModul}{UModul}{}
\Eintrag{}{\AoA}{AoA}{}
\Eintrag{}{\AnzA}{AnzA}{}
\Eintrag{}{\AnzM}{AnzM}{}
\Eintrag{}{\AnzUM}{AnzUM}{}
\end{array}
}